\documentclass[a4paper, 10pt]{article}
\pdfoutput=1

\usepackage[utf8]{luainputenc}
\usepackage[USenglish]{babel}
\usepackage{csquotes}

\usepackage[a4paper, top=0.91in, bottom=1.29in, left=0.599in, right=0.59in]{geometry}  

\usepackage[plainpages=false,pdfpagelabels,hidelinks,unicode]{hyperref}
\makeatletter
\hypersetup{pdfauthor={\@author}}
\hypersetup{pdftitle={\@title}}
\makeatother

\usepackage[%
  backend=bibtex,
  bibencoding=ascii,
  firstinits=true, uniquename=init, 
  natbib=true,
  url=false,
  doi=false,
  isbn=false,
  backref=false,
  maxnames=99,
  ]{biblatex}
\usepackage{amsmath}
\allowdisplaybreaks
\usepackage{amssymb}
\usepackage{commath}
\usepackage{mathtools}
\usepackage{bbm}
\usepackage{nicefrac}

\usepackage{siunitx} 
\usepackage{amsthm}
\theoremstyle{plain}
  \newtheorem{theorem}{Theorem}[section]
  
\theoremstyle{definition}
  \newtheorem{definition}[theorem]{Definition}
  \newtheorem{remark}[theorem]{Remark}
  \newtheorem{example}[theorem]{Example}

\usepackage{ifluatex}
\ifluatex
  \usepackage[no-math]{fontspec}
\else
  \usepackage[T1]{fontenc}
\fi
\usepackage{newpxtext,newpxmath}

\usepackage{color}
\usepackage{graphicx}
\usepackage[small]{caption}
\usepackage{subcaption}

\makeatletter
\DeclareOldFontCommand{\rm}{\normalfont\rmfamily}{\mathrm}
\DeclareOldFontCommand{\sf}{\normalfont\sffamily}{\mathsf}
\DeclareOldFontCommand{\tt}{\normalfont\ttfamily}{\mathtt}
\DeclareOldFontCommand{\bf}{\normalfont\bfseries}{\mathbf}
\DeclareOldFontCommand{\it}{\normalfont\itshape}{\mathit}
\DeclareOldFontCommand{\sl}{\normalfont\slshape}{\@nomath\sl}
\DeclareOldFontCommand{\sc}{\normalfont\scshape}{\@nomath\sc}
\makeatother

\begingroup\expandafter\expandafter\expandafter\endgroup
\expandafter\ifx\csname pdfsuppresswarningpagegroup\endcsname\relax
\else
\fi

\usepackage{booktabs}
\usepackage{rotating}
\usepackage{multirow}

\usepackage{multicol}
\usepackage{enumitem}

\usepackage{calc}
\usepackage{xparse}

\renewcommand{\vec}[1]{\underline{#1}}
\NewDocumentCommand{\mat}{mo}{%
  \IfValueTF{#2}{%
    \underline{\underline{#1}}{#2}
  }{%
    \underline{\underline{#1}}\,
  }%
}
\newcommand{\diag}[1]{\operatorname{diag}\left(#1\right)}

\AtBeginDocument{%
  
}

\DeclarePairedDelimiterX\newset[1]\lbrace\rbrace{\setaux #1||\endsetaux}
\def\setaux#1|#2|#3\endsetaux{\if\relax\detokenize{#2}\relax #1 \else #1 \;\delimsize\vert\; #2 \fi}
\renewcommand{\set}[1]{\newset*{#1}}

\newcommand{\interior}[1]{\mathring{#1}}

\newcommand{\vect}[1]{\begin{pmatrix} #1 \end{pmatrix}}

\newcommand{\xmin}{{x_L}}
\newcommand{\xmax}{{x_R}}

\newcommand{\I}{\operatorname{I}}
\newcommand{\fnum}{f^{\mathrm{num}}}

\newcommand{\vecfnum}{\vec{f}^{\mathrm{num}}}

\newcommand{\fvol}{f^{\mathrm{vol}}}

\newcommand{\dt}{\Delta t}

\newcommand{\VOL}{\vec{\mathrm{VOL}}}

\newcommand{\SAT}{\vec{\mathrm{SAT}}}
\newcommand{\DISS}{\vec{\mathrm{DISS}}}

\NewDocumentCommand{\ENO}{o}{%
  \IfValueTF{#1}{%
    ENO$(#1)$
  }{%
    ENO
  }%
}
\NewDocumentCommand{\mENO}{o}{%
  \IfValueTF{#1}{%
    mENO$(#1)$
  }{%
    mENO
  }%
}

\AtBeginDocument{%
  \let\rho\varrho
  \let\phi\varphi
  \let\epsilon\varepsilon
}

\newcommand{\R}{\mathbb{R}}

\renewcommand{\emptyset}{\varnothing}

\newsavebox{\DelimiterBox}
\newlength{\DelimiterHeight}
\newlength{\DelimiterDepth}
\newsavebox{\ArgumentBox}
\newlength{\ArgumentHeight}
\newlength{\ArgumentDepth}
\newlength{\ResizedDelimiterHeight}
\newlength{\ResizedDelimiterDepth}

\newcommand{\encloseby}[3]{%
  \savebox{\ArgumentBox}{$\displaystyle #1$}%
  \settoheight{\ArgumentHeight}{\usebox{\ArgumentBox}}%
  \settodepth{\ArgumentDepth}{\usebox{\ArgumentBox}}%
  \savebox{\DelimiterBox}{#2}%
  \settoheight{\DelimiterHeight}{\usebox{\DelimiterBox}}%
  \settodepth{\DelimiterDepth}{\usebox{\DelimiterBox}}%
  \setlength{\ResizedDelimiterHeight}{%
    \maxof{1.2\ArgumentHeight}{\DelimiterHeight}%
  }
  \setlength{\ResizedDelimiterDepth}{%
    \maxof{1.2\ArgumentDepth}{\DelimiterDepth}%
  }
  \raisebox{-\ResizedDelimiterDepth}{%
    \resizebox{\width}{\ResizedDelimiterHeight+\ResizedDelimiterDepth}{%
      \raisebox{\DelimiterDepth}{#2}%
    }%
  }
  #1
  \raisebox{-\ResizedDelimiterDepth}{%
    \resizebox{\width}{\ResizedDelimiterHeight+\ResizedDelimiterDepth}{%
      \raisebox{\DelimiterDepth}{#3}%
    }%
  }
}

\ifluatex
  \newcommand{\mean}[1]{\encloseby{#1}{$\{\mkern-5mu\{$}{$\}\mkern-5mu\}$}}
  \newcommand{\jump}[1]{\encloseby{#1}{$[\mkern-4mu[$}{$]\mkern-4mu]$}}
\else
  \newcommand{\mean}[1]{\encloseby{#1}{$\{\mkern-6mu\{$}{$\}\mkern-6mu\}$}}
  \newcommand{\jump}[1]{\encloseby{#1}{$[\mkern-3mu[$}{$]\mkern-3mu]$}}
\fi

\newcommand{\logmean}[1]{\mean{#1}_\mathrm{log}}

\begin{document}

\title{\sl Kinetic functions for nonclassical shocks, 
\\
entropy stability, and discrete summation by parts}
\author{Philippe G. LeFloch\footnote{Laboratoire Jacques-Louis Lions and Centre National de la Recherche Scientifique,
Sorbonne Universit\'e, 
4 Place Jussieu, 75252 Paris, France. Email: {\sl contact@philippelefloch.org}.}
\, 
and Hendrik Ranocha\footnote{King Abdullah University of Science and Technology (KAUST), Computer Electrical and Mathematical Science and Engineering Division (CEMSE), Thuwal, 23955-6900, Saudi Arabia.
Present address:
Applied Mathematics M\"unster,
University of M\"unster,
M\"unster, Germany.
Email: {\sl mail@ranocha.de.}}
}
\date{December 13, 2020}

\maketitle

\begin{abstract}
We study nonlinear hyperbolic conservation laws with non-convex flux
in one space dimension and, for a
broad class of numerical methods based on summation by parts operators, we
compute numerically the kinetic functions associated with each scheme.
As established by LeFloch and collaborators, kinetic functions (for continuous
or discrete models) uniquely characterize the macro-scale dynamics of small-scale
dependent, undercompressive, nonclassical shock waves. We show here that various
entropy-dissipative numerical schemes can yield nonclassical solutions containing
classical shocks, including Fourier methods with (super-) spectral viscosity,
finite difference schemes with artificial dissipation, discontinuous Galerkin
schemes with or without modal filtering, and TeCNO schemes.  We demonstrate
numerically that entropy stability does not imply uniqueness of the limiting
numerical solutions for scalar conservation laws in one space dimension,
and we compute the associated kinetic functions in order to distinguish between
these schemes.
In addition, we design entropy-dissipative schemes for the Keyfitz-Kranzer system
whose solutions are measures with delta shocks. This system illustrates the fact
that entropy stability does not imply boundedness under grid refinement.

\end{abstract}

\section{Introduction}
\label{sec:introduction}

\subsection{Objective and background}

We present and analyze here several classes of entropy-stable and semi-discrete schemes for nonlinear hyperbolic problems,
next investigate numerically the behavior of weak solutions, and demonstrate certain important features or limitations of these schemes. In particular, we observe that entropy-stable schemes can converge to {\sl different weak solutions.}
Hence, we are interested in qualitative properties of weak solutions $u=u(t,x)$ to nonlinear hyperbolic conservation laws
\begin{equation}
\label{eq:HCL-intro}
  \partial_t u + \partial_x f(u) = 0,
\qquad
  u|_{t=0} = u_0,
\end{equation}
posed on a bounded domain $\Omega \subset \R$ (subjected to suitable boundary conditions).
Here, $u_0$ is a prescribed initial data defined on $\Omega$, while the flux $f=f(u)$ is a prescribed nonlinear function of the unknown $u$.
In general, weak solutions to \eqref{eq:HCL-intro} (understood in the sense of distributions) contain shock waves, and
a central issue in the theory of nonlinear hyperbolic equations is formulating
suitable admissibility criteria for the selection of shocks.  For scalar
conservation laws with convex flux (such as Burgers' equation) a single entropy
inequality
\begin{equation}
\label{eq:entropy-intro}
  \partial_t U(u) + \partial_x F(u) \leq 0,
\end{equation}
associated with a strictly convex entropy pair $(U,F)$, suffices
to single out a unique weak solution \cite{delellis2004minimal,panov1994uniqueness} to the initial value problem \eqref{eq:HCL-intro}.

However, this is not true for conservation laws with non-convex flux
\cite[Remark~2]{panov1994uniqueness}, for instance for the cubic conservation law
$\partial_t u + \partial_x u^3 = 0$ \cite[Chapter~II]{lefloch2002hyperbolic}.
While it is well-known that weak solutions to conservation laws can be generated as vanishing
viscosity limits, that is, as limits (when $\epsilon \to 0$) of solutions to
\begin{equation}
  \partial_t u + \partial_x f(u) = \epsilon \partial_x^2 u,
\end{equation}
other regularization operators are equally relevant in physical applications
and
generate weak solutions that may not satisfy the standard selection criteria.
In fact, the single entropy inequality \eqref{eq:entropy-intro}
permits also {\sl nonclassical shocks} of undercompressive type, which in turn should be
characterized via the notion of a {\sl kinetic function}; for an overview of the theory, see LeFloch
\cite{lefloch1999introduction,lefloch2002hyperbolic,lefloch2010kinetic}.
The role of small-scale effects in weak solutions (for instance when capillarity effect is included) and the numerical approximation of nonclassical solutions were extensively investigated in the past two decades; see the pioneering papers \cite{hou1994nonconservative,hayes1997nonclassical,hayes1998nonclassical}, as well as the advances
in \cite{lefloch2008many,castro2008many,lefloch2012many}
and the references therein. More recently, the class of well-controlled dissipation (WCD) schemes which capture diffusive-dispersive shocks at any arbitrary order of accuracy, was proposed in \cite{lefloch2014numerical,leflochtesdall-two,leflochtesdall-one}.

\subsection{Main contributions in this paper}

We proceed here by presenting first a broad class of numerical methods which are based on summation-by-parts (SBP) operators and, importantly, are entropy-satisfying in the sense that they satisfy a discrete form of the entropy inequality.
Recall that entropy-conservative schemes were constructed in a pioneering work by Tadmor~\cite{tadmor1987numerical,tadmor2003entropy} (second-order accuracy) and LeFloch and Rohde \cite{lefloch2000high} (third-order accuracy), and later extended
in \cite{lefloch2002fully} (high-order accuracy in a periodic domain) and \cite{fisher2013high,ranocha2018comparison,chen2017entropy} (bounded domain and non-uniform grids).
The entropy inequality is essential since it ensures a fundamental $L^2$ nonlinear stability property, but does not guarantee  convergence to classical entropy solutions.
The entropy-satisfying numerical methods developed and analyzed in the present
paper are based on the notion of SBP operators in \cite{kreiss1974finite},
which recently have gained a lot of interest. Nowadays, many classes of
numerical methods can be formulated within a unifying SBP framework, including
finite difference \cite{strand1994summation},
finite volume \cite{nordstrom2001finite},
finite element \cite{gassner2013skew,hicken2020entropy}, and flux reconstruction
schemes \cite{ranocha2016summation}. Importantly, these semi-discretizations
can be made entropy-satisfying and the main stability estimates can also be
transferred to fully discrete numerical methods by applying a relaxation approach
\cite{ketcheson2019relaxation,ranocha2020relaxation,ranocha2020general,ranocha2020relaxationHamiltonian,ranocha2020fully,ranocha2020broad,mitsotakis2020conservative}.

It is precisely our purpose here to demonstrate that a {\sl
broad class of entropy-satisfying schemes generate nonclassical shocks} and,
in addition, to compute the associated kinetic functions.  We focus on the
spatial part and apply sufficiently accurate time integration schemes in order
to eliminate any significant errors from that part.
We restrict convergence studies of the numerical methods to numerical experiments
and grid refinement studies. In particular, we observe convergence of numerical
solutions of nonlinear scalar conservation laws in one space dimension,
in contrast to the behavior of nonlinear systems in multiple space dimensions
\cite{fjordholm2016computation,fjordholm2017construction}. However, even when
the numerical methods satisfy a single entropy inequality, the numerically
converged solutions can still approximate nonclassical shocks.

\subsection{Outline of this paper}

In Section~\ref{sec:SBP}, entropy-stable discretizations based on SBP operators are recalled and discussed.
Next, in Sections~\ref{sec:cubic} and \ref{sec:cubic-TeCNO}, we revisit the
uniqueness issue for nonlinear conservation laws and,
in Section~\ref{sec:kinetic-function-cubic}, for a variety of schemes we numerically compute the corresponding
kinetic function associated with the cubic conservation law.
This study is then extended
to a quartic conservation law in Section~\ref{sec:kinetic-function-quartic}.
Finally, in Section~\ref{section7} we turn our attention toward the Keyfitz-Kranzer system, for which we develop and apply a broad class of entropy-stable schemes.


\section{Summation-by-parts operators and entropy stability}
\label{sec:SBP}

\subsection{Notation}

In this section, some general notions about SBP operators are reviewed and a
notation to be used throughout the following sections is introduced.
For more information on SBP operators, we refer to the review articles
\cite{svard2014review,fernandez2014review} and references cited therein.
We consider the nonlinear hyperbolic system of conservation laws
\begin{equation}
\label{eq:HCL}
\begin{aligned}
  \partial_t u(t,x) + \partial_x f\bigl( u(t,x) \bigr) &= 0,
    && t \in (0, T),\, x \in \Omega,
  \\
  u(0,x) &= u_0(x),
    && x \in \Omega,
\end{aligned}
\end{equation}
posed
on the spatial domain $\Omega \subseteq \R$ in one space dimension, supplemented with appropriate boundary data or periodic boundary conditions.
Here, $u\colon (0,T) \times \Omega \to \Upsilon \subseteq \R^m$ are the conserved
variables and $f\colon \Upsilon \to \R^m$ is referred to as the flux.
Using the method of lines, a semi-discretization is introduced at first and a
suitable time integration scheme is applied to the resulting set of ordinary
differential equations, e.g.\ a Runge-Kutta method.
For the semi-discretization, the spatial domain $\Omega \subseteq \R$ is divided
into non-overlapping elements $\Omega_l$, i.e.~we have $\bigcup_l \Omega_l = \Omega$
and $\interior\Omega_l \cap \interior\Omega_k = \emptyset$ if $l \neq k$, where
$\interior\Omega_l$ is the interior of the element $\Omega_l$. In each element
the numerical solution is represented by its values $\vec{u} = ( \vec{u}_1,
\dots, \vec{u}_K )^T$ on a grid with nodes $x_1, \dots, x_K $, i.e.\ we write $\vec{u}_i = u(x_i)$.
All nonlinear operations of interest are then defined pointwise, for instance $\vec{f}_i = f(\vec{u}_i)$.
The discrete derivative will be represented by a matrix $\mat{D} \in \R^{K \times K}$, referred to as the derivative matrix, and the discrete scalar product approximating the standard $L^2$ scalar product will be represented by a symmetric and positive definite matrix, denoted by
$\mat{M} \in \R^{K \times K}$ (mass/norm matrix).
In the following, the notation in \cite{ranocha2016summation,ranocha2017extended} is used. A table with translation rules to other common notations for finite difference and spectral element methods can be found in \cite{ranocha2018generalised}.

\subsection{Periodic setting}

Consider the domain $\Omega = (\xmin, \xmax)$ and absolutely continuous functions
$u,v\colon \R \to \R$ that are periodic with period $\abs{\Omega} = \xmax - \xmin$.
Then, integration by parts gives
\begin{equation}
  \int_\xmin^\xmax u \, (\partial_x v)
  =
\big[  u v \big]_\xmin^\xmax
  - \int_\xmin^\xmax (\partial_x u) \, v,
\end{equation}
where the boundary term vanishes due to periodicity. This is mimicked at the discrete level by writing
$
  \vec{u}^T \mat{M} \mat{D} \vec{v}
  =
  - \vec{u}^T \mat{D}[^T] \mat{M} \vec{v}.
$
Requiring that this relation holds for all $\vec{u},\vec{v} \in \R^K$ results in
the following definition.
With some abuse of notation, the derivative matrix $\mat{D}$ will also be called
  SBP (derivative) operator if the corresponding mass/norm matrix can be deduced
  from the context.

\begin{definition}
\label{def:SBP-periodic}
  A \emph{periodic SBP operator} (or summation-by-parts operator) consists of a derivative matrix $\mat{D}$ and
  a symmetric and positive definite mass/norm matrix $\mat{M}$ such that
$
    \mat{M} \mat{D} + \mat{D}[^T] \mat{M} = \mat{0}.
$
\end{definition}

\begin{example}
\label{ex:periodic-FD}
  Periodic central finite difference approximations of the first derivative yield
  SBP operators with mass matrix $\mat{M} = \Delta x \mat{\I}$, where $\mat{\I}$ is the
  identity matrix and $\Delta x$ the grid spacing. In other words, periodic
  central finite difference approximations of the first derivative result
  in skew-symmetric derivative matrices $\mat{D}$.
\end{example}

\begin{example}
\label{ex:Fourier}
  Fourier (pseudo-) spectral methods computing the derivative via the discrete
  Fourier transform (via FFT) also yield SBP operators with a multiple of
  the identity matrix as mass matrix.
\end{example}

\subsection{Non-periodic setting}

Consider again the domain $\Omega = (\xmin, \xmax)$ and absolutely continuous functions
$u,v\colon \R \to \R$. Without requiring periodicity, integration by parts gives
\begin{equation}
  \int_\xmin^\xmax u \, (\partial_x v) + \int_\xmin^\xmax (\partial_x u) \, v
  =
\big[  u v \big]_\xmin^\xmax.
\end{equation}
Similarly to the periodic case, this is mimicked at the discrete level by
\begin{equation}
\label{eq:SBP-motivation}
  \vec{u}^T \mat{M} \mat{D} \vec{v} + \vec{u}^T \mat{D}[^T] \mat{M} \vec{v}
  =
  \vec{u}^T \mat{R}[^T] \mat{B} \mat{R} \vec{v},
\end{equation}
where a restriction matrix $\mat{R}$ performing an interpolation to the boundary
nodes $\xmin, \xmax$ and a boundary matrix $\mat{B} = \diag{-1,1}$ have been
introduced. Requiring that the relation above holds for all $\vec{u},\vec{v} \in \R^K$
results in the following definition.

\begin{definition}
\label{def:SBP}
  A \emph{(non-periodic) SBP operator} (i.e.\ summation-by-parts operator) consists of a derivative matrix $\mat{D}$,
  a symmetric and positive definite mass/norm matrix $\mat{M}$, a restriction
  matrix $\mat{R}$, and the boundary matrix $\mat{B} = \diag{-1,1}$ such that
$
    \mat{M} \mat{D} + \mat{D}[^T] \mat{M} = \mat{R}[^T] \mat{B} \mat{R}.
$
  The order of accuracy of the approximation
$
    \vec{u}^T \mat{R}[^T] \mat{B} \mat{R} \vec{v}
    \approx
    [ u v ]_\xmin^\xmax
$
  should be at least the order of accuracy of the derivative matrix $\mat{D}$.
\end{definition}

SBP operators can be used both in single element discretizations and in multiple
element discretizations where the domain $\Omega$ is divided into smaller elements
as described at the beginning of this section. In the latter case, SBP operators
are used on each element. Typically, the operators are developed on a reference
element and a coordinate transformation is used to map all quantities between
the physical and the reference element. If not stated otherwise, affine-linear
coordinate mappings will be used in the following.

\begin{example}
\label{ex:classical-SBP}
\begin{enumerate}
  \item Classical finite difference SBP operators have been proposed by many researchers and many examples can be found in the review articles
  \cite{svard2014review,fernandez2014review} and references cited therein.

  \item
  Polynomial collocation methods based on Legendre polynomials of degree $p = K-1$
  yield SBP operators if the derivative matrix $\mat{D}$ and the restriction matrix
  $\mat{R}$ are exact for polynomials of degree $\leq p$ and the mass matrix $\mat{M}$
  is chosen such that it is exact for polynomials of degree $\leq p-1$, since
  all integrals in \eqref{eq:SBP-motivation} can be evaluated exactly in this case,
  cf.\ \cite{kopriva2010quadrature}. In particular, polynomial collocation methods
  based on the Gauss, Radau, or Lobatto Legendre nodes yield SBP operators,
  see also \cite{gassner2013skew,fernandez2014generalized,hicken2016multidimensional}.
\end{enumerate}
\end{example}

\subsection{Boundary procedures}

If periodic SBP operators are used (for periodic problems), no boundary conditions
have to be enforced. If non-periodic boundaries or multiple elements are considered,
exterior (at the boundary of the domain $\Omega$) or interior (between two elements)
boundary conditions have to be enforced. Here, the prescription of boundary
conditions via simultaneous approximation terms using numerical fluxes as in
finite volume methods will be used.

\begin{definition}
\label{def:fnum}
  A numerical flux is a Lipschitz continuous mapping $\fnum\colon \Upsilon^2 \to \R^m$
  that is consistent with the flux $f$ of the conservation law \eqref{eq:HCL},
  i.e.\ $\text{ for all } u \in \Upsilon\colon \fnum(u,u) = f(u)$.
\end{definition}

\begin{figure}[th]
  \centering
  \includegraphics{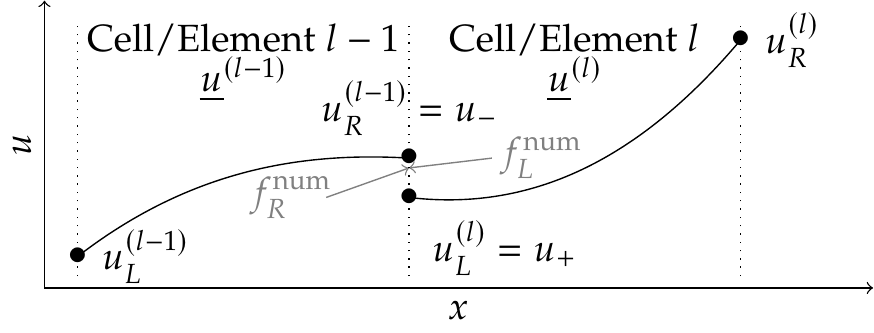}
  \caption{Visualization of the notation used for multiple element discretizations
           using SATs with numerical fluxes.
           The numerical flux at the boundary between the elements $l-1$ and $l$
           is computed as $\fnum(u_-,u_+)$, where $u_-$ ($u_+$) is the value of
           the numerical solution to the left (right) of the interface. In element
           $l-1$ ($l$), this numerical flux is at the right (left) hand side of
           the element and called $\fnum_R$ ($\fnum_L$).}
  \label{fig:notation-numerical-fluxes}
\end{figure}

A semi-discretization of the conservation law $\partial_t u + \partial_x f(u) = 0$
on an element will usually be of the form
\begin{equation}
  \partial_t \vec{u} + \VOL = \SAT,
\end{equation}
where $\VOL$ are volume terms approximating the divergence $\partial_x f(u)$
in the interior and $\SAT$ is a simultaneous approximation term
\cite{carpenter1994time,carpenter1999stable}
enforcing boundary conditions in a suitably weak sense. Typical forms might be
\begin{equation}
\VOL = \mat{D} \vec{f},
\qquad
\SAT = - \mat{M}[^{-1}] \mat{R}[^T] \mat{B} \left( \vecfnum - \mat{R} \vec{f} \right),
\end{equation}
where $\mat{R} \vec{f} = (f_L, f_R)^T$ is the interpolation of the flux $\vec{f}$
to the boundary and $\vecfnum = (\fnum_L, \fnum_R)^T$ is the numerical flux
at the boundaries. At a given boundary between the elements $l$ and $l+1$, two values
of the numerical solution are given via interpolation: The value $u_- = u^{(l)}_R$
of the numerical solution at the right-hand side of cell $l$ and the value
$u_+ = u^{(l+1)}_L$ of the numerical solution at the left-hand side of cell $l+1$.
Then, the unique numerical flux between the elements $l$ and $l+1$ is computed as
$\fnum(u_-,u_+)$; see the illustration in \autoref{fig:notation-numerical-fluxes}.

\subsection{Entropy stability}

Throughout, we assume that our system of conservation laws \eqref{eq:HCL} is
equipped with an entropy function.
%
Recall that a convex function $U\colon \Upsilon \to \R$ is an \emph{entropy} for the
  conservation law \eqref{eq:HCL}, if there is an entropy flux $F\colon \Upsilon \to \R$
  fulfilling $\partial_u U(u) \cdot \partial_u f(u) = \partial_u F(u)$.
  The \emph{entropy variables} are $w(u) = U'(u)$ and the \emph{flux potential}
  is $\psi = w \cdot f - F$.
Thus, if $U$ is an entropy and $u$ a smooth solution of the conservation law
\eqref{eq:HCL}, the entropy generates the conservation law
$\partial_t U(u) + \partial_x F(u) = 0$. As an admissibility criterion,
the entropy inequality
\begin{equation}
\label{eq:entropy-ineq}
  \partial_t U(u) + \partial_x F(u) \leq 0.
\end{equation}
is imposed for weak solutions. The flux potential $\psi$ is the potential
of the flux $f$ with respect to the entropy variables $w$, i.e.\
$\partial_w \psi(w) = f(w)$, where $f(w)$ should be read as $f\bigl( u(w) \bigr)$.

\begin{example}
\label{ex:L2-entropy}
  For scalar conservation laws, the $L^2$ entropy $U(u) = \frac{1}{2} u^2$ can
  be used. It is strictly convex and the entropy variables are $w(u) = U'(u) = u$.
  The flux potential is given by $\psi(u) = \int^u f(v) \dif v$.
\end{example}

Building on the seminal work of Tadmor~\cite{tadmor1987numerical,tadmor2003entropy}
on second-order schemes, extended to arbitrary order of accuracy in
\cite{lefloch2000high},
semi-discrete entropy-stable schemes can be constructed as entropy-conservative
schemes and additional dissipation mechanisms. The basic ingredient are entropy
conservative two-point numerical fluxes used in finite volume methods.

\begin{definition}
\label{def:fnum-EC-ES}
  A numerical flux is \emph{entropy-conservative} (EC) for an entropy $U$ with
  entropy variables $w$ and flux potential $\psi$, if
  \begin{equation}
  \label{eq:fnum-EC}
    \text{ for all } u_-,u_+ \in \Upsilon\colon \quad
    \bigl( w(u_+) - w(u_-) \bigr) \cdot \fnum(u_-,u_+) = \psi(u_+) - \psi(u_-).
  \end{equation}
  The numerical flux is \emph{entropy-stable} (ES), if
  \begin{equation}
  \label{eq:fnum-ES}
    \text{ for all } u_-,u_+ \in \Upsilon\colon \quad
    \bigl( w(u_+) - w(u_-) \bigr) \cdot \fnum(u_-,u_+) \leq \psi(u_+) - \psi(u_-).
  \end{equation}
\end{definition}
Using the jump operator $\jump{a} := a_+ - a_-$ and the notation $a_\pm = a(u_\pm)$,
these equations can be written as $\jump{w} \cdot \fnum = \jump{\psi}$ and
$\jump{w} \cdot \fnum \leq \jump{\psi}$, respectively.

\begin{remark}
\label{rem:entropy-dissipative}
  Entropy stability of a numerical method does not imply general stability of
  the scheme, since further robustness properties might be necessary to guarantee
  that the numerical solution does not blow up. Moreover, for linear equations,
  (strong) stability is the discrete analogue of (strong) well-posedness of the
  continuous problem, cf.\ \cite{svard2014review,nordstrom2017roadmap}. This is
  in general not the case for nonlinear equations. Thus, it might be better to
  speak about \emph{entropy-conservative} and \emph{entropy-dissipative} schemes.
  Nevertheless, since it is more common in the literature to speak about
  \emph{entropy-stable} schemes, this term will be used here.
\end{remark}

Entropy-conservative numerical fluxes can be used to construct high-order entropy
conservative semi-discretizations using SBP operators via the volume terms
\begin{equation}
\label{eq:volume-terms-flux-diff}
  \VOL_i
  =
  \sum_{k=1}^K {\mat{D}}_{i,k} 2 \fvol(\vec{u}_i, \vec{u}_k),
\end{equation}
where $\fvol$ is an entropy-conservative numerical flux called \emph{volume flux}
(since it it used for the volume terms)
and $\vec{u}_i, \vec{u}_k$ are the values of the discrete solution at the grid
nodes. The following result can be found in
\cite{ranocha2018comparison,chen2017entropy} and is a generalisation of
\cite[Theorem 3.1]{fisher2013high}.
\begin{theorem}
\label{thm:flux-diff-order}
  If the volume flux $\fvol$ is smooth and symmetric, the flux differencing form
  \eqref{eq:volume-terms-flux-diff} is an approximation of the same order of
  accuracy as the derivative matrices $\mat{D}$.
  Here, the derivative matrix does not need to be given by SBP operators.
\end{theorem}

If the corresponding mass matrix $\mat{M}$ of an SBP derivative $\mat{D}$ is
diagonal, the volume terms \eqref{eq:volume-terms-flux-diff} can be written in a
locally conservative form. Moreover, if the volume flux is entropy-conservative
and the boundary nodes are included, high-order entropy-stable semi-discretizations
can be constructed, cf.\ \cite{fisher2013high,ranocha2018comparison,chen2017entropy}.
\begin{theorem}
\label{thm:flux-diff-ES}
  Consider the semi-discretization $\partial_t \vec{u} + \VOL = \SAT$ with
  volume terms given by \eqref{eq:volume-terms-flux-diff} and the surface terms
  \begin{equation}
  \label{eq:surface-terms-flux-diff}
    \SAT
    =
    - \mat{M}[^{-1}] \mat{R}[^T] \mat{B} \left( \vecfnum- \mat{R} \vec{f} \right).
  \end{equation}
  If the numerical volume flux $\fvol$ is consistent with $f$, symmetric, and
  entropy-conservative, and both the mass matrix $\mat{M}$ and the boundary
  operator $\mat{R}[^T] \mat{B} \mat{R}$ are diagonal, the semi-discretization is
  entropy-conservative/stable across elements, if the numerical surface flux
  $\fnum$ is entropy-conservative/stable.
  Moreover, there is a locally conservative form for the semi-discrete entropy
  equation.
\end{theorem}

\begin{remark}
  Requiring that $\mat{R}[^T] \mat{B} \mat{R}$ be diagonal seems to be
  necessary for general conservation laws. Basically, it states that the boundary
  nodes have to be included in the computational grid.
  For some conservation laws, the surface terms \eqref{eq:surface-terms-flux-diff}
  can be adapted to allow general SBP operators not including the boundary nodes,
  e.g.\ polynomial collocation methods on Gauss-Legendre nodes, cf.\
  \cite{ranocha2016summation,ranocha2017shallow,ranocha2018generalised,offner2019stability,ranocha2018thesis}.
\end{remark}

\subsection{Dissipation operators}

A semi-discretization $\partial_t \vec{u} + \VOL = \SAT$ can be enhanced by
artificial dissipation terms $\DISS$, resulting in
\begin{equation}
  \partial_t \vec{u} + \VOL = \SAT + \DISS.
\end{equation}
Since entropy stability is investigated by multiplying the semi-discretization
by $\vec{w}^{T} \mat{M}$, the dissipation term $\DISS$ should fulfil
$\vec{w}^{T} \mat{M} \DISS \leq 0$. Moreover, in order not to influence conservation
across elements, the dissipation term should satisfy $\vec{1}^T \mat{M} \DISS = \vec{0}$,
where $\vec{1}$ is a vector with entries $1$, i.e.\ the discrete version of the function
$x \mapsto 1$.

\begin{example}
\label{ex:diss-L2-FD}
  Considering the $L^2$ entropy $U(u) = \frac{1}{2} u^2$ as in Example~\ref{ex:L2-entropy},
  dissipation operators can be constructed as $\DISS = -\mat{M}[^{-1}] \mat{S}$,
  where $\mat{S}$ is (symmetric and) positive semidefinite and satisfies
  $\vec{1}^T \mat{S} = \vec{0}$. Such dissipation operators approximating weighted
  derivatives of even degree have been proposed in \cite{mattsson2004stable}.
\end{example}

\begin{example}
\label{ex:diss-Legendre}
  For collocation schemes using Legendre polynomials, suitable discretizations
  of the Legendre derivative operator $u \mapsto \partial_x ( a \partial_x u)$
  with $a(x) = 1-x^2$ can be used. The Legendre polynomials $\phi_n$ are eigenvectors
  of this operator with eigenvalues $\lambda_n = - n(n+1)$. An investigation using
  this kind of dissipation and SBP operators can be found inter alia in
   \cite{ranocha2018stability}.
\end{example}

\begin{example}
\label{ex:diss-spectral}
  The (super-) spectral viscosity operators investigated in
  \cite{maday1989analysis,tadmor1989convergence,schochet1990rate,maday1993legendre,tadmor1993super,guo2001spectral,tadmor2012adaptive}
  are suitable dissipation operators for the $L^2$ entropy $U(u) = \frac{1}{2} u^2$,
  cf.\ Example~\ref{ex:L2-entropy}.
\end{example}

\begin{example}
\label{ex:diss-TeCNO}
  The so-called TeCNO schemes presented in \cite{fjordholm2012arbitrarily,fjordholm2013high}
  are designed for periodic boundary conditions. The volume
  terms are exactly the entropy-conservative ones \eqref{eq:volume-terms-flux-diff}.
  The dissipation operators are constructed using the ENO procedure and are based
  on a recent stability result of this reconstruction \cite{fjordholm2013eno}.
\end{example}

\subsection{Filtering}

Another possibility introducing dissipation is given by filtering. Here, the
baseline scheme is used to compute one time step (or only one stage of a
Runge-Kutta method) and a dissipative filter is applied afterwards, reducing the
total amount of entropy without modifying the total mass.
Classical modal filters can also be motivated by the application of a splitting
procedure to a semi-discretization enhanced with artificial dissipation. If this
dissipation corresponds to a modal basis such as a Legendre basis, the action of
the dissipation operator can be realized as modal filtering. Splitting a complete
time step into a step using the baseline scheme and an exact integration of the
dissipation operator results in this modal filtering approach. However, other filter
functions can be used as well.

\begin{example}
\label{ex:Legendre-filtering}
  Applying an exponential filter of the form
  $\exp(- \epsilon \, \Delta t \, n (n+1) )$
  corresponds to the exact solution of the equation
  $
  \partial_t u(t,x) = \partial_x \bigl( (1-x^2) \partial_x u(t,x) \bigr),
  $
  as in Example~\ref{ex:diss-Legendre}. Similar to the super-spectral viscosity
  approach, the Legendre dissipation operator can be applied $s$ times,
  resulting in the filter function
  $\exp(- \epsilon \, \Delta t \, (n (n+1))^s )$.
  Results for modal filters can be found in
  \cite{vandeven1991family,hesthaven2008filtering,ranocha2018stability}.
\end{example}


\section{Revisiting the uniqueness issue for conservation laws}
\label{sec:cubic}

\subsection{Periodic boundary conditions}

In the following numerical experiments, the fourth-order, ten stage, strong stability
preserving explicit Runge-Kutta method SSPRK(10,4) of \cite{ketcheson2008highly}
will be used to advance the numerical solutions up to the final time $T=1$.
We treat the cubic conservation law
\begin{equation}
\label{eq:cubic}
\begin{aligned}
  \partial_t u(t,x) + \partial_x u(t,x)^3 &= 0,
    \quad && t \in (0, T),\, x \in (\xmin,\xmax),
  \\
  u(0,x) &= u_0(x),
    \quad && x \in (\xmin,\xmax),
\end{aligned}
\end{equation}
supplemented with periodic or other appropriate boundary conditions as described
in the following. The initial condition is chosen as $u_0(x) = -\sin(\pi x)$ in
the domain $(\xmin,\xmax) = (-1,1)$.
Using periodic boundary conditions and the $L^2$ entropy $U(u) = \frac{1}{2} u^2$,
the total entropy $\int U$ is conserved for smooth solutions and bounded
from above by its initial value for entropy weak solutions. Using the unique
entropy-conservative flux (cf.\ Definition~\ref{def:fnum-EC-ES})
\begin{equation}
\label{eq:cubic-fnum-EC}
  \fvol(u_-,u_+)
  =
  \frac{1}{4} \frac{u_+^4 - u_-^4}{u_+ - u_-}
  =
  \frac{1}{4} \left( u_+^3 + u_+^2 u_- + u_+ u_-^2 + u_-^3 \right)
\end{equation}
in the flux difference discretization described in Theorem~\ref{thm:flux-diff-ES}
results in the volume terms
\begin{equation}
\label{eq:cubic-split-form}
  \VOL
  =
  \frac{1}{2} \mat{D} \mat{u}^2 \vec{u}
  + \frac{1}{2} \mat{u} \mat{D} \mat{u} \vec{u}
  + \frac{1}{2} \mat{u}^2 \mat{D} \vec{u},
\end{equation}
where $\mat{u} = \diag{\vec{u}}$ are diagonal multiplication operators, performing
pointwise multiplication with the values of $\vec{u}$ at the nodes of the grid.
These volume terms can be interpreted as approximations to the split form
$\frac{1}{2} \bigl( \partial_x u^3 + u \partial_x u^2 + u^2 \partial_x u \bigr)$.
Traditionally, one could also use the following unsplit form,
without obtaining an entropy estimate,
\begin{equation}
\label{eq:cubic-unsplit-form}
  \VOL
  =
  \mat{D} \mat{u}[^2] \vec{u}.
\end{equation}

The following semi-discretizations will be used for periodic boundary conditions.
For these schemes, the time step is chosen as $\Delta t = \frac{1}{5N}$, where
$N$ is the number of grid points.

\begin{itemize}
  \item {\bf Periodic (central) finite difference methods.}
  The interior schemes of the dissipation operators proposed in \cite{mattsson2004stable}
  multiplied by a fixed strength $\epsilon$ will be used. These operators
  approximate derivatives of even degree.

  \item {\bf Fourier methods.}
  The spectral viscosity operators described in \cite{tadmor2012adaptive} will
  be used. These operators are described by a cutoff frequency $m = \sqrt{N}$,
  a strength $\epsilon \sim \frac{1}{N}$, and the dissipation coefficients $Q_k$.
  The ``standard'' choice of these coefficients is \cite[eq.~(1.7)]{tadmor2012adaptive}
  \begin{equation}
  \label{eq:tadmor2012adaptive-standard}
    \widehat{Q}_k
    =
    \begin{cases}
      0, & \abs{k} \leq m,
      \\
      \exp\Bigl( -\frac{(N-k)^2}{(k-m)^2} \Bigr), & m < \abs{k}.
    \end{cases}
  \end{equation}
  The ``convergent'' scheme inspired by results of \cite{schochet1990rate} is
  \cite[eq.~(4.3)]{tadmor2012adaptive}
  \begin{equation}
  \label{eq:tadmor2012adaptive-convergent}
    \widehat{Q}_k
    =
    \begin{cases}
      0, & \abs{k} \leq m,
      \\
      \exp\Bigl( -\frac{(2m-k)^2}{(k-m)^2} \Bigr), & m < \abs{k} < 2m,
      \\
      1, & \abs{k} \geq 2m.
    \end{cases}
  \end{equation}
\end{itemize}

\begin{figure}
\centering
  \begin{subfigure}{0.45\textwidth}
    \centering
    \includegraphics[width=\textwidth]{%
      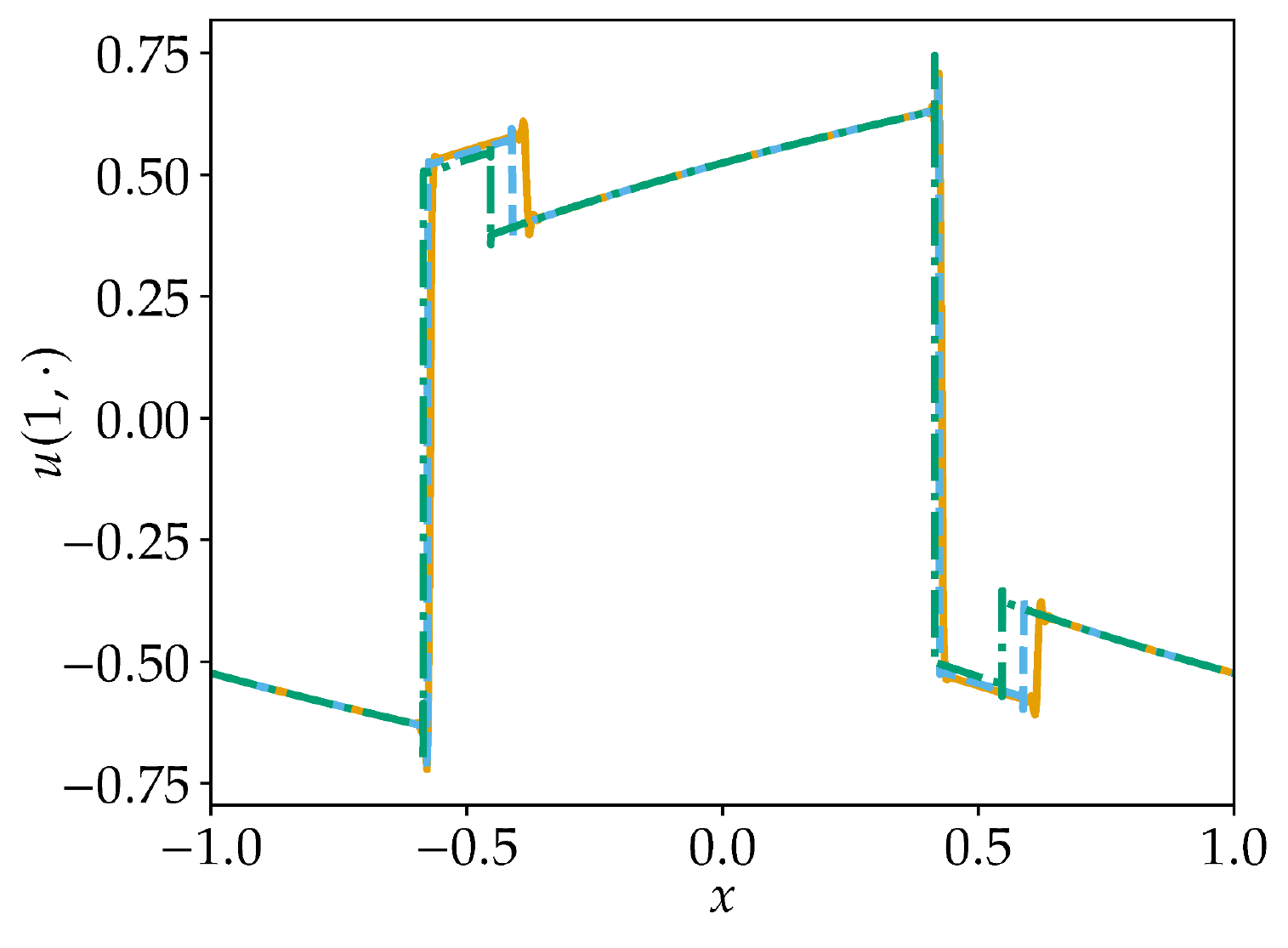}
    \caption{Before postprocessing.}
  \end{subfigure}%
  \hspace*{\fill}
  \begin{subfigure}{0.45\textwidth}
    \centering
    \includegraphics[width=\textwidth]{%
      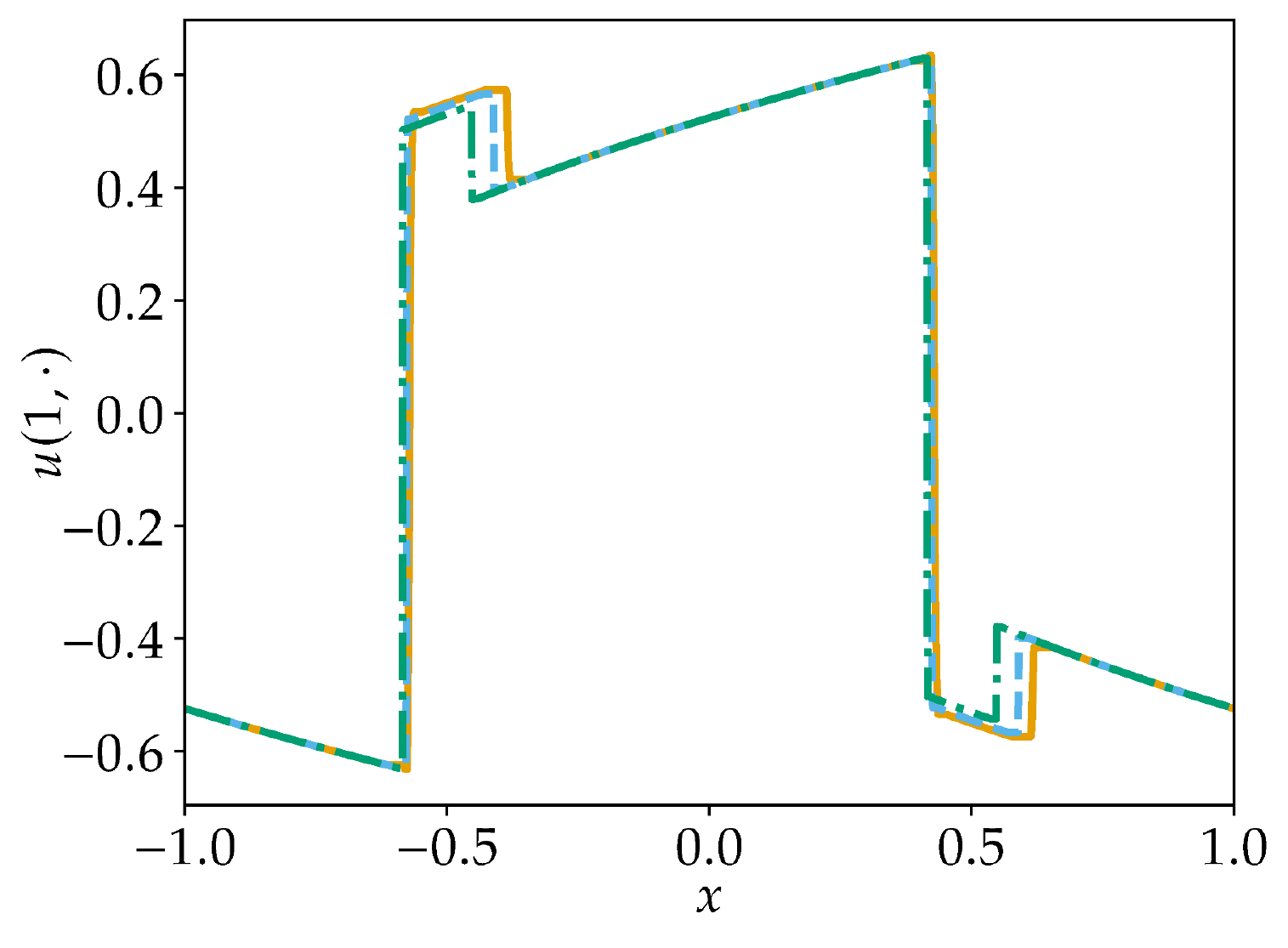}
    \caption{After postprocessing.}
  \end{subfigure}%
  \caption{Numerical solutions shown in
           \autoref{fig:Cubic_PeriodicFD_Split_1_DissOrder_2}
           before and after postprocessing.}
  \label{fig:postprocessing}
\end{figure}

\begin{figure}
\centering
  \begin{subfigure}{0.5\textwidth}
    \centering
    \includegraphics[width=\textwidth]{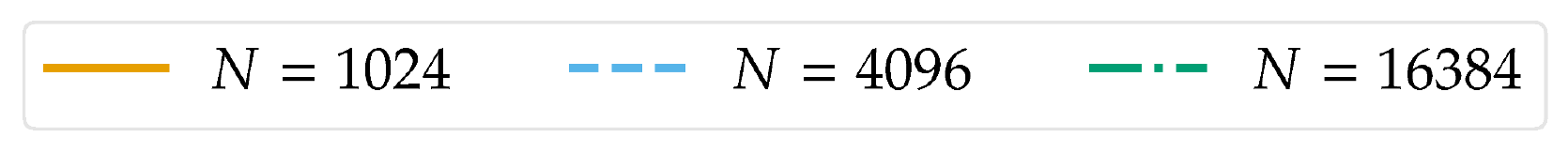}
  \end{subfigure}%
  \\
  \begin{subfigure}{0.45\textwidth}
    \centering
    \includegraphics[width=\textwidth]{%
      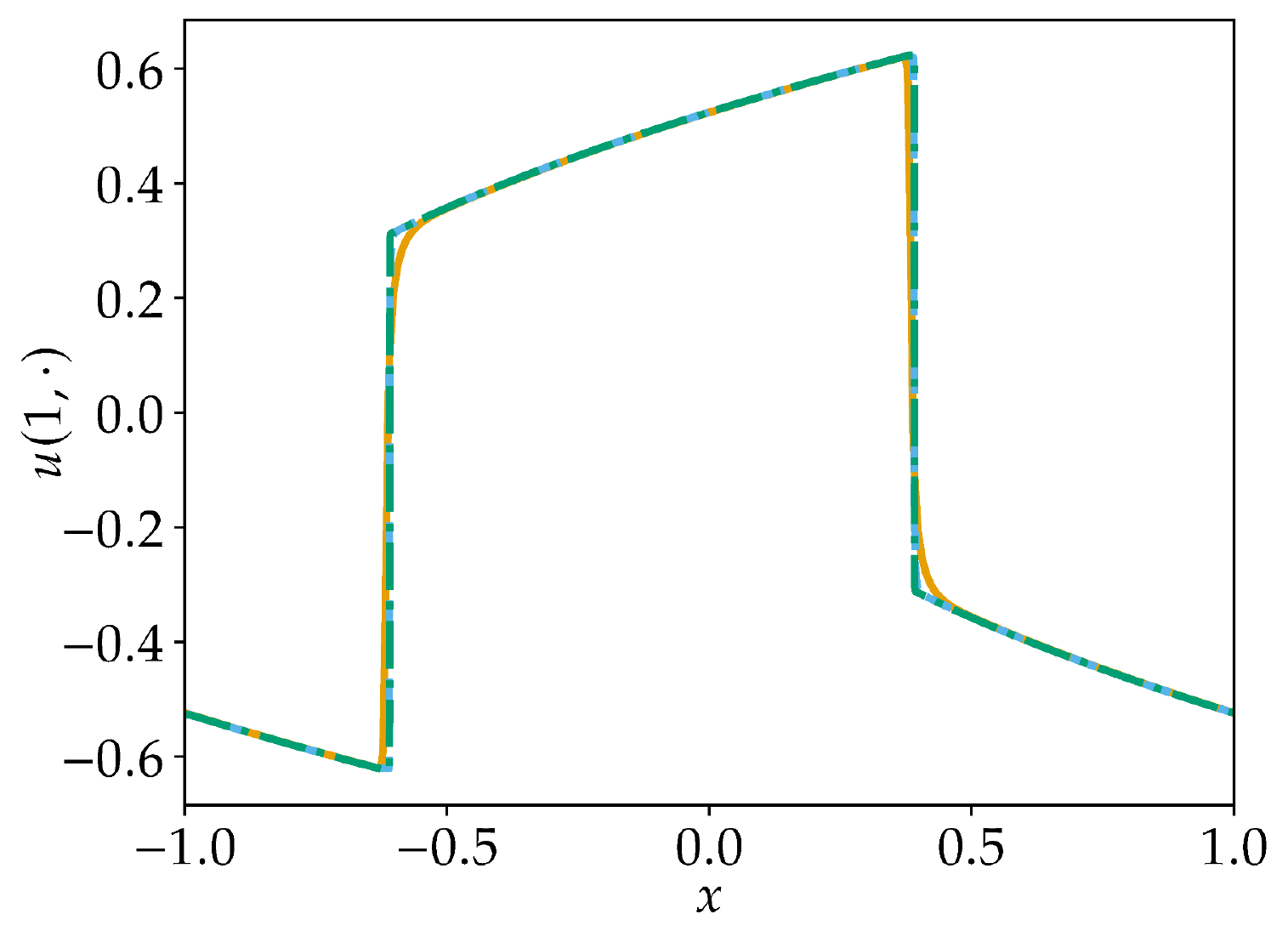}
    \caption{Split form, second-order dissipation, $\epsilon = 400$;
             numerical solutions at $t=1$.}
  \label{fig:Cubic_PeriodicFD_Split_1_DissOrder_2}
  \end{subfigure}%
  \hspace*{\fill}
  \begin{subfigure}{0.45\textwidth}
    \centering
    \includegraphics[width=\textwidth]{%
      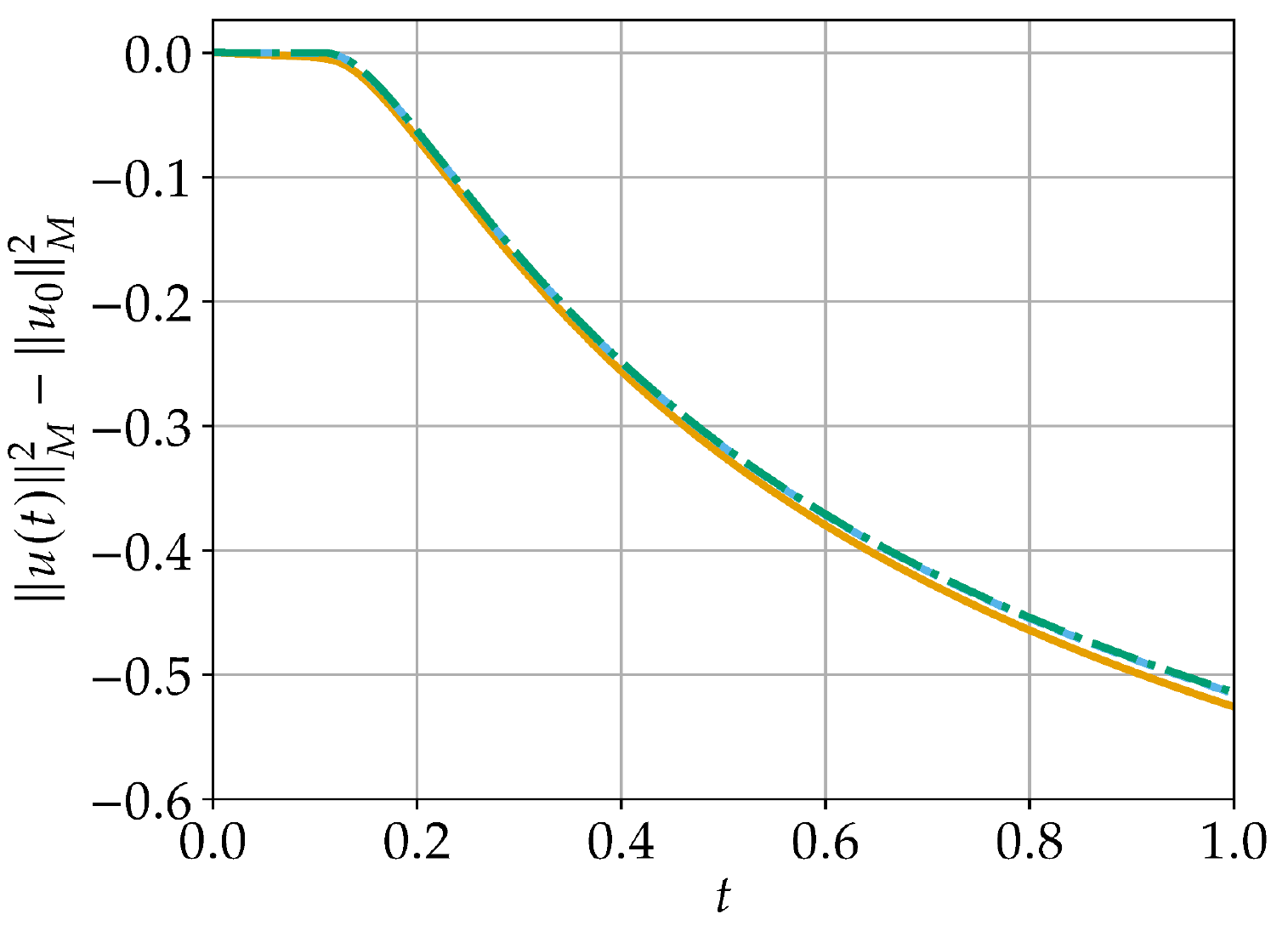}
    \caption{Split form, second-order dissipation, $\epsilon = 400$;
             evolution of the $L^2$ entropy.}
  \end{subfigure}%
  \\
  \begin{subfigure}{0.45\textwidth}
    \centering
    \includegraphics[width=\textwidth]{%
      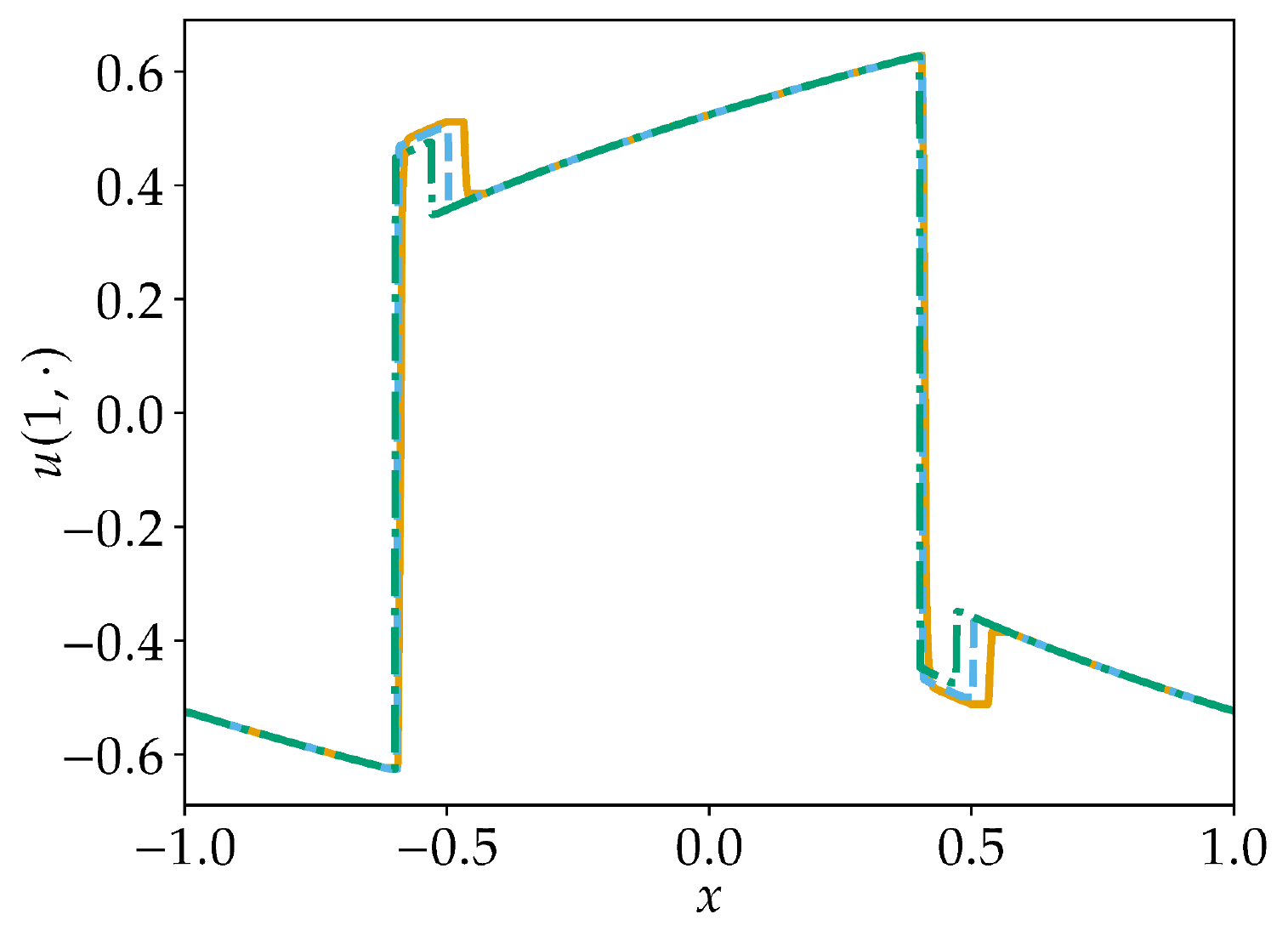}
    \caption{Split form, fourth-order dissipation, $\epsilon = 400$;
             numerical solutions at $t=1$.}
  \label{fig:Cubic_PeriodicFD_Split_1_DissOrder_4}
  \end{subfigure}%
  \hspace*{\fill}
  \begin{subfigure}{0.45\textwidth}
    \centering
    \includegraphics[width=\textwidth]{%
      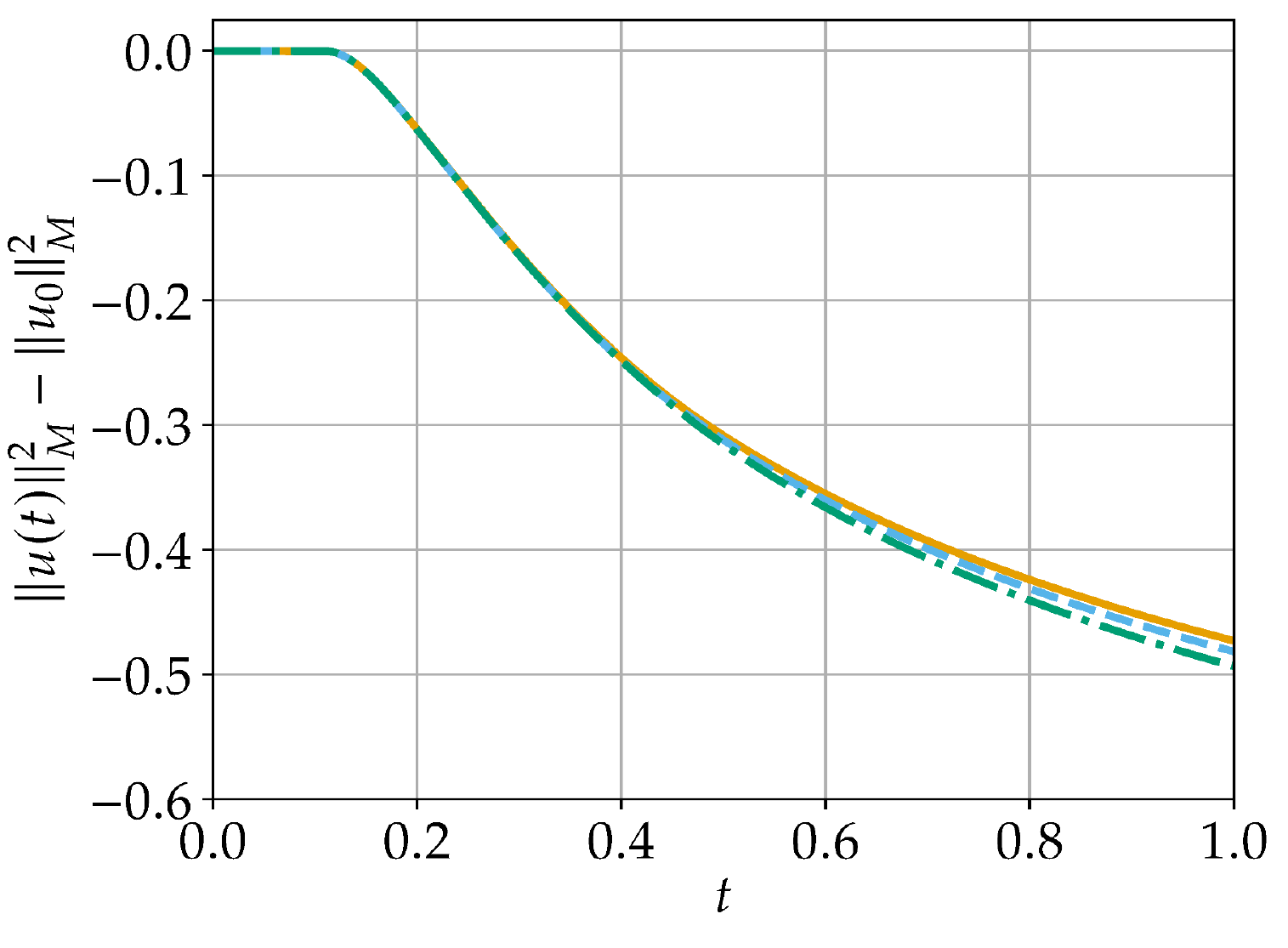}
    \caption{Split form, fourth-order dissipation, $\epsilon = 400$;
             evolution of the $L^2$ entropy.}
  \end{subfigure}%
  \\
  \begin{subfigure}{0.45\textwidth}
    \centering
    \includegraphics[width=\textwidth]{%
      figures/Cubic_PeriodicFD_AccOrder_6_DissOrder_6_Split_1__solution}
    \caption{Split form, sixth-order dissipation, $\epsilon = 400$;
             numerical solutions at $t=1$.}
  \label{fig:Cubic_PeriodicFD_Split_1_DissOrder_6}
  \end{subfigure}%
  \hspace*{\fill}
  \begin{subfigure}{0.45\textwidth}
    \centering
    \includegraphics[width=\textwidth]{%
      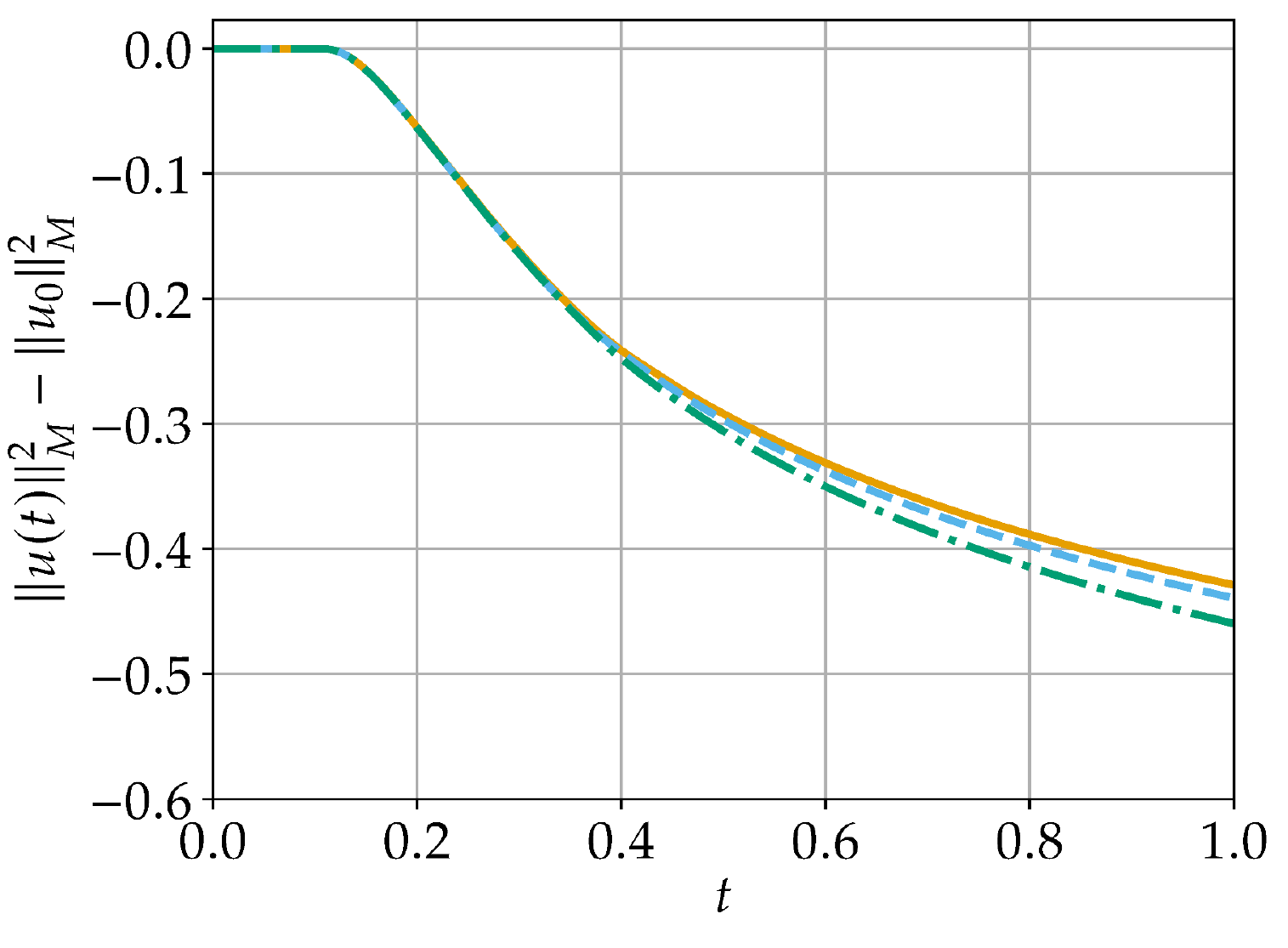}
    \caption{Split form, sixth-order dissipation, $\epsilon = 400$;
             evolution of the $L^2$ entropy.}
  \end{subfigure}%
  \caption{Numerical results for periodic finite differences using the split form
           \eqref{eq:cubic-split-form} in order to approximate solutions of the
           cubic conservation law \eqref{eq:cubic} with dissipation operators of
           different degrees for $N=1024$ (solid), $N=4096$ (dashed), and
           $N=\num{16384}$ (dotted) grid points.}
  \label{fig:Cubic_PeriodicFD_Split_1}
\end{figure}

\begin{figure}
\centering
  \begin{subfigure}{0.5\textwidth}
    \centering
    \includegraphics[width=\textwidth]{figures/Cubic__legend}
  \end{subfigure}%
  \\
  \begin{subfigure}{0.45\textwidth}
    \centering
    \includegraphics[width=\textwidth]{%
      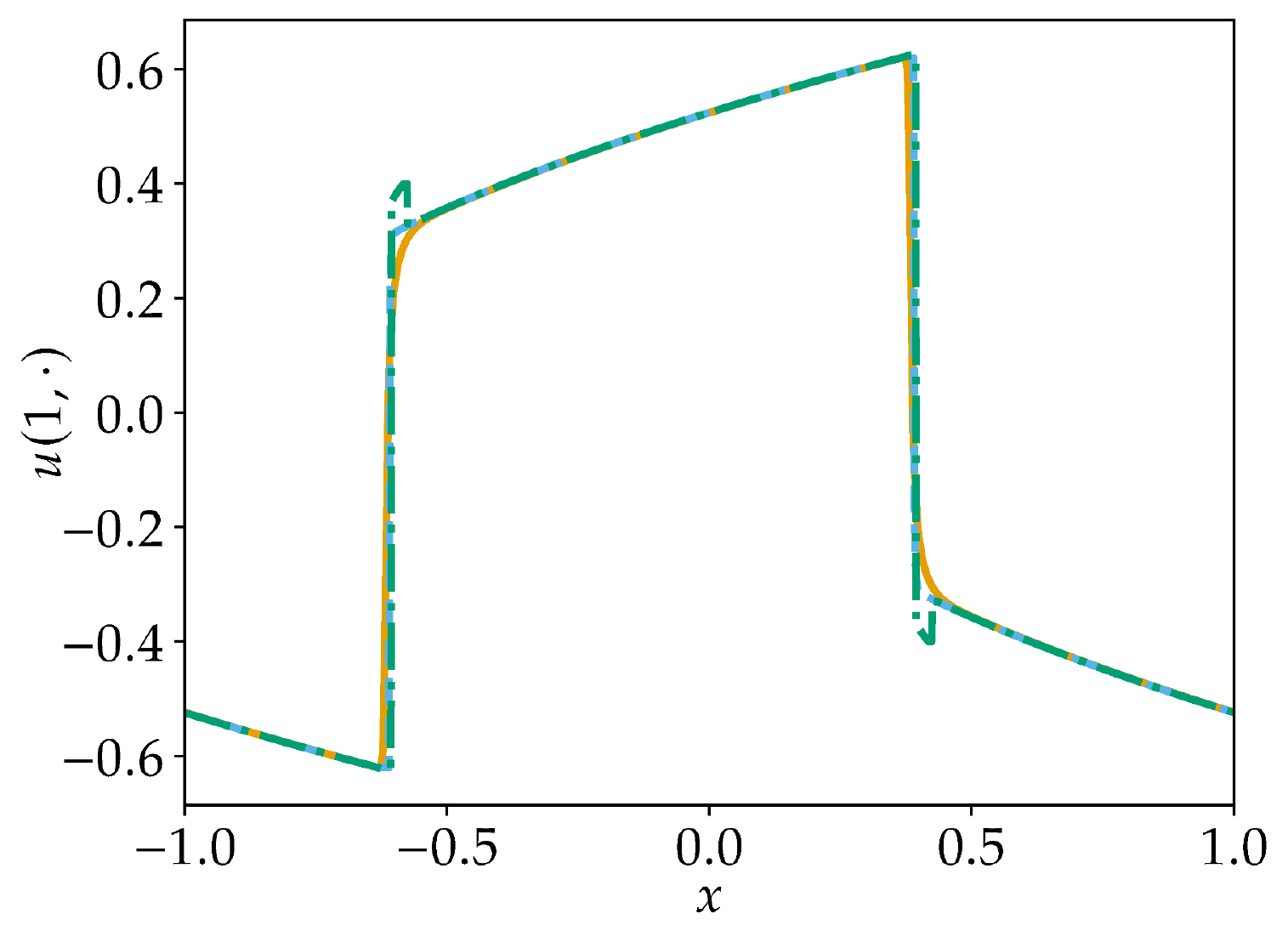}
    \caption{Unsplit form, second-order dissipation, $\epsilon = 400$;
             numerical solutions at $t=1$.}
  \label{fig:Cubic_PeriodicFD_Split_0_DissOrder_2}
  \end{subfigure}%
  \hspace*{\fill}
  \begin{subfigure}{0.45\textwidth}
    \centering
    \includegraphics[width=\textwidth]{%
      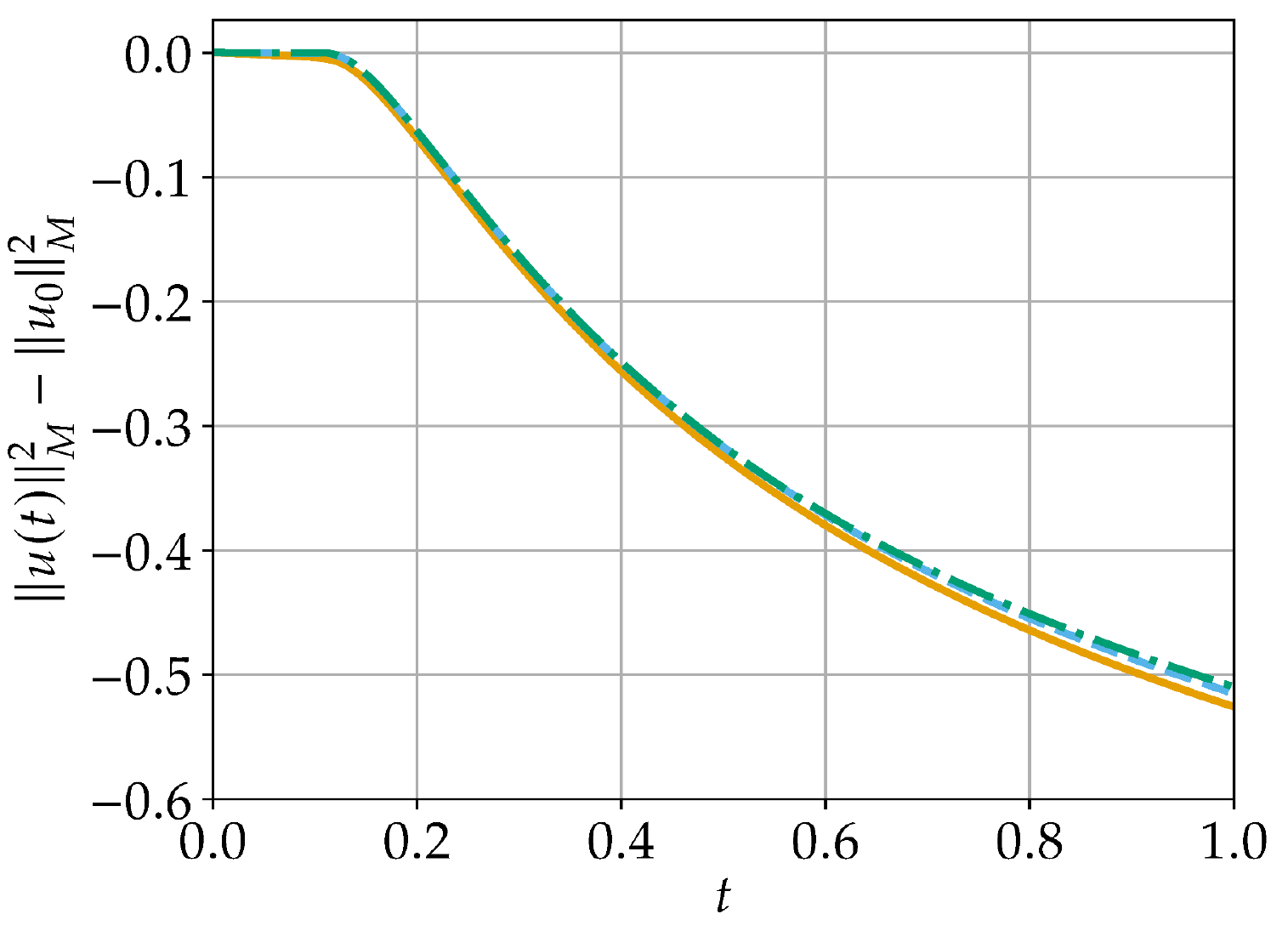}
    \caption{Unsplit form, second-order dissipation, $\epsilon = 400$;
             evolution of the $L^2$ entropy.}
  \end{subfigure}%
  \\
  \begin{subfigure}{0.45\textwidth}
    \centering
    \includegraphics[width=\textwidth]{%
      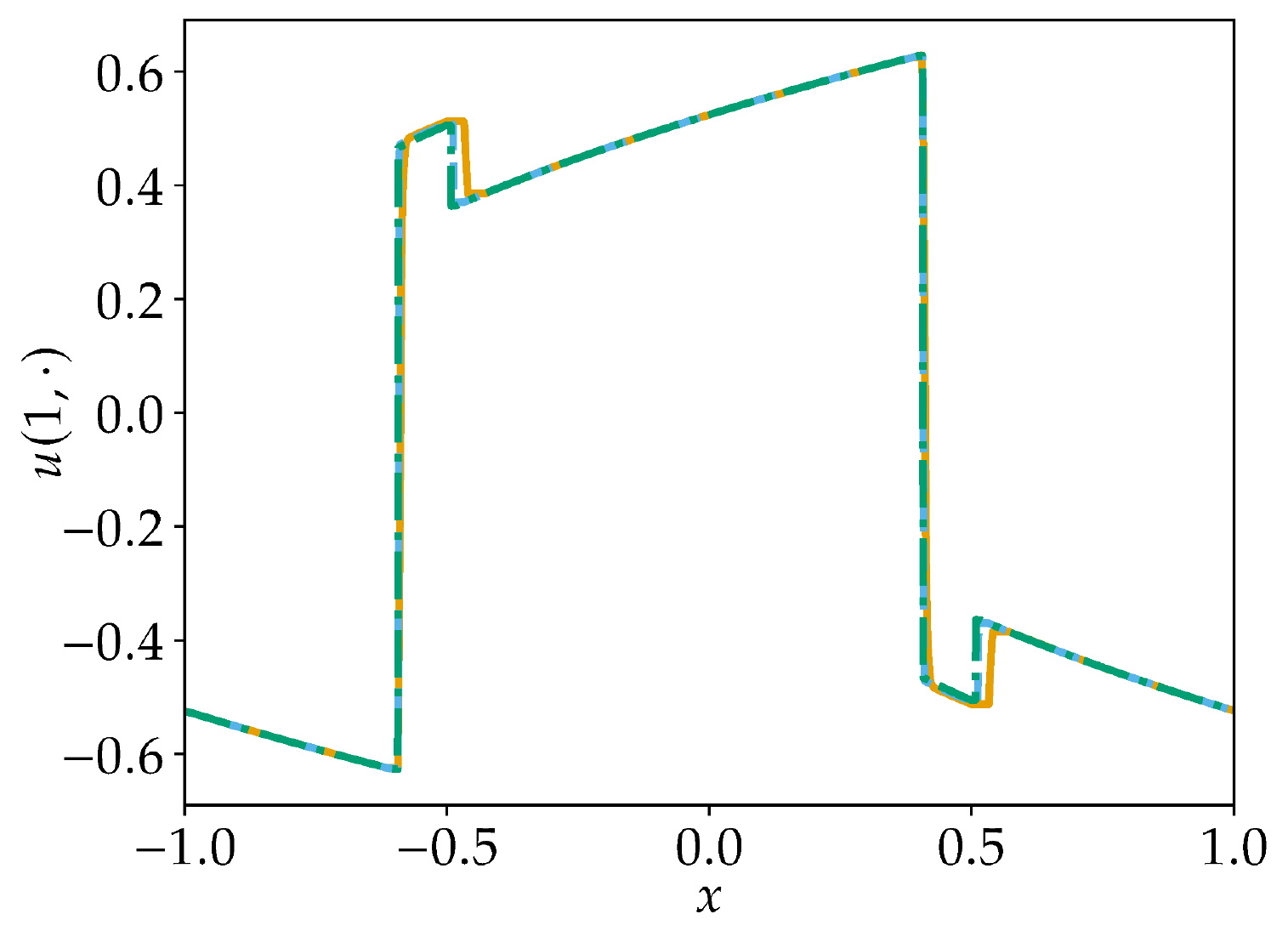}
    \caption{Unsplit form, fourth-order dissipation, $\epsilon = 400$;
             numerical solutions at $t=1$.}
  \label{fig:Cubic_PeriodicFD_Split_0_DissOrder_4}
  \end{subfigure}%
  \hspace*{\fill}
  \begin{subfigure}{0.45\textwidth}
    \centering
    \includegraphics[width=\textwidth]{%
      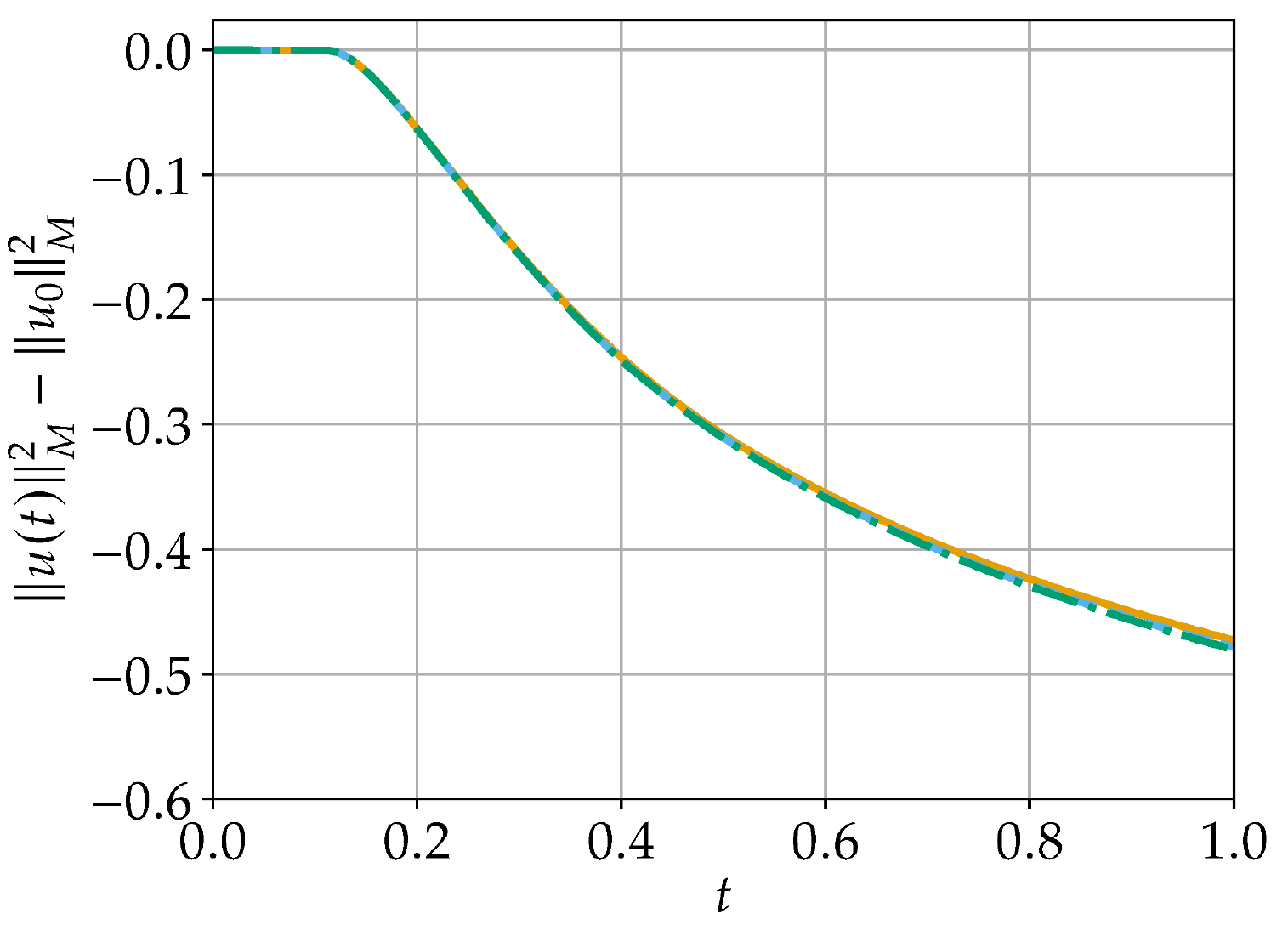}
    \caption{Unsplit form, fourth-order dissipation, $\epsilon = 400$;
             evolution of the $L^2$ entropy.}
  \end{subfigure}%
  \\
  \begin{subfigure}{0.45\textwidth}
    \centering
    \includegraphics[width=\textwidth]{%
      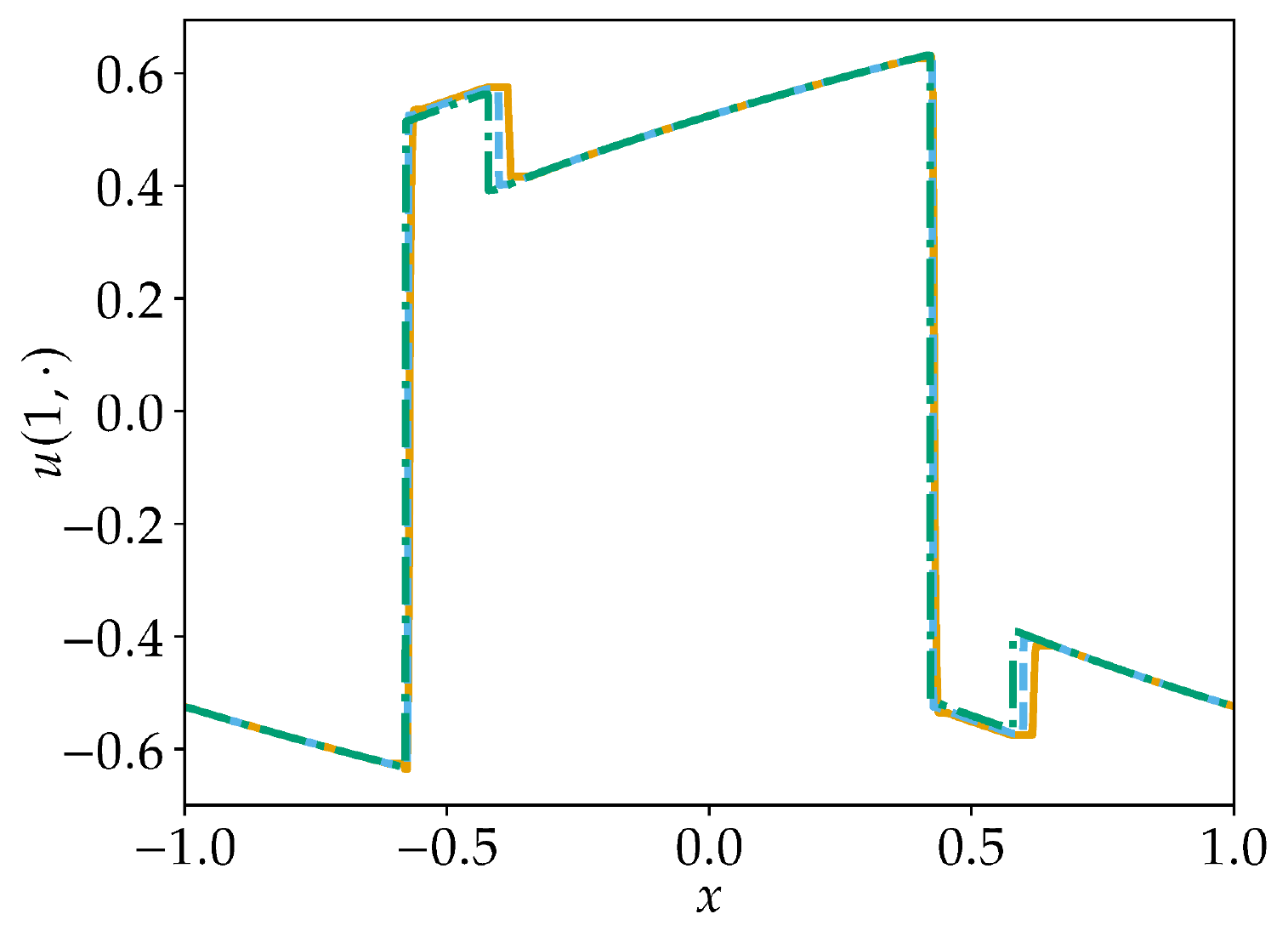}
    \caption{Unsplit form, sixth-order dissipation, $\epsilon = 400$;
             numerical solutions at $t=1$.}
  \label{fig:Cubic_PeriodicFD_Split_0_DissOrder_6}
  \end{subfigure}%
  \hspace*{\fill}
  \begin{subfigure}{0.45\textwidth}
    \centering
    \includegraphics[width=\textwidth]{%
      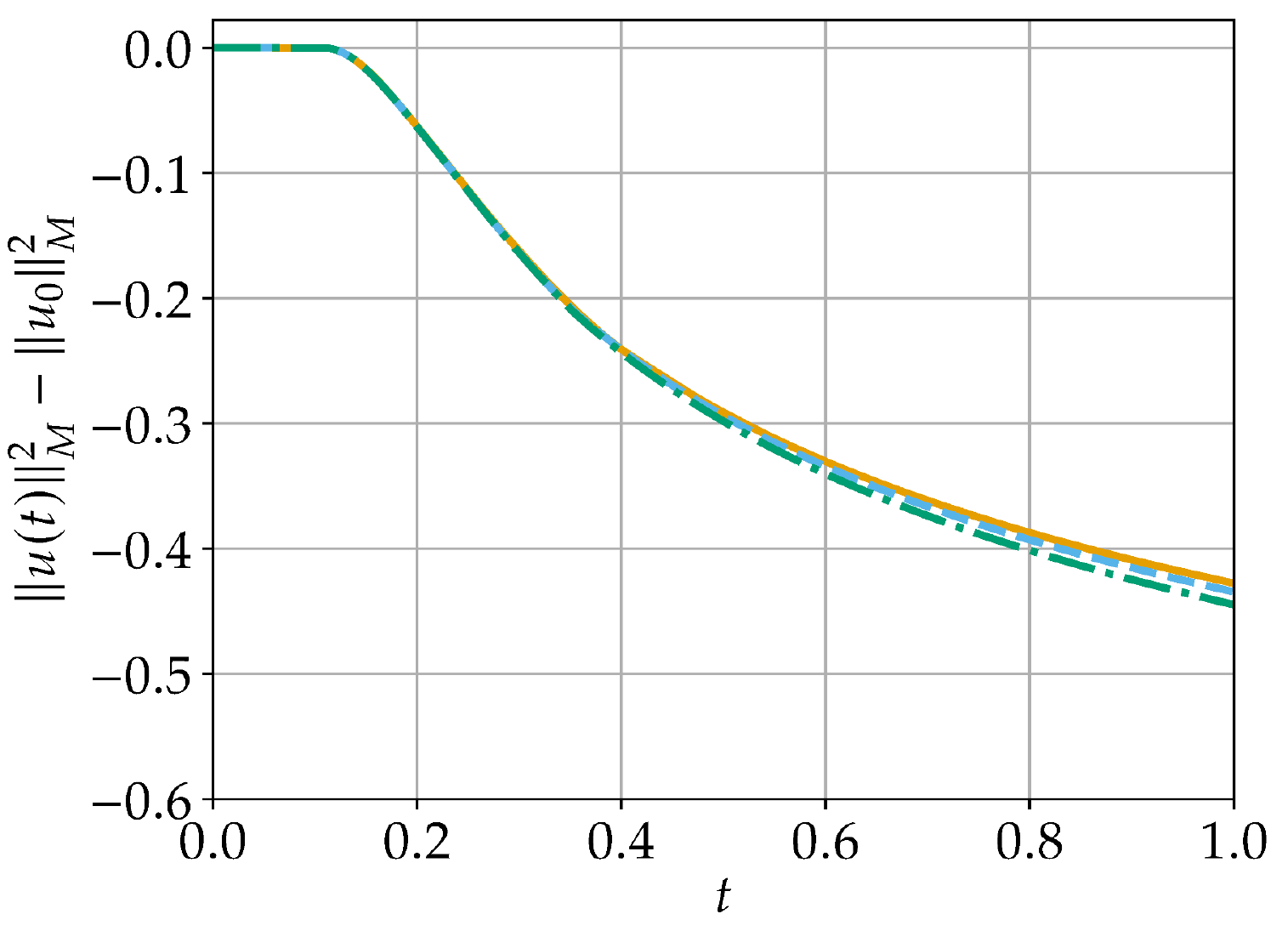}
    \caption{Unsplit form, sixth-order dissipation, $\epsilon = 400$;
             evolution of the $L^2$ entropy.}
  \end{subfigure}%
  \caption{Numerical results for periodic finite differences using the unsplit form
           \eqref{eq:cubic-unsplit-form} in order to approximate solutions of the
           cubic conservation law \eqref{eq:cubic} with dissipation operators of
           different degrees for $N=1024$ (solid), $N=4096$ (dashed), and
           $N=\num{16384}$ (dotted) grid points.}
  \label{fig:Cubic_PeriodicFD_Split_0}
\end{figure}

Typical numerical results with nonclassical solutions show some oscillations
located near discontinuities. Since they are not essential (as far the {\sl limiting solutions} are concerned) and distract from
the main observations; these have been removed using a simple total variation
denoising algorithm \cite{condat2013direct}. Exemplary numerical results before
and after this postprocessing are shown in \autoref{fig:postprocessing}.
Results of the numerical experiments using sixth-order periodic central finite difference
methods are visualized in \autoref{fig:Cubic_PeriodicFD_Split_1} using the split form
\eqref{eq:cubic-split-form} and in \autoref{fig:Cubic_PeriodicFD_Split_0} using
the unsplit form \eqref{eq:cubic-unsplit-form}. On the left-hand sides, numerical
solutions using different parameters are shown at the final time $t = 1$. On the
right-hand sides, the corresponding evolution of the discrete total entropy
$\int_M u^2 = \norm{u}_M^2 = \vec{u}^T \mat{M} \vec{u}$ can be found.

Using the split form and second-order artificial dissipation, the numerical
solutions seem to converge to the classical entropy solution
(\autoref{fig:Cubic_PeriodicFD_Split_1_DissOrder_2}). However, if artificial
dissipation operators approximating higher-order-derivatives are used, nonclassical
shocks appear in the numerical solutions (\autoref{fig:Cubic_PeriodicFD_Split_1_DissOrder_4}
and \autoref{fig:Cubic_PeriodicFD_Split_1_DissOrder_6}). The nonclassical parts
of the numerical solutions become smaller for increased number of grid points
but are still clearly visible even for $N = 2^{14} = \num{16384}$ nodes, which
seems to be pretty much for a relatively simple problem in one space dimension.
Increasing the strength $\epsilon$ of the artificial dissipation operators of
higher order does not yield substantially better results, the nonclassical parts
remain (not plotted). In accordance with the theory of nonclassical shocks of
LeFloch~\cite{lefloch1999introduction,lefloch2002hyperbolic}
and the entropy rate admissibility criterion of Dafermos~\cite{dafermos1973entropy},
the final entropy of the classical solutions is
smaller than the final entropy of the numerical solutions containing nonclassical
shocks.

Using the unsplit form, similar numerical results are obtained
(\autoref{fig:Cubic_PeriodicFD_Split_0}). Nevertheless, there are some important
differences. Firstly, while the numerical solutions with $N = 2^{10}$
and $N = 2^{12}$ nodes and second-order artificial dissipation seem
to be classical, small regions containing nonclassical shocks appear for
$N = 2^{14} = \num{16384}$ nodes (\autoref{fig:Cubic_PeriodicFD_Split_0_DissOrder_2}).
Moreover, the nonclassical parts of the numerical solutions for higher-order
artificial dissipation do not seem to become smaller at the same rate as for the
split form \eqref{eq:cubic-split-form}. In fact, using fourth-order dissipation
operators, the numerical solutions for $N = 2^{12} = \num{4096}$ and
$N = 2^{14} = \num{16384}$ nodes are visually nearly indistinguishable
(\autoref{fig:Cubic_PeriodicFD_Split_0_DissOrder_4}).

\begin{figure}
\centering
  \begin{subfigure}{0.5\textwidth}
    \centering
    \includegraphics[width=\textwidth]{figures/Cubic__legend}
  \end{subfigure}%
  \\
  \begin{subfigure}{0.45\textwidth}
    \centering
    \includegraphics[width=\textwidth]{%
      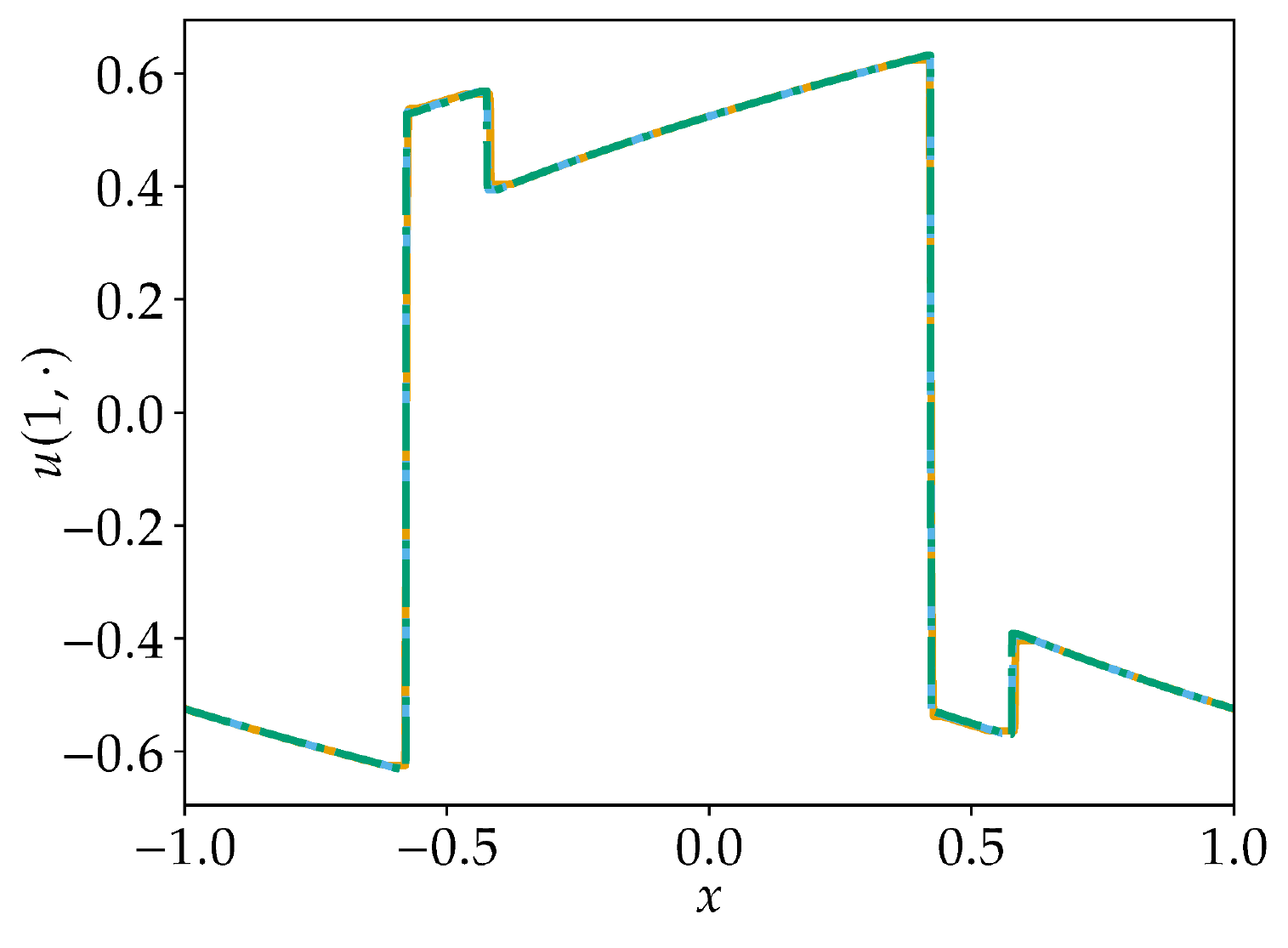}
    \caption{``Standard'' dissipation \eqref{eq:tadmor2012adaptive-convergent},
             $\epsilon = \frac{1}{N}$; numerical solutions at $t=1$.}
    \label{fig:Cubic_Fourier_Standard10}
  \end{subfigure}%
  \hspace*{\fill}
  \begin{subfigure}{0.45\textwidth}
    \centering
    \includegraphics[width=\textwidth]{%
      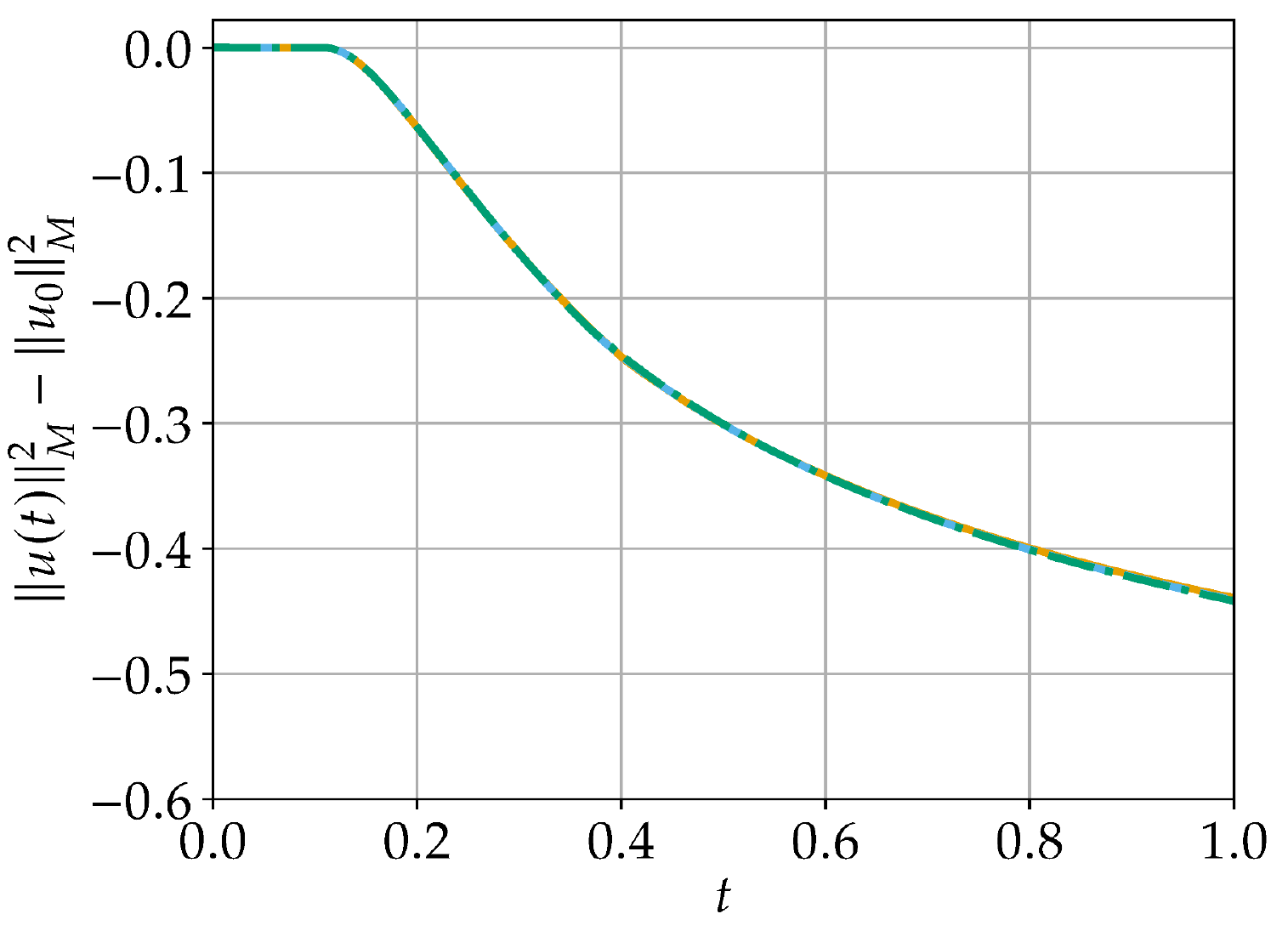}
    \caption{``Standard'' dissipation \eqref{eq:tadmor2012adaptive-convergent},
             $\epsilon = \frac{1}{N}$; evolution of the $L^2$ entropy.}
  \end{subfigure}%
  \\
  \begin{subfigure}{0.45\textwidth}
    \centering
    \includegraphics[width=\textwidth]{%
      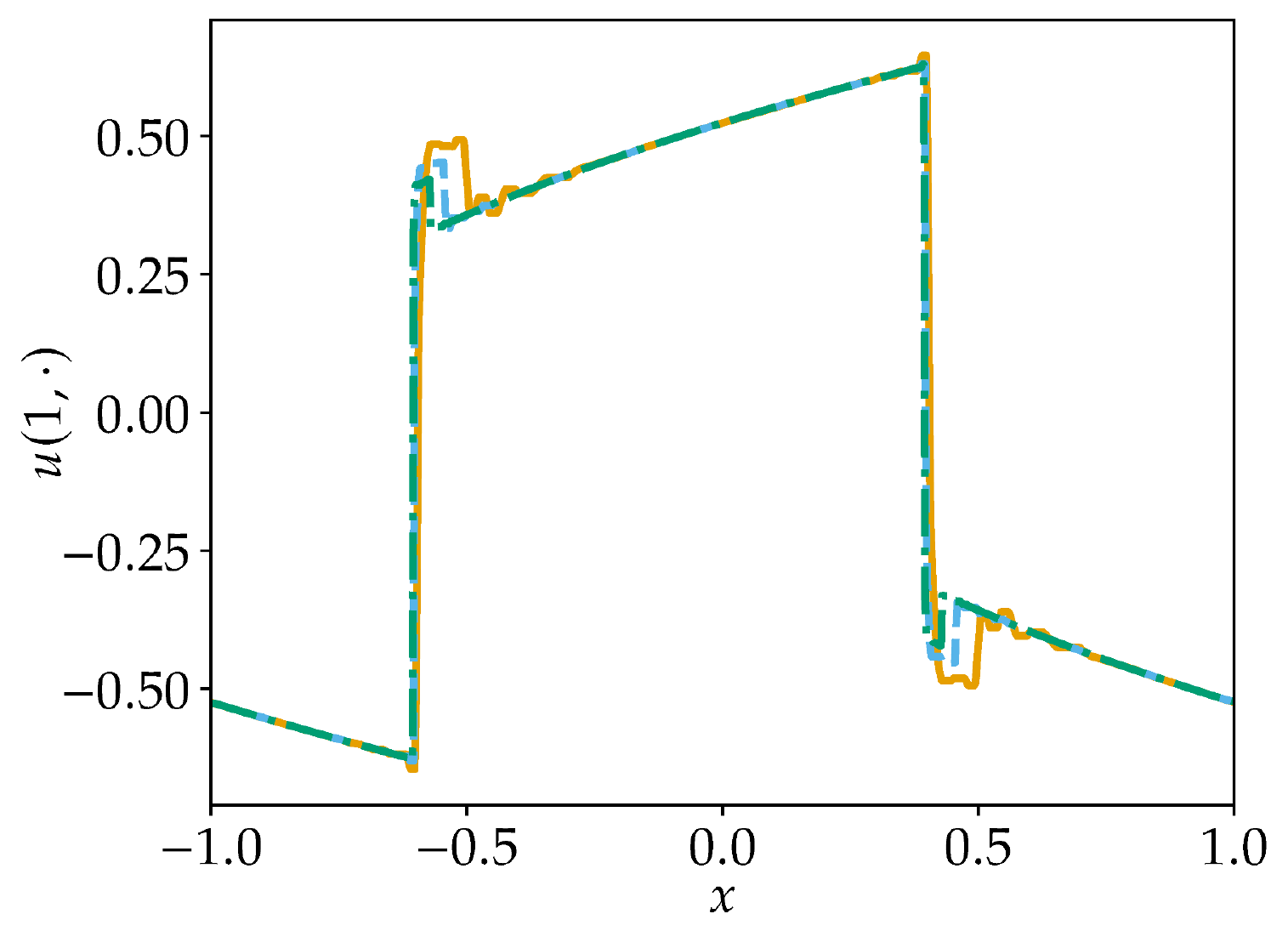}
    \caption{``Convergent'' dissipation \eqref{eq:tadmor2012adaptive-convergent},
             $\epsilon = \frac{1}{N}$; numerical solutions at $t=1$.}
    \label{fig:Cubic_Fourier_Convergent10}
  \end{subfigure}%
  \hspace*{\fill}
  \begin{subfigure}{0.45\textwidth}
    \centering
    \includegraphics[width=\textwidth]{%
      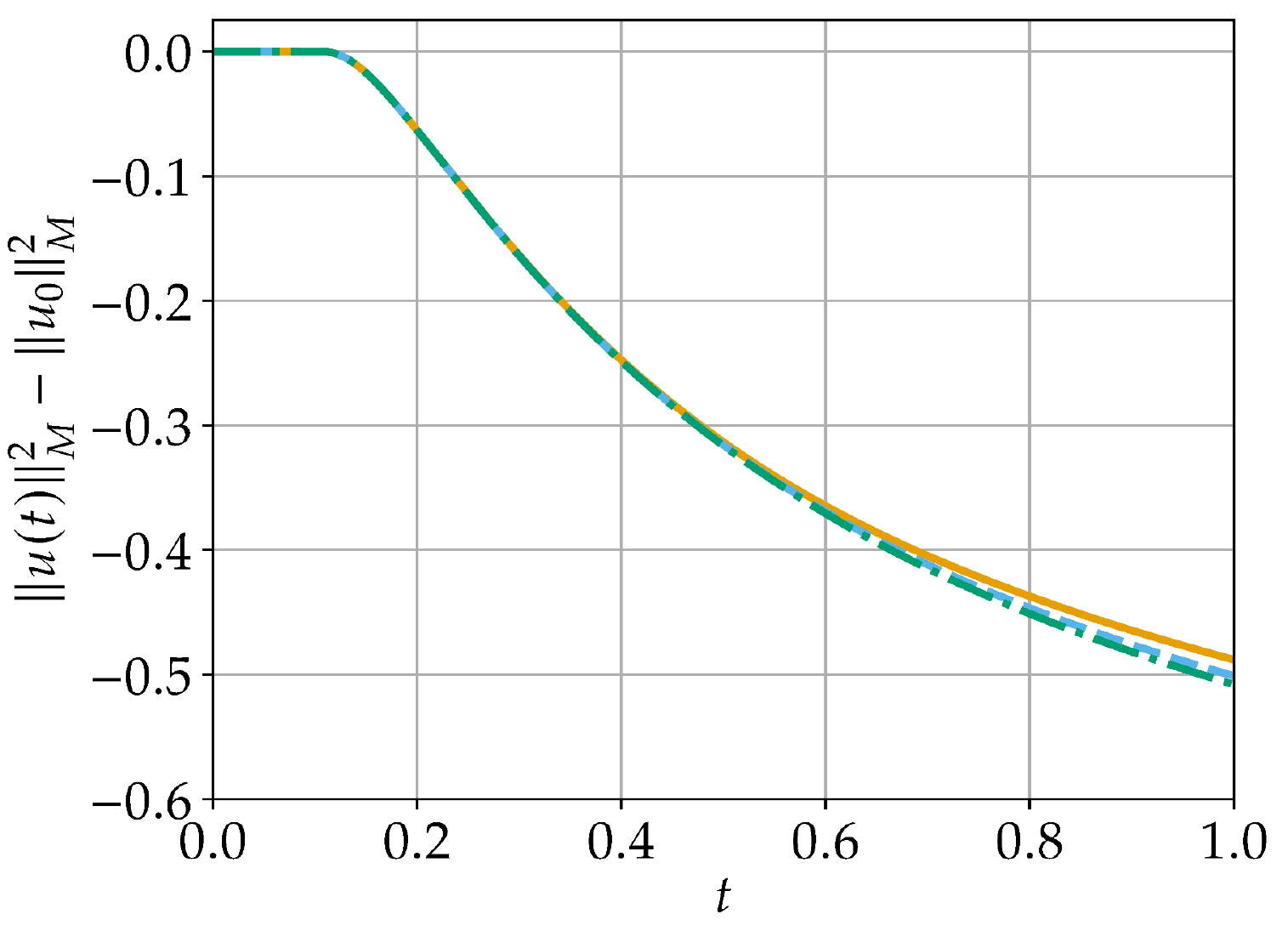}
    \caption{``Convergent'' dissipation \eqref{eq:tadmor2012adaptive-convergent},
             $\epsilon = \frac{1}{N}$; evolution of the $L^2$ entropy.}
  \end{subfigure}%
  \\
  \begin{subfigure}{0.45\textwidth}
    \centering
    \includegraphics[width=\textwidth]{%
      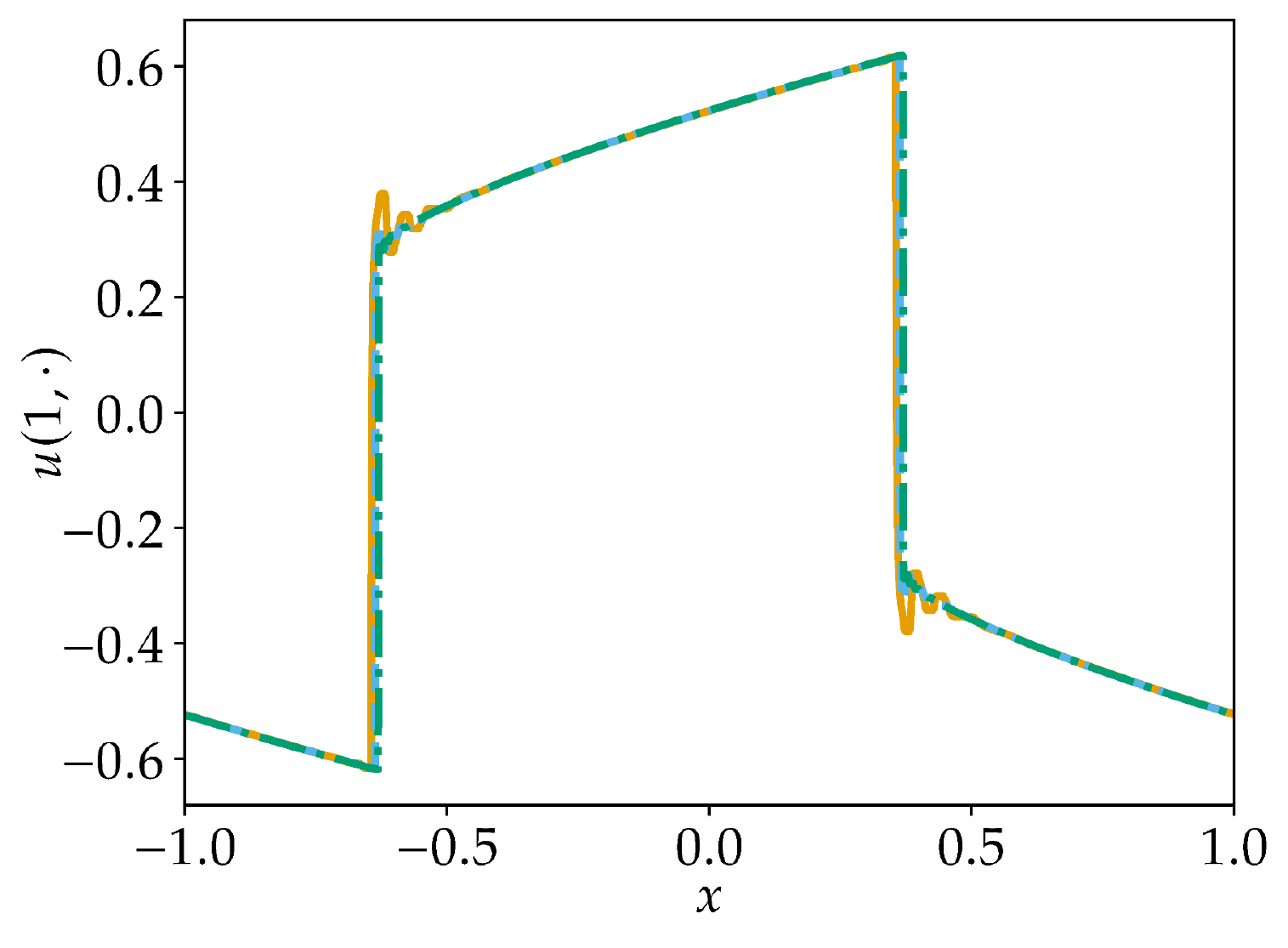}
    \caption{``Convergent'' dissipation \eqref{eq:tadmor2012adaptive-convergent},
             $\epsilon = \frac{1}{5N}$; numerical solutions at $t=1$.}
    \label{fig:Cubic_Fourier_Convergent2}
  \end{subfigure}%
  \hspace*{\fill}
  \begin{subfigure}{0.45\textwidth}
    \centering
    \includegraphics[width=\textwidth]{%
      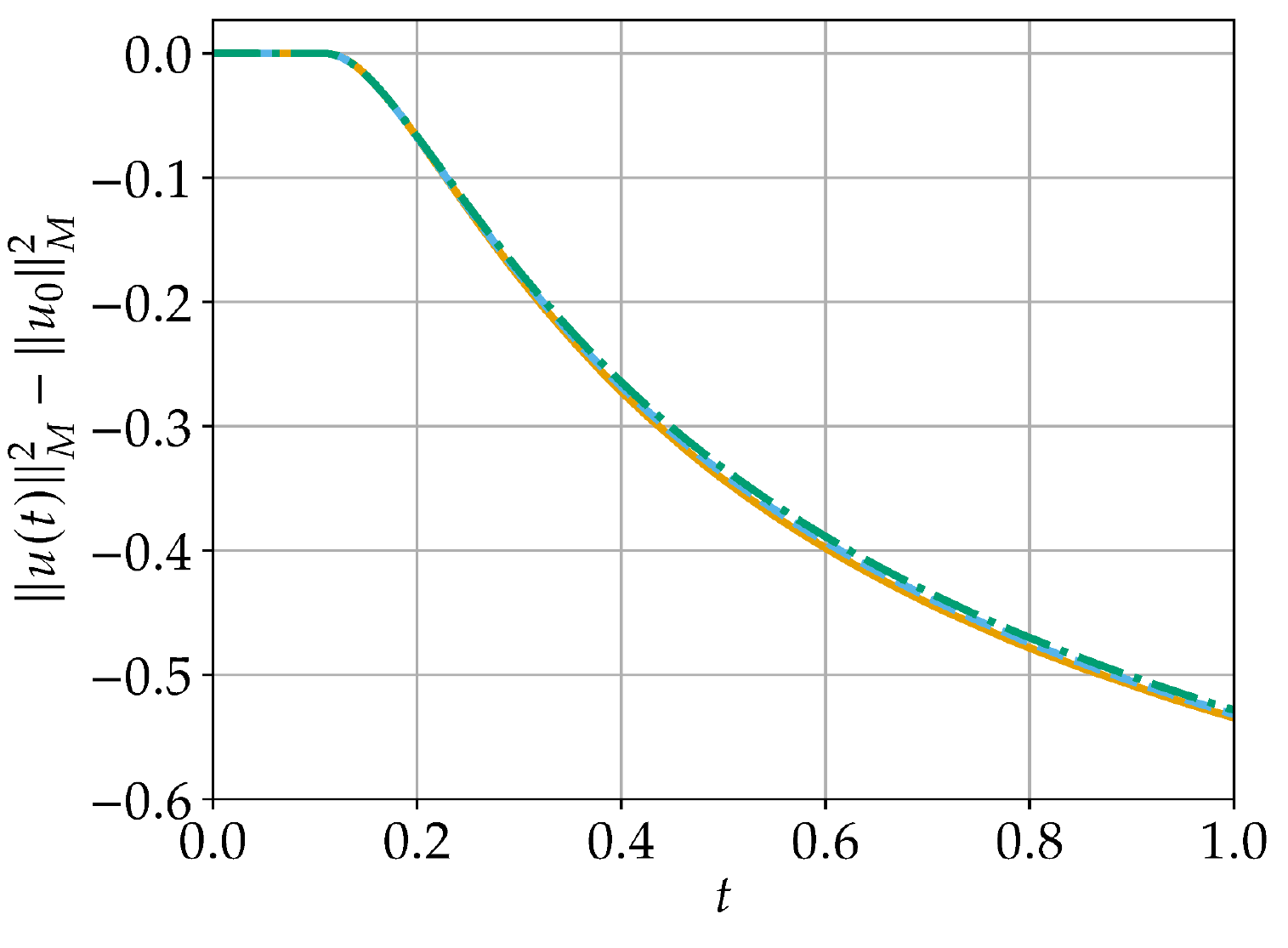}
    \caption{``Convergent'' dissipation \eqref{eq:tadmor2012adaptive-convergent},
             $\epsilon = \frac{1}{5N}$; evolution of the $L^2$ entropy.}
  \end{subfigure}%
  \caption{Numerical results for Fourier methods approximating solutions
           of the cubic conservation law \eqref{eq:cubic} with viscosity operators
           of different forms and strengths for $N=1024$ (solid), $N=4096$ (dashed),
           and $N=\num{16384}$ (dotted) grid points.}
  \label{fig:Cubic_Fourier}
\end{figure}

Results of the numerical experiments using entropy-stable Fourier methods are
shown in \autoref{fig:Cubic_Fourier}. Using the ``standard'' dissipation
\eqref{eq:tadmor2012adaptive-standard}, the numerical solutions for $N \in
\{2^{10}, 2^{12}, 2^{14}\}$ are visually indistinguishable and include nonclassical
shocks (\autoref{fig:Cubic_Fourier_Standard10}). Using the ``convergent'' choice
\eqref{eq:tadmor2012adaptive-convergent} of the dissipation coefficients with
strength $\epsilon = \frac{1}{N}$, the numerical solutions still contain nonclassical
shocks but these nonclassical regions become smaller for increasing numbers of
grid points $N$ (\autoref{fig:Cubic_Fourier_Convergent10}). However, nonclassical
shocks are still visible even for $N = 2^{14} = \num{16384}$ nodes, which seems
to be a very high resolution for this one-dimensional problem. Again, the final
total entropy becomes smaller as the nonclassical regions become smaller.
Surprisingly, reducing the strength of the spectral viscosity operators to
$\epsilon = \frac{1}{5N}$, the convergence to the classical entropy solution
becomes faster and more clearly visible (\autoref{fig:Cubic_Fourier_Convergent2}).
Thus, \emph{reducing} the strength $\epsilon$ of the artificial dissipation
\emph{increases} the total dissipation, since the classical solution is approached
and the final value of the total entropy is smaller than for nonclassical
solutions.

\subsection{Non-periodic boundary conditions}

Considering non-periodic boundary conditions, it might be expected that a boundary
condition has to be provided at the left-hand side (i.e.\ at $\xmin$), since the
advection speed $f'(u) = 3 u^2$ is non-negative. Similarly, no boundary data
should be given at the right-hand side, i.e.\ at $\xmax$. Indeed, smooth solutions
of the initial boundary value problem
\begin{equation}
\label{eq:cubic-IBVP}
\begin{aligned}
  \partial_t u(t,x) + \partial_x u(t,x)^3 &= 0,
    && t \in (0, T),\, x \in (\xmin,\xmax),
  \\
  u(0,x) &= u_0(x),
    && x \in (\xmin,\xmax),
  \\
  u(t,\xmin) &= g_L(t),
    && t \in (0,T),
\end{aligned}
\end{equation}
with compatible initial and boundary data fulfil
\begin{equation}
\begin{aligned}
  \frac{1}{2} \od{}{t} \int_\xmin^\xmax u(t,x)^2 \dif x
  &=
  \int_\xmin^\xmax u(t,x) \, \partial_t u(t,x) \dif x
  =
  - \int_\xmin^\xmax u(t,x) \, \partial_x u(t,x)^3 \dif x
  \\
  &=
  - \frac{3}{4} u(t,x)^4 \big|_\xmin^\xmax
  =
  \frac{3}{4} g_L(t)^4 - \frac{3}{4} u(t,\xmax)^4.
\end{aligned}
\end{equation}
Thus, the total entropy is bounded by initial and boundary data.

Considering the semi-discretization described in Theorem~\ref{thm:flux-diff-ES}
results in the semi-discrete entropy balance
\begin{equation}
\begin{aligned}
  \od{}{t} \norm{\vec{u}}_M^2
  &=
  \vec{u}^T \mat{M} \partial_t \vec{u}
  =
  \bigl( u_L \fnum_L -
  \frac{1}{4} u_L^4
  \bigr)
  - \bigl( u_R \fnum_R -
  \frac{1}{4} u_R^4
 \bigr),
\end{aligned}
\end{equation}
in which $\frac{1}{4} u_L^4 = \psi_L$ and $\frac{1}{4} u_R^4=\psi_R$.
Using Godunov's flux $\fnum(u_-,u_+) = u_-^3$ yields
\begin{multline}
  \od{}{t} \norm{\vec{u}}_M^2
  =
  \biggl( u_L g_L^3 - \frac{1}{4} u_L^4 \biggr)
  - \biggl( u_R^4 - \frac{1}{4} u_R^4 \biggr)
  \\
  =
  \frac{3}{4} g_L^4 - \frac{3}{4} u_R^4
  - \frac{1}{4} \underbrace{(u_L - g_L)^2 \bigl( 3 g_L^2 + 2 u_L g_L + u_L^2 \bigr)}_{\geq 0},
\end{multline}
\vskip-.4cm
\noindent since
\begin{equation}
\begin{aligned}
  0
  &\leq
  (u_L - g_L)^2 \bigl( 3 g_L^2 + 2 u_L g_L + u_L^2 \bigr)
  \\
  &=
  3 u_L^2 g_L^2 + 2 u_L^3 g_L + u_L^4
  - 6 u_L g_L^3 - 4 u_L^2 g_L^2 - 2 u_L^3 g_L
  + 3 g_L^4 + 2 u_L g_L^3 + u_L^2 g_L^2
  \\
  &=
  u_L^4 - 4 u_L g_L^3 + 3 g_L^4.
\end{aligned}
\end{equation}
Thus, the entropy rate of the numerical solution is bounded from above by
the corresponding analytical entropy rate. Hence, using Godunov's flux at
the (exterior) boundaries results in an entropy-stable scheme.
If multiple elements are used, the numerical flux at inter-element boundaries
can be any entropy-stable flux (cf.\ Definition~\ref{def:fnum-EC-ES}),
resulting in entropy-stable semi-discretizations.

In the following, homogeneous boundary data $g_L(t) \equiv 0$ and the initial
condition $u_0(x) = -\sin(\pi x)$ will be used. The domain $(\xmin,\xmax) = (-1,3)$
is divided uniformly into $N$ elements. A DG approach on Lobatto nodes is used,
i.e.\ the solution is represented on each element as a polynomial of degree $\leq p$,
represented in a nodal basis. Godunov's flux is used both at the exterior and the
interior boundaries. Moreover, the entropy-conservative flux
\eqref{eq:cubic-fnum-EC} is used for the volume terms as described in
Theorem~\ref{thm:flux-diff-ES}. Again, SSPRK(10,4) is used to advance the numerical
solutions in time, up to the final time $T=1.5$. The time step is chosen as
$\Delta t = \frac{1}{N (p^2+1)}$.

\begin{figure}
\centering
  \begin{subfigure}{0.5\textwidth}
    \centering
    \includegraphics[width=\textwidth]{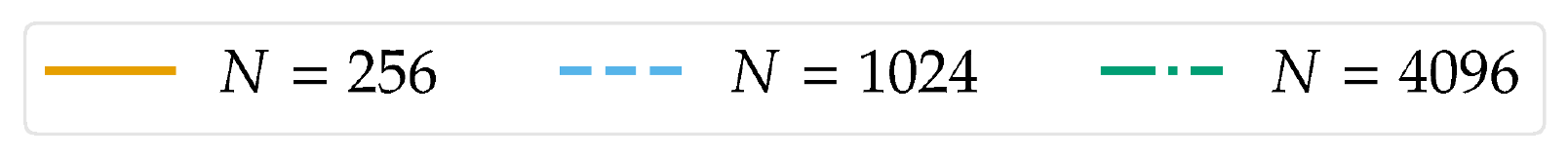}
  \end{subfigure}%
  \\
  \begin{subfigure}{0.45\textwidth}
    \centering
    \includegraphics[width=\textwidth]{%
      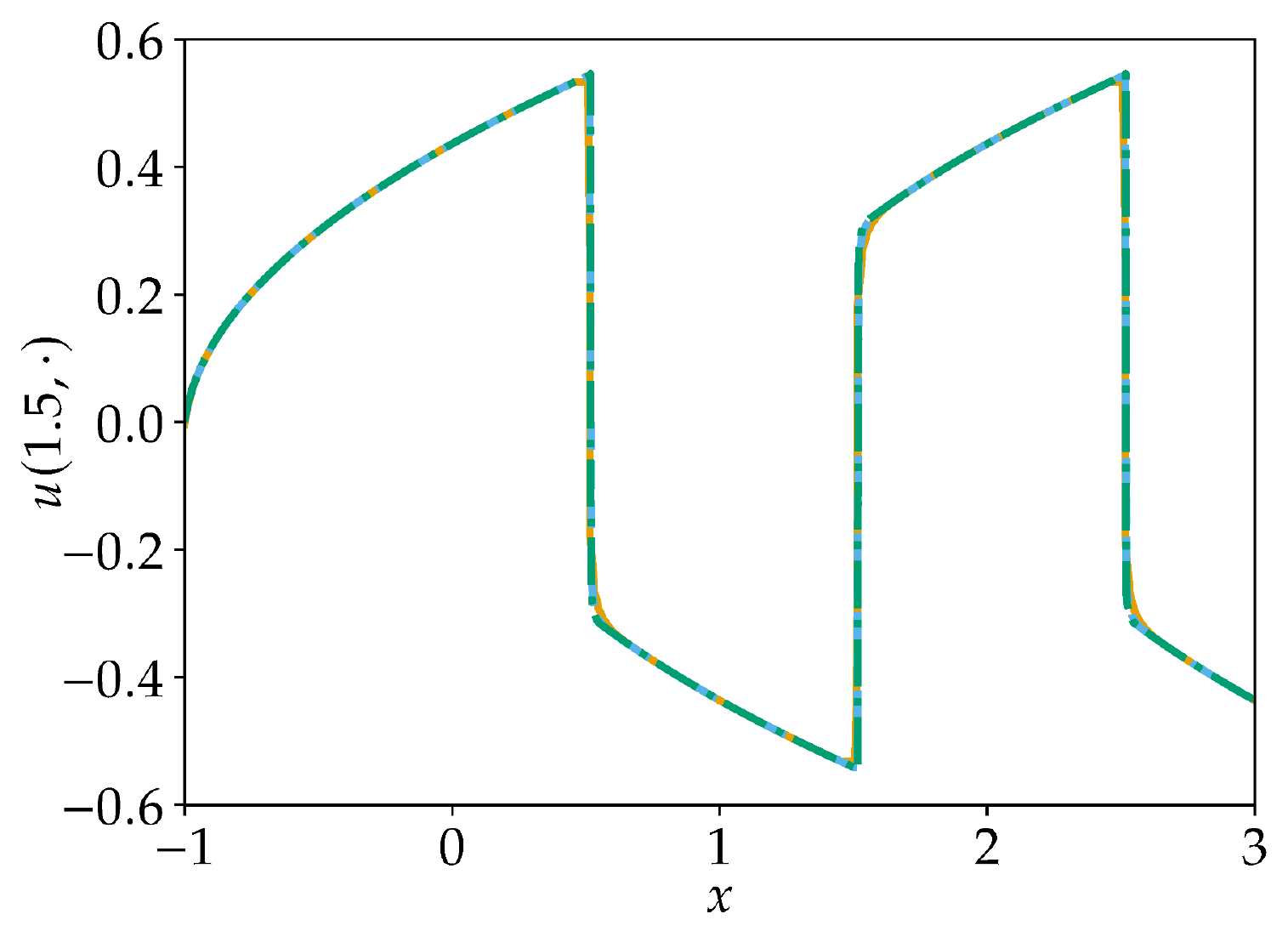}
    \caption{Polynomial degree $p=1$;
             numerical solutions at $t=1.5$.}
  \label{fig:Cubic_DG_p1}
  \end{subfigure}%
  \hspace*{\fill}
  \begin{subfigure}{0.45\textwidth}
    \centering
    \includegraphics[width=\textwidth]{%
      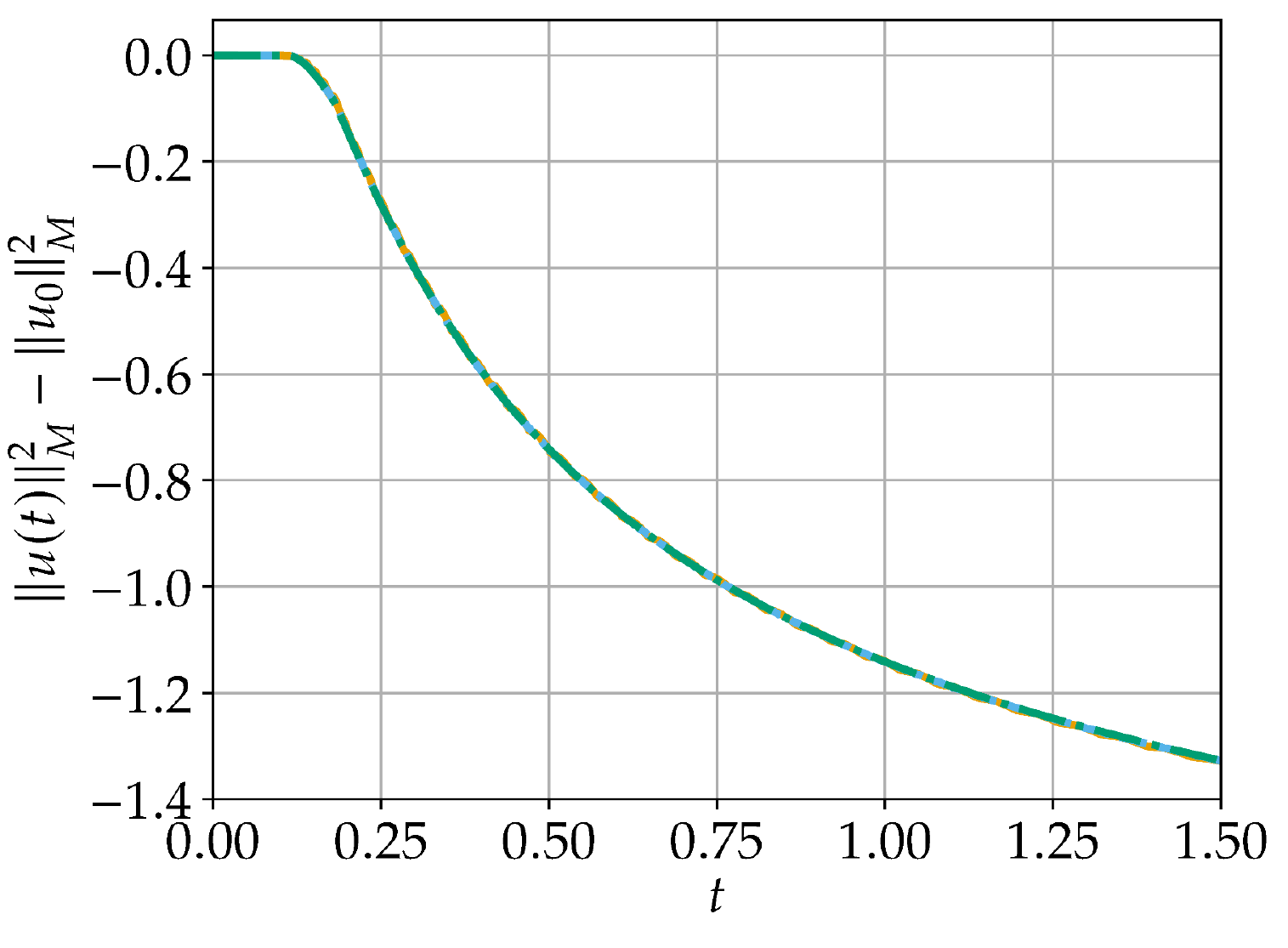}
    \caption{Polynomial degree $p=1$;
             evolution of the $L^2$ entropy.}
  \end{subfigure}%
  \\
  \begin{subfigure}{0.45\textwidth}
    \centering
    \includegraphics[width=\textwidth]{%
      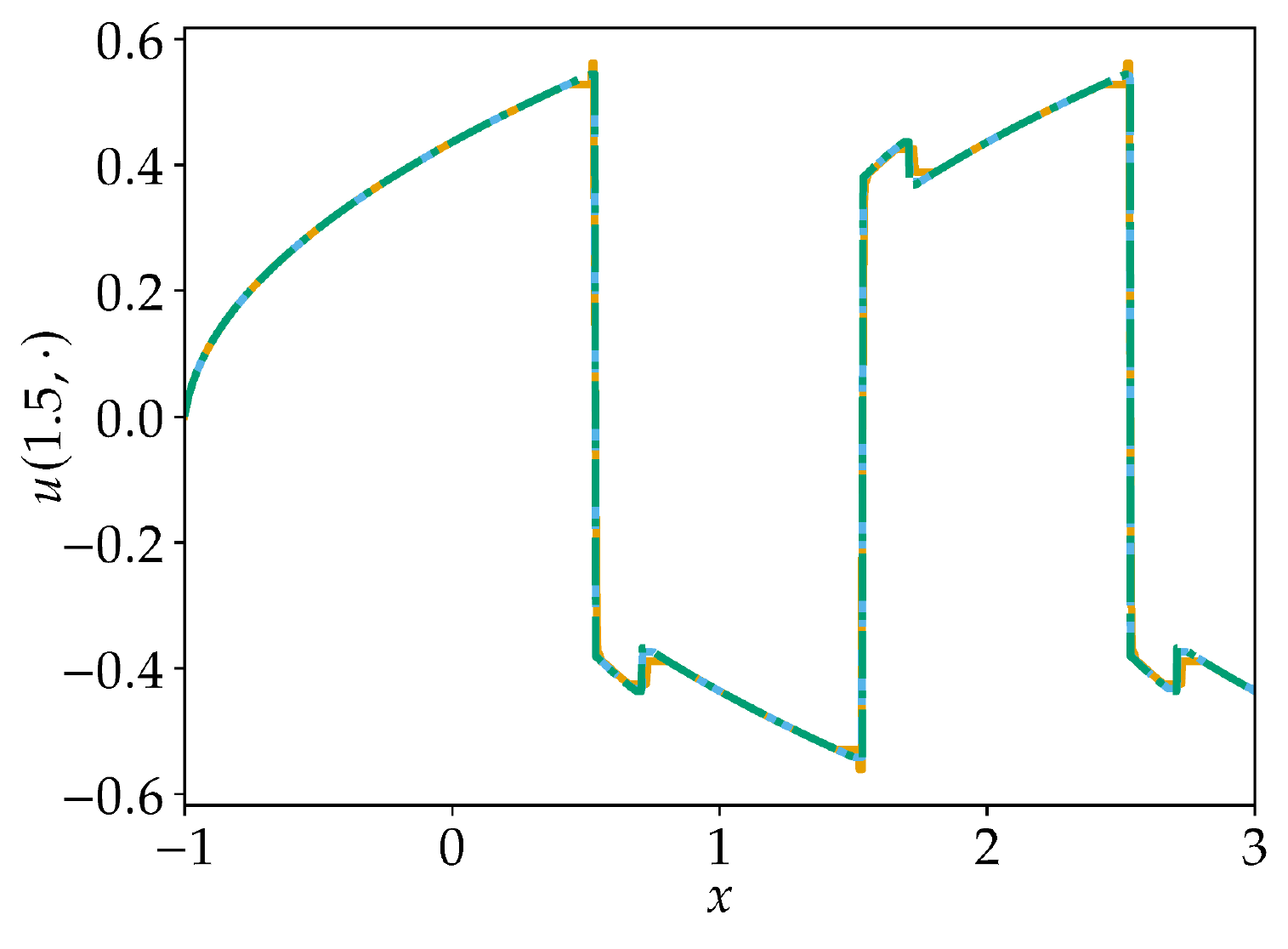}
    \caption{Polynomial degree $p=2$;
             numerical solutions at $t=1.5$.}
  \label{fig:Cubic_DG_p2}
  \end{subfigure}%
  \hspace*{\fill}
  \begin{subfigure}{0.45\textwidth}
    \centering
    \includegraphics[width=\textwidth]{%
      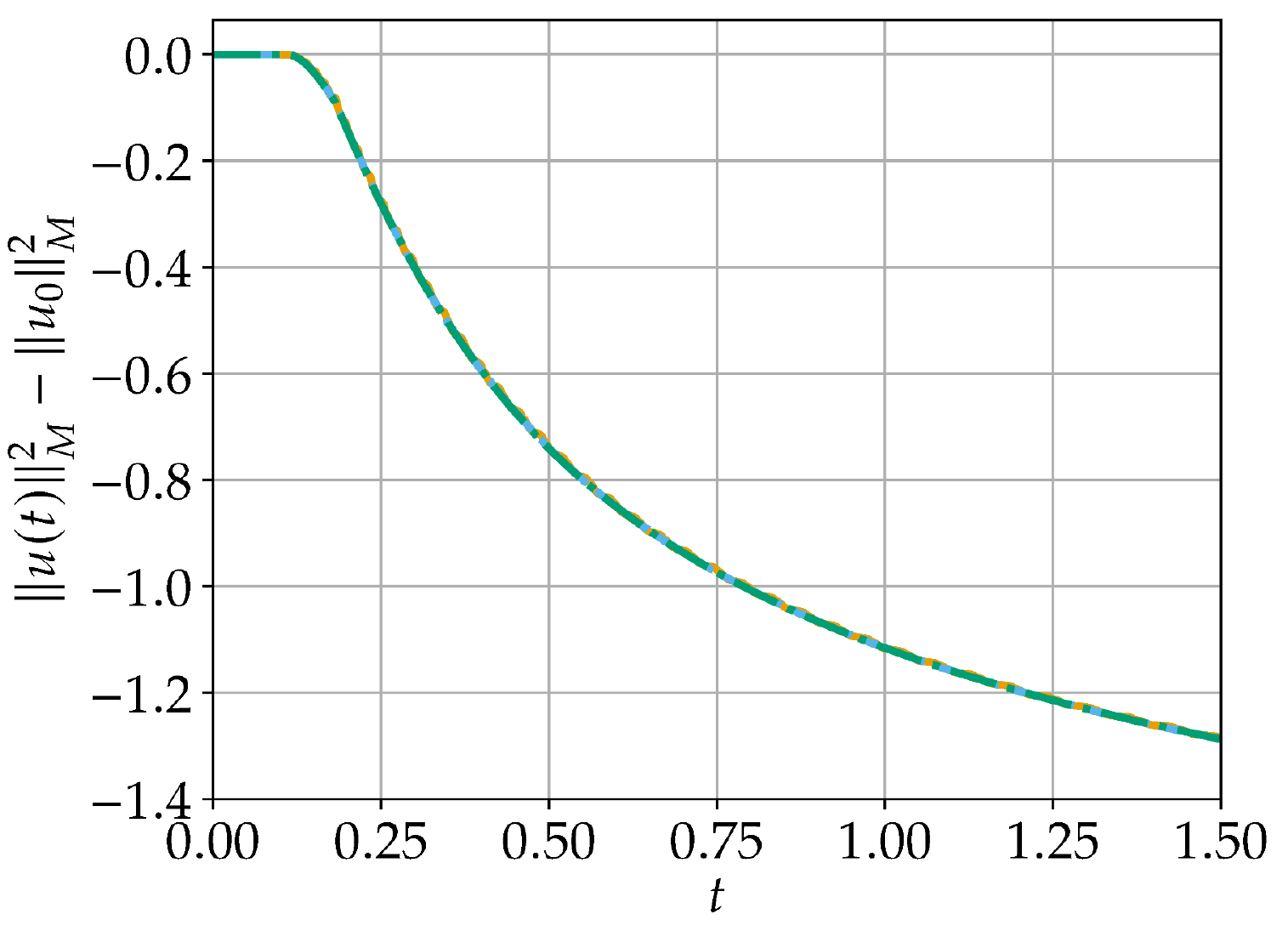}
    \caption{Polynomial degree $p=2$;
             evolution of the $L^2$ entropy.}
  \end{subfigure}%
  \\
  \begin{subfigure}{0.45\textwidth}
    \centering
    \includegraphics[width=\textwidth]{%
      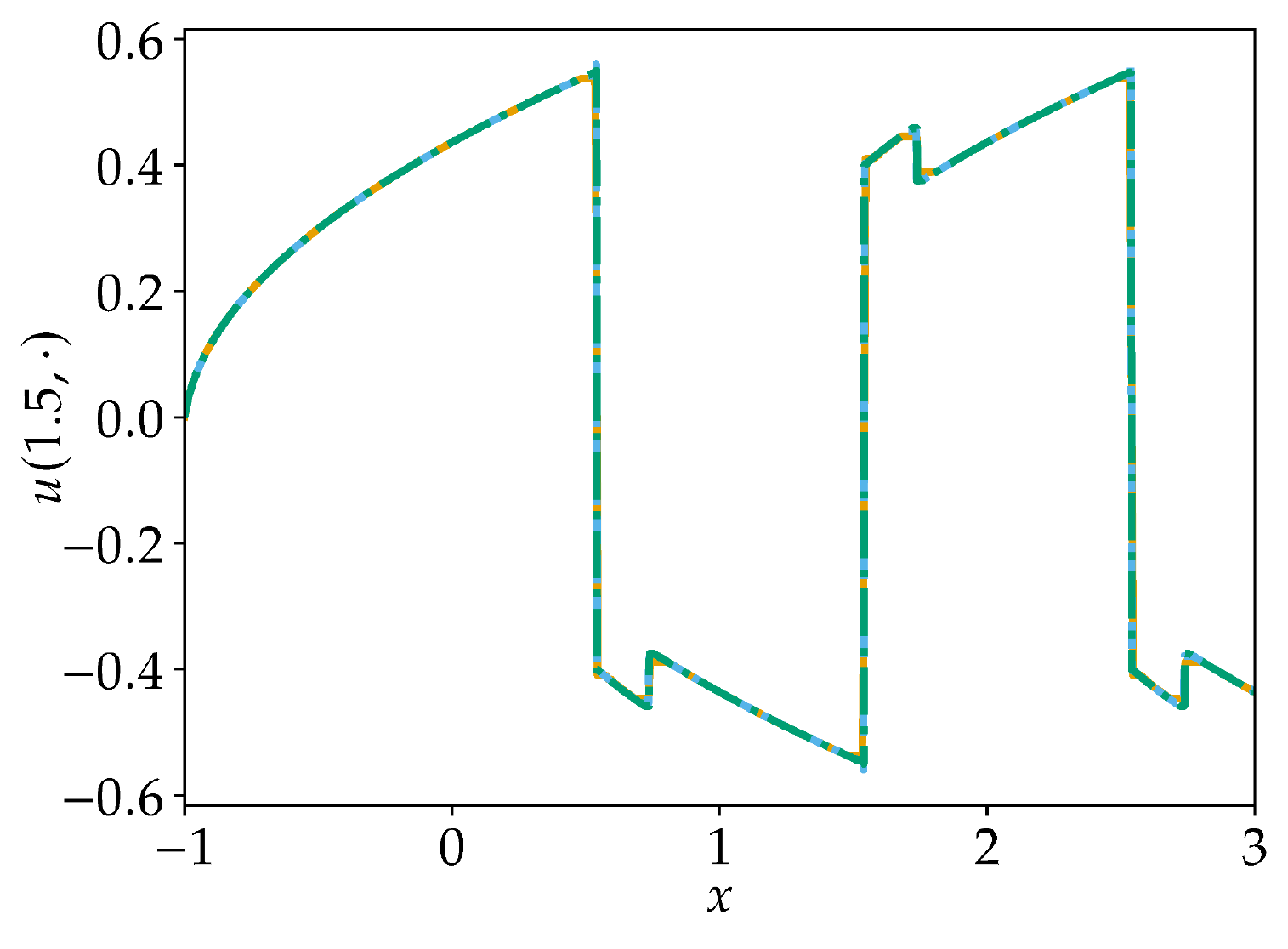}
    \caption{Polynomial degree $p=3$;
             numerical solutions at $t=1.5$.}
  \label{fig:Cubic_DG_p3}
  \end{subfigure}%
  \hspace*{\fill}
  \begin{subfigure}{0.45\textwidth}
    \centering
    \includegraphics[width=\textwidth]{%
      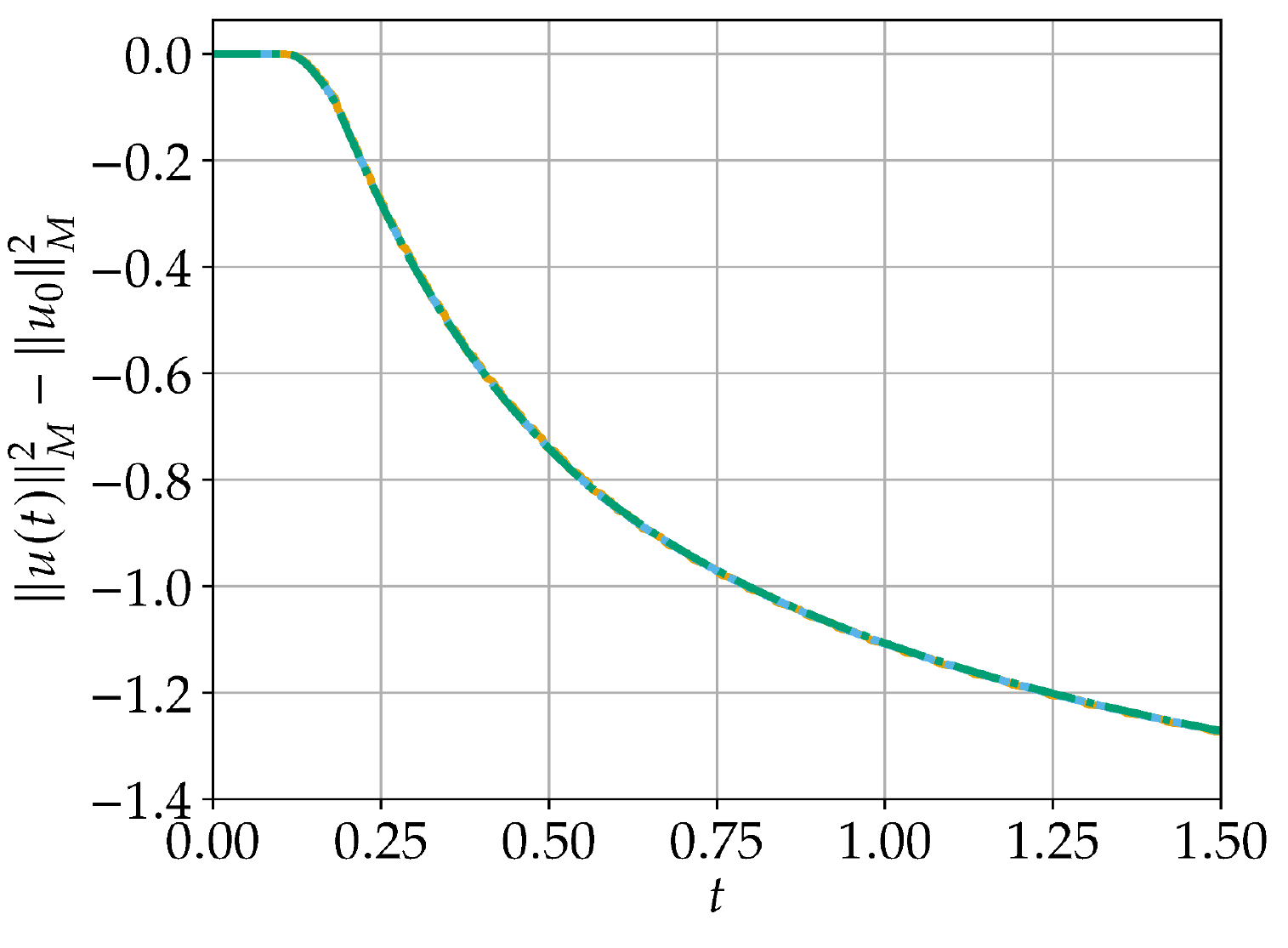}
    \caption{Polynomial degree $p=3$;
             evolution of the $L^2$ entropy.}
  \end{subfigure}%
  \caption{Numerical results for DG methods using Godunov's flux and the entropy
           stable semi-discretization of Theorem~\ref{thm:flux-diff-ES} in order
           to approximate solutions of the cubic conservation law \eqref{eq:cubic-IBVP}
           with different polynomial degrees $p$ for $N=256$ (solid), $N=1024$ (dashed),
           and $N=4096$ (dotted) elements.}
  \label{fig:Cubic_DG_p}
\end{figure}

Results of the numerical experiments are visualized in \autoref{fig:Cubic_DG_p},
where numerical solutions at the final time $T$ are shown at the left-hand side.
On the right-hand side, the evolution of the discrete total entropy
$\norm{\vec{u}}_M^2$ (summed over all elements) is visualized. Here, the mass
matrix $\mat{M}$ is a diagonal matrix containing the Lobatto Legendre quadrature
weights on the diagonal, scaled by an appropriate factor to account for the
width of each cell.

The numerical solutions for the same polynomial degree $p$ are visually nearly
indistinguishable.
Using polynomials of degree $p=1$, the numerical solutions converge to the classical
entropy solution (\autoref{fig:Cubic_DG_p1}). If polynomials of higher degree $p \geq 2$
are used, nonclassical shocks develop and remain stable and unchanged if the number
of elements is increased. As in the periodic case, the final value of the total
entropy is higher if nonclassical shocks occur.

Introducing another kind of dissipation (besides the dissipation provided by the
numerical fluxes as in the previous examples), modal filtering is applied after
every complete time step of the Runge-Kutta method. This modal filtering can be
seen as a discretization of the $s$th power of the Legendre dissipation operator,
cf.\ Example~\ref{ex:Legendre-filtering}. Here, the strength $\epsilon$ is chosen as
$\epsilon = -\log(\mathtt{eps}) / (\Delta t \, (n (n+1))^s)$, where $\mathtt{eps}$
is the machine accuracy, i.e.\ $\mathtt{eps} = \num{2.220446049250313e-16}$ for
$64$ bit floating point numbers (\texttt{Float64} in Julia \cite{bezanson2017julia})
used in the calculations.

\begin{figure}
\centering
  \begin{subfigure}{0.5\textwidth}
    \centering
    \includegraphics[width=\textwidth]{figures/Cubic_DG_p__legend}
  \end{subfigure}%
  \\
  \begin{subfigure}{0.45\textwidth}
    \centering
    \includegraphics[width=\textwidth]{%
      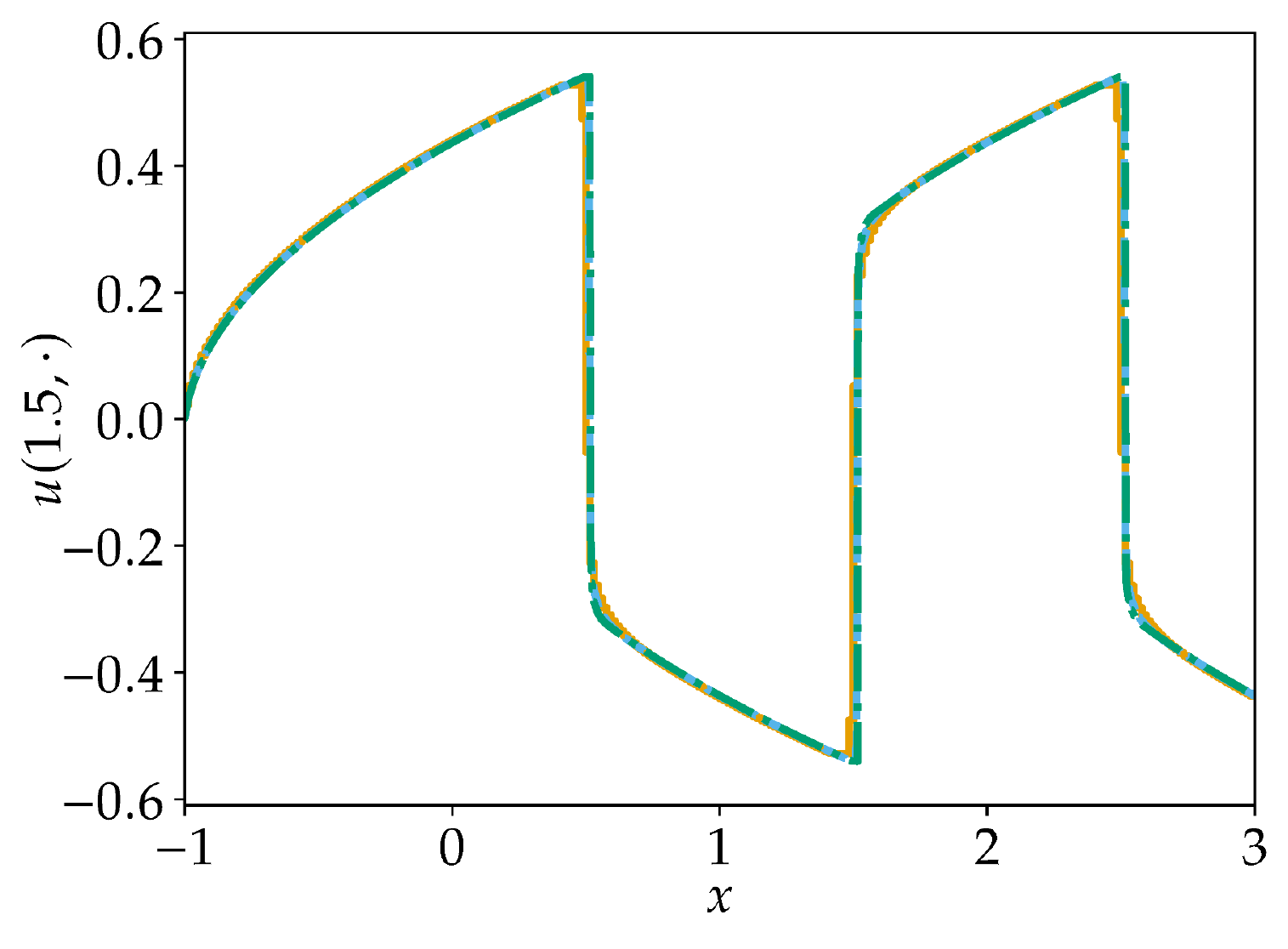}
    \caption{Filter order $s=1$;
             numerical solutions at $t=1.5$.}
  \label{fig:Cubic_DG_s1}
  \end{subfigure}%
  \hspace*{\fill}
  \begin{subfigure}{0.45\textwidth}
    \centering
    \includegraphics[width=\textwidth]{%
      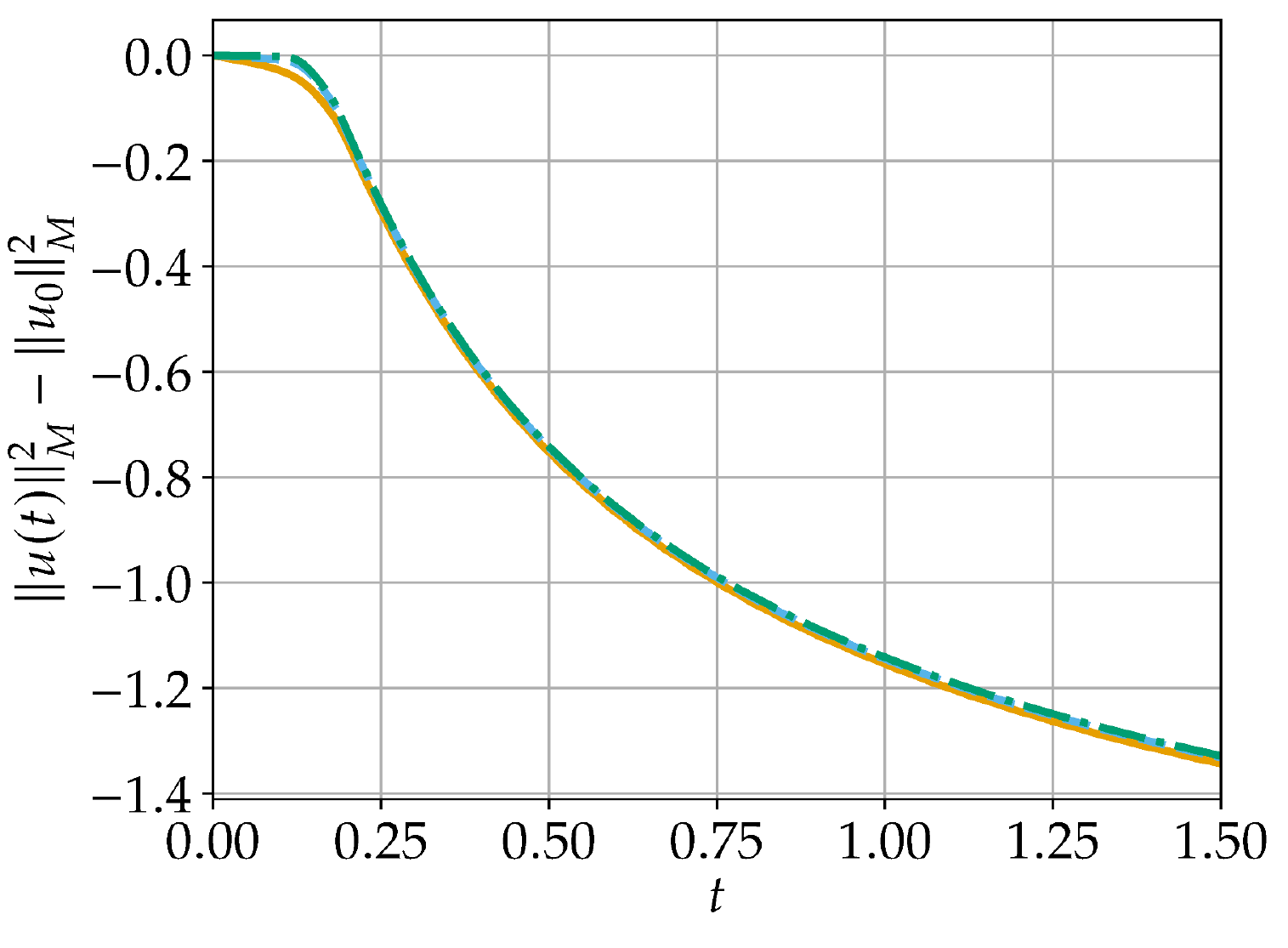}
    \caption{Filter order $s=1$;
             evolution of the $L^2$ entropy.}
  \end{subfigure}%
  \\
  \begin{subfigure}{0.45\textwidth}
    \centering
    \includegraphics[width=\textwidth]{%
      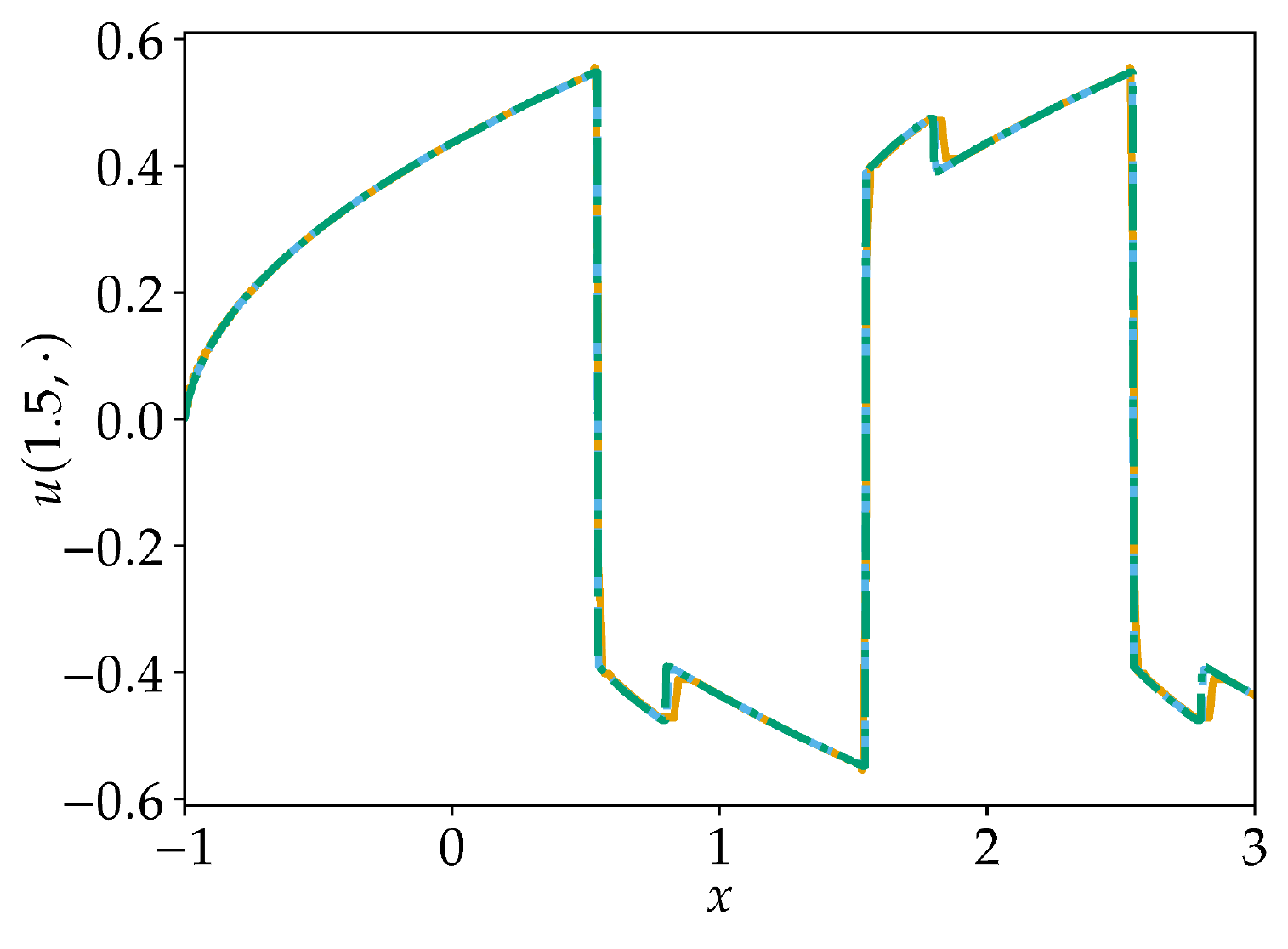}
    \caption{Filter order $s=4$;
             numerical solutions at $t=1.5$.}
  \label{fig:Cubic_DG_s4}
  \end{subfigure}%
  \hspace*{\fill}
  \begin{subfigure}{0.45\textwidth}
    \centering
    \includegraphics[width=\textwidth]{%
      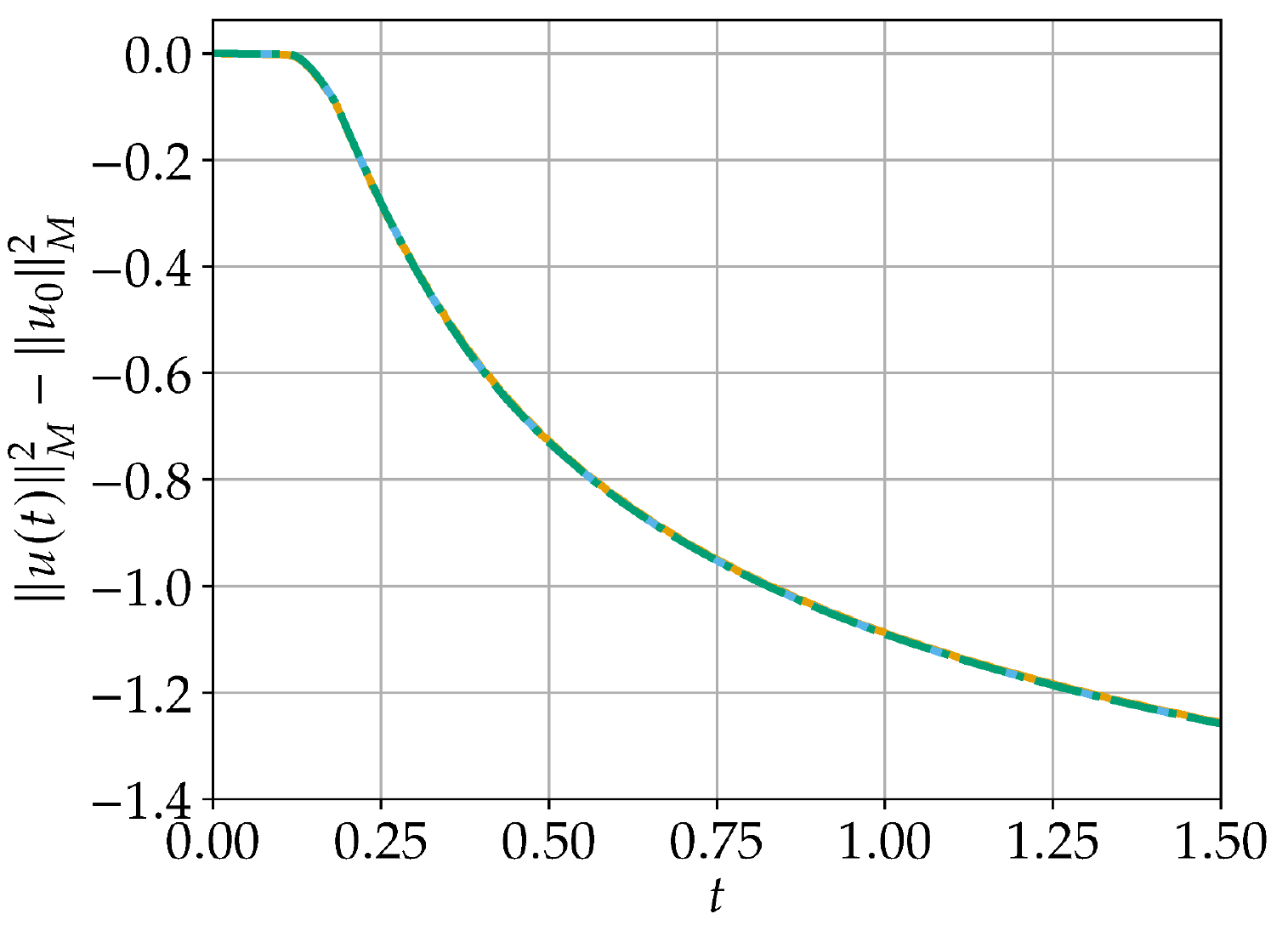}
    \caption{Filter order $s=4$;
             evolution of the $L^2$ entropy.}
  \end{subfigure}%
  \\
  \begin{subfigure}{0.45\textwidth}
    \centering
    \includegraphics[width=\textwidth]{%
      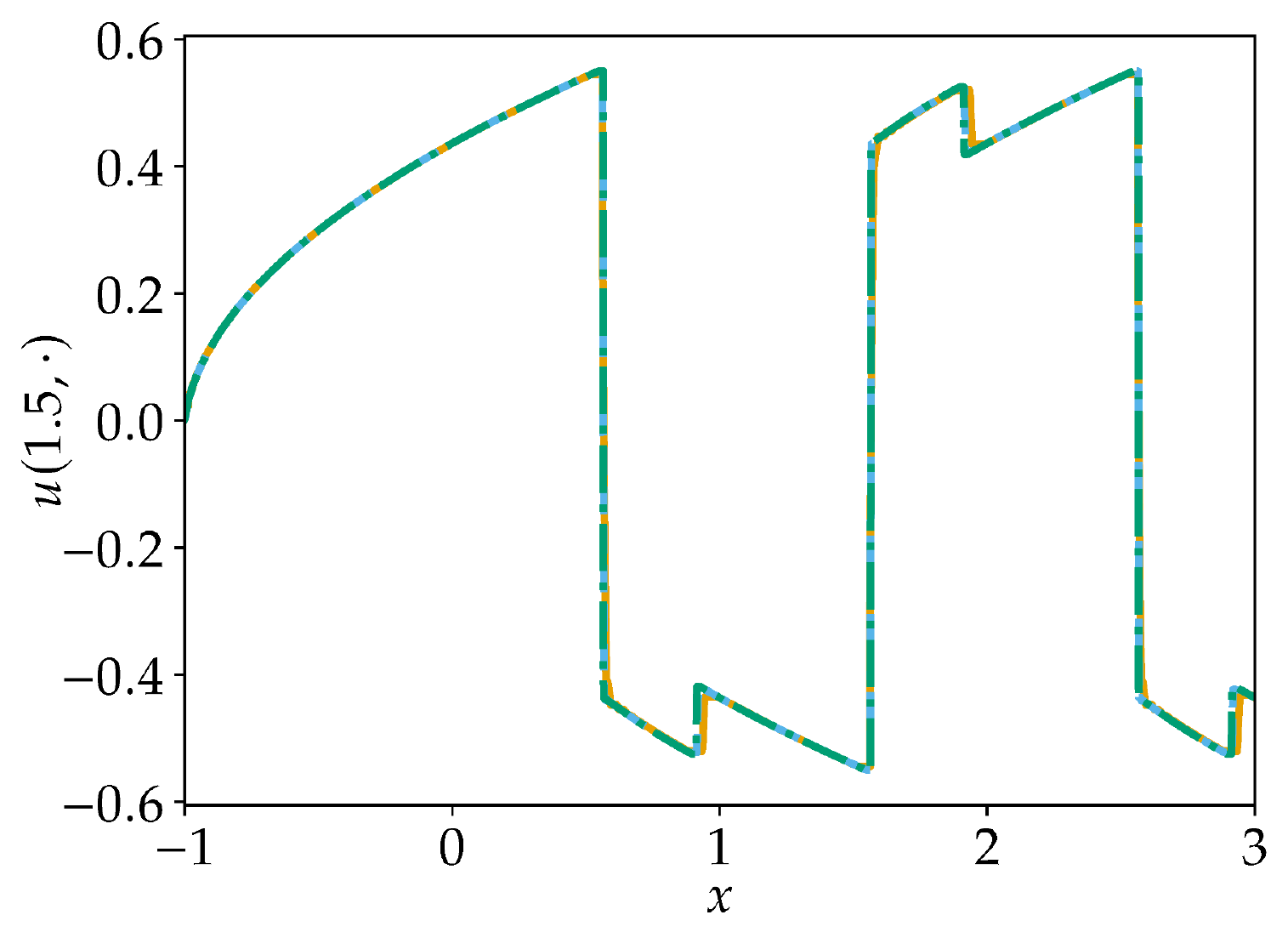}
    \caption{Filter order $s=5$;
             numerical solutions at $t=1.5$.}
  \label{fig:Cubic_DG_s5}
  \end{subfigure}%
  \hspace*{\fill}
  \begin{subfigure}{0.45\textwidth}
    \centering
    \includegraphics[width=\textwidth]{%
      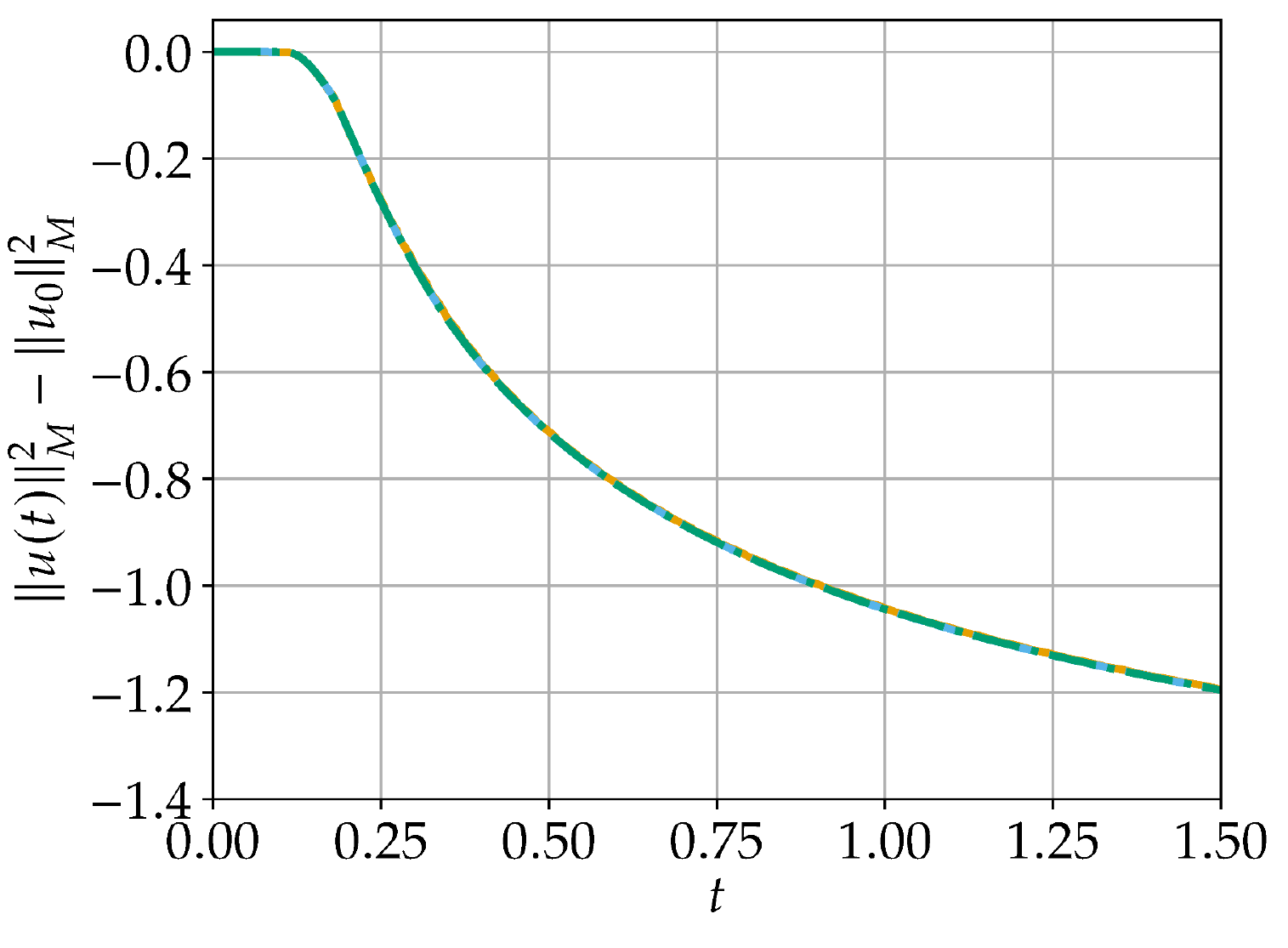}
    \caption{Filter order $s=5$;
             evolution of the $L^2$ entropy.}
  \end{subfigure}%
  \caption{Numerical results for DG methods using Godunov's flux and the entropy
           stable semi-discretization of Theorem~\ref{thm:flux-diff-ES} in order
           to approximate solutions of the cubic conservation law \eqref{eq:cubic-IBVP}
           with polynomial degree $p=4$ and different filter orders $s$ for
           $N=256$ (solid), $N=1024$ (dashed), and $N=4096$ (dotted) elements.}
  \label{fig:Cubic_DG_s}
\end{figure}

Results of the numerical experiments using DG methods with modal filtering of
different orders $s$ are shown in \autoref{fig:Cubic_DG_s}. The
numerical solutions using different numbers $N$ of elements are visually nearly
indistinguishable. Using the filter order $s=1$, the numerical solutions converge
to the classical entropy weak solution (\autoref{fig:Cubic_DG_s1}). However, if
the filter order is increased, the numerical solutions converge to nonclassical
solutions (\autoref{fig:Cubic_DG_s4} and \autoref{fig:Cubic_DG_s5}). As before,
the appearance of nonclassical shocks is linked with less entropy dissipation.

\begin{figure}
\centering
  \begin{subfigure}{0.5\textwidth}
    \centering
    \includegraphics[width=\textwidth]{figures/Cubic__legend}
  \end{subfigure}%
  \\
  \begin{subfigure}{0.45\textwidth}
    \centering
    \includegraphics[width=\textwidth]{%
      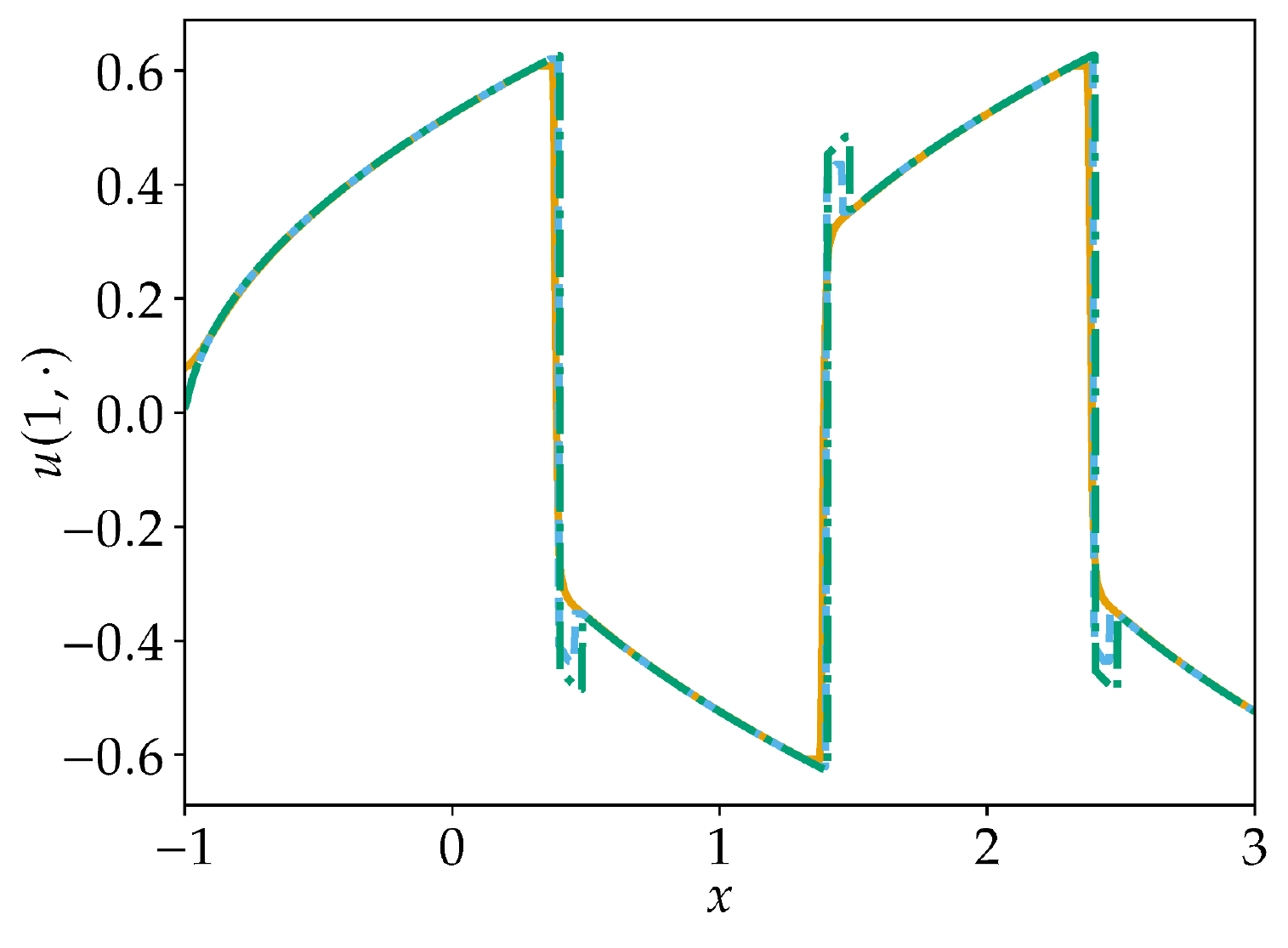}
    \caption{Dissipation strengths $\epsilon_2 = 100$, $\epsilon_4 = 100$,
             $\epsilon_6 = 100$;
             numerical solutions at $t=1.5$.}
  \label{fig:Cubic_FD_Strength2_100_Strength4_100_Strength6_100}
  \end{subfigure}%
  \hspace*{\fill}
  \begin{subfigure}{0.45\textwidth}
    \centering
    \includegraphics[width=\textwidth]{%
      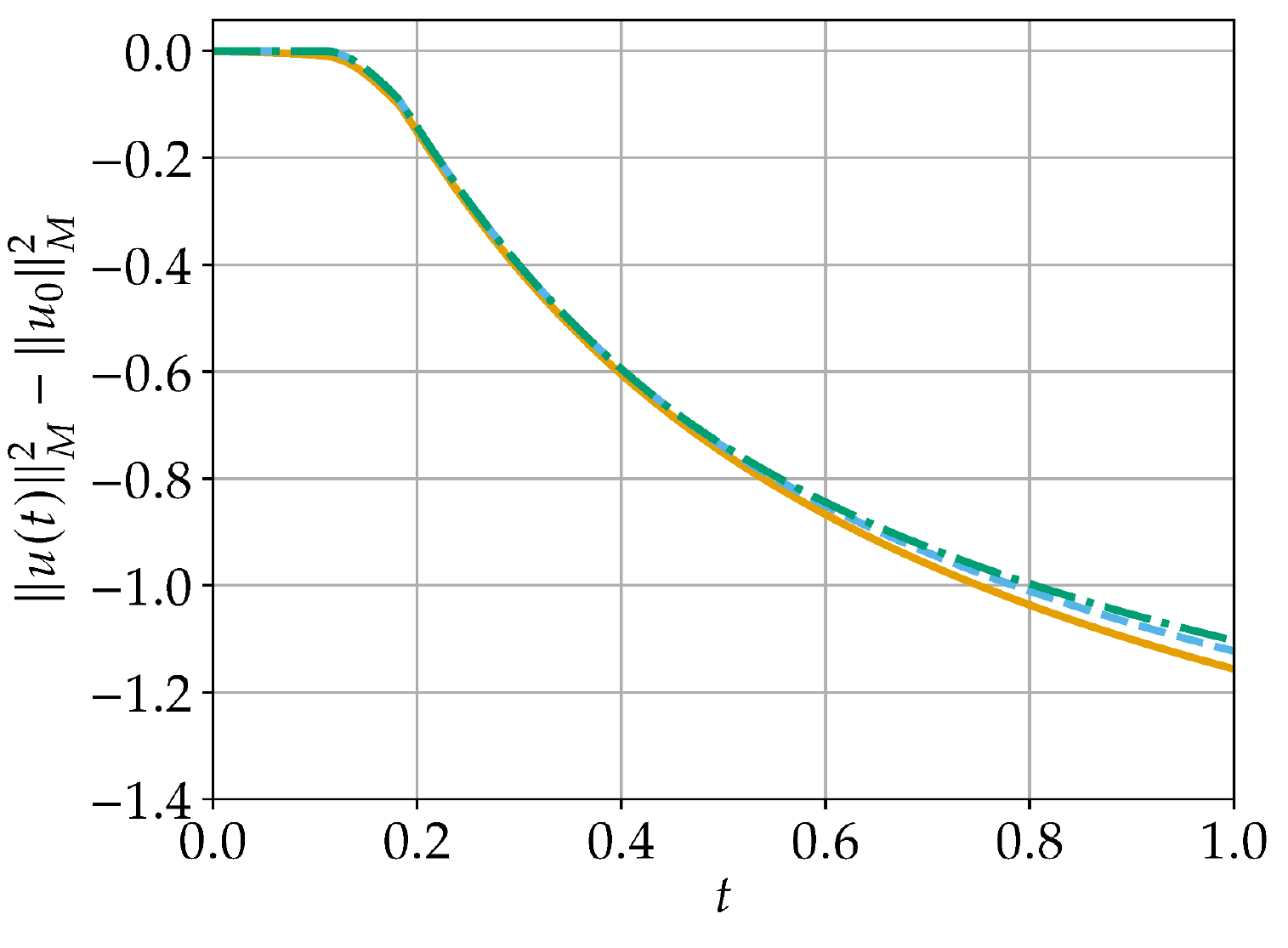}
    \caption{Dissipation strengths $\epsilon_2 = 100$, $\epsilon_4 = 100$,
             $\epsilon_6 = 100$;
             evolution of the $L^2$ entropy.}
  \end{subfigure}%
  \\
  \begin{subfigure}{0.45\textwidth}
    \centering
    \includegraphics[width=\textwidth]{%
      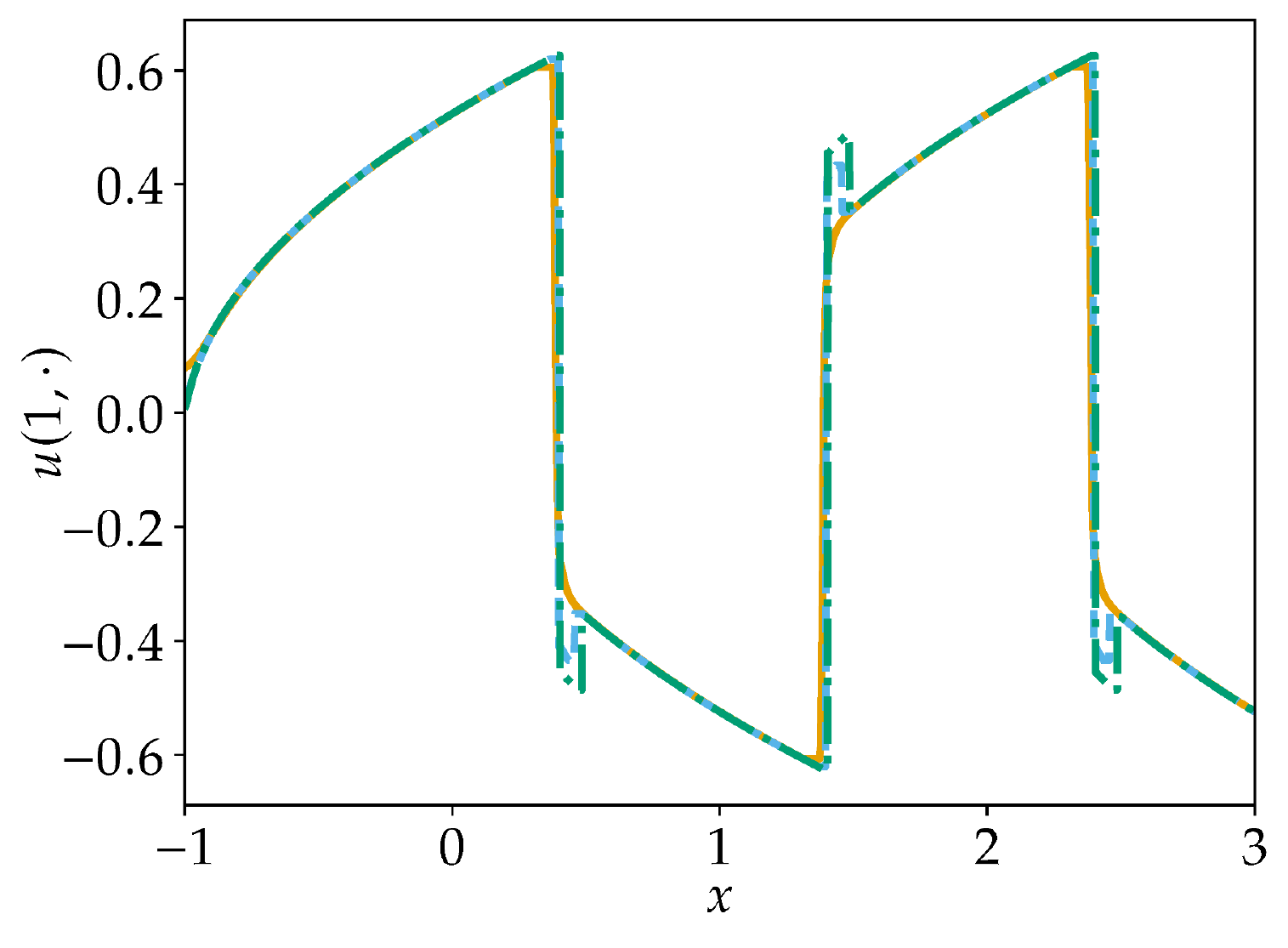}
    \caption{Dissipation strengths $\epsilon_2 = 0$, $\epsilon_4 = 100$,
             $\epsilon_6 = 100$;
             numerical solutions at $t=1.5$.}
  \label{fig:Cubic_FD_Strength2_100_Strength4_0_Strength6_100}
  \end{subfigure}%
  \hspace*{\fill}
  \begin{subfigure}{0.45\textwidth}
    \centering
    \includegraphics[width=\textwidth]{%
      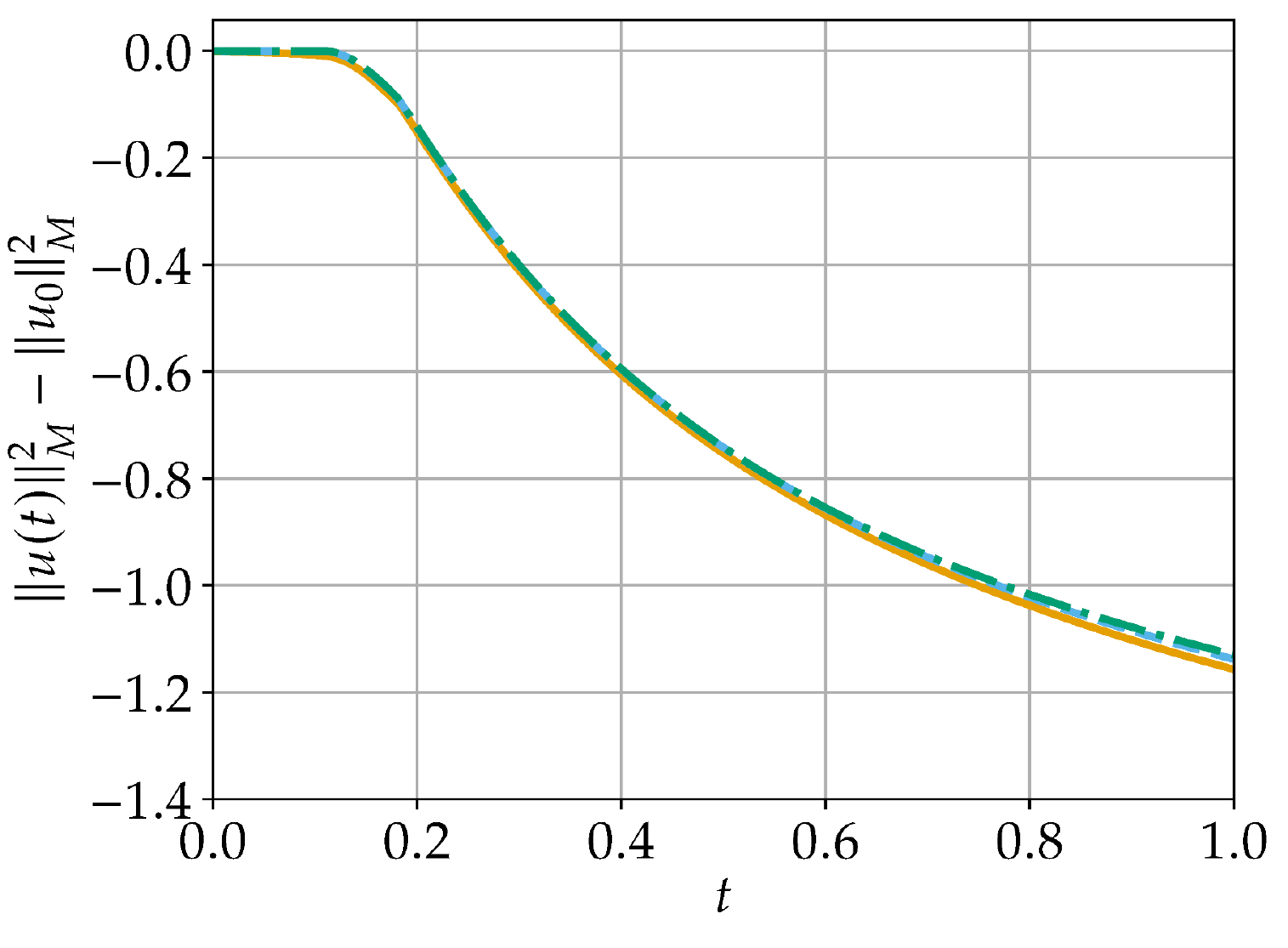}
    \caption{Dissipation strengths $\epsilon_2 = 0$, $\epsilon_4 = 100$,
             $\epsilon_6 = 100$;
             evolution of the $L^2$ entropy.}
  \end{subfigure}%
  \\
  \begin{subfigure}{0.45\textwidth}
    \centering
    \includegraphics[width=\textwidth]{%
      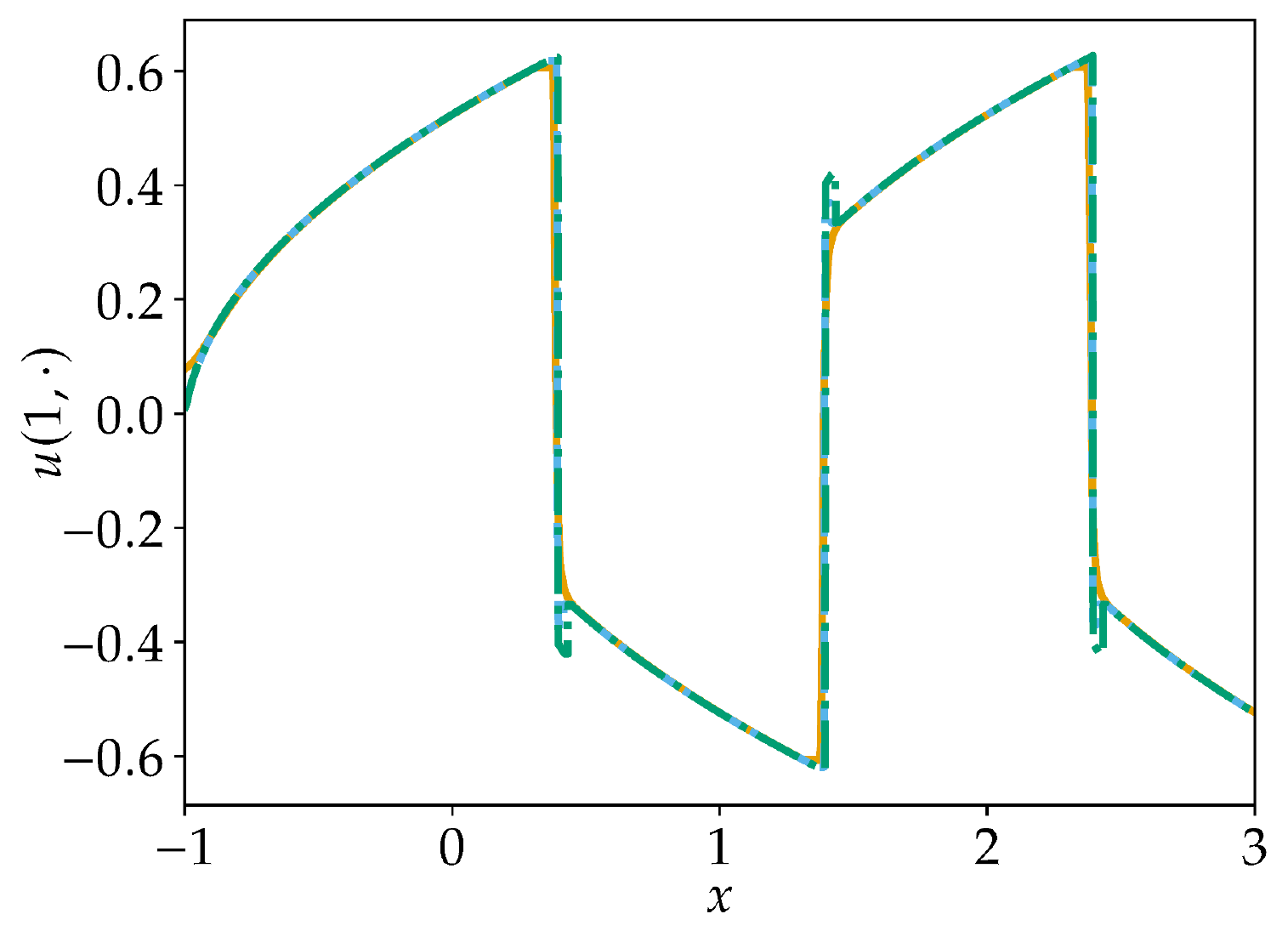}
    \caption{Dissipation strengths $\epsilon_2 = 100$, $\epsilon_4 = 100$,
             $\epsilon_6 = 0$;
             numerical solutions at $t=1.5$.}
  \label{fig:Cubic_FD_Strength2_100_Strength4_100_Strength6_0}
  \end{subfigure}%
  \hspace*{\fill}
  \begin{subfigure}{0.45\textwidth}
    \centering
    \includegraphics[width=\textwidth]{%
      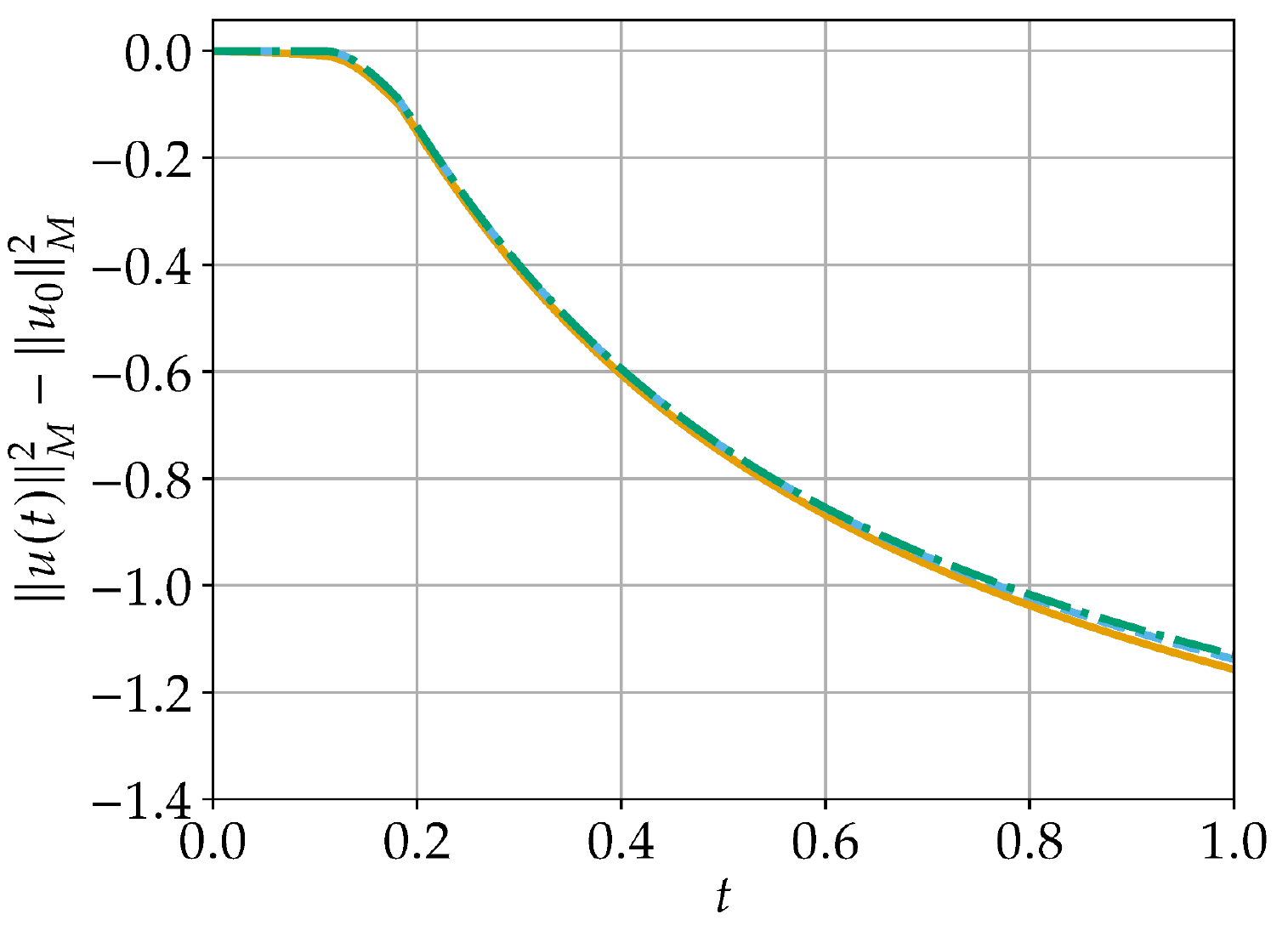}
    \caption{Dissipation strengths $\epsilon_2 = 100$, $\epsilon_4 = 100$,
             $\epsilon_6 = 0$;
             evolution of the $L^2$ entropy.}
  \end{subfigure}%
  \caption{Numerical results for SBP FD methods with interior order of accuracy
           six using Godunov's flux at the boundaries and the entropy-stable
           semi-discretization of Theorem~\ref{thm:flux-diff-ES} in order to
           approximate solutions of the cubic conservation law \eqref{eq:cubic-IBVP}
           with different artificial dissipation operators for $N=1024$ (solid),
           $N=4096$ (dashed), and $N=\num{16384}$ (dotted) grid points.}
  \label{fig:Cubic_FD}
\end{figure}

Furthermore, finite difference SBP methods with the same
boundary procedures as the DG schemes are tested as well. The results using only
a single kind of artificial dissipation operator are similar to the ones in the
periodic case and thus not shown here. Instead, the artificial dissipation operators
of \cite{mattsson2004stable} are weighted with strengths $\epsilon_2$ (second-order
dissipation), $\epsilon_4$ (fourth-order dissipation), and $\epsilon_6$
(sixth-order dissipation) and added to the semi-discretization. The time step is
chosen as $\Delta t = \frac{1}{N}$.
The results of these numerical experiments are presented in \autoref{fig:Cubic_FD}.
As can be seen there, for a fixed choice of the strengths $\epsilon_i$, nonclassical
shocks occur under grid refinement. Thus, the influence of the higher-order
dissipation operators can destroy the convergence to the classical solution induced
by the second-order dissipation operator.

\section{Uniqueness and entropy properties for the cubic conservation law}
\label{sec:cubic-TeCNO}

\subsection{Preliminaries}

Here, the TeCNO schemes in \cite{fjordholm2012arbitrarily,fjordholm2013high}
mentioned in Example~\ref{ex:diss-TeCNO} are now used to compute numerical
solutions of the cubic conservation law \eqref{eq:cubic} with periodic boundary
conditions. As in the Section~\ref{sec:cubic}, the initial condition is chosen as
$u_0(x) = - \sin(\pi x)$ in the domain $(\xmin,\xmax) = (-1,1)$ and the numerical
solutions are evolved up to the final time $T=1$ using the fourth-order, ten stage,
strong stability preserving explicit Runge-Kutta method SSPRK(10,4) of
\cite{ketcheson2008highly} with time steps $\Delta t = 1 / N$. The following
entropy functions will be considered.
\begin{itemize}
  \item
  The $L^2$ entropy $U(u) = \frac{1}{2} u^2$.
  As in Section~\ref{sec:cubic}, the entropy variables are $w(u) = U'(u) = u$
  and the flux potential is $\psi(u) = \frac{1}{4} u^4$. Thus, the corresponding
  entropy-conservative numerical flux is \eqref{eq:cubic-fnum-EC}.

  \item
  The $L^4$ entropy $U(u) = \frac{1}{4} u^4$.
  In this case, $w(u) = U'(u) = u^3$ and the flux potential is $\psi(u) =
  \frac{1}{2} u^6$. Hence, the entropy-conservative numerical flux is the central
  flux
$
    \fnum(u_-,u_+)
    =
    \frac{\psi(u_+) - \psi(u_-)}{w(u_+) - w(u_-)}
    =
    \frac{u_+^3 + u_-^3}{2}.
$

  \item
  The $L^2 \cap L^4$ entropy $U(u) = \frac{1}{4} u^4 + \frac{\alpha}{2} u^2$,
  $\alpha > 0$.
  For this strictly convex entropy, $w(u) = U'(u) = u^3 + \alpha u$ and $\psi(u)
  = \frac{1}{2} u^6 + \frac{\alpha}{4} u^4$. The corresponding entropy-conservative
  numerical flux is
  \begin{equation}
  \begin{aligned}
    &\quad
    \fnum(u_-,u_+)
    =
    \frac{1}{2} \frac{u_+^6 - u_-^6 + \frac{\alpha}{2} u_+^4 - \frac{\alpha}{2} u_-^4}
                     {u_+^3 - u_-^3 + \alpha u_+ - \alpha u_-}
    \\
    &=
    \frac{1}{2} \frac{u_+^5 + u_+^4 u_- + u_+^3 u_-^2 + u_+^2 u_-^3 + u_+ u_-^4 + u_-^5
    + \frac{\alpha}{2} \bigl( u_+^3 + u_+^2 u_- + u_+ u_-^2 + u_-^3 \bigr)}
    {u_+^2 + u_+ u_- + u_-^2 + \alpha},
  \end{aligned}
  \end{equation}
  where the fraction has been reduced by $(u_+ - u_-)$. In the numerical experiments
  presented in the following, $\alpha = \frac{1}{100}$ has been chosen.
\end{itemize}

\subsection{Numerical results}

\begin{figure}
\centering
  \begin{subfigure}{0.55\textwidth}
    \centering
    \includegraphics[width=\textwidth]{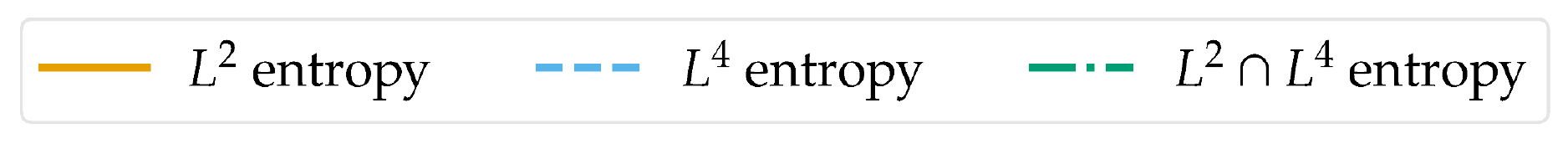}
  \end{subfigure}%
  \\
  \begin{subfigure}{0.45\textwidth}
    \centering
    \includegraphics[width=\textwidth]{%
      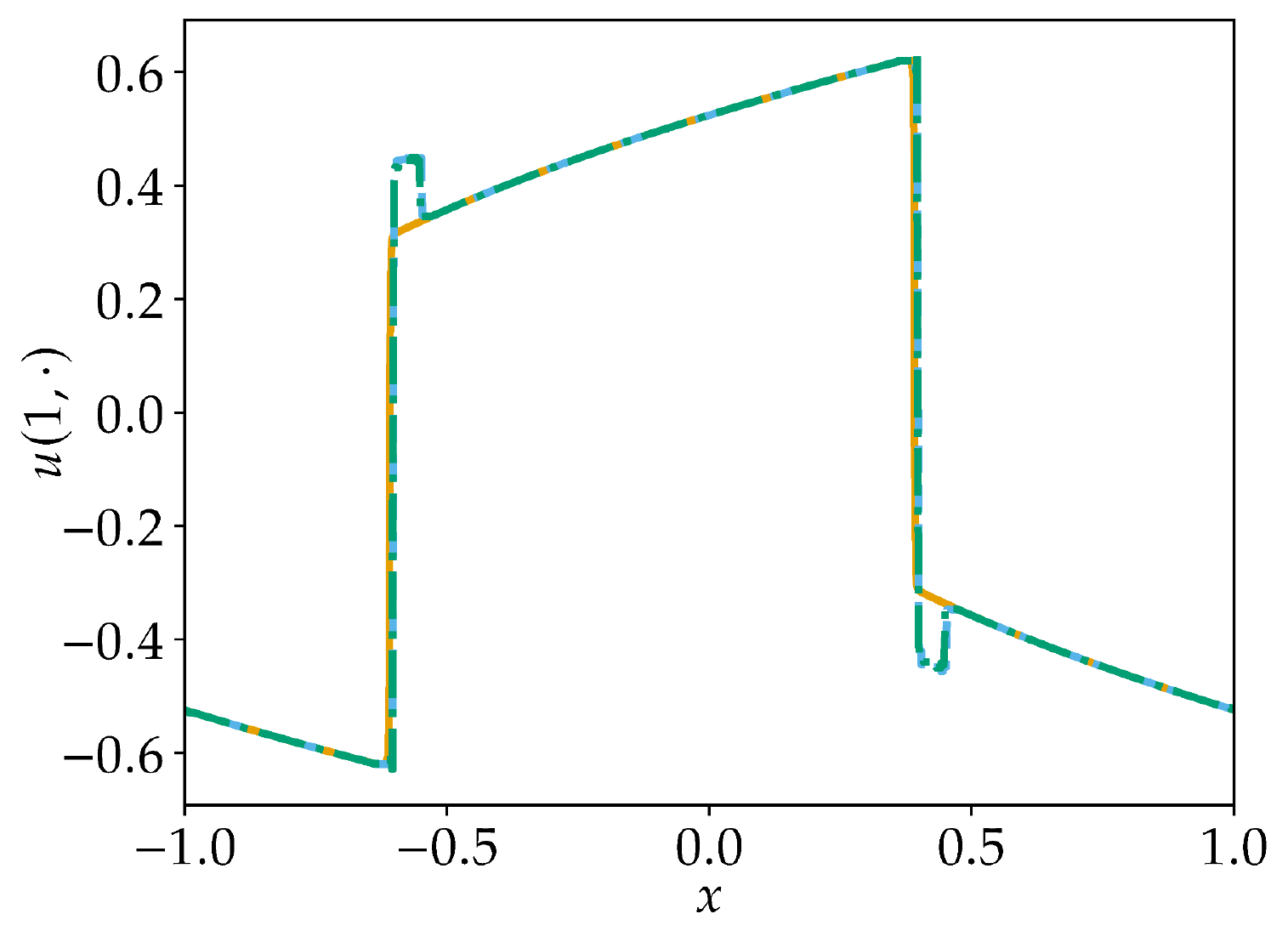}
    \caption{$N = 2^{10} = \num{1024}$; numerical solutions at $t=1$.}
  \end{subfigure}%
  \hspace*{\fill}
  \begin{subfigure}{0.45\textwidth}
    \centering
    \includegraphics[width=\textwidth]{%
      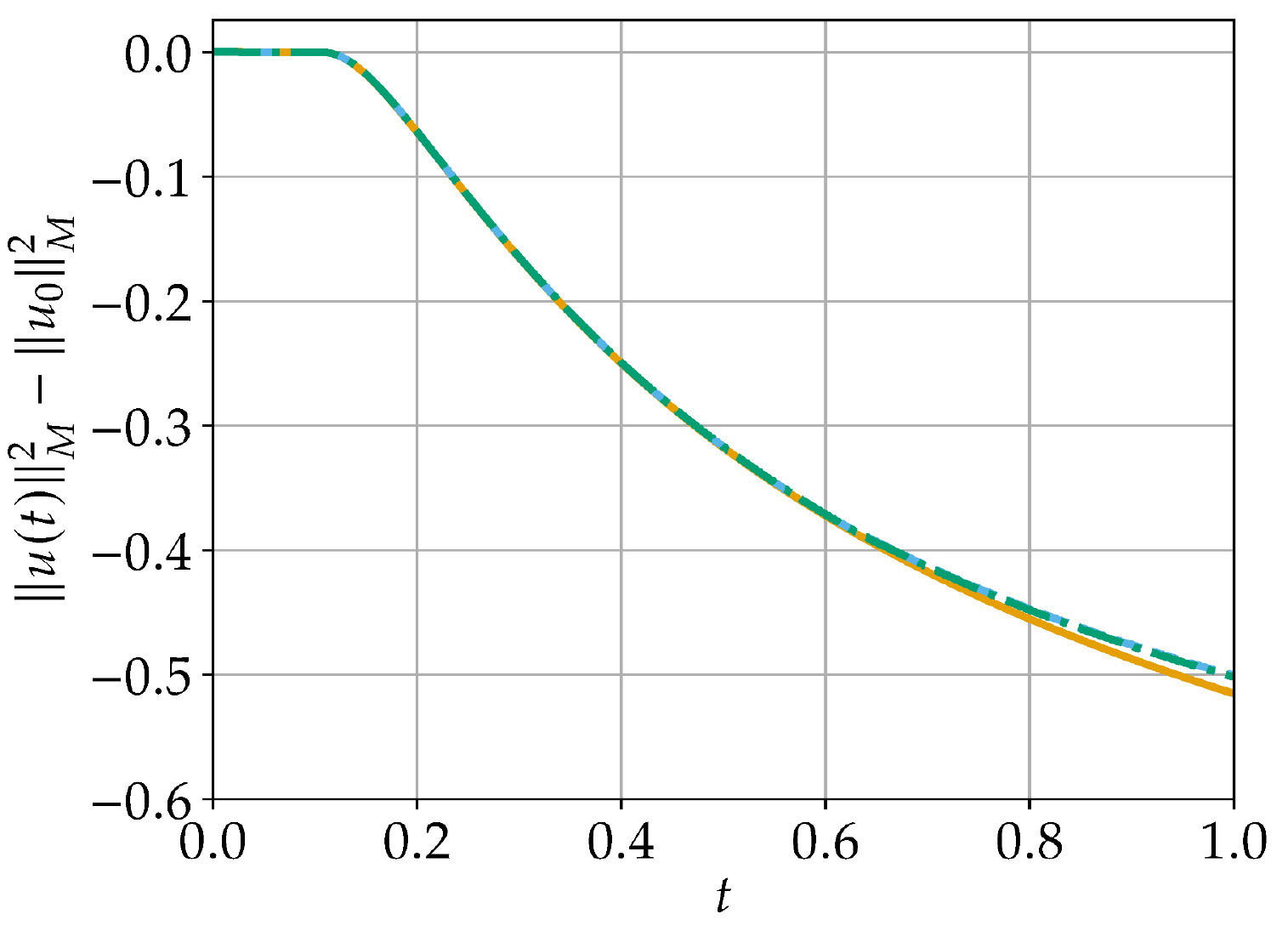}
    \caption{$N = 2^{10} = \num{1024}$; evolution of the $L^2$ entropy.}
  \end{subfigure}%
  \\
  \begin{subfigure}{0.45\textwidth}
    \centering
    \includegraphics[width=\textwidth]{%
      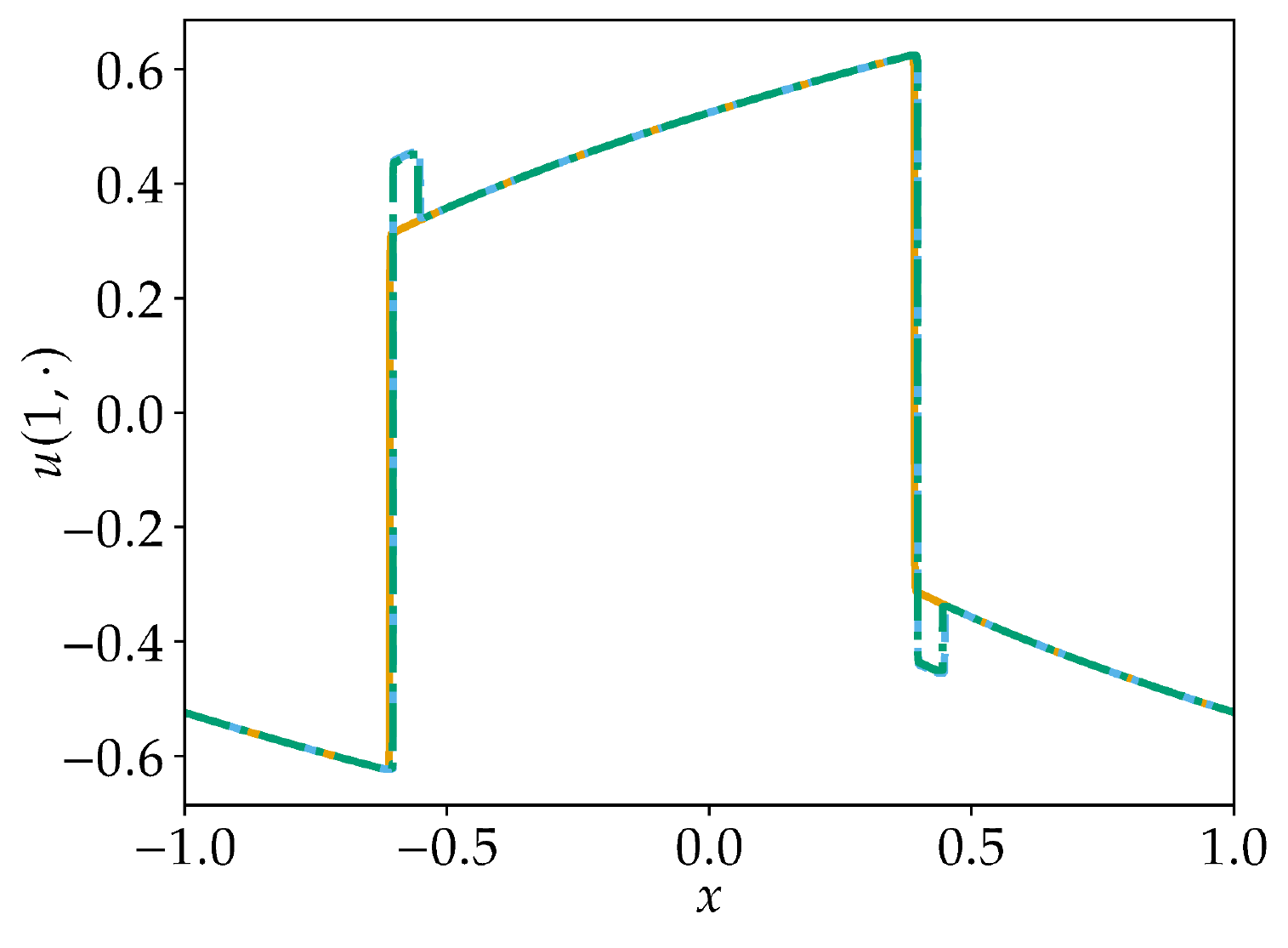}
    \caption{$N = 2^{12} = \num{4096}$; numerical solutions at $t=1$.}
  \end{subfigure}%
  \hspace*{\fill}
  \begin{subfigure}{0.45\textwidth}
    \centering
    \includegraphics[width=\textwidth]{%
      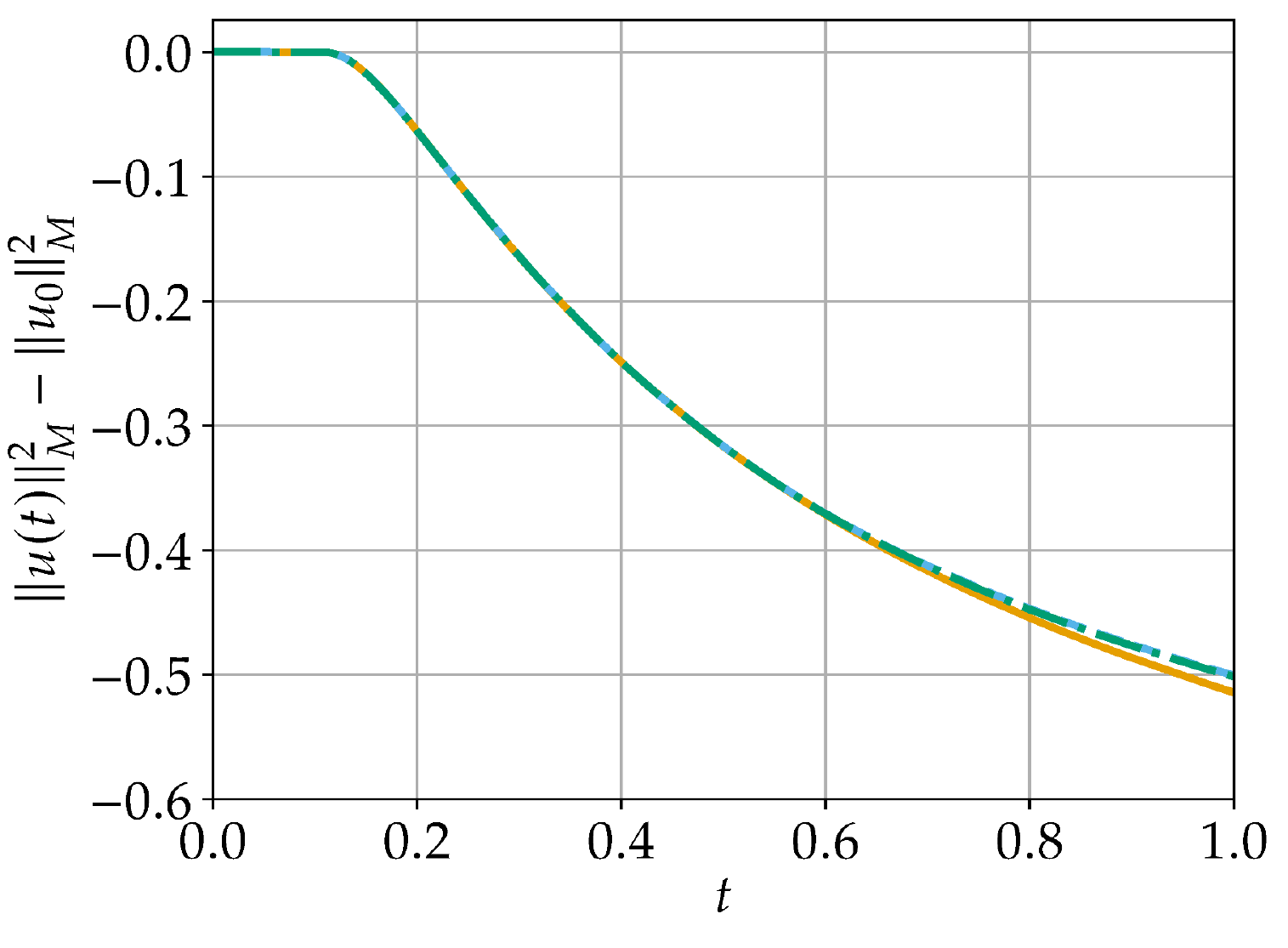}
    \caption{$N = 2^{12} = \num{4096}$; evolution of the $L^2$ entropy.}
  \end{subfigure}%
  \\
  \begin{subfigure}{0.45\textwidth}
    \centering
    \includegraphics[width=\textwidth]{%
      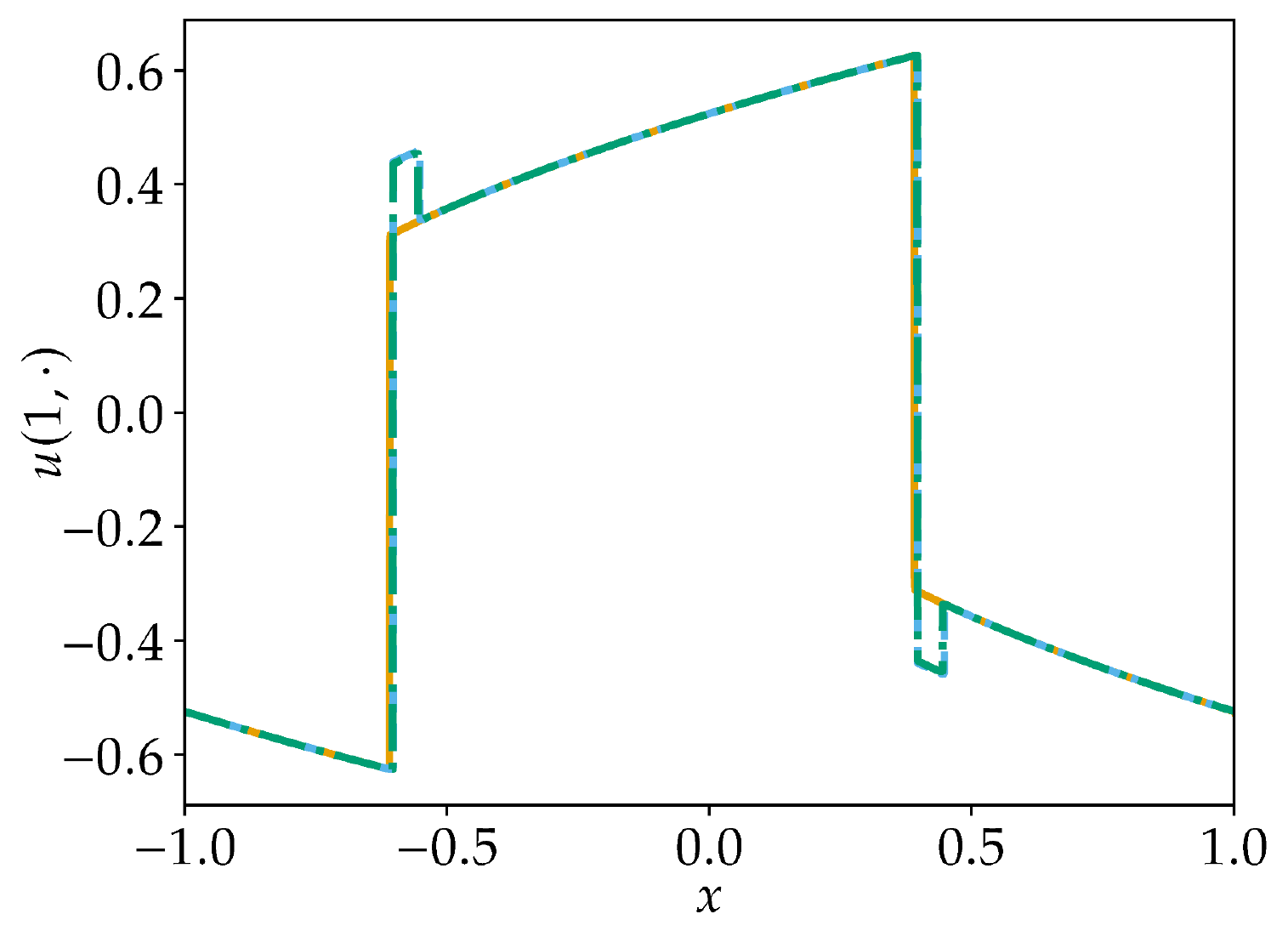}
    \caption{$N = 2^{14} = \num{16384}$; numerical solutions at $t=1$.}
  \end{subfigure}%
  \hspace*{\fill}
  \begin{subfigure}{0.45\textwidth}
    \centering
    \includegraphics[width=\textwidth]{%
      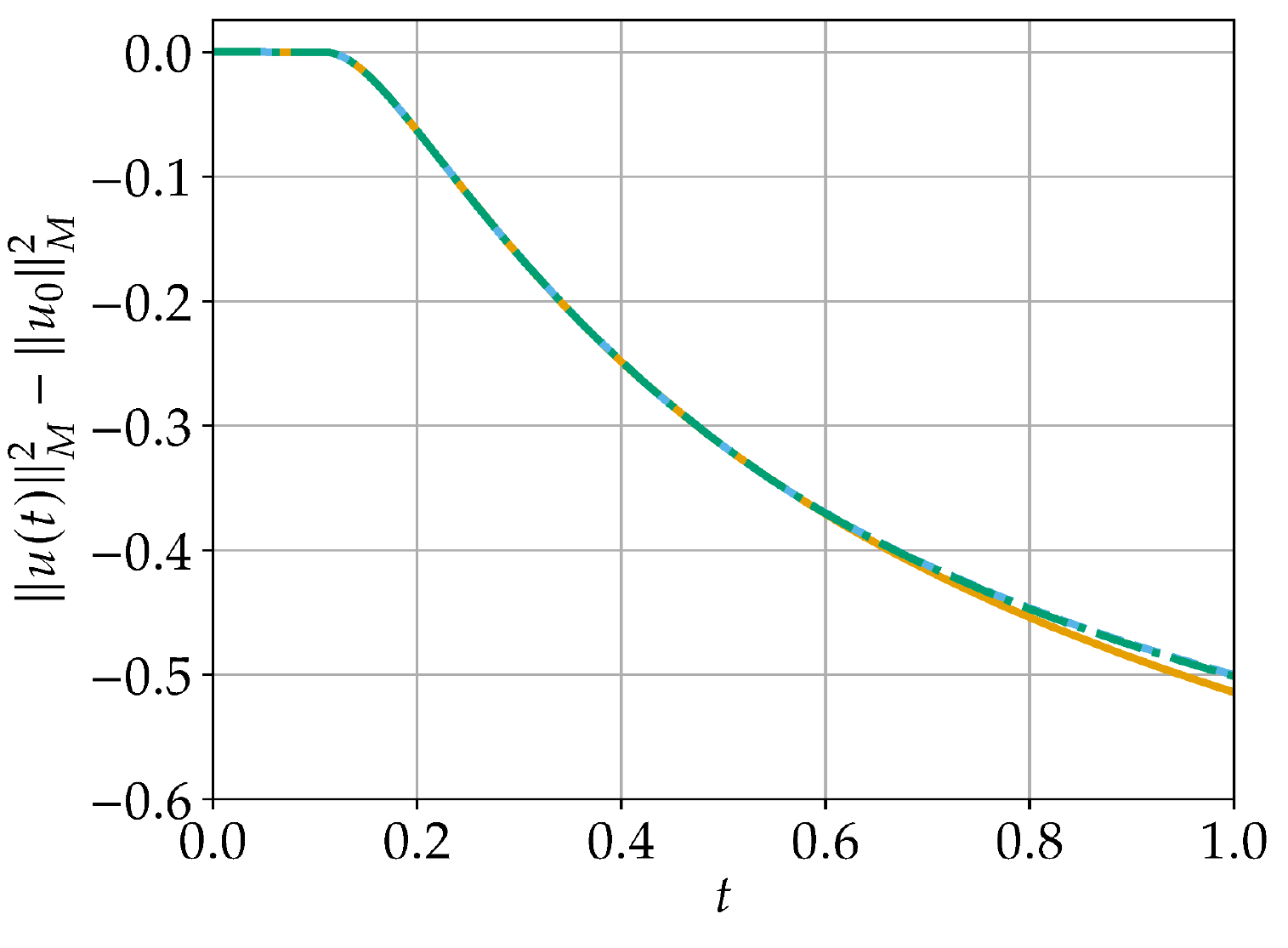}
    \caption{$N = 2^{14} = \num{16384}$; evolution of the $L^2$ entropy.}
  \end{subfigure}%
  \caption{Numerical results for TeCNO(3) methods approximating solutions
           of the cubic conservation law \eqref{eq:cubic} based on different
           entropy functions for $N$ of cells.}
  \label{fig:Cubic_TeCNO}
\end{figure}

Results of numerical experiments with TeCNO(3) schemes can be seen in
\autoref{fig:Cubic_TeCNO}. For the schemes based on the $L^2$ entropy, the numerical
solutions are visually indistinguishable from the entropy weak solution. In contrast,
the schemes based on the $L^4$ entropy yield numerical solutions with overshoots
at the discontinuities that do not vanish as the grid is refined. Results using
TeCNO(2) schemes are similar; the overshoot regions are a bit smaller but clearly
present for all investigated numbers $N$ of cells.
It might be conjectured that the ``better'' behavior of the schemes based on the $L^2$
entropy $U(u) = \frac{1}{2} u^2$ compared to the $L^4$ entropy $U(u) = \frac{1}{4} u^4$
is influenced by the fact that the former entropy is strictly convex while the latter
is only convex. However, the schemes based on the strictly convex $L^2 \cap L^4$
entropy show the same behavior as the ones based on the $L^4$ entropy, i.e.\
the overshoots persist.

\section{Computing kinetic functions numerically}
\label{sec:kinetic-function-cubic}

\subsection{Preliminaries}

To analyze the behavior of the provably entropy-dissipative numerical methods
described above in more detail, the corresponding kinetic functions will be
computed.
The basic motivation is as follows \cite[Chapter~II]{lefloch2002hyperbolic}.
The locally smooth parts of a weak solution of a nonlinear scalar conservation
law are unique but non-uniqueness can arise from discontinuities if the flux
is non-convex and only a single entropy inequality is required. Hence, an approach
to single out one specific weak solutions among all weak solutions satisfying
a single entropy inequality is given by prescribing the allowed forms of
discontinuities. This is exactly the purpose of kinetic functions
\cite[Chapter~II]{lefloch2002hyperbolic}.

Consider a scalar nonlinear conservation law $\partial_t u + \partial_x f(u) = 0$
in one space dimension and a weak solution of an associated Riemann problem
with left- and right-hand states $u_L$ and $u_R$. This weak solution is a
combination of rarefaction waves, classical shock waves (satisfying all entropy
inequalities locally), and nonclassical shock waves (which are only required
to satisfy a single entropy inequality locally).
In this context, the kinetic function $\phi^\flat$ is the mapping of the
left state $u_L$ to the middle state $u_M$ if a nonclassical solution appears.
We use the following definition of a kinetic function in the context of
numerical solutions.
\begin{definition}
  Given a numerical solution of a Riemann problem with left- and right-hand
  states $u_L$ and $u_R$, define $\Upsilon_L \subseteq \R$ as the set of all
  left-hand states $u_L$ such that the numerical solution results contains an
  (approximately) constant part with value $u_M$ that increases the
  total variation and is (approximately) connected to the left-hand state
  via a single discontinuity.
  The \emph{kinetic function} associated to the numerical method is the mapping
  $\phi^\flat\colon \Upsilon_L \to \R$, $\phi^\flat(u_L) = u_M$.
\end{definition}

Here, a series of Riemann problems with right-hand state $u_R = -2$ and varying
left-hand state $u_L$ has been solved for each scheme.
The FD SBP and DG methods use a domain $[-1, 3]$ with initial discontinuity
located at $x = -0.5$. The solution is computed until $t = 5 / \max\{3 u_0^2\}$
with a time step of $\dt = ((p^2+1) N \max\{3 u_0^2\})^{-1}$ for DG methods
and $\dt = (N \max\{3 u_0^2\})^{-1}$ for FD methods.

\begin{figure}[htp]
\centering
  \includegraphics[width=0.7\textwidth]{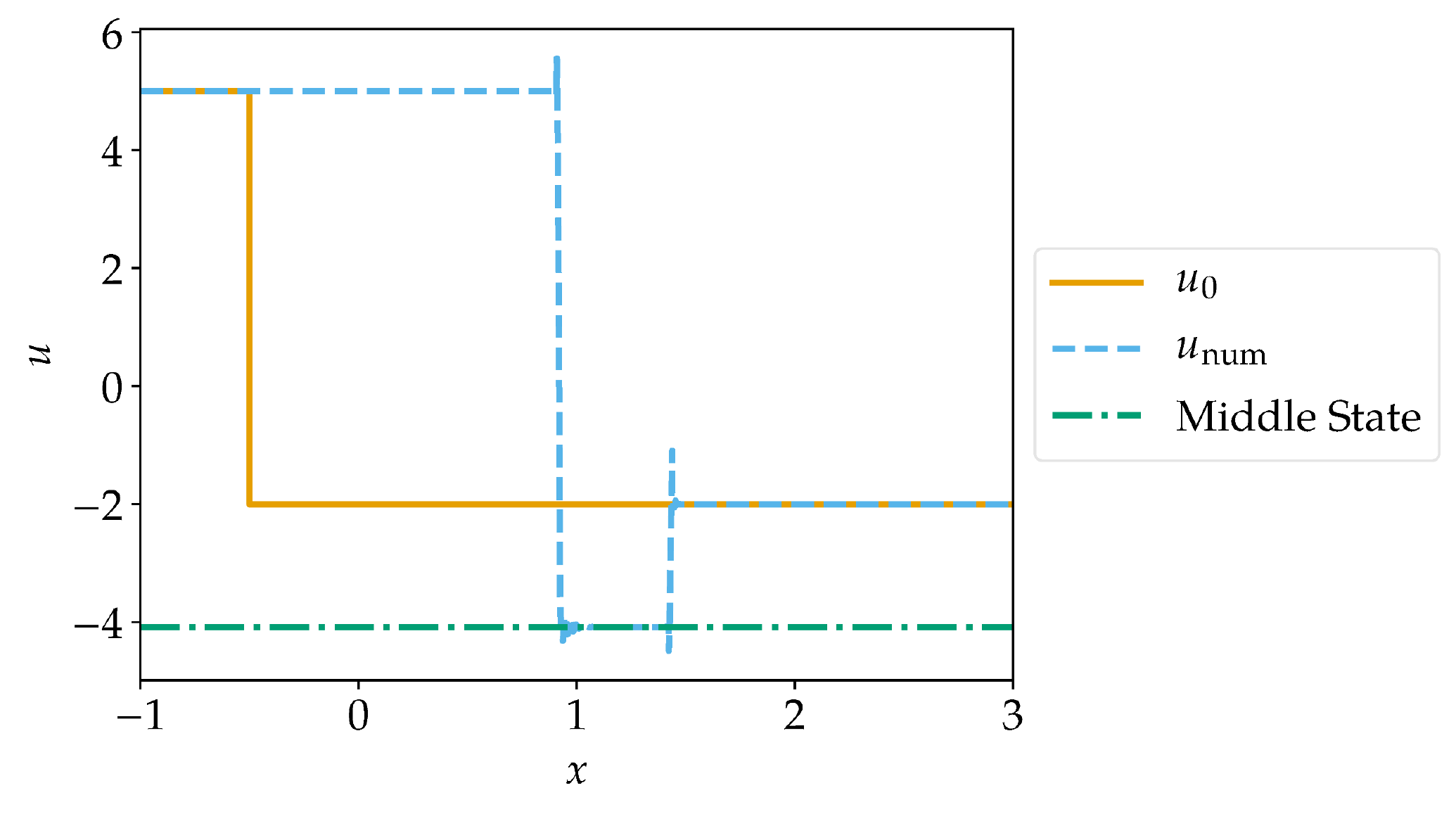}
  \caption{Numerical Riemann solution $u_\mathrm{num}$ (without postprocessing)
with $u_L = 5$ and $u_R = -2$, and a DG method with parameters:
           polynomial degree $p = 5$,
           filter order $s = 5$,
           number of elements $N = 256$.}
  \label{fig:Cubic_DG_p_5_s_5_N_00256__Riemann}
\end{figure}

The Fourier methods use a domain $[-6, 6]$ with initial value
\begin{equation}
  u_0(x) =
    u_L  \, \text{ for }  x \in [-4.5, 0],
    \,
 \text{ while }     u_0(x) =  u_R \,  \text{otherwise}.
\end{equation}
Again, the time step is $\dt = (N \max\{3 u_0^2\})^{-1}$.
The other parameters are the same as described above.
A typical numerical solution obtained using a DG scheme is shown in
Figure~\ref{fig:Cubic_DG_p_5_s_5_N_00256__Riemann}.
The middle state has been computed as follows. At first, the discontinuities
are detected in a simple way by averaging the solution locally and computing
the standard deviation. If there are two discontinuities of the form allowing
a nonclassical middle state, its value is computed as the median of the values
of the numerical solution between the two discontinuities. The general theory of
nonclassical shocks predicts bounds of the kinetic function $\phi^\flat$.
Here, as established in \cite{lefloch2002hyperbolic}, we have
\begin{equation}
\label{eq:bounds-kinetic-function}
  \phi^\flat_0 \leq \phi^\flat \leq \phi^\sharp,
  \qquad
  \text{where }
  \phi^\flat_0(u_L) = - u_L
  \text{ and }
  \phi^\sharp(u_L) = - u_L/2.
\end{equation}

\subsection{Nodal DG methods}
\label{sec:kinetic-function-cubic-DG}

For nodal DG methods, polynomial degrees $p \in \{0, 1, \dots, 5\}$,
filter order $s \in \{0, 1, \dots, 5\}$, and numbers of elements
$N \in \{2^6, 2^7, \dots, 2^{11}\}$ have been used.
The following observations have been made.
\begin{enumerate}
  \item
  The schemes with filter order $s \in \{1, 2, 3\}$ did not result
  in nonclassical solutions.

  \item
  The schemes with polynomial degrees $p \in \{0, 1\}$ did not result
  in nonclassical solutions. The schemes with polynomial degree $p = 2$
  resulted in nonclassical solutions only if no filtering was applied
  ($s = 0$).

  \item
  If the strength of the initial discontinuity is too big, no nonclassical
  solutions occur. For the investigated range of parameters, nonclassical
  solutions occurred only for $u_L < 10$ (as a necessary criterion). Depending
  on the other parameters, the maximal value of $u_L$ for which nonclassical
  solution occurred can be smaller.

  \item
  The numerically obtained kinetic functions $\phi^\flat$ are affine linear.
  These affine linear functions remain visually indistinguishable under
  grid refinement by increasing the number of elements $N$. This is shown
  for $p = 5$ and $s = 5$ in Figure~\ref{fig:Cubic_DG_p_5_s_5__kinetic_functions}.
  \begin{figure}[htb]
  \centering
    \includegraphics[width=0.7\textwidth]{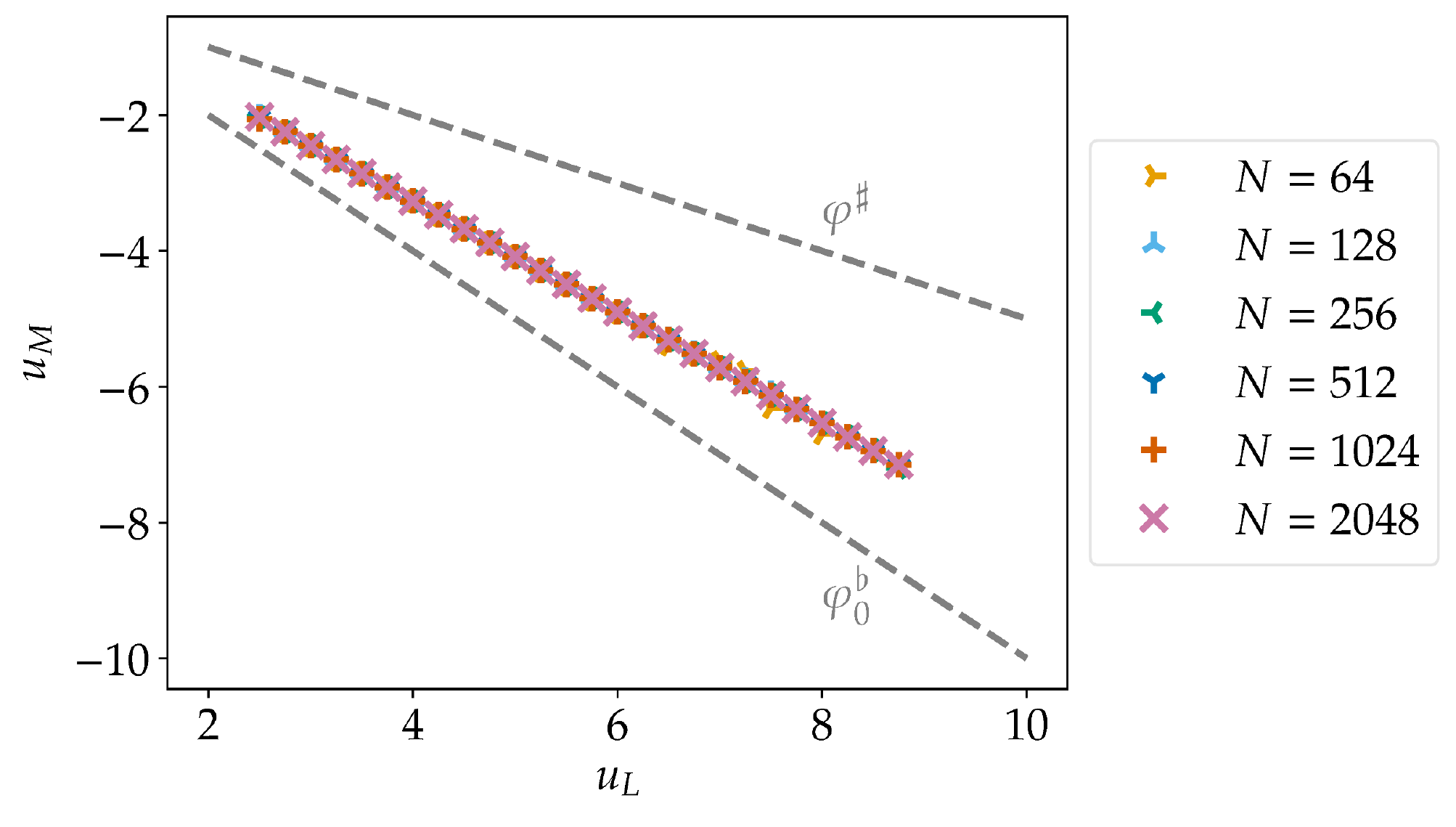}
    \caption{Kinetic function for a DG method with polynomial degree $p = 5$
            and filter order $s = 5$.}
    \label{fig:Cubic_DG_p_5_s_5__kinetic_functions}
  \end{figure}

  \item
  The kinetic functions depend on the polynomial degree $p$. For increasing
  $p$, the slope becomes steeper, i.e.\ smaller, since it is negative.
  This is shown for $s = 0$ and $N = \num{1024}$ in
  Figure~\ref{fig:Cubic_DG_s_0_N_01024__kinetic_functions}.
  \begin{figure}[htb]
  \centering
    \includegraphics[width=\textwidth]{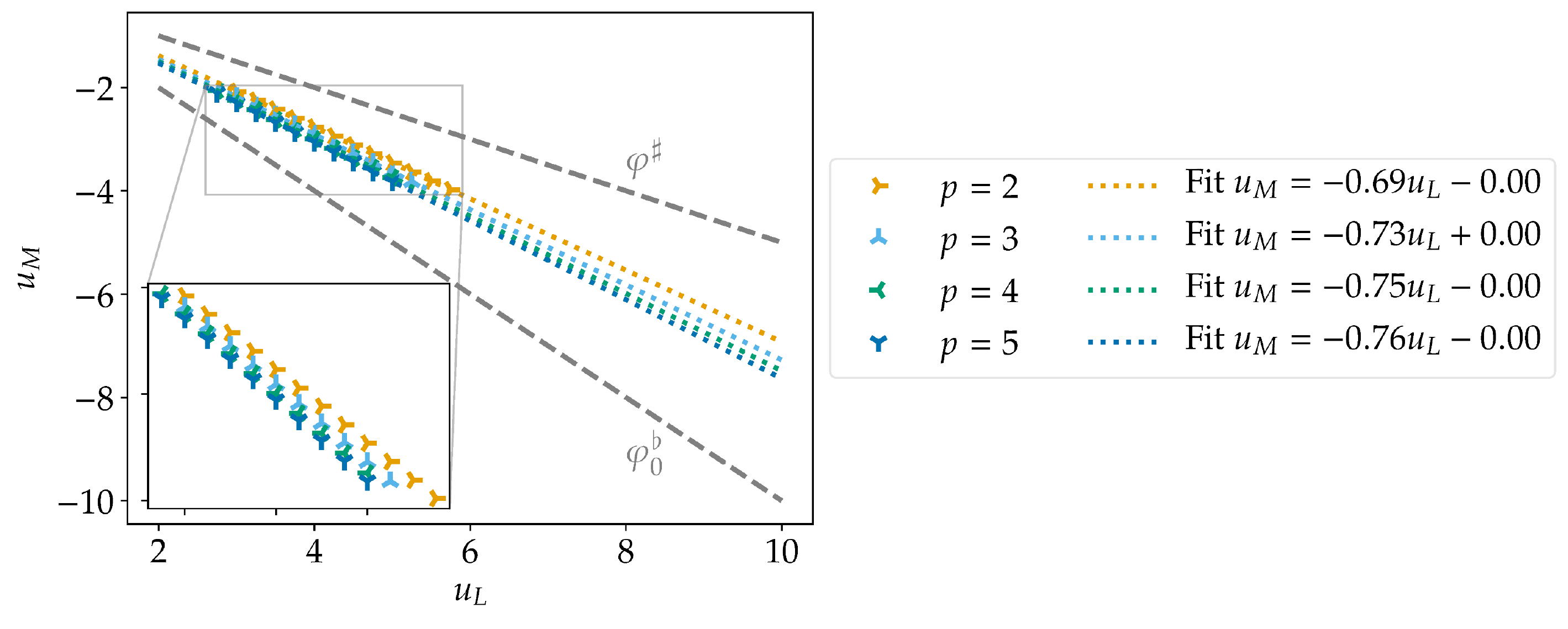}
    \caption{Kinetic function for DG methods without filtering ($s = 0$)
             and $N = \num{1024}$ elements.}
    \label{fig:Cubic_DG_s_0_N_01024__kinetic_functions}
  \end{figure}

  \item
  The kinetic functions depend on the filter order $s$.
  For $s = 5$, the slope is steeper than for $s < 5$. However, there is
  no clear relation for the other values of $s$. For $p = 5$, the slopes
  for $s = 0$ and $s = 4$ are visually indistinguishable while the slope
  for $s = 0$ is steeper than for $s = 4$ for $p = 4$.
  This is shown for $N = \num{1024}$ in
  Figure~\ref{fig:Cubic_DG_p_45_N_01024__kinetic_functions}.
  \begin{figure}[htb]
  \centering
    \begin{subfigure}{0.5\textwidth}
    \centering
      \includegraphics[width=\textwidth]{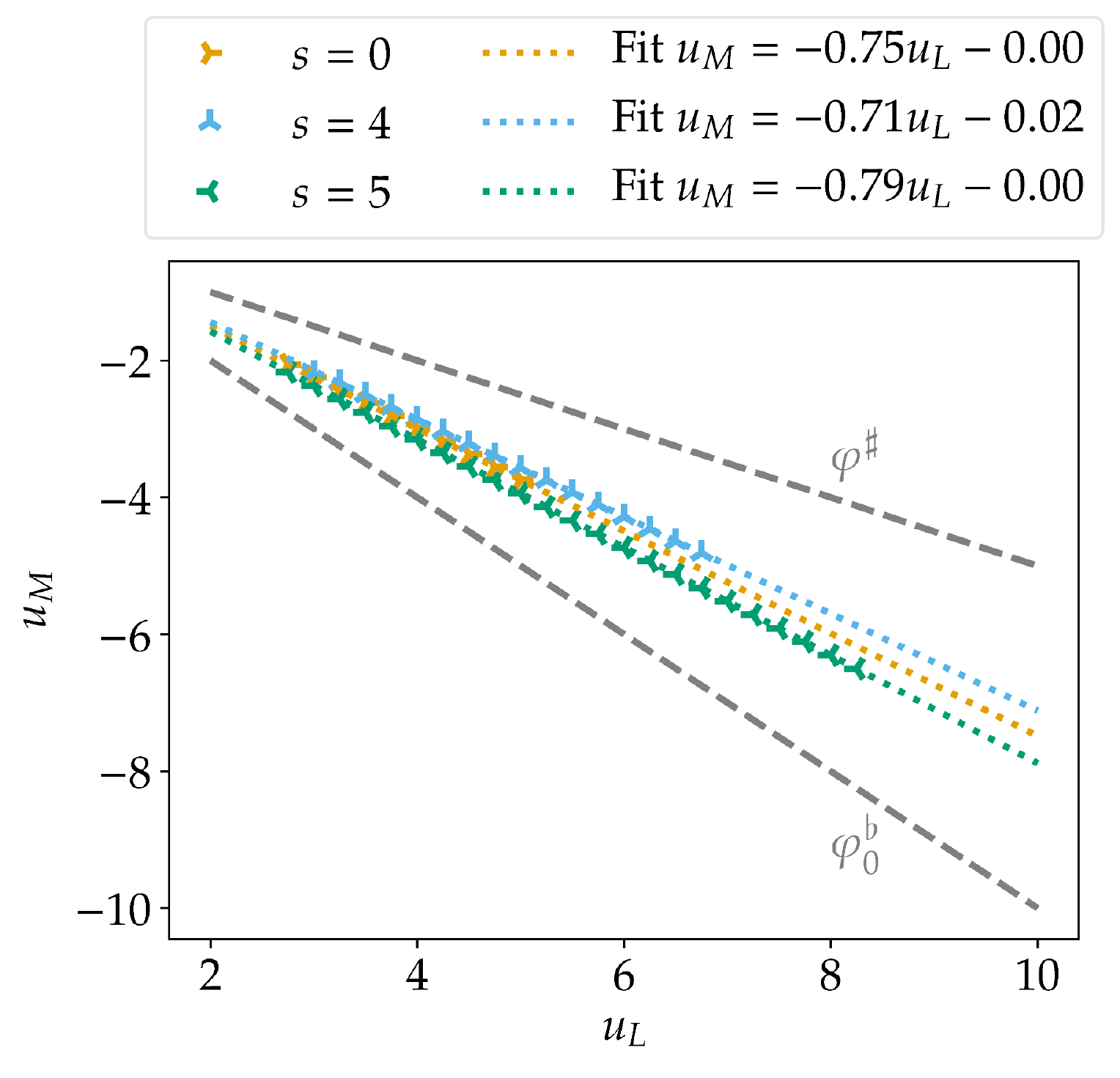}
      \caption{$p = 4$.}
    \end{subfigure}%
    \begin{subfigure}{0.5\textwidth}
    \centering
      \includegraphics[width=\textwidth]{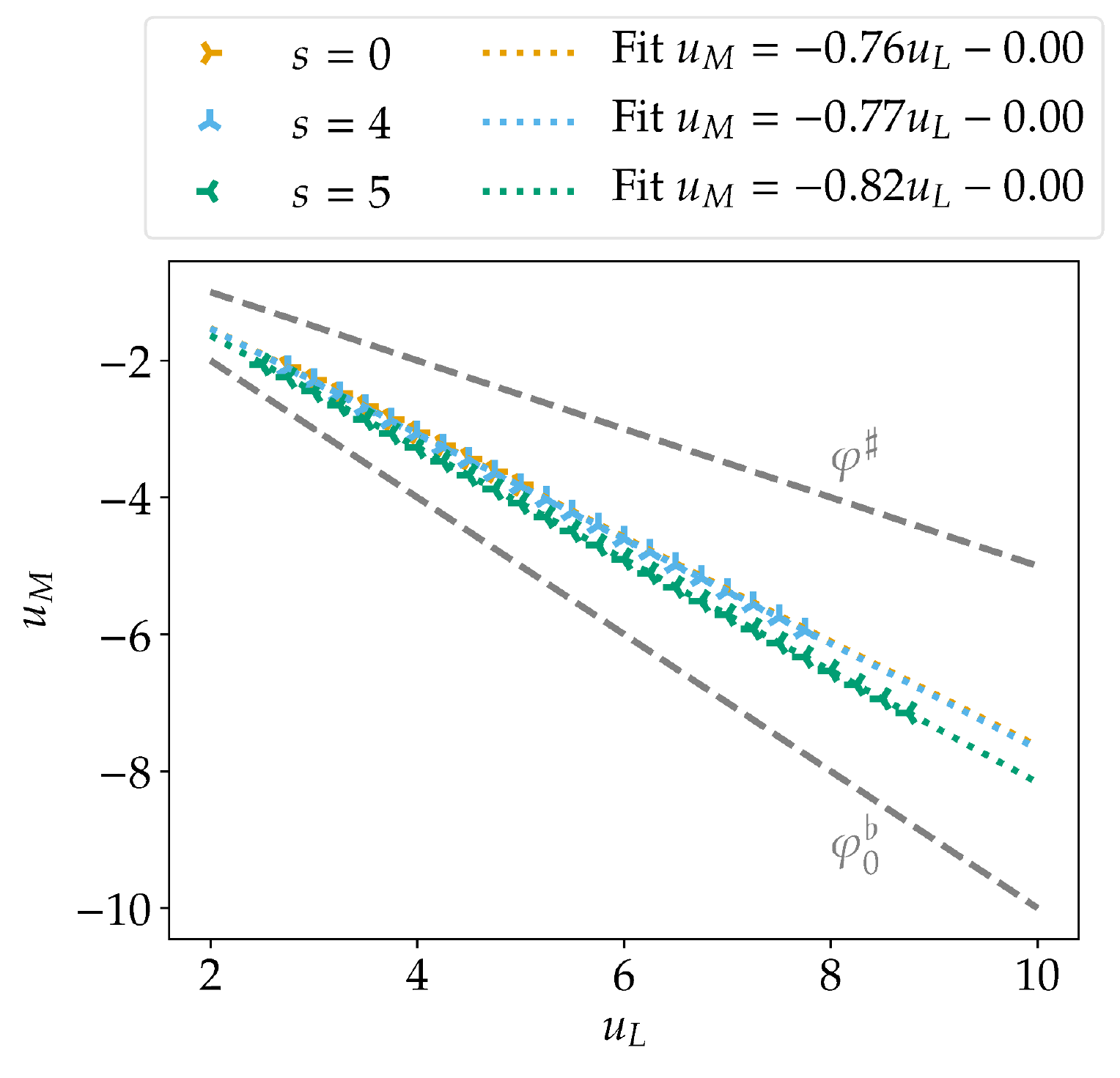}
      \caption{$p = 5$.}
    \end{subfigure}%
    \caption{Kinetic function for DG methods with $N = \num{1024}$ elements.}
    \label{fig:Cubic_DG_p_45_N_01024__kinetic_functions}
  \end{figure}

  \item
  The offset of the affine linear kinetic functions is nearly zero,
  i.e.\ the kinetic functions are approximately linear.

  \item
  The kinetic functions satisfy the bounds \eqref{eq:bounds-kinetic-function}.
\end{enumerate}


\subsection{Finite difference methods}
\label{sec:kinetic-function-cubic-FD}

For FD methods, accuracy order $p \in \{2, 4, 6\}$,
artificial dissipation strength $\epsilon_i \in \{0, 100, \dots, 400\}$,
and $N \in \{2^8, 2^9, \dots, 2^{14}\}$ grid nodes have been used.
The following observations have been made.
\begin{enumerate}
  \item
  If artificial dissipation was applied, nonclassical solutions occurred
  at least for some values of $N$, even if only the second-order artificial
  dissipation was used ($\epsilon_2 \neq 0, \epsilon_4 = \epsilon_6 = 0$).
  The only exception is given by the second-order method ($p = 2$) with
  second-order artificial dissipation ($\epsilon_2 \neq 0$), where no
  nonclassical solutions occurred.
  This is in agreement with the results of \cite{ranocha2019mimetic}:
  Discretizations of the second derivative can only be entropy-dissipative
  for all entropies if the order of accuracy is at most two. Additionally,
  the second-order discretizations applied here are dissipative for all
  entropies.
  Some examples are shown in Figure~\ref{fig:Cubic_FD_Strength4_0_Strength6_0_N_04096__kinetic_functions}.
  \begin{figure}[!htb]
  \centering
    \begin{subfigure}{0.5\textwidth}
    \centering
      \includegraphics[width=\textwidth]{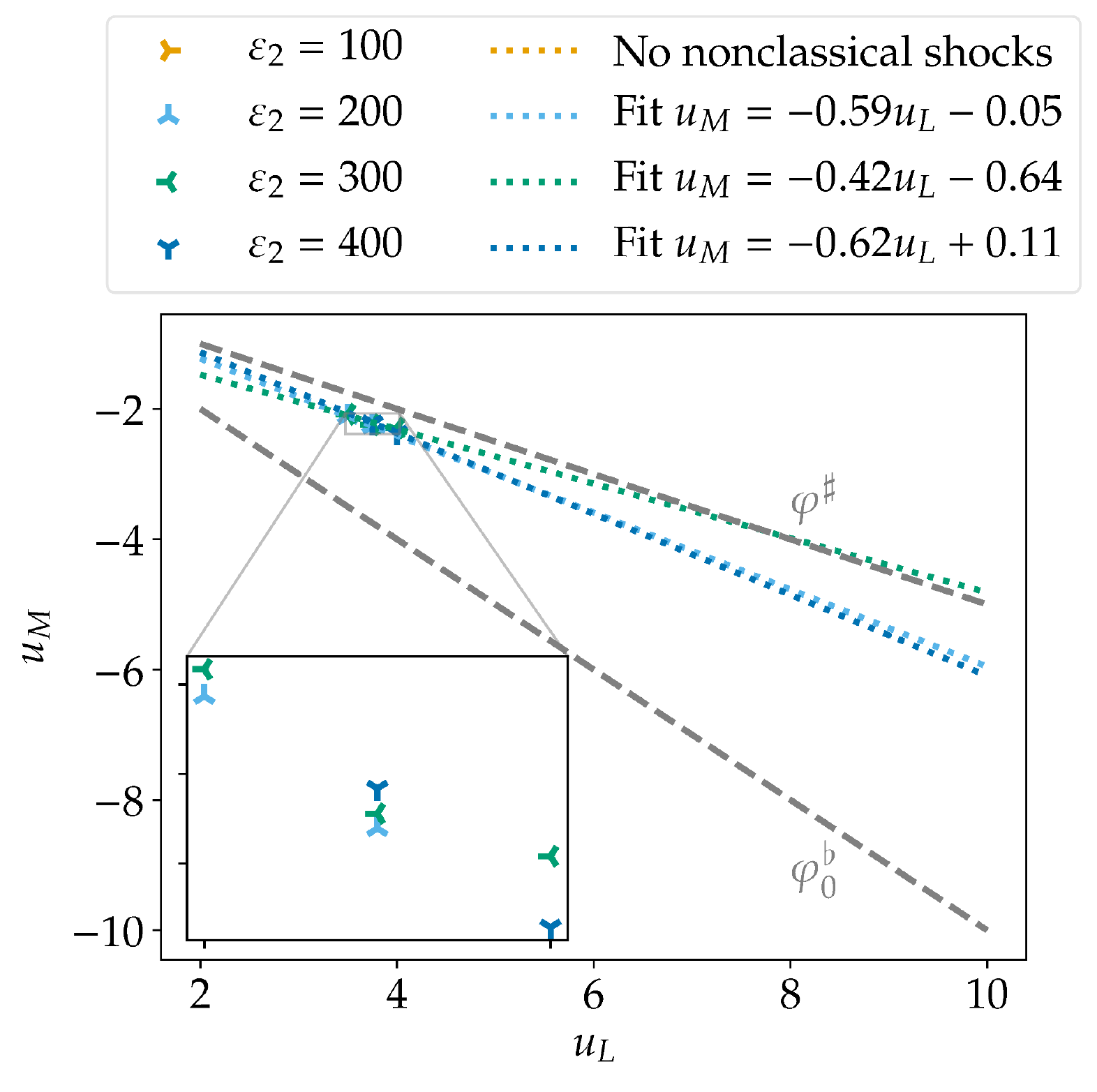}
      \caption{$p = 4$, $\epsilon_4 = 0$, $\epsilon_6 = 0$.}
    \end{subfigure}%
    \begin{subfigure}{0.5\textwidth}
    \centering
      \includegraphics[width=\textwidth]{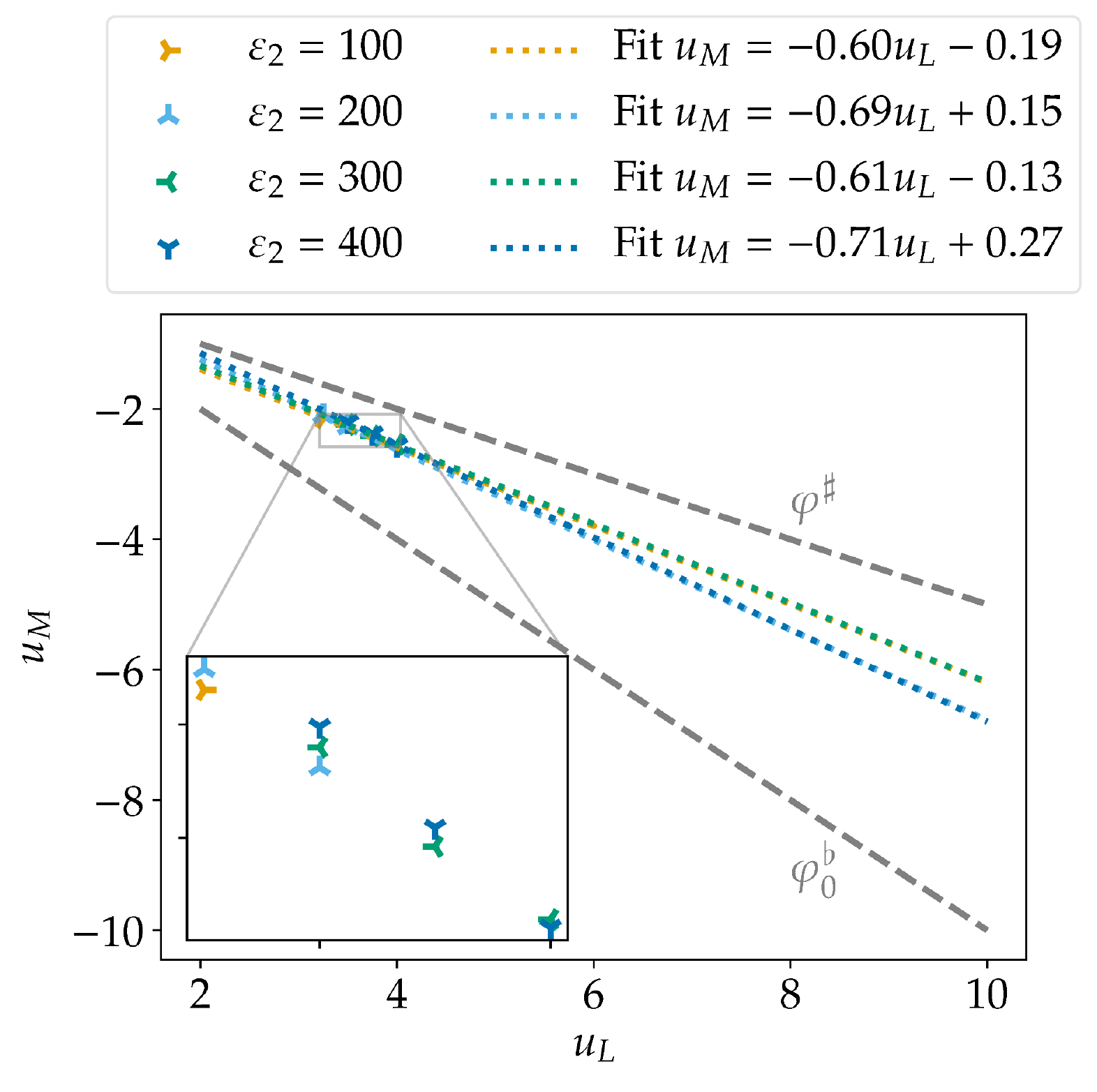}
      \caption{$p = 6$, $\epsilon_4 = 0$, $\epsilon_6 = 0$.}
    \end{subfigure}%
    \caption{Kinetic function for FD methods with order of accuracy $p \in \{4, 6\}$,
             $N = \num{4096}$ grid nodes, and different strengths of the artificial
             dissipation.}
    \label{fig:Cubic_FD_Strength4_0_Strength6_0_N_04096__kinetic_functions}
  \end{figure}

  \item
  As for DG methods, if the strength of the initial discontinuity is too big,
  no nonclassical solutions occur. For the investigated range of parameters,
  nonclassical solutions occurred only for $u_L < 10$ (as a necessary criterion).
  Depending   on the other parameters, the maximal value of $u_L$ for which
  nonclassical solution occurred can be smaller.

  \item
  The numerically obtained kinetic functions $\phi^\flat$ are affine linear.
  These affine linear functions vary slightly under grid refinement but seem
  to converge. However, the maximal value of $u_L$ leading to nonclassical
  solutions depends on $N$ and typically decreases when $N$ is increased.
  This is shown in Figure~\ref{fig:Cubic_FD_p_6_Strength2_0_Strength4_0400_Strength6_0400__kinetic_functions}.
  \begin{figure}[!htb]
  \centering
    \begin{subfigure}{\textwidth}
    \centering
      \includegraphics[width=\textwidth]{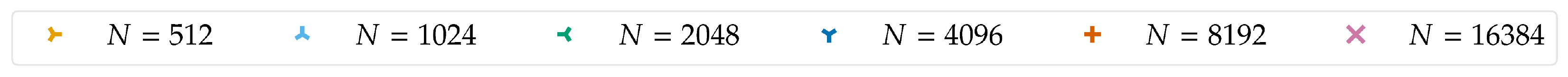}
    \end{subfigure}%
    \\
    \begin{subfigure}{0.5\textwidth}
    \centering
      \includegraphics[width=\textwidth]{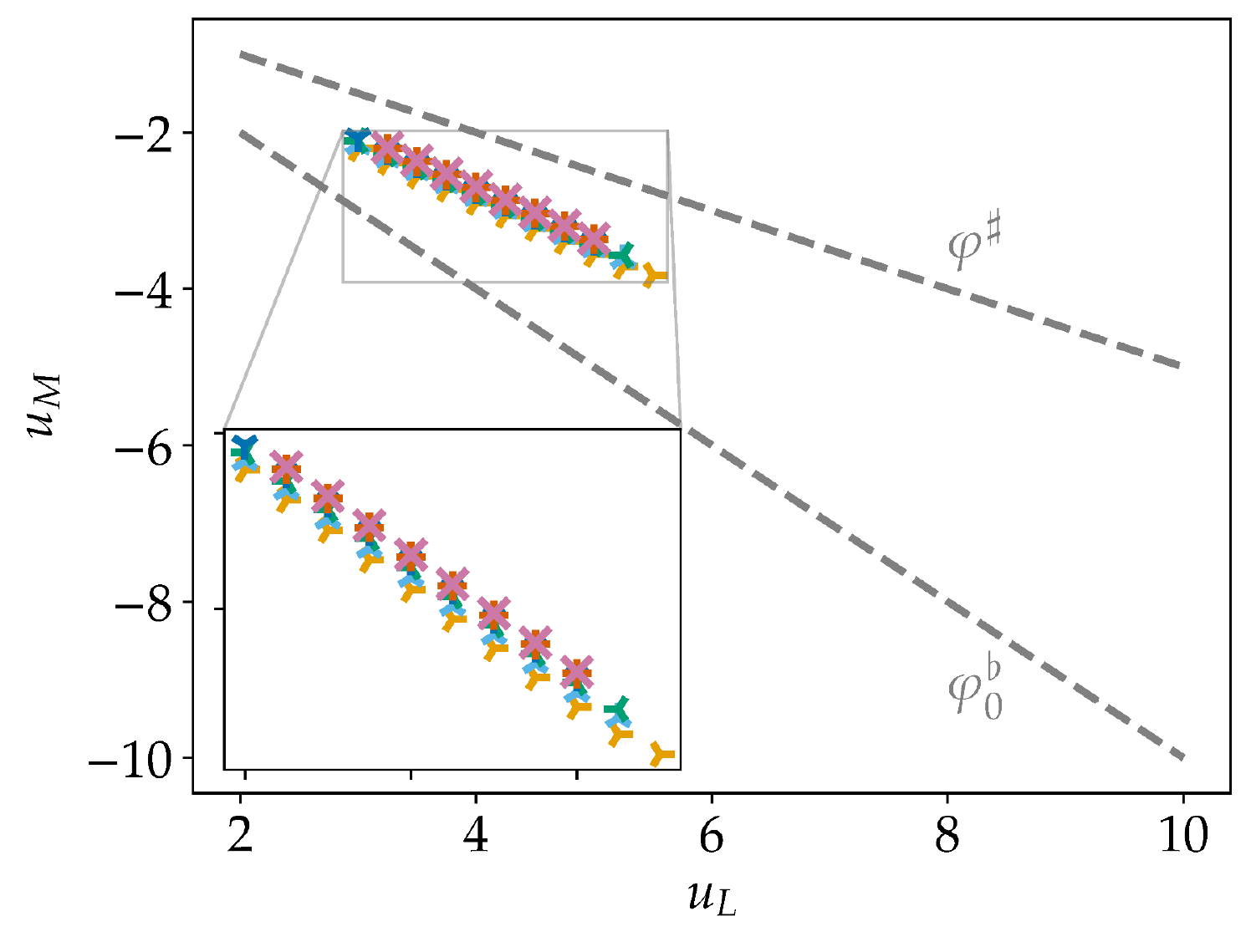}
      \caption{$\epsilon_2 = 0$, $\epsilon_4 = 400$, $\epsilon_6 = 0$.}
    \end{subfigure}%
    \begin{subfigure}{0.5\textwidth}
    \centering
      \includegraphics[width=\textwidth]{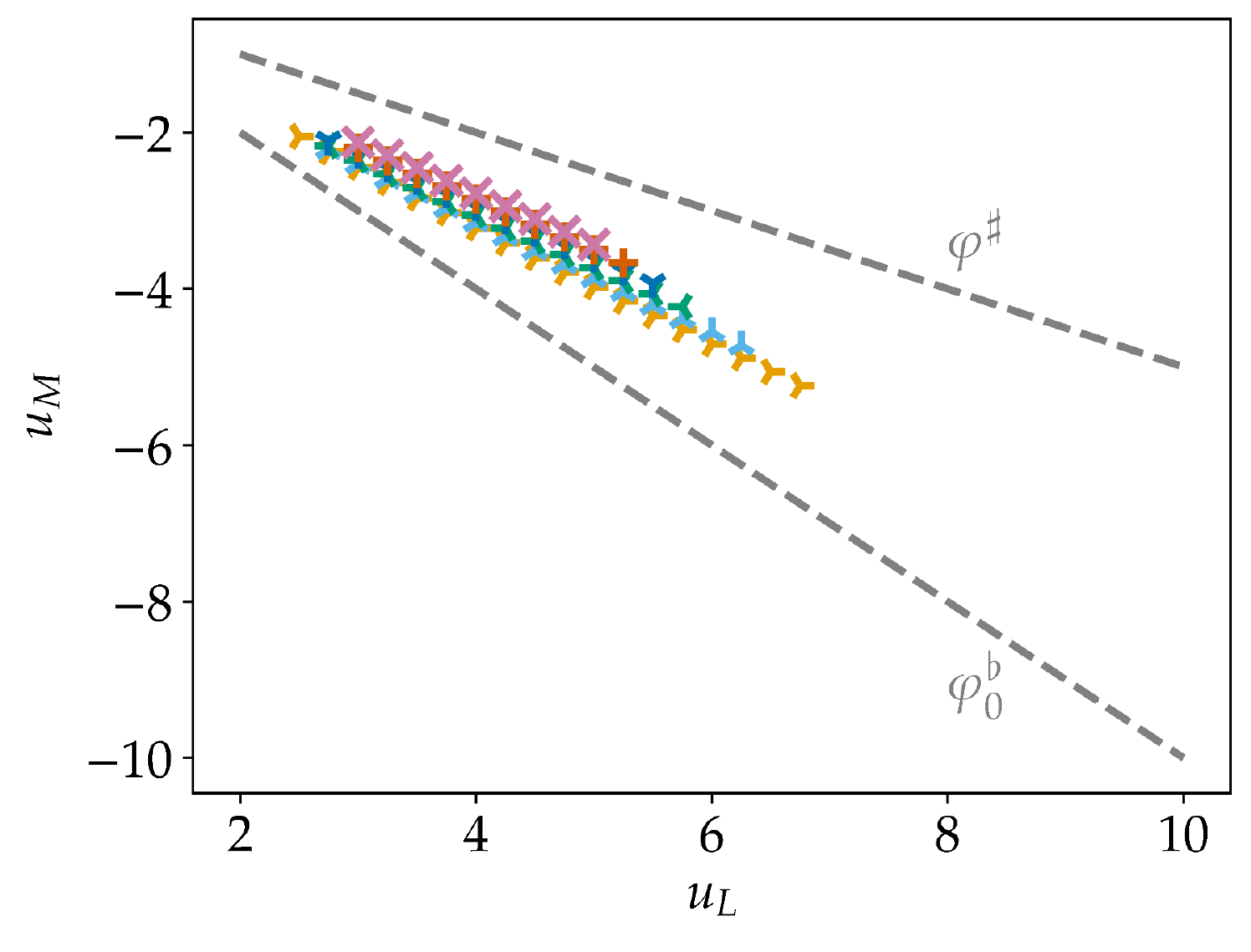}
      \caption{$\epsilon_2 = 0$, $\epsilon_4 = 0$, $\epsilon_6 = 400$.}
    \end{subfigure}%
    \caption{Kinetic function for FD methods with order of accuracy $p = 6$
             and varying number of grid nodes $N$ for different strengths
             of the artificial dissipation.}
    \label{fig:Cubic_FD_p_6_Strength2_0_Strength4_0400_Strength6_0400__kinetic_functions}
  \end{figure}

  \item
  Choosing a fixed order of the artificial dissipation, i.e.\ $\epsilon_i \neq 0$
  and $\epsilon_j = 0$ for $j \neq i$, the kinetic functions are nearly
  indistinguishable if the strength $\epsilon_i \in \{100, 200, 300, 400\}$
  is varied.
  This is shown in Figure~\ref{fig:Cubic_FD_p_6_Strength2_0_N_04096__kinetic_functions}.
  \begin{figure}[!htb]
  \centering
    \begin{subfigure}{0.5\textwidth}
    \centering
      \includegraphics[width=\textwidth]{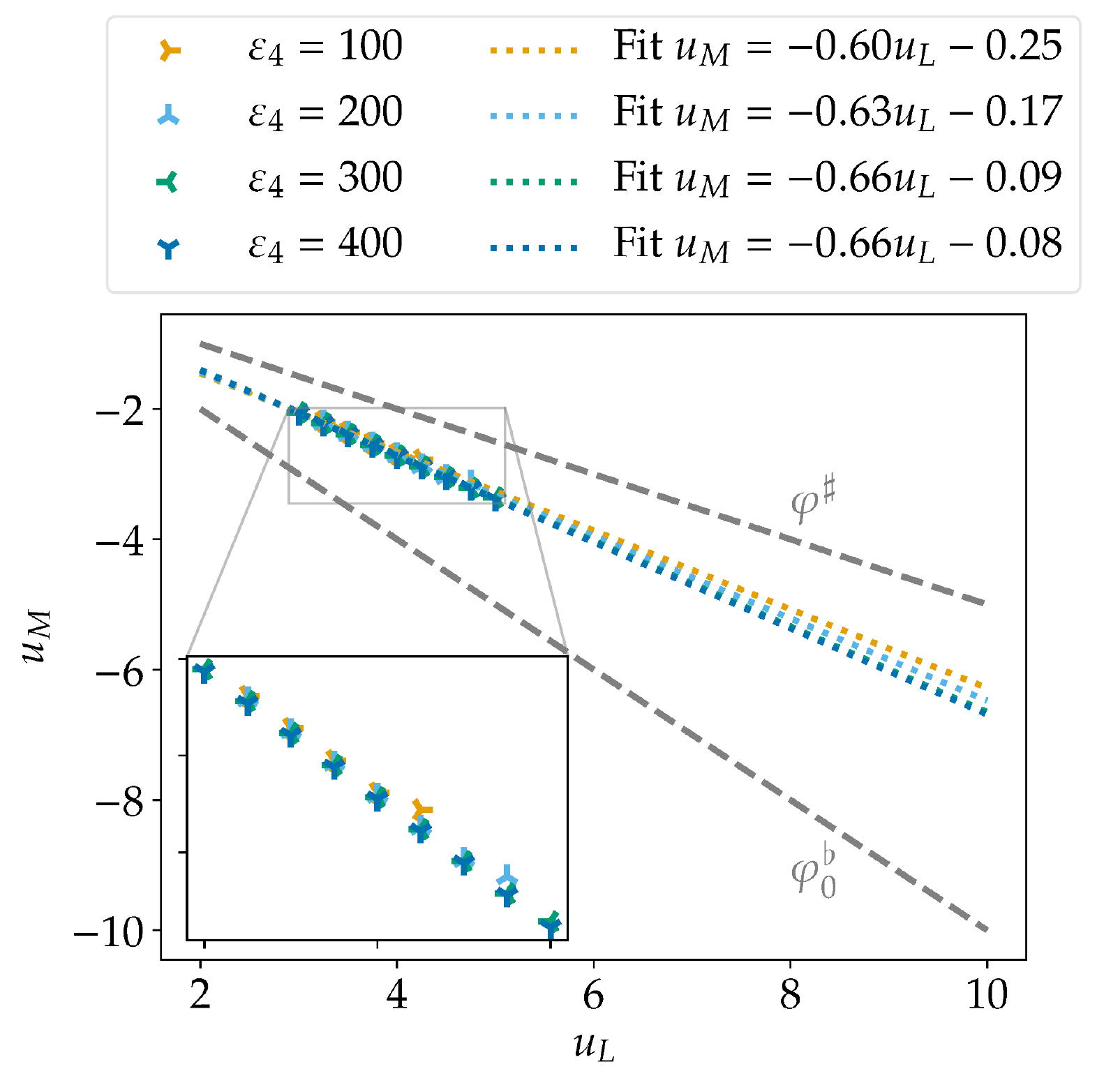}
      \caption{$\epsilon_2 = 0$, $\epsilon_6 = 0$.}
    \end{subfigure}%
    \begin{subfigure}{0.5\textwidth}
    \centering
      \includegraphics[width=\textwidth]{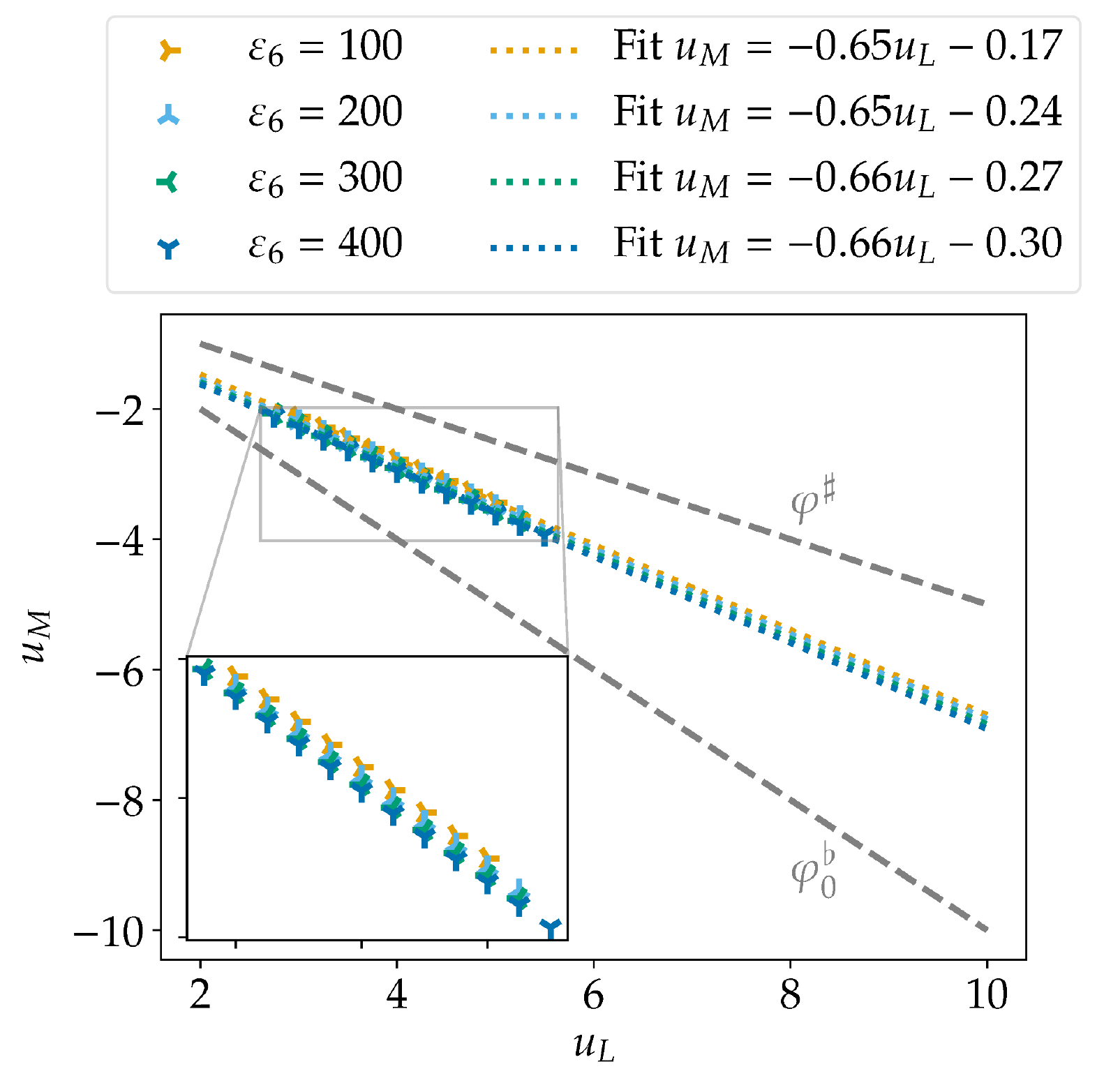}
      \caption{$\epsilon_2 = 0$, $\epsilon_4 = 0$.}
    \end{subfigure}%
    \caption{Kinetic function for FD methods with order of accuracy $p = 6$,
             $N = \num{4096}$ grid nodes, and different strengths of the artificial
             dissipation.}
    \label{fig:Cubic_FD_p_6_Strength2_0_N_04096__kinetic_functions}
  \end{figure}

  \item
  The kinetic function varies with the order of accuracy $p$. Typically, it
  becomes steeper (more negative) for higher values of $p$.
  This is shown in Figure~\ref{fig:Cubic_FD_Strength2_0_Strength4_0400_Strength6_0400_N_04096__kinetic_functions}.
  \begin{figure}[!htb]
  \centering
    \begin{subfigure}[b]{0.5\textwidth}
    \centering
      \includegraphics[width=\textwidth]{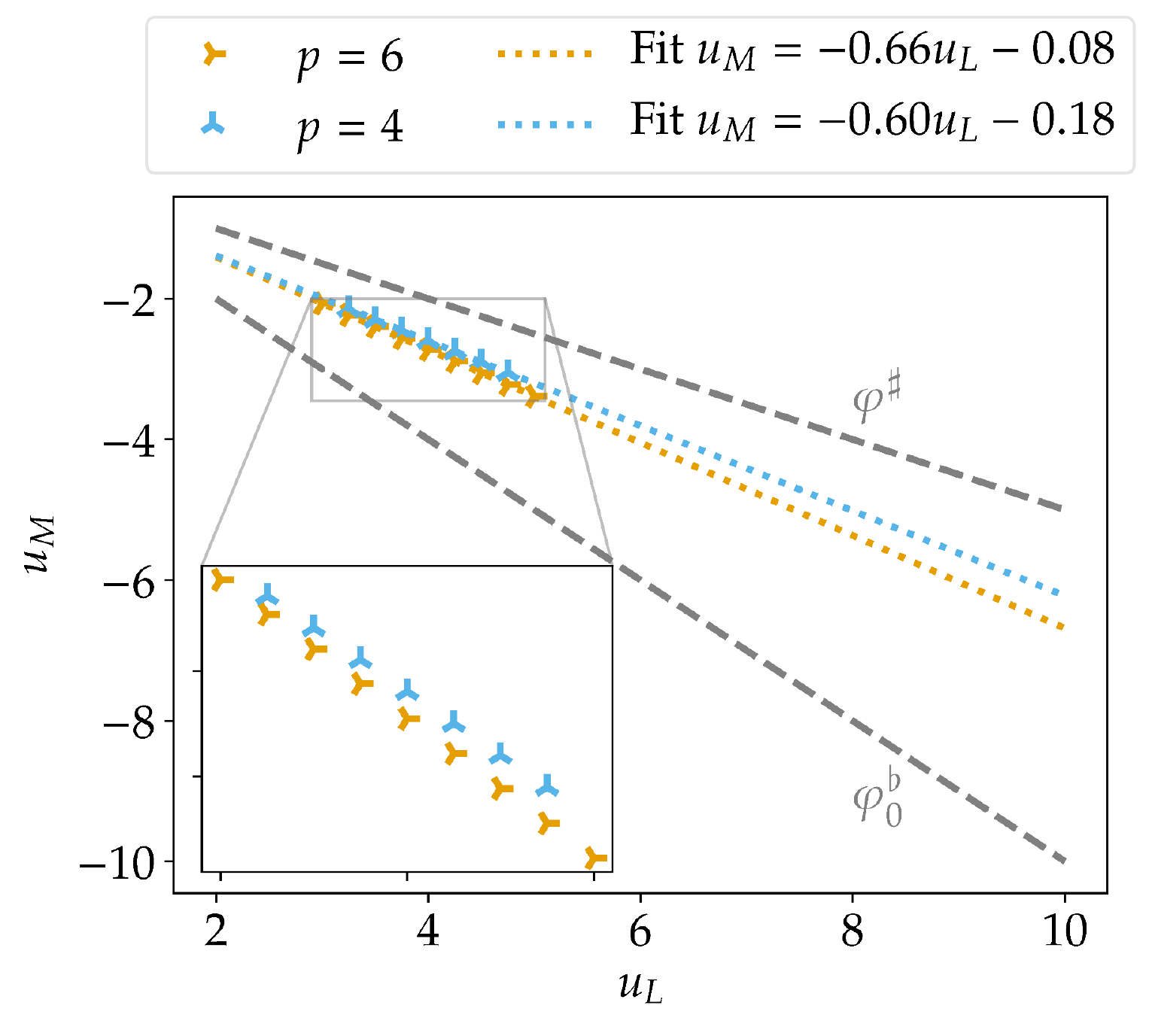}
      \caption{$\epsilon_2 = 0$, $\epsilon_4 = 400$, $\epsilon_6 = 0$.}
      \label{fig:Cubic_FD_Strength2_0_Strength4_400_Strength6_0_N_04096__kinetic_functions}
    \end{subfigure}%
    \begin{subfigure}[b]{0.5\textwidth}
    \centering
      \includegraphics[width=\textwidth]{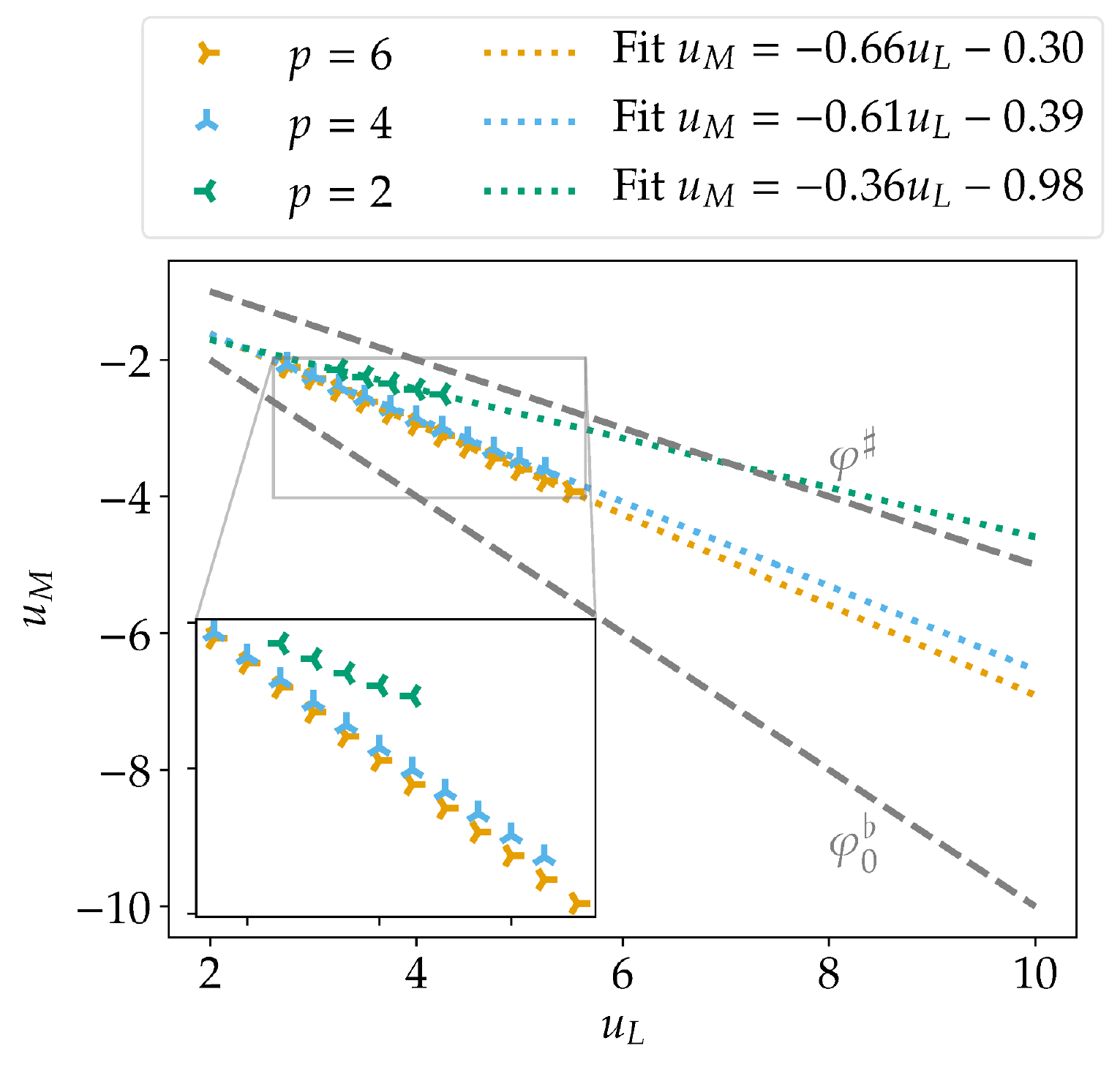}
      \caption{$\epsilon_2 = 0$, $\epsilon_4 = 0$, $\epsilon_6 = 400$.}
    \end{subfigure}%
    \caption{Kinetic function for FD methods with different orders of accuracy $p$,
             $N = \num{4096}$ grid nodes, and different strengths of the artificial
             dissipation. For $p = 2$ and $\epsilon_2 = 0$, $\epsilon_4 = 400$,
             $\epsilon_6 = 0$, no nonclassical solutions occur.}
    \label{fig:Cubic_FD_Strength2_0_Strength4_0400_Strength6_0400_N_04096__kinetic_functions}
  \end{figure}

  \item
  In contrast to DG methods, the offset of the affine linear kinetic functions
  is in general not zero.

  \item
  The kinetic functions satisfy the bounds \eqref{eq:bounds-kinetic-function}
  in the region where nonclassical solutions occur. For $p = 2$, the kinetic
  functions can violate the bounds \eqref{eq:bounds-kinetic-function} if they
  are extrapolated to bigger values of $u_L$, cf.\ Figure~\ref{fig:Cubic_FD_Strength2_0_Strength4_400_Strength6_0_N_04096__kinetic_functions}.
  For $p \in \{4, 6\}$ such a behavior did not occur.
\end{enumerate}


\subsection{Fourier collocation methods}
\label{sec:kinetic-function-cubic-Fourier}

For Fourier methods, the viscosity strengths $\epsilon \in \{10 / N,
50 / N, 100 / N\}$ for the standard and convergent choices of
\cite{tadmor2012adaptive} and $N \in \{2^{10}, 2^{11}, \dots, 2^{14}\}$
grid nodes have been used.
The following observations have been made.
\begin{enumerate}
  \item
  As for DG and FD methods, if the strength of the initial discontinuity
  is too big, no nonclassical solutions occurred. For the investigated range
  of parameters, nonclassical solutions occurred only for $u_L < 10$
  (as a necessary criterion).
  Depending   on the other parameters, the maximal value of $u_L$ for which
  nonclassical solution occurred can be smaller.

  \item
  The numerically obtained kinetic functions $\phi^\flat$ are affine linear
  and seem to converge under grid refinement (or are already visually
  indistinguishable). This can be seen in Figure~\ref{fig:Cubic_Fourier_grid_refinement__kinetic_function}.
  For lower resolution such as $N = \num{1024}$ and strength $\epsilon = 10 / N$,
  the nonclassical part is highly oscillatory, resulting in some errors of
  the measurement of the middle state. For higher resolutions, these problems
  are less severe or disappear completely.
  \begin{figure}[!htb]
  \centering
    \begin{subfigure}{\textwidth}
    \centering
      \includegraphics[width=\textwidth]{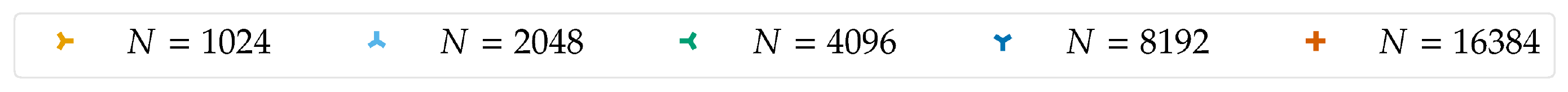}
    \end{subfigure}%
    \\
    \begin{subfigure}[b]{0.5\textwidth}
    \centering
      \includegraphics[width=\textwidth]{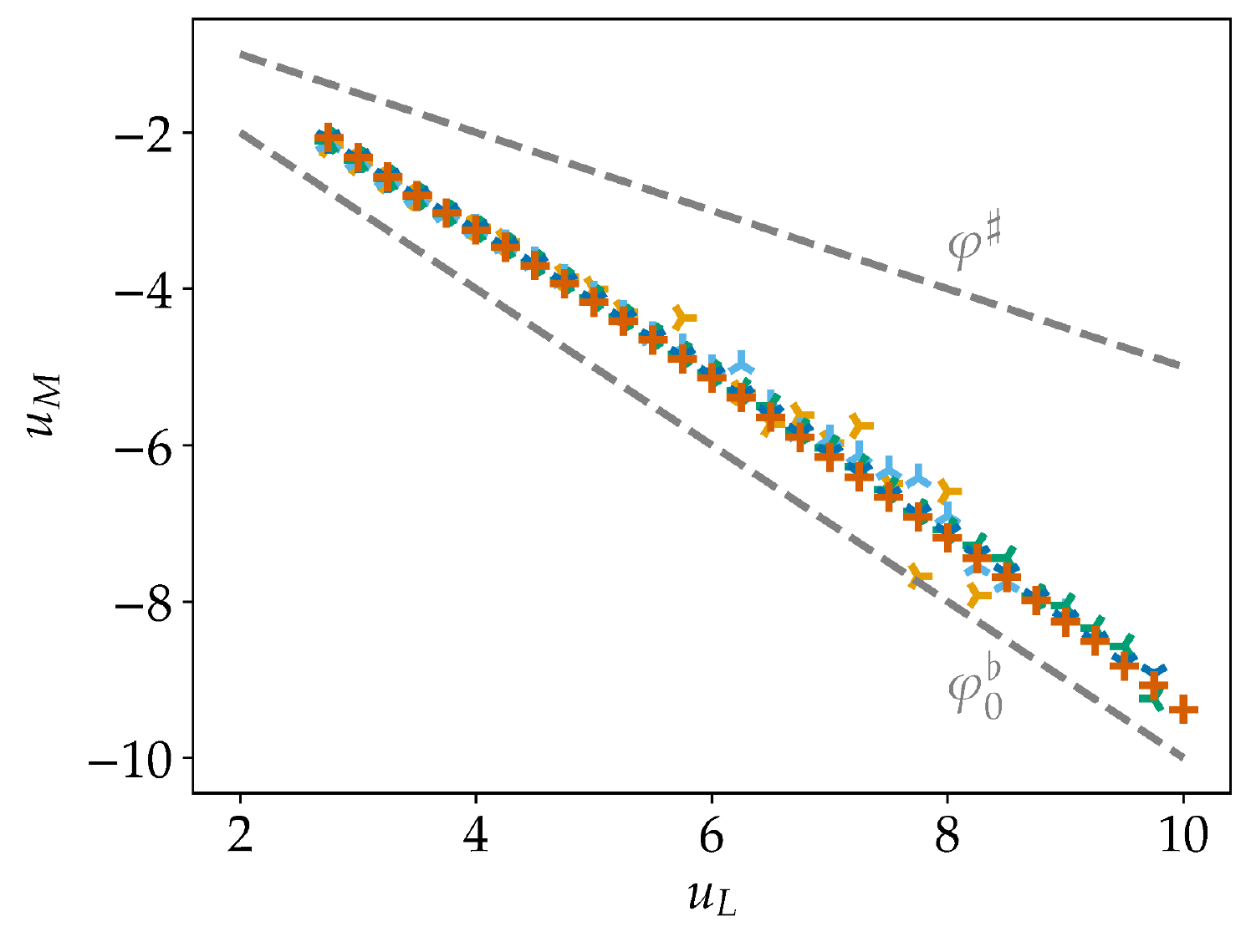}
      \caption{Standard choice of \cite{tadmor2012adaptive}, $\epsilon = 10 / N$.}
    \end{subfigure}%
    \begin{subfigure}[b]{0.5\textwidth}
    \centering
      \includegraphics[width=\textwidth]{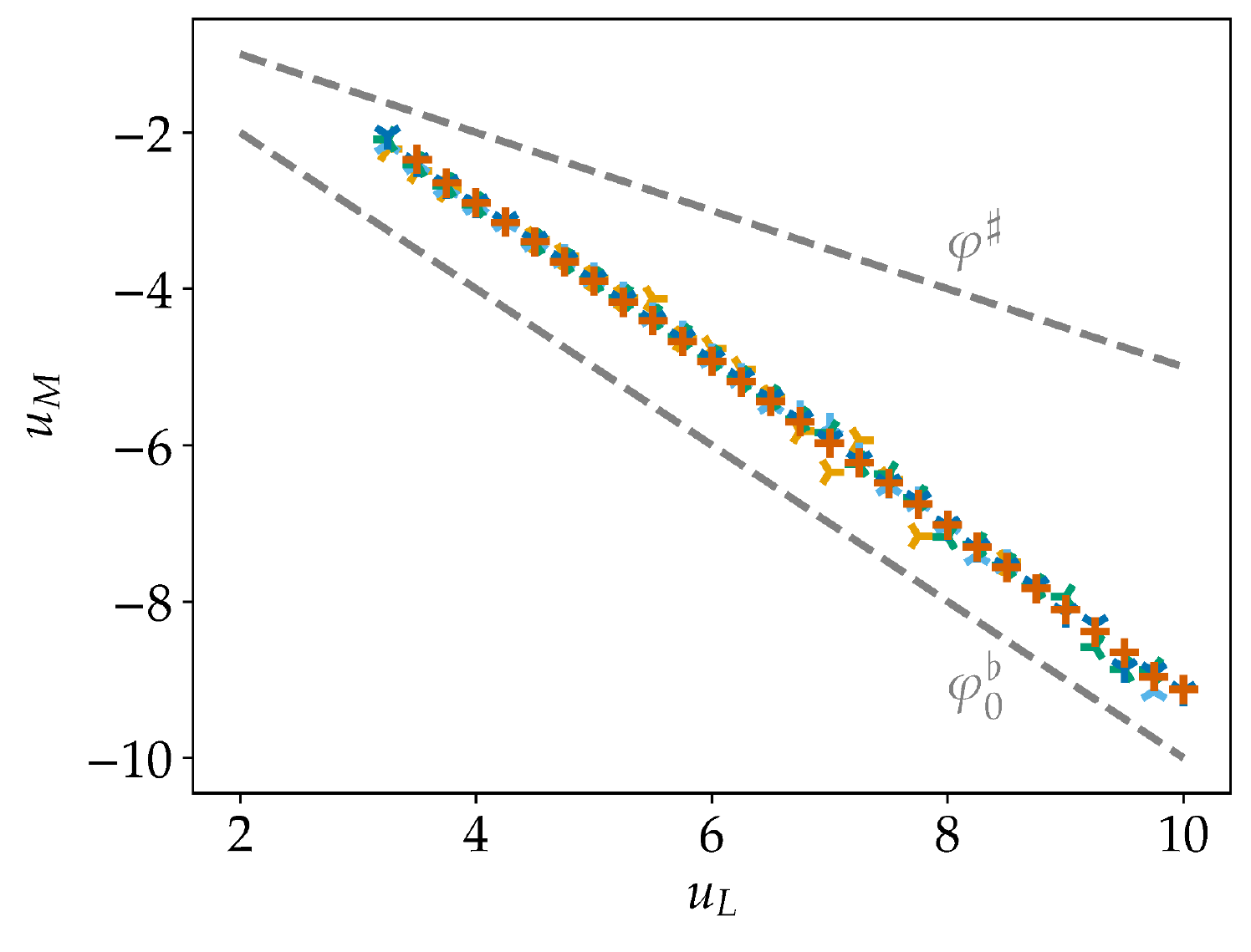}
      \caption{Convergent choice of \cite{tadmor2012adaptive}, $\epsilon = 10 / N$.}
    \end{subfigure}%
    \caption{Kinetic function for Fourier methods with different choices of
             the spectral viscosity and numbers $N$ of grid nodes.}
    \label{fig:Cubic_Fourier_grid_refinement__kinetic_function}
  \end{figure}

  \item
  The kinetic functions depend on the strength of the spectral viscosity.
  For the convergent choice of \cite{tadmor2012adaptive}, nonclassical shocks
  can occur or not, depending on the strength $\epsilon$.
  The strength $\epsilon = 100 / N$ results in steeper kinetic function than
  the choice $\epsilon = 10 / N$. For $\epsilon = 50 / N$ and the convergent
  choice of \cite{tadmor2012adaptive}, no nonclassical solutions appear,
  in contrast to the classical choice of the spectral viscosity which always
  results in nonclassical solutions. However, for smaller or bigger strengths
  of the spectral viscosity, nonclassical solutions occur even for the convergent
  choice of \cite{tadmor2012adaptive} and $N = \num{16384}$ grid nodes.
  This can be seen in Figure~\ref{fig:Cubic_Fourier_N_16384_TadmorWaagan2012_Strength_100_Split_1__kinetic_function}.
  \begin{figure}[!htb]
  \centering
    \begin{subfigure}[b]{0.5\textwidth}
    \centering
      \includegraphics[width=\textwidth]{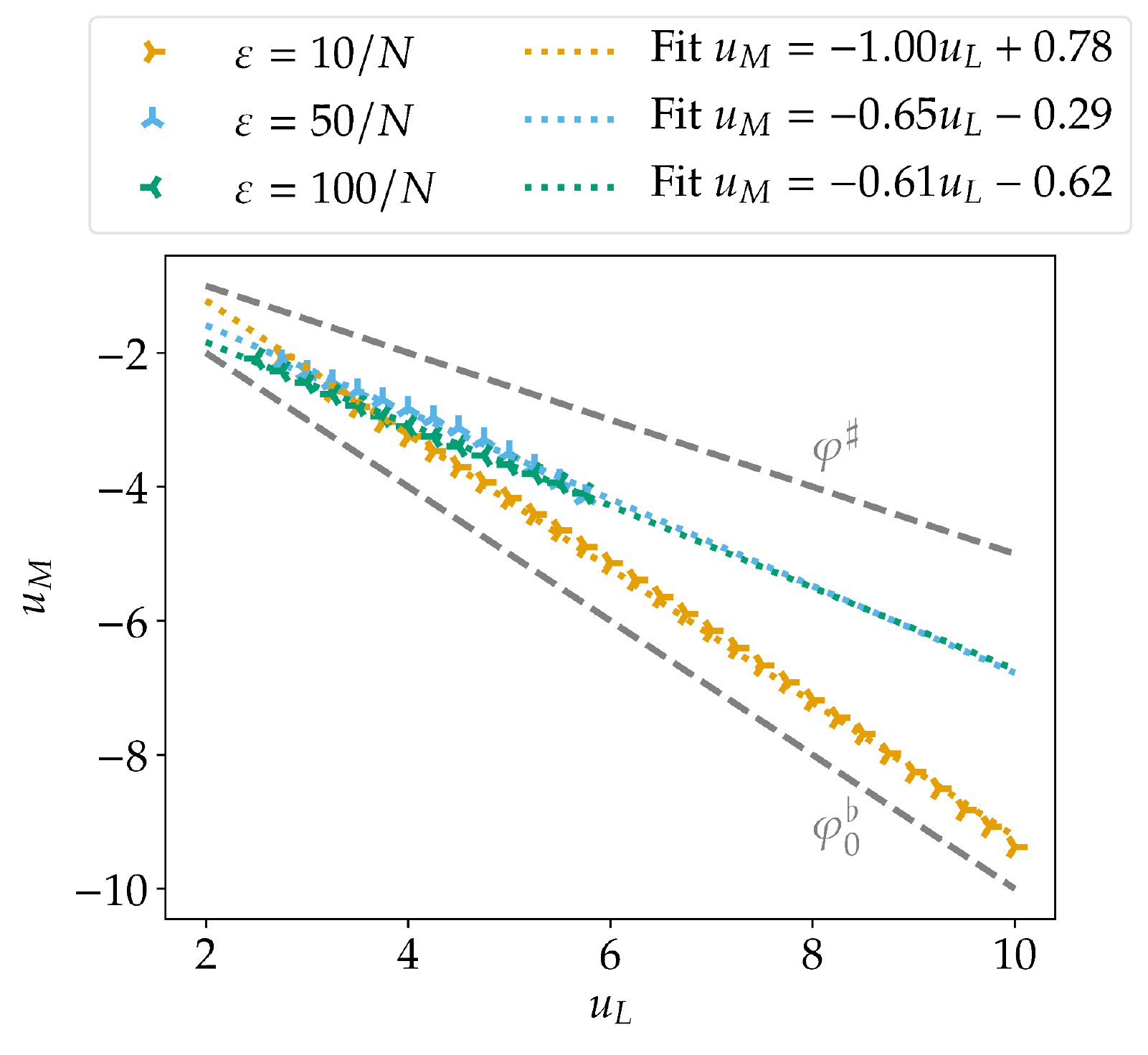}
      \caption{Standard choice of \cite{tadmor2012adaptive}.}
    \end{subfigure}%
    \begin{subfigure}[b]{0.5\textwidth}
    \centering
      \includegraphics[width=\textwidth]{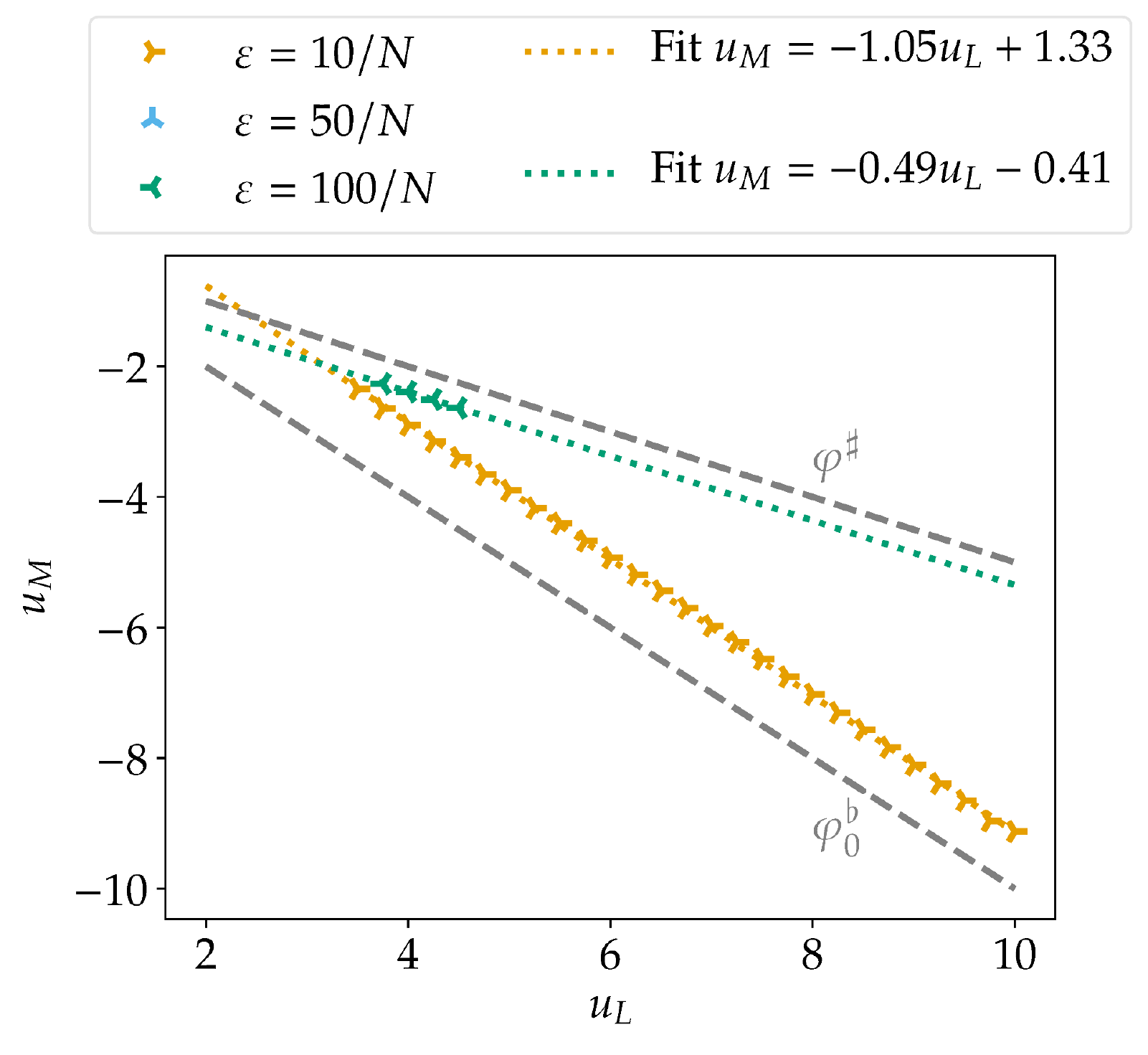}
      \caption{Convergent choice of \cite{tadmor2012adaptive}.}
    \end{subfigure}%
    \caption{Kinetic function for Fourier methods with different choices of
             the spectral viscosity and $N = \num{16384}$ grid nodes.}
    \label{fig:Cubic_Fourier_N_16384_TadmorWaagan2012_Strength_100_Split_1__kinetic_function}
  \end{figure}

  \item
  In contrast to DG methods but similar to FD methods, the offset of the
  affine linear kinetic functions is in general not zero.

  \item
  The kinetic functions satisfy the bounds \eqref{eq:bounds-kinetic-function}.
\end{enumerate}


\subsection{WENO methods}
\label{sec:kinetic-function-cubic-WENO}

In addition to the provably entropy-stable methods tested above, some
standard shock capturing have been studied. Specifically, high-order WENO
methods implemented in Clawpack
\cite{clawpack561,mandli2016clawpack,ketcheson2013high,ketcheson2012pyclaw}
have been investigated.
The discretizations used WENO orders $p \in \{13, 15, 17\}$ and
$N \in \{2^8, 2^9, \dots, 2^{14} \}$ cells with a CFL number of $0.5$.
The other parameters of the problem are the same as for the SBP FD methods.
The following observations have been made.
\begin{enumerate}
  \item
  Nonclassical solutions have only been observed for very high-order WENO
  methods with $p \in \{13, 15, 17\}$. For $p \leq 11$, no nonclassical
  shocks occurred.

  \item
  As for DG and FD methods, if the strength of the initial discontinuity
  is too big, no nonclassical solutions occurred and the critical strength
  of the initial discontinuity depends on the WENO order $p$.

  \item
  The numerically obtained kinetic functions $\phi^\flat$ are affine linear
  and seem to converge under grid refinement (or are already visually
  indistinguishable). This can be seen in Figure~\ref{fig:Cubic_WENO_p_17__kinetic_function}.

  \item
  The kinetic functions do not seem to depend strongly on the WENO order $p$
  when nonclassical solutions occur. For higher-order WENO methods, nonclassical
  shocks occurred also for bigger values of $u_L$.
  This can be seen in Figure~\ref{fig:Cubic_WENO_N_16384__kinetic_function}.
  \begin{figure}[!htb]
  \centering
    \begin{subfigure}[b]{0.5\textwidth}
    \centering
      \includegraphics[width=\textwidth]{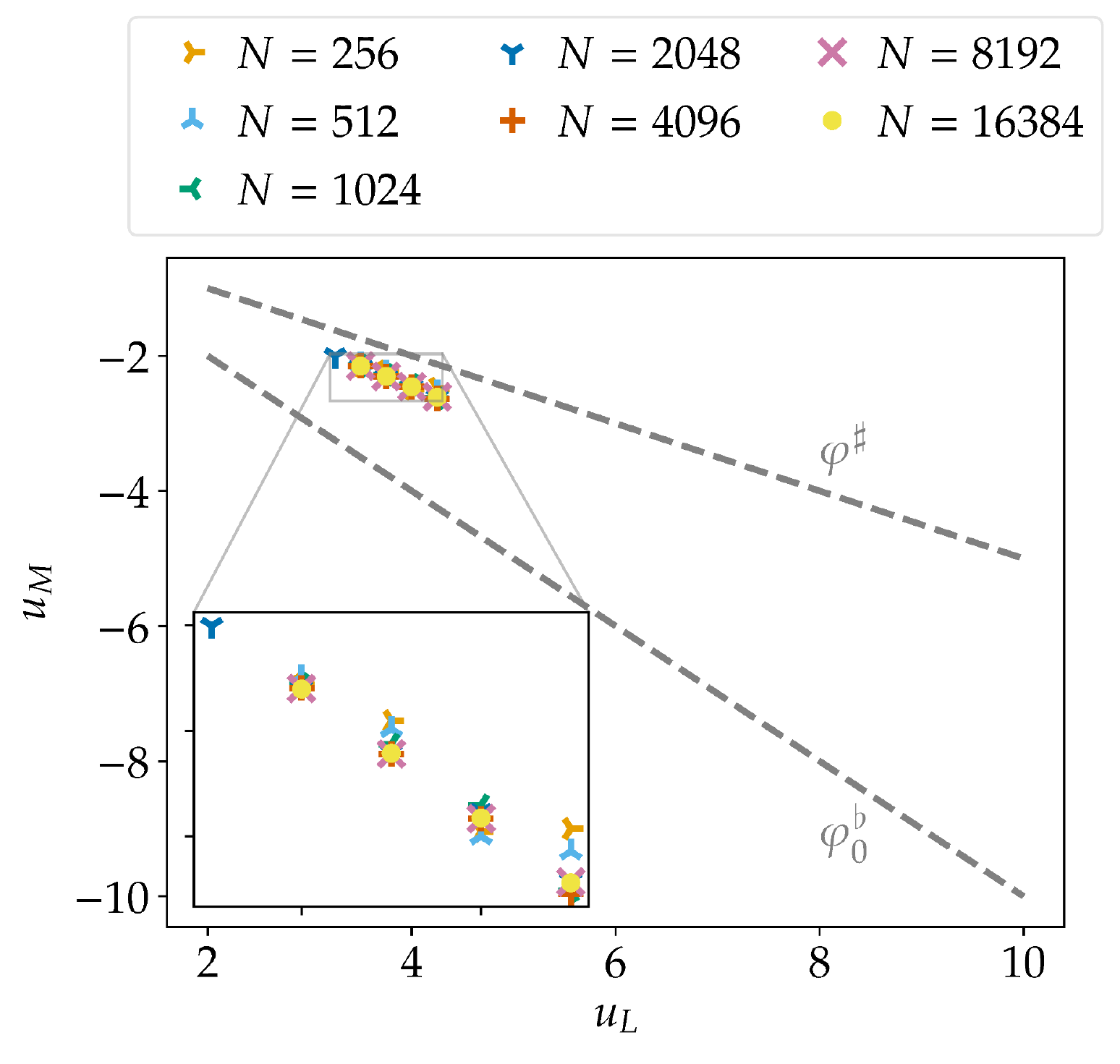}
      \caption{Grid refinement study for $p = 17$.}
      \label{fig:Cubic_WENO_p_17__kinetic_function}
    \end{subfigure}%
    \begin{subfigure}[b]{0.5\textwidth}
    \centering
      \includegraphics[width=\textwidth]{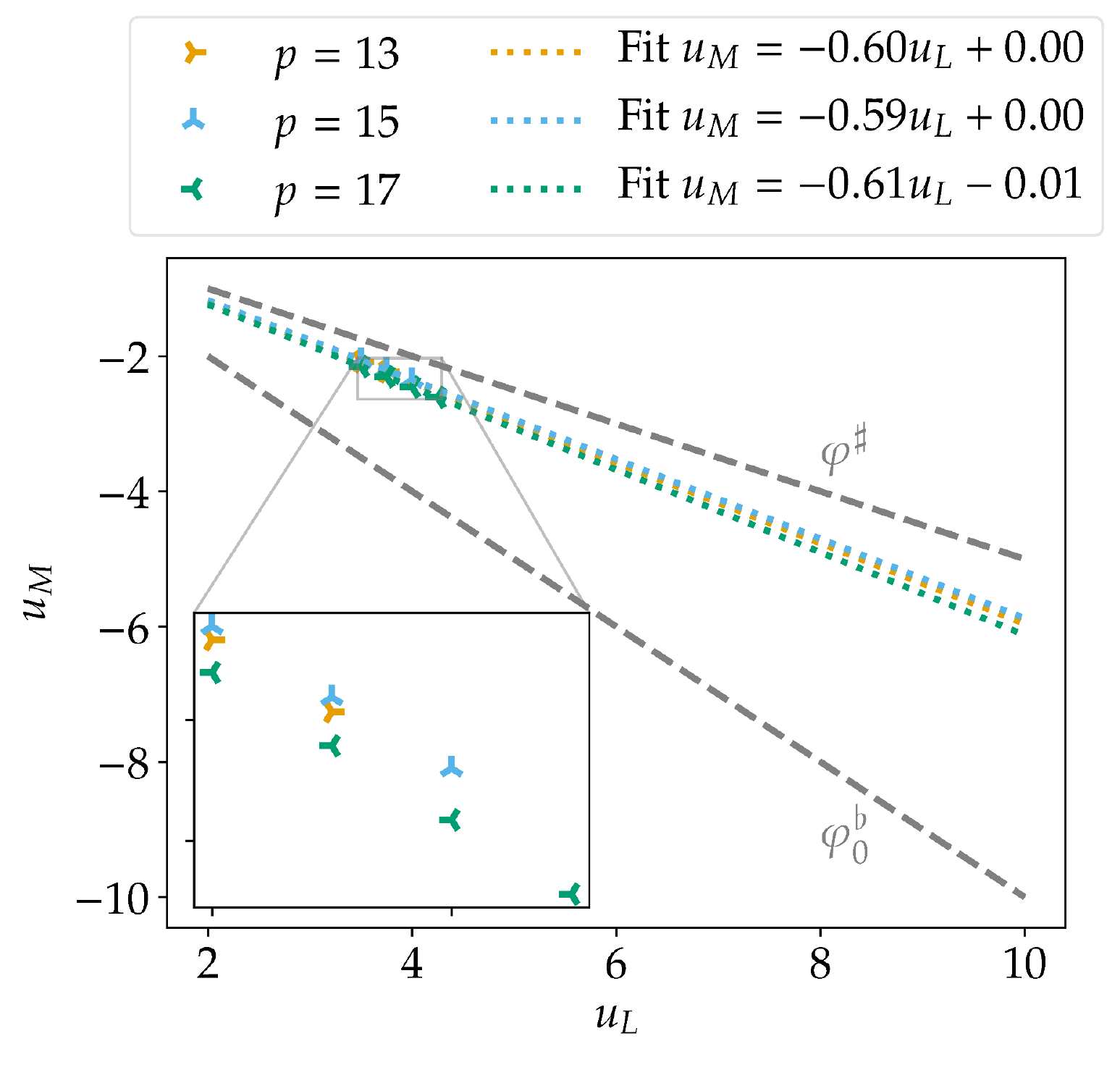}
      \caption{Different order $p$ for $N = \num{16384}$.}
      \label{fig:Cubic_WENO_N_16384__kinetic_function}
    \end{subfigure}%
    \caption{Kinetic functions of WENO methods.}
    \label{fig:Cubic_WENO__kinetic_function}
  \end{figure}

  \item
  If few cells are used, the offset of the affine linear kinetic functions
  is not necessarily zero, similar to FD methods. For increased grid
  resolutions, the offset becomes (nearly) zero, similar to DG methods.

  \item
  The kinetic functions satisfy the bounds \eqref{eq:bounds-kinetic-function}.
\end{enumerate}


\section{Generalization to a quartic conservation law}
\label{sec:kinetic-function-quartic}

\subsection{Preliminaries}

\begin{figure}[!htb]
\centering
  \includegraphics[width=.4\textwidth]{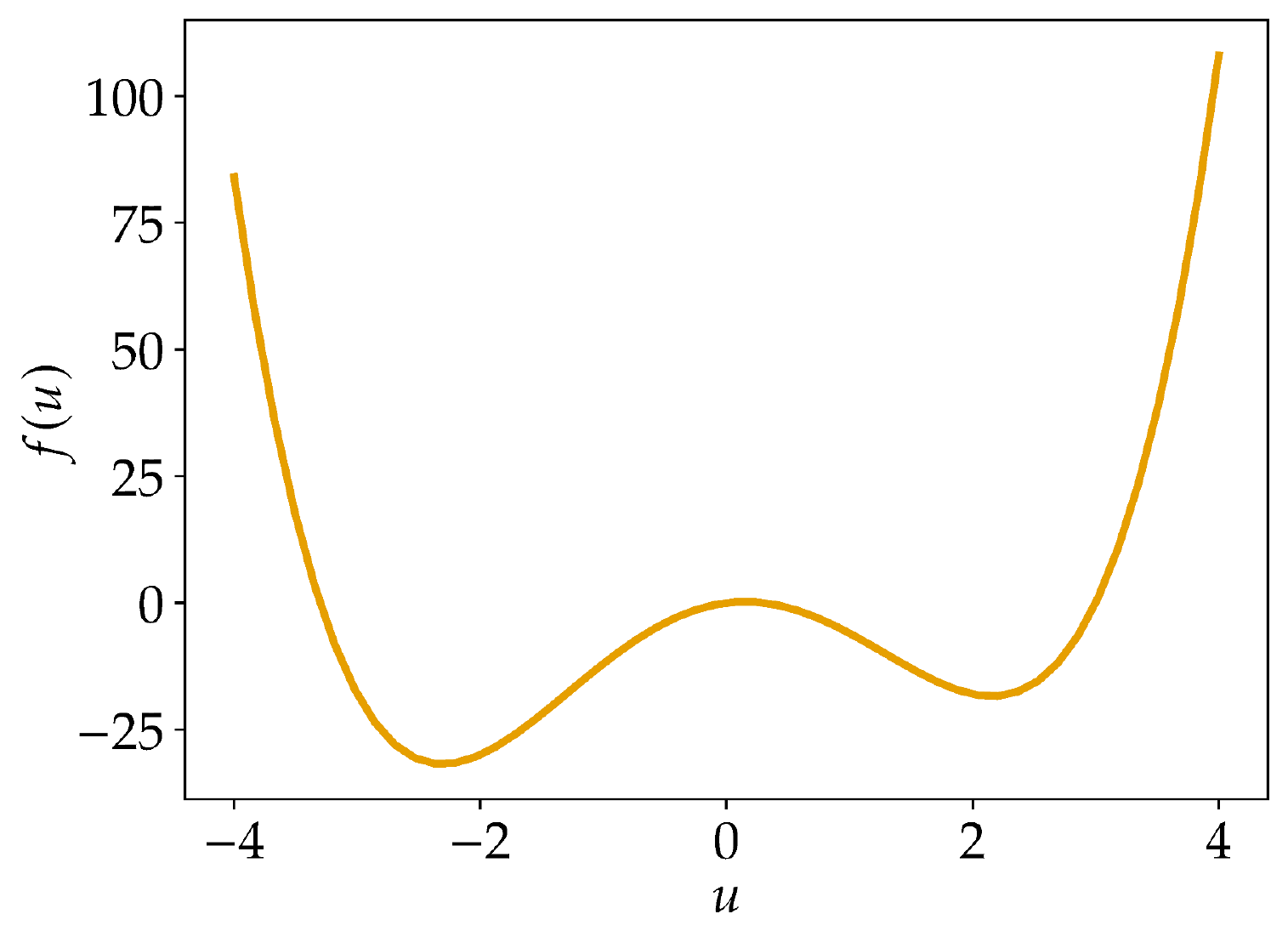}
  \caption{Non-convex flux function of the quartic conservation law
           \eqref{eq:quartic}.}
  \label{fig:quartic-flux}
\end{figure}

Entropy-conservative semi-discretizations of the scalar conservation law
\begin{equation}
\label{eq:quartic}
\begin{aligned}
  \partial_t u(t,x) + \partial_x f(u(t,x)) &= 0,
    \quad && t \in (0, T),\, x \in (\xmin,\xmax),
    &&&
  \\
  u(0,x) &= u_0(x),
    \quad && x \in (\xmin,\xmax),
  \\
  f(u) &= u^2 (u^2 - 10) + 3 u,
\end{aligned}
\end{equation}
can be constructed using the numerical flux
\begin{equation}
\label{eq:quartic-fnum-EC}
  \fnum(u_-, u_+)
  =
  \frac{u_+^4 + u_- u_+^3 + u_-^2 u_+^2 + u_-^3 u_+ + u_-^4}{5}
  - 10 \frac{u_+^2 + u_- u_+ + u_-^2}{3}
  + 3 \frac{u_+ + u_-}{2}
\end{equation}
with corresponding split form
\begin{equation}
  \frac{2}{5} \left(
    \mat{D} \vec{u}^4
    + \mat{u} \mat{D} \vec{u}^3
    + \mat{u}^2 \mat{D} \vec{u}^2
    + \mat{u}^3 \mat{D} \vec{u}
  \right)
  - \frac{20}{3} \left(
    \mat{D} \vec{u}^2
    + \mat{u} \mat{D} \vec{u}
  \right)
  + 3 \mat{D} \vec{u},
\end{equation}
cf.\ \cite[Section~4.5]{ranocha2018thesis}. The flux is visualized
in \autoref{fig:quartic-flux}.
To compute the kinetic function $\phi^\flat$,
Riemann problems with an initial condition
\begin{equation}
  u_0(x) =
  \begin{cases}
    u_R, & x \in [0, 4.5], \\
    u_L, & \text{otherwise},
  \end{cases}
\end{equation}
have been solved in the periodic domain $[-7, 7]$ till the final time
$t = 3 / \max\set{\abs{f'(u)} | u_L \leq u \leq u_R}$ for the fixed right
state $u_R = 2$.

\begin{figure}[!htp]
\centering
  \begin{subfigure}{0.45\textwidth}
  \centering
    \includegraphics[width=\textwidth]{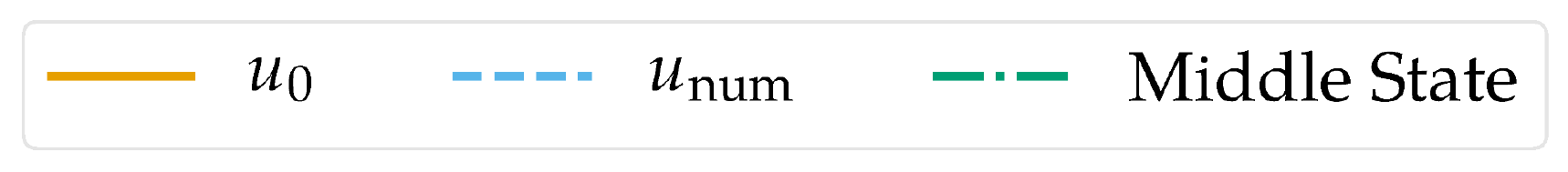}
  \end{subfigure}%
  \\
  \begin{subfigure}{0.33\textwidth}
  \centering
    \includegraphics[width=\textwidth]{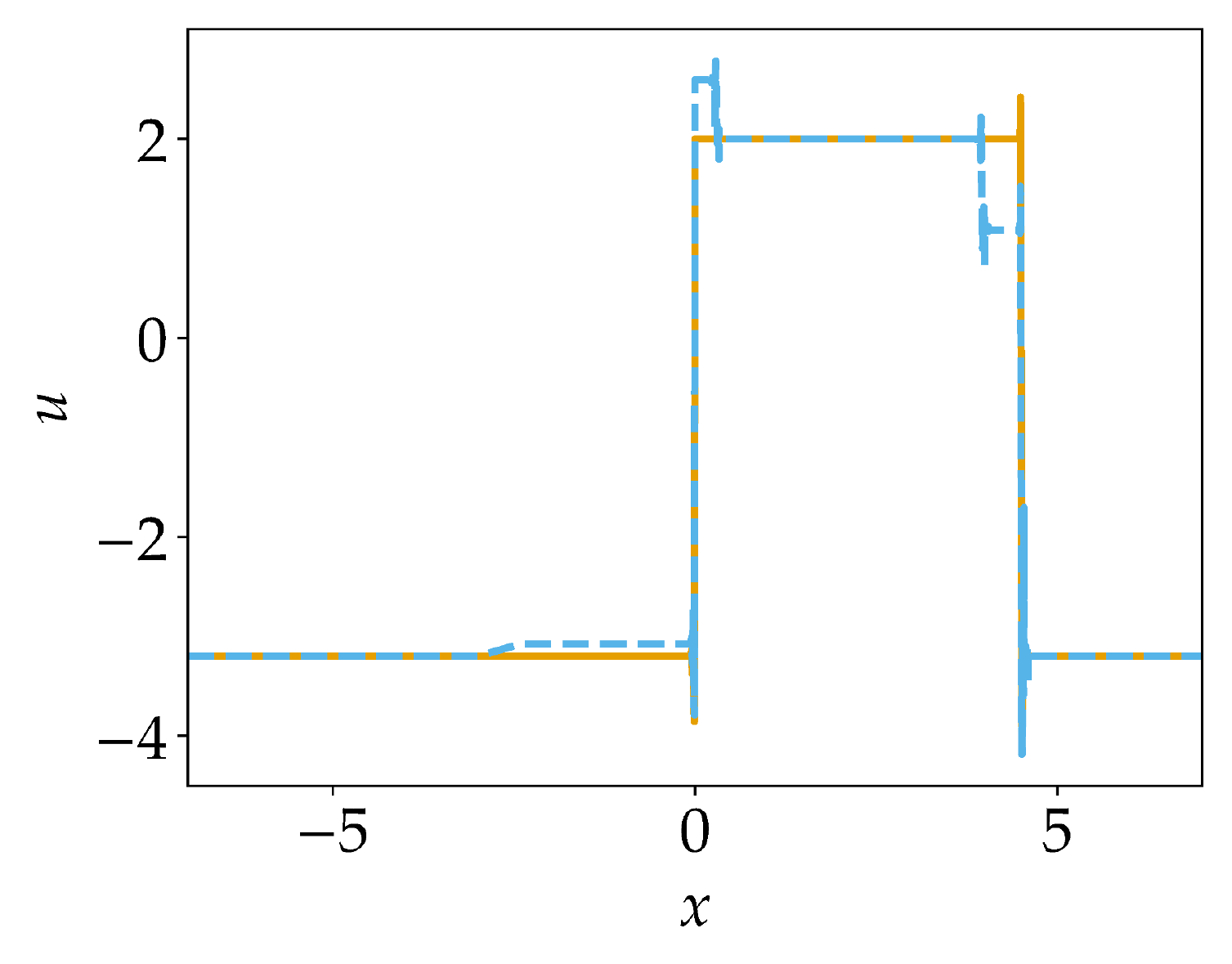}
    \caption{$u_L = -3.2$.}
  \end{subfigure}%
  \hspace*{\fill}
  \begin{subfigure}{0.33\textwidth}
  \centering
    \includegraphics[width=\textwidth]{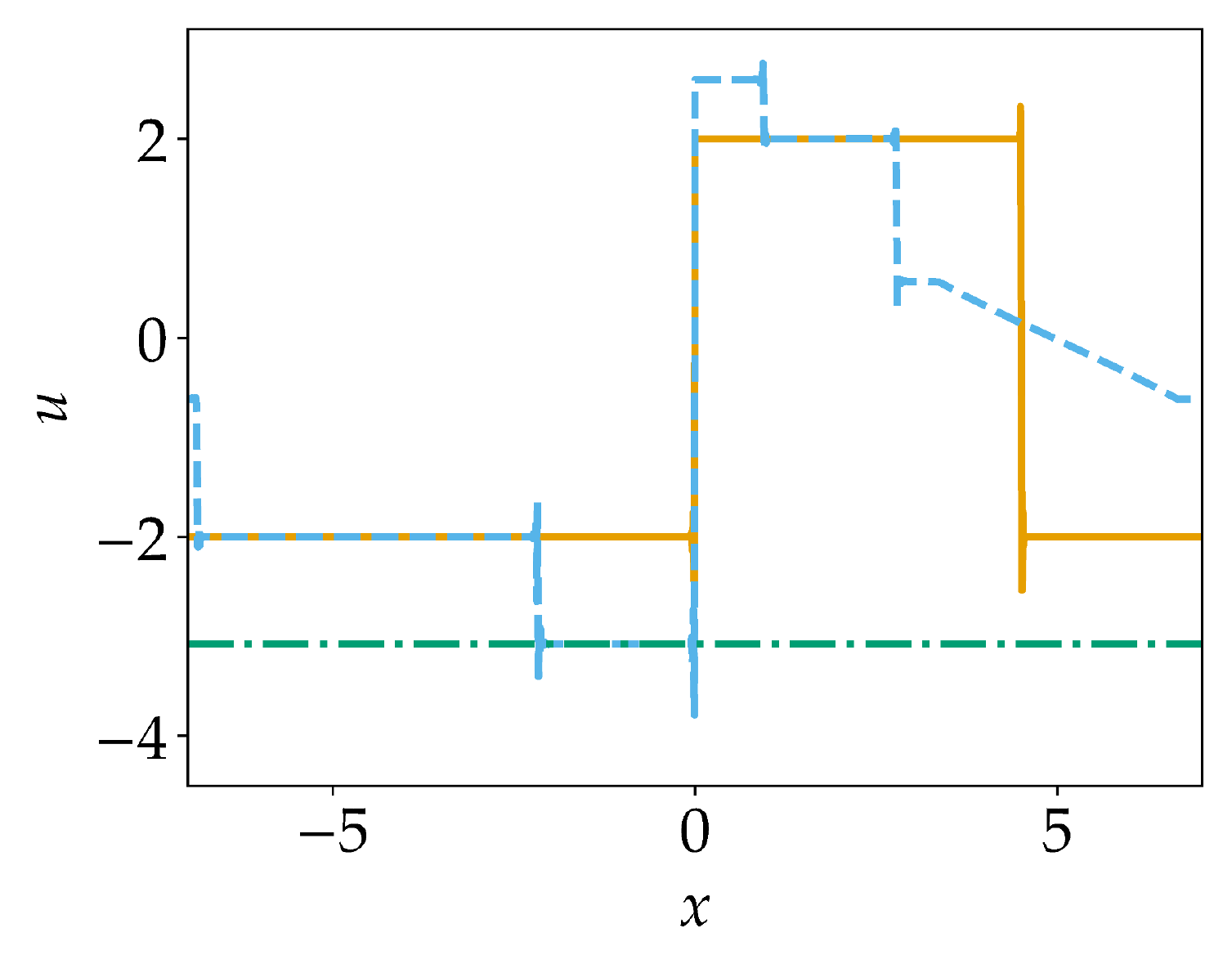}
    \caption{$u_L = -2$.}
  \end{subfigure}%
  \hspace*{\fill}
  \begin{subfigure}{0.33\textwidth}
  \centering
    \includegraphics[width=\textwidth]{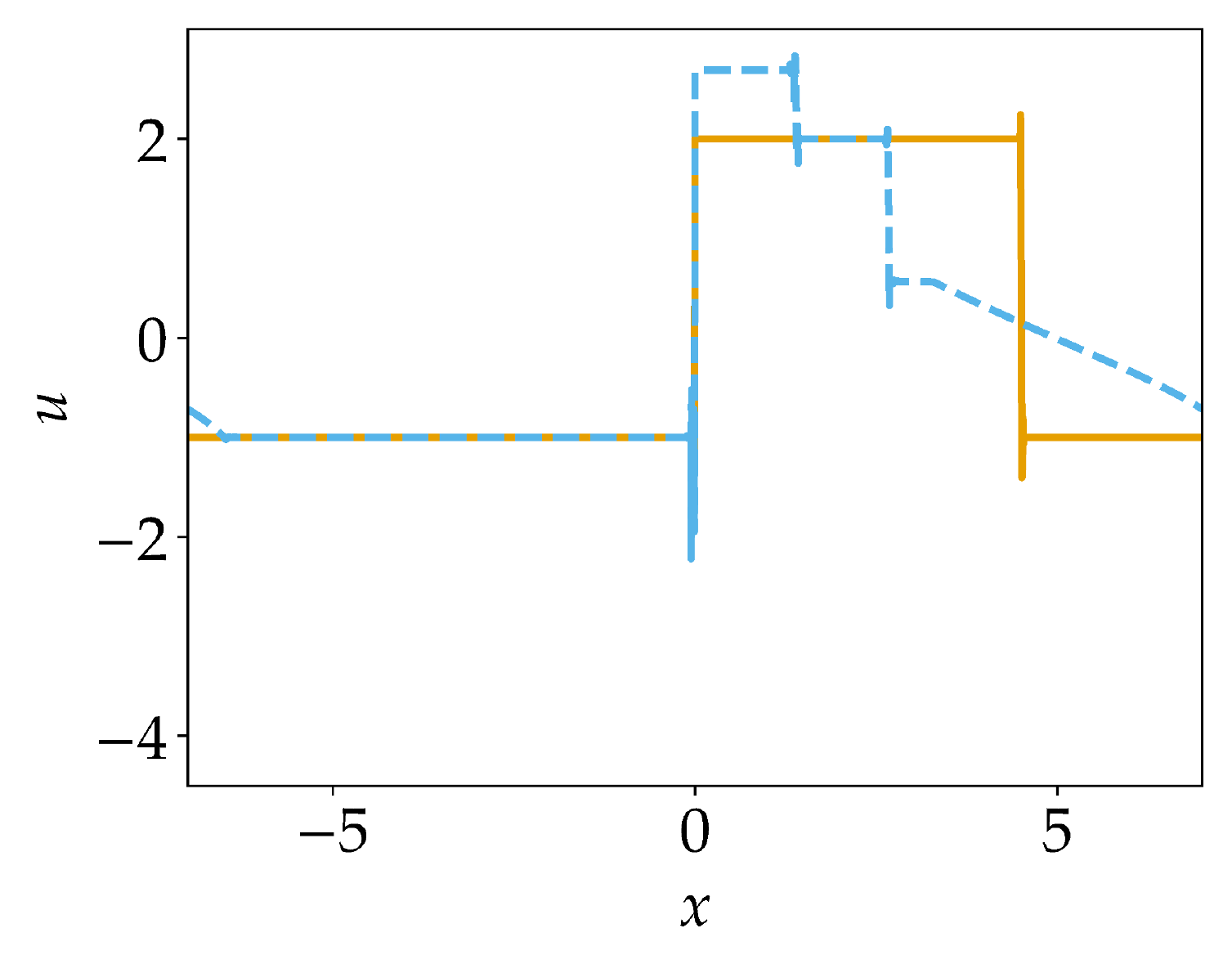}
    \caption{$u_L = -1$.}
  \end{subfigure}%
  \caption{Numerical solutions $u_\mathrm{num}$ (without postprocessing)
           of the Riemann problems
           with left-hand state $u_L$, right-hand state $u_R = 2$, and a DG method
           with the following parameters:
           polynomial degree $p = 5$,
           filter order $s = 0$,
           number of elements $N = 256$.}
  \label{fig:Quartic_DG_p_5_s_0_N_00256__Riemann}
\end{figure}

Typical solutions obtained by a DG method are presented in
Figure~\ref{fig:Quartic_DG_p_5_s_0_N_00256__Riemann}.
For small left-hand states such as $u_L = -3.2$, the left-hand state is
connected to classical wave. However, the right-hand state can be connected
to a nonclassical state above $u_R$.
For intermediate left-hand states such as $u_L = -2.0$, two nonclassical
states occur, one below $u_L$ and one above $u_R$.
For larger left-hand states such as $u_L = -1$, both the left and the right
state are connected to a nonclassical state above $u_R$.
In the following, the kinetic function $\phi^\flat$ is computed
as the mapping from a left-hand state $u_L$ to a nonclassical state below
the left-hand state.

\subsection{Nodal DG methods}
\label{sec:kinetic-function-quartic-DG}

For nodal DG methods, polynomial degrees $p \in \{0, 1, \dots, 5\}$,
filter order $s \in \{0, 1, \dots, 5\}$, and numbers of elements
$N \in \{2^6, 2^7, \dots, 2^{10}\}$ have been used. The numerical
flux between elements is the entropy-dissipative flux obtained by
adding local Lax-Friedrichs/Rusanov dissipation to the entropy
conservative numerical flux \eqref{eq:quartic-fnum-EC}.
The following observations have been made.
\begin{enumerate}
  \item
  The schemes with filter order $s \in \{1, 2, 3\}$ did not result
  in nonclassical solutions.

  \item
  The finite volume schemes ($p = 0$) did not result
  in nonclassical solutions. The schemes with polynomial degree
  $p \in \{1, 2\}$
  resulted in nonclassical solutions only if no filtering was applied
  ($s = 0$).

  \item
  For the investigated range of parameters, nonclassical
  solutions occurred only for $u_L \in [-4, -0.5]$
  (as a necessary criterion). Depending on the other parameters,
  the extremal values of $u_L$ for which nonclassical
  solution occurred can be different.

  \item
  The numerically obtained kinetic functions $\phi^\flat$ are
  approximately constant.
  They remain visually indistinguishable under
  grid refinement by increasing the number of elements $N$.
  This is shown for $p = 5$ and $s = 5$ in
  Figure~\ref{fig:Quartic_DG_p_5_s_5__kinetic_functions}.
  The slight deviations from the constant value near the
  extremal values of $u_L$ are caused by oscillations of the
  numerical solutions.
  \begin{figure}[!htb]
  \centering
    \includegraphics[width=0.6\textwidth]{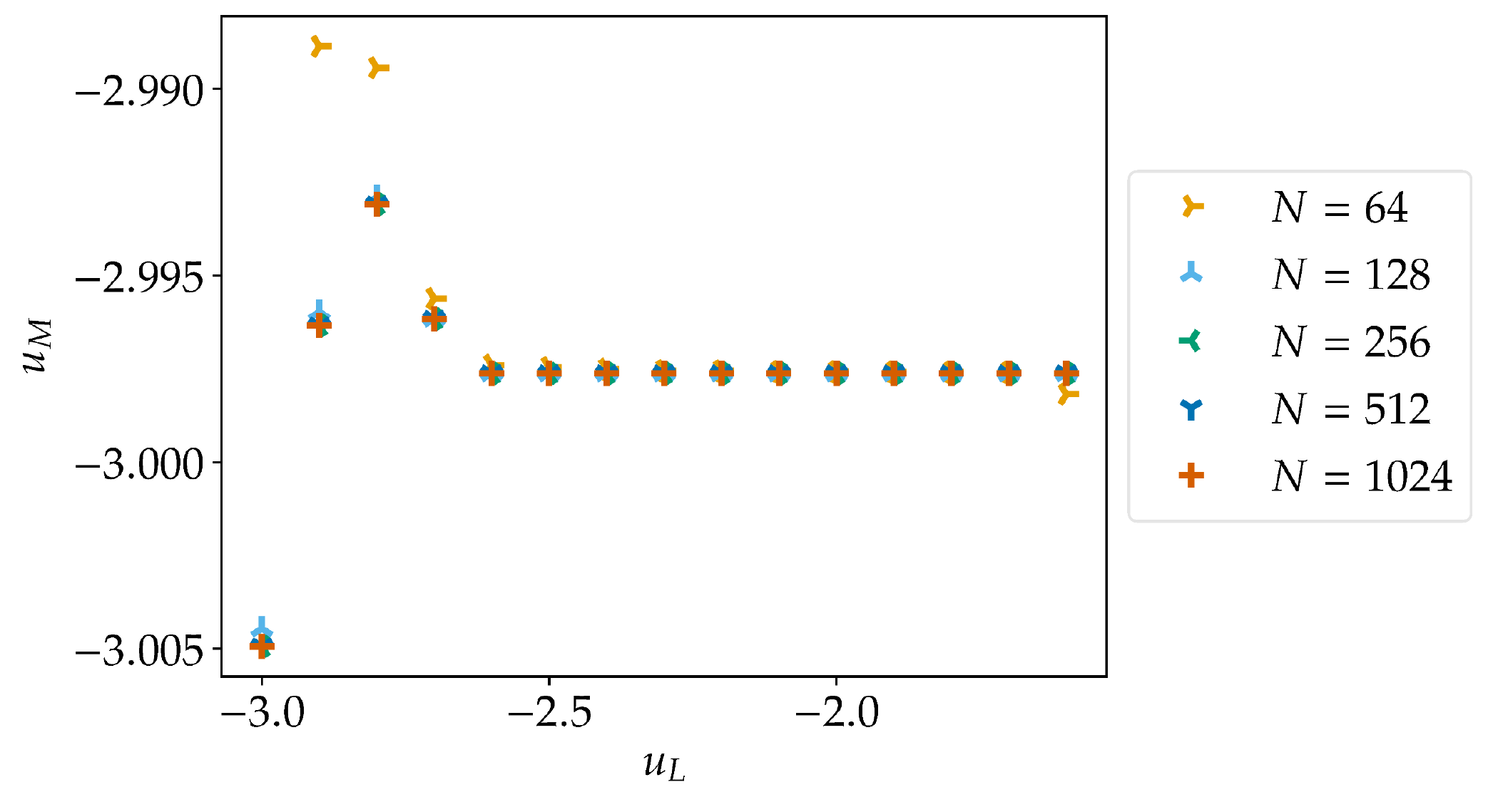}
    \caption{Kinetic function for a DG method with polynomial degree $p = 5$
            and filter order $s = 5$.}
    \label{fig:Quartic_DG_p_5_s_5__kinetic_functions}
  \end{figure}

  \item
  The kinetic functions do not depend on the polynomial degree $p$
  if no filtering is applied, i.e.\ $s = 0$.
  Otherwise, they depend on $p$.
  For $s = 4$, $u_M$ is bigger for $p = 4$ than for $p = 5$ while
  it is the other way round for $s = 5$.
  This is shown in
  Figure~\ref{fig:Quartic_DG_N_01024__kinetic_functions}.
  \begin{figure}[!htb]
  \centering
    \begin{subfigure}{0.6\textwidth}
    \centering
      \includegraphics[width=\textwidth]{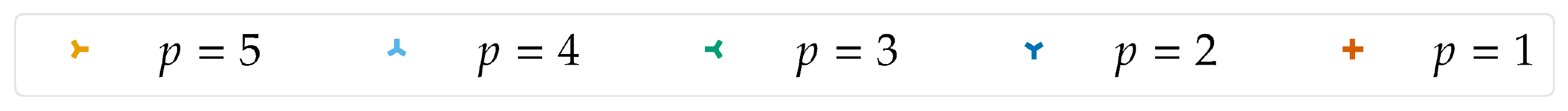}
    \end{subfigure}%
    \\
    \begin{subfigure}{0.33\textwidth}
    \centering
      \includegraphics[width=\textwidth]{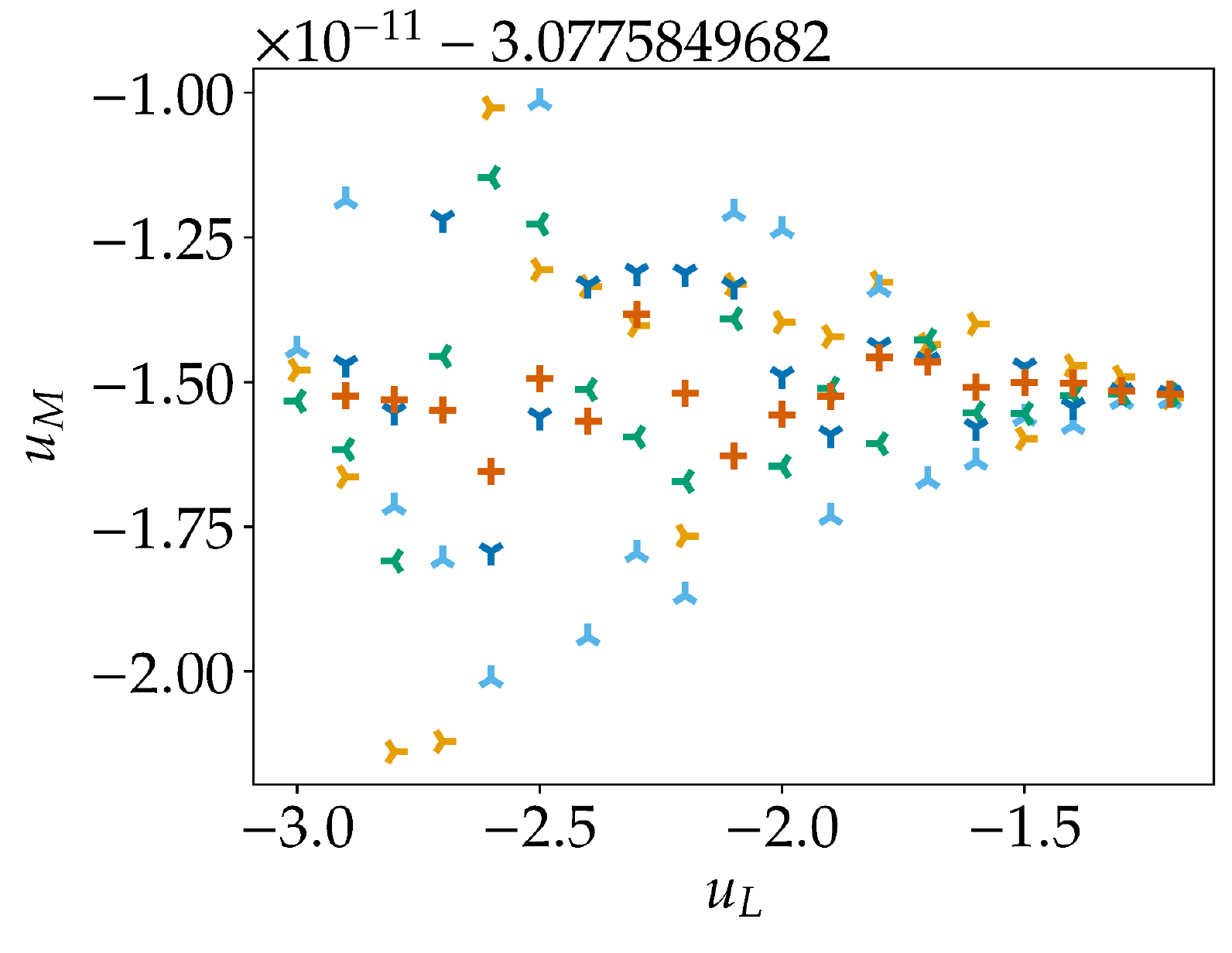}
      \caption{$s = 0$.}
    \end{subfigure}%
    \hspace*{\fill}
    \begin{subfigure}{0.33\textwidth}
    \centering
      \includegraphics[width=\textwidth]{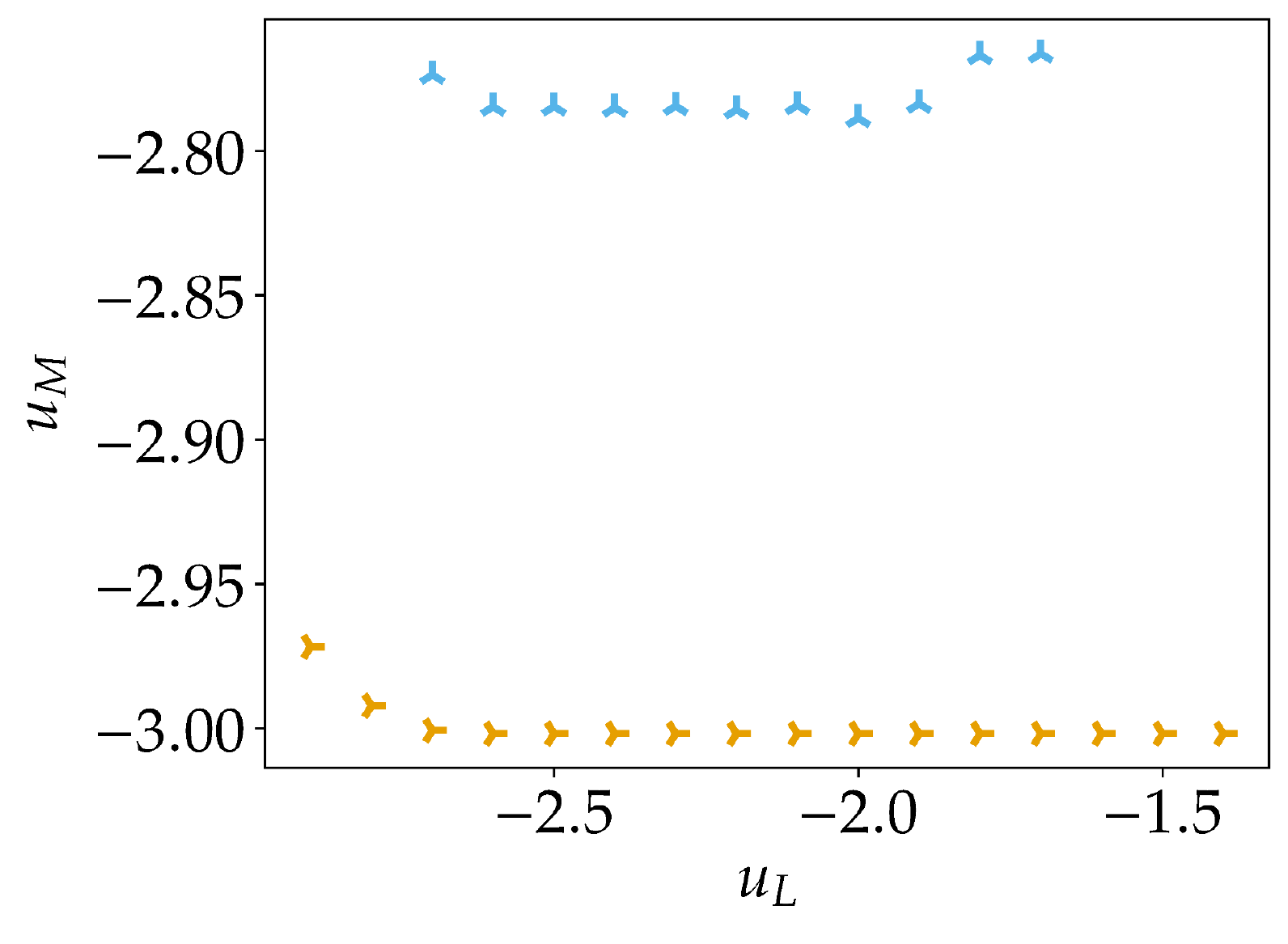}
      \caption{$s = 4$.}
    \end{subfigure}%
    \hspace*{\fill}
    \begin{subfigure}{0.33\textwidth}
    \centering
      \includegraphics[width=\textwidth]{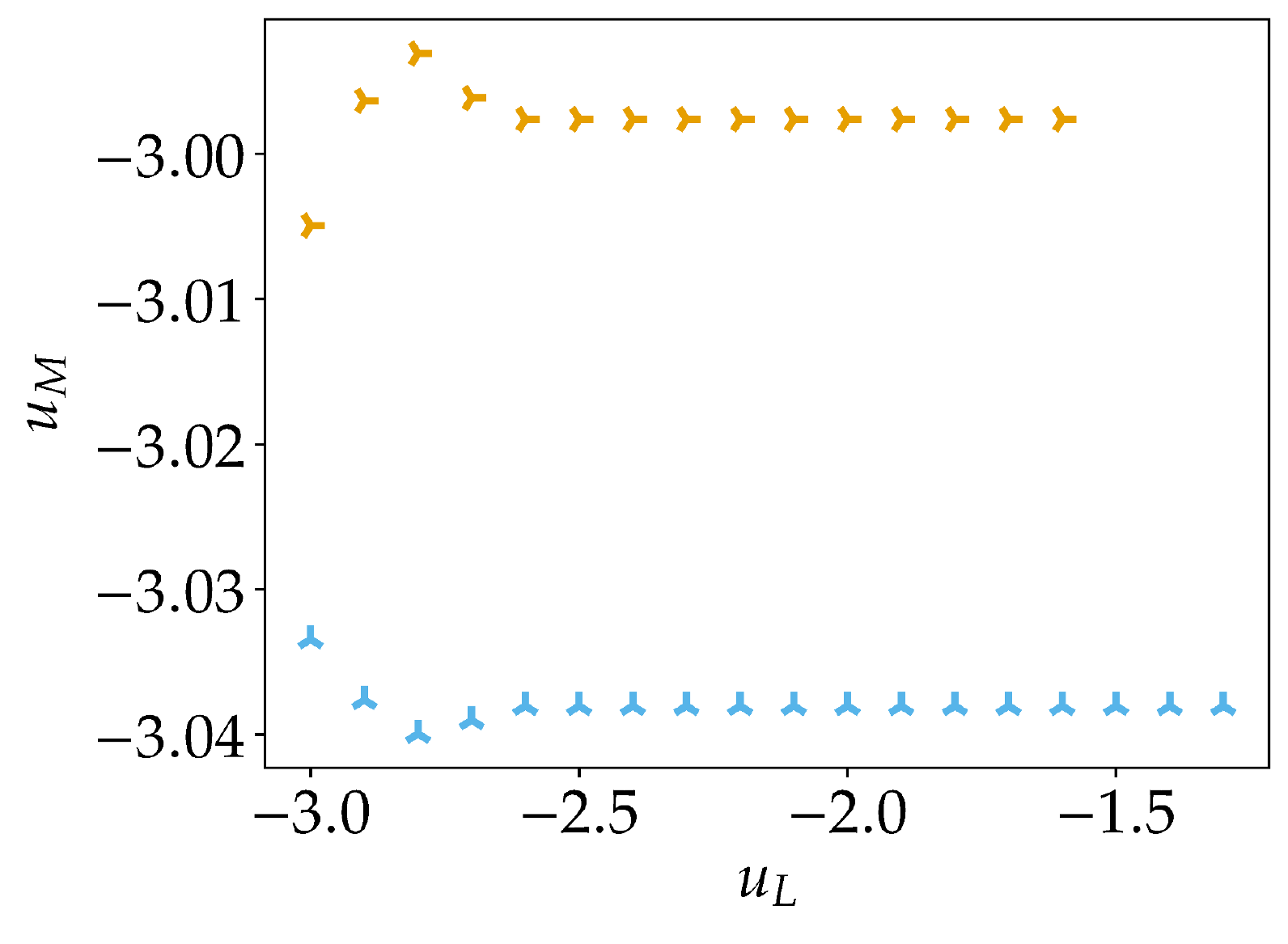}
      \caption{$s = 5$.}
    \end{subfigure}%
    \caption{Kinetic function for DG methods with $N = \num{1024}$ elements.}
    \label{fig:Quartic_DG_N_01024__kinetic_functions}
  \end{figure}

  \item
  The kinetic functions depend on the filter order $s$.
  However, there is no clear relation as shown in
  Figure~\ref{fig:Quartic_DG_N_01024__kinetic_functions}.
\end{enumerate}

\subsection{Finite difference methods}
\label{sec:kinetic-function-quartic-FD}

For FD methods, accuracy order $p \in \{2, 4, 6\}$,
artificial dissipation strength $\epsilon_i \in \{0, 100, \dots, 400\}$,
and $N \in \{2^9, 2^{10}, \dots, 2^{14}\}$ grid nodes have been used.
\begin{enumerate}
  \item
  If artificial dissipation was applied, nonclassical solutions
  occurred at least for some values of $N$, even if only the
  second-order artificial dissipation was used
  ($\epsilon_2 \neq 0, \epsilon_4 = \epsilon_6 = 0$).
  The only exception is given by the second-order method ($p = 2$) with
  second-order artificial dissipation if $\epsilon_2 \neq 0$ is
  sufficiently big.
  This is again in agreement with  \cite{ranocha2019mimetic}.

  \item
  As for DG methods, nonclassical solutions occured only for
  $u_L \in [-4, -0.5]$ for the investigated range of parameters.

  \item
  The numerically obtained kinetic functions $\phi^\flat$ are
  again approximately constant.
  They vary under grid refinement and the range of left-hand states
  $u_L$ for which nonclassical solutions with middle state
  $u_M < u_L$ occur typically increases with the number of grid
  points. This is shown in Figure~\ref{fig:Quartic_FD_p_6_Strength2_0_Strength4_0100_Strength6_0100__kinetic_functions}.
  \begin{figure}[!htb]
  \centering
    \begin{subfigure}{0.9\textwidth}
    \centering
      \includegraphics[width=\textwidth]{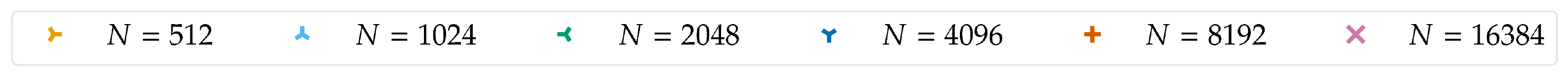}
    \end{subfigure}%
    \\
    \begin{subfigure}{0.5\textwidth}
    \centering
      \includegraphics[width=\textwidth]{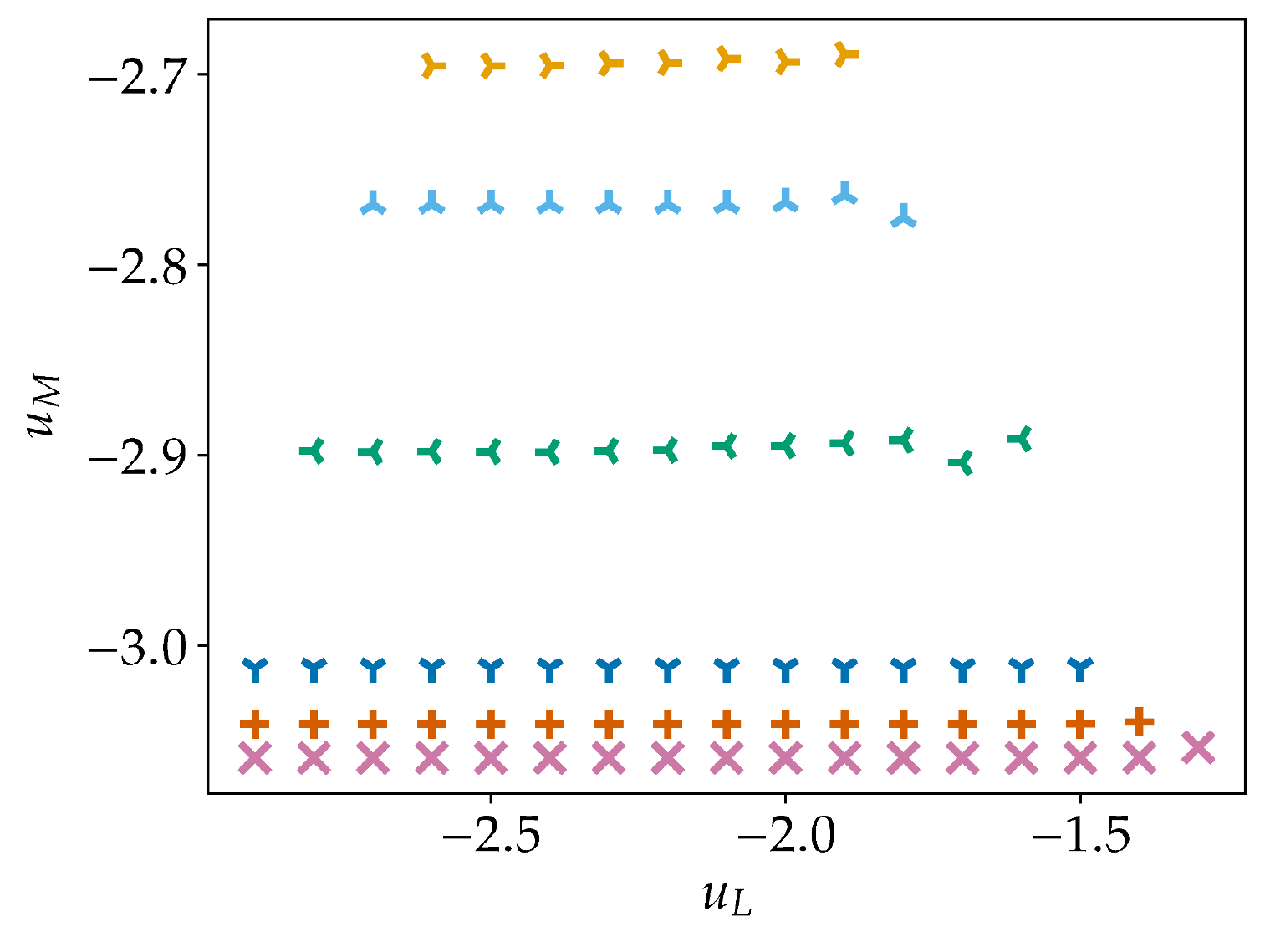}
      \caption{$\epsilon_2 = 0$, $\epsilon_4 = 100$, $\epsilon_6 = 0$.}
    \end{subfigure}%
    \begin{subfigure}{0.5\textwidth}
    \centering
      \includegraphics[width=\textwidth]{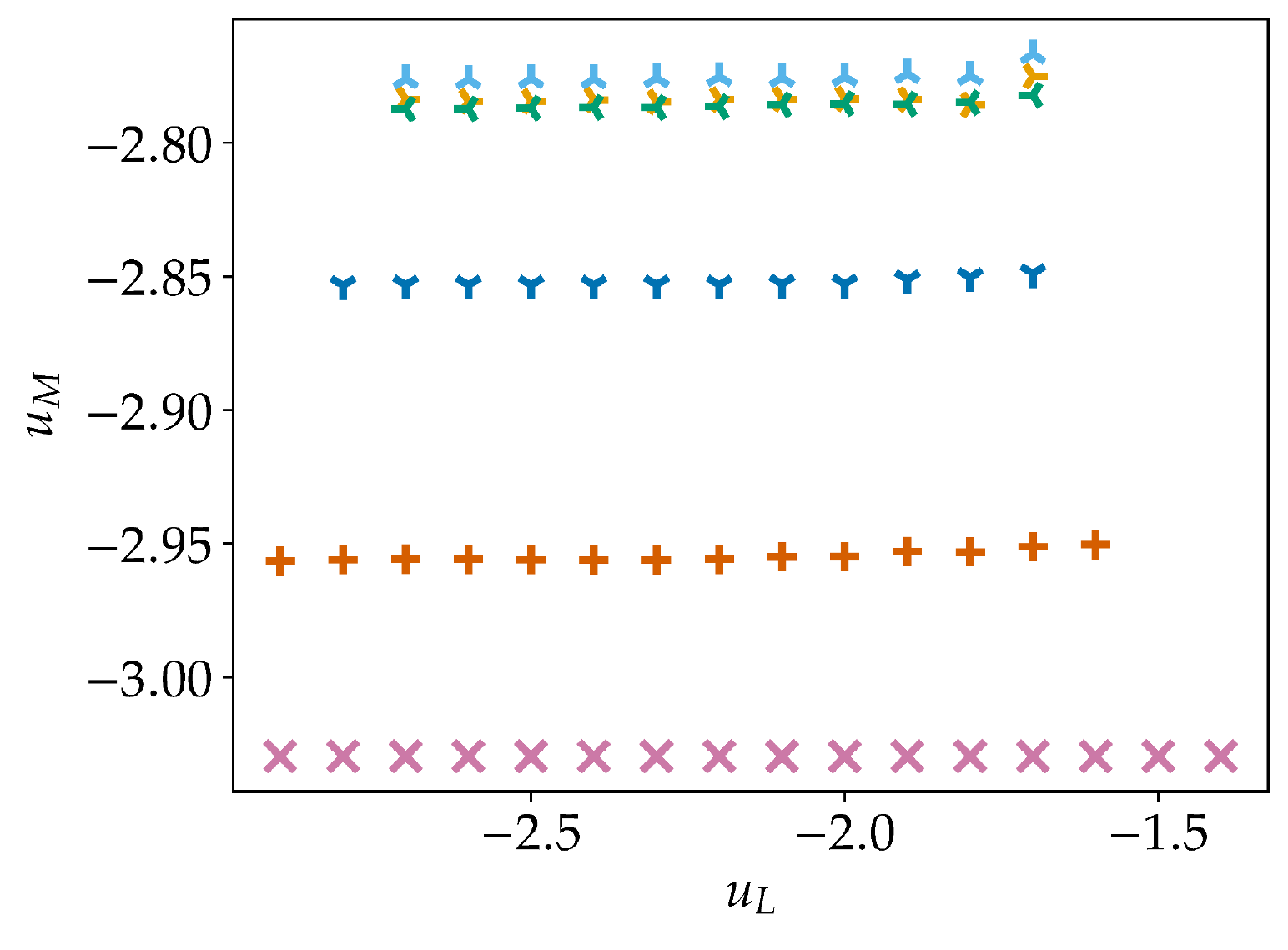}
      \caption{$\epsilon_2 = 0$, $\epsilon_4 = 0$, $\epsilon_6 = 100$.}
    \end{subfigure}%
    \caption{Kinetic function for FD methods with order of accuracy
             $p = 6$ and varying number of grid nodes $N$ for different
             orders of the artificial dissipation.}
    \label{fig:Quartic_FD_p_6_Strength2_0_Strength4_0100_Strength6_0100__kinetic_functions}
  \end{figure}

  \item
  Choosing a fixed order of the artificial dissipation, i.e.\
  $\epsilon_i \in \{100, 200, 300, 400\}$ and $\epsilon_j = 0$
  for $j \neq i$, the kinetic functions depend on the strength
  $\epsilon_i$. For high resolutions, the value of the kinetic
  function typically increases with the strength $\epsilon_i$.
  This is shown in
  Figure~\ref{fig:Quartic_FD_p_6_Strengthi_N_16384__kinetic_functions}.
  \begin{figure}[!htb]
  \centering
    \begin{subfigure}{0.55\textwidth}
    \centering
      \includegraphics[width=\textwidth]{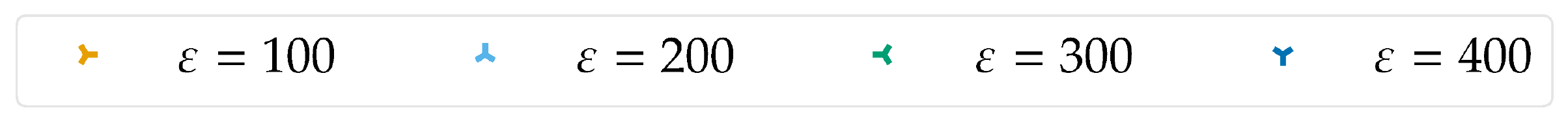}
    \end{subfigure}%
    \\
    \begin{subfigure}{0.33\textwidth}
    \centering
      \includegraphics[width=\textwidth]{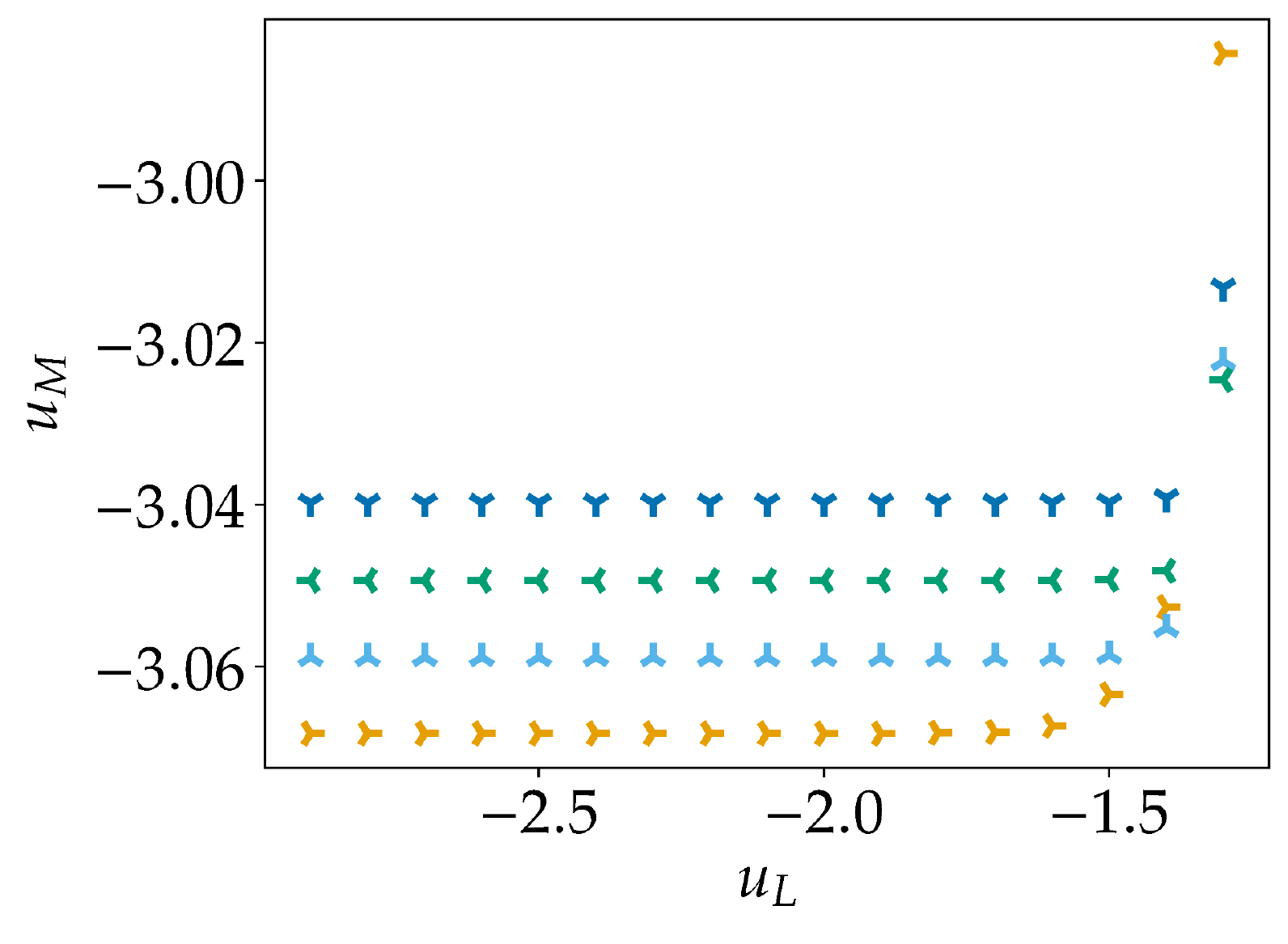}
      \caption{$\epsilon_2 = \epsilon$, $\epsilon_4 = \epsilon_6 = 0$.}
    \end{subfigure}%
    \hspace*{\fill}
    \begin{subfigure}{0.33\textwidth}
    \centering
      \includegraphics[width=\textwidth]{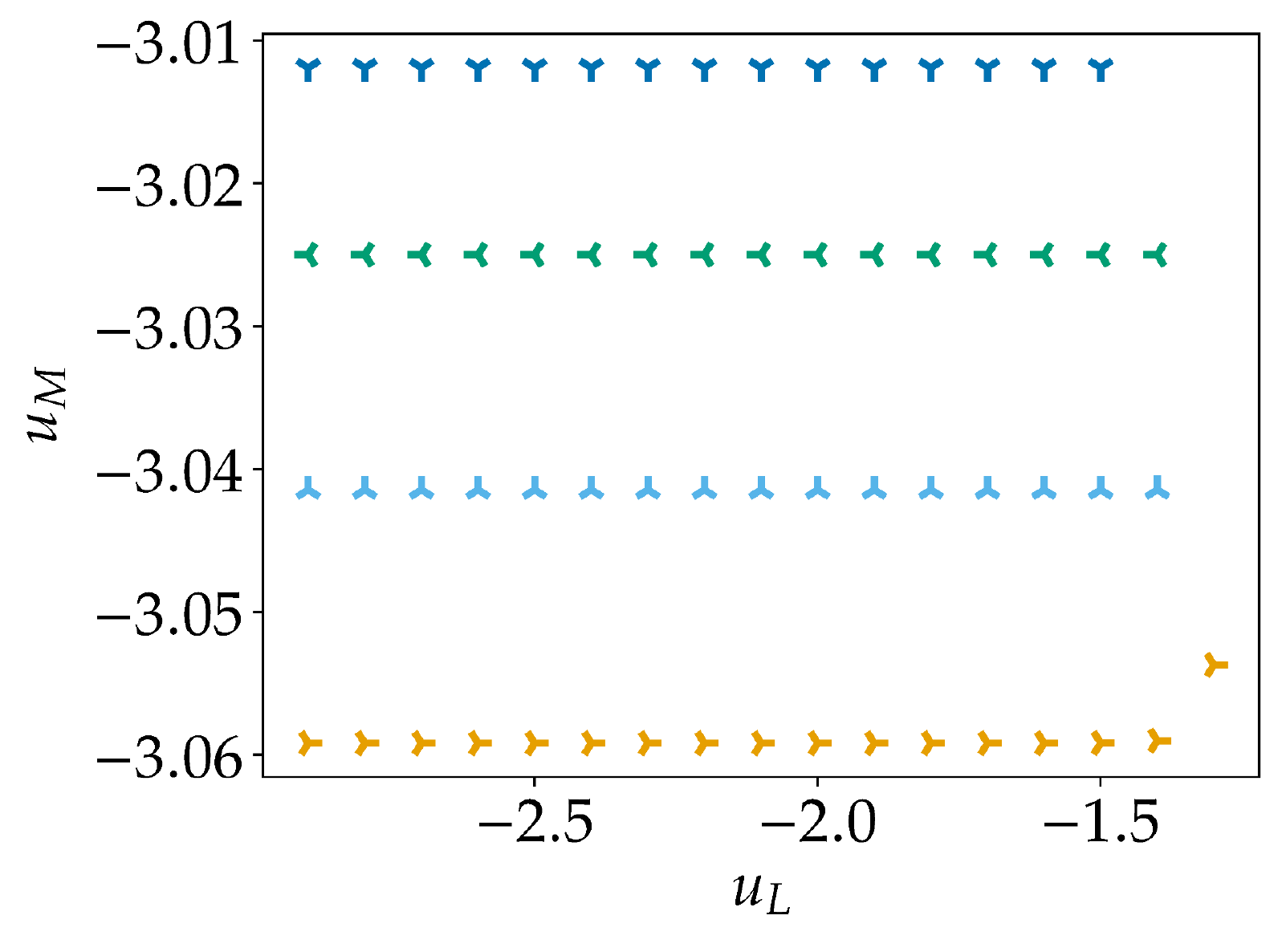}
      \caption{$\epsilon_4 = \epsilon$, $\epsilon_2 = \epsilon_6 = 0$.}
    \end{subfigure}%
    \hspace*{\fill}
    \begin{subfigure}{0.33\textwidth}
    \centering
      \includegraphics[width=\textwidth]{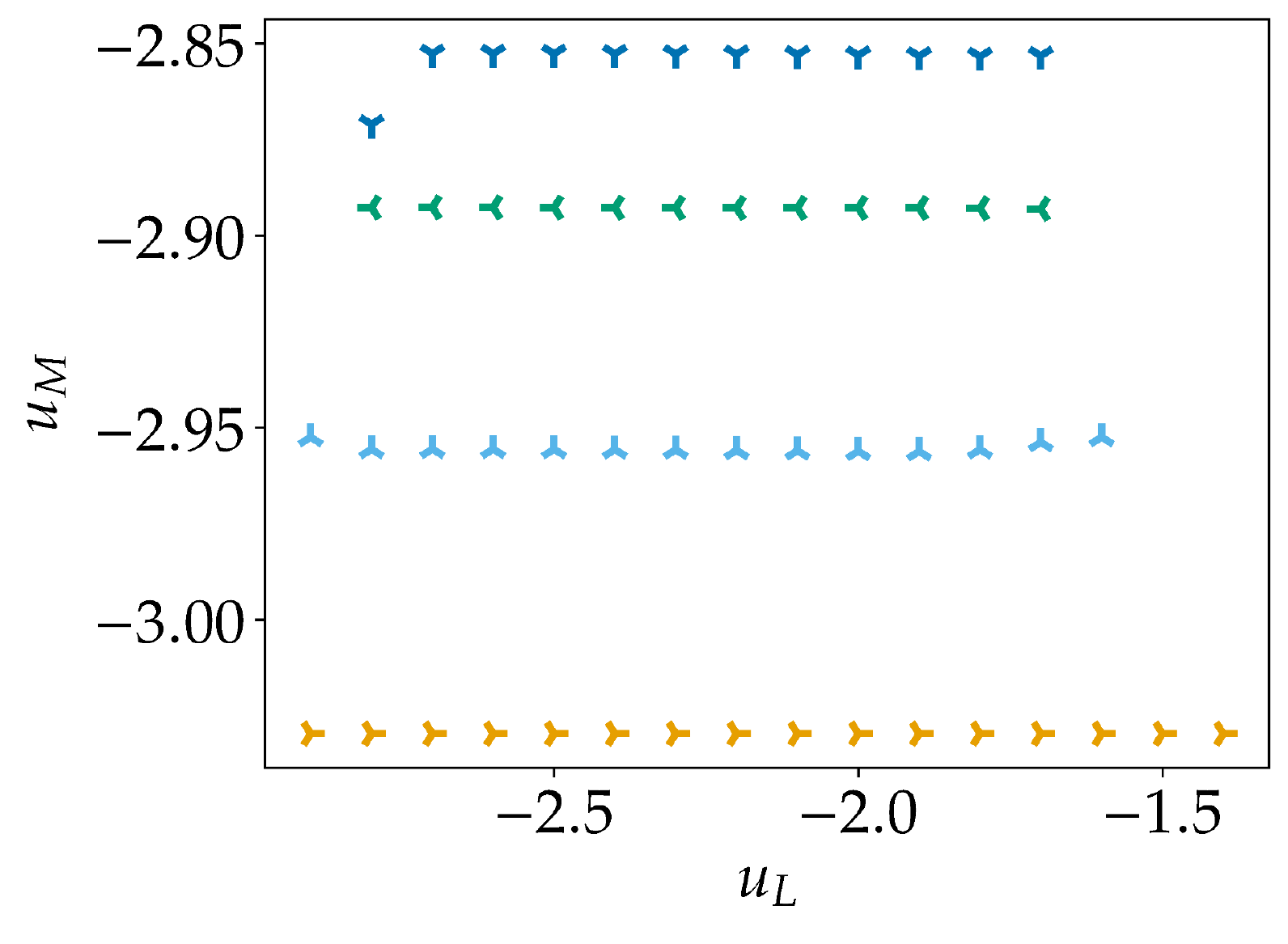}
      \caption{$\epsilon_6 = \epsilon$, $\epsilon_2 = \epsilon_4 = 0$.}
    \end{subfigure}%
    \caption{Kinetic function for FD methods with order of accuracy
             $p = 6$, $N = \num{16384}$ grid nodes, and different strengths
             of the artificial dissipation.}
    \label{fig:Quartic_FD_p_6_Strengthi_N_16384__kinetic_functions}
  \end{figure}

  \item
  The kinetic function varies slightly with the order of accuracy $p$.
  Typically, it increases (more negative) for higher values of $p$.
  This is shown in
  Figure~\ref{fig:Quartic_FD_Strength2_0_Strength4_0300_Strength6_0400_N_16384__kinetic_functions}.
  \begin{figure}[!htb]
  \centering
    \begin{subfigure}[b]{0.35\textwidth}
    \centering
      \includegraphics[width=\textwidth]{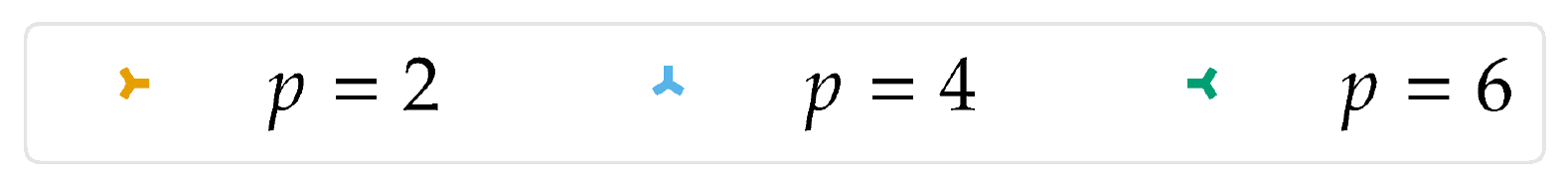}
    \end{subfigure}%
    \\
    \begin{subfigure}[b]{0.45\textwidth}
    \centering
      \includegraphics[width=\textwidth]{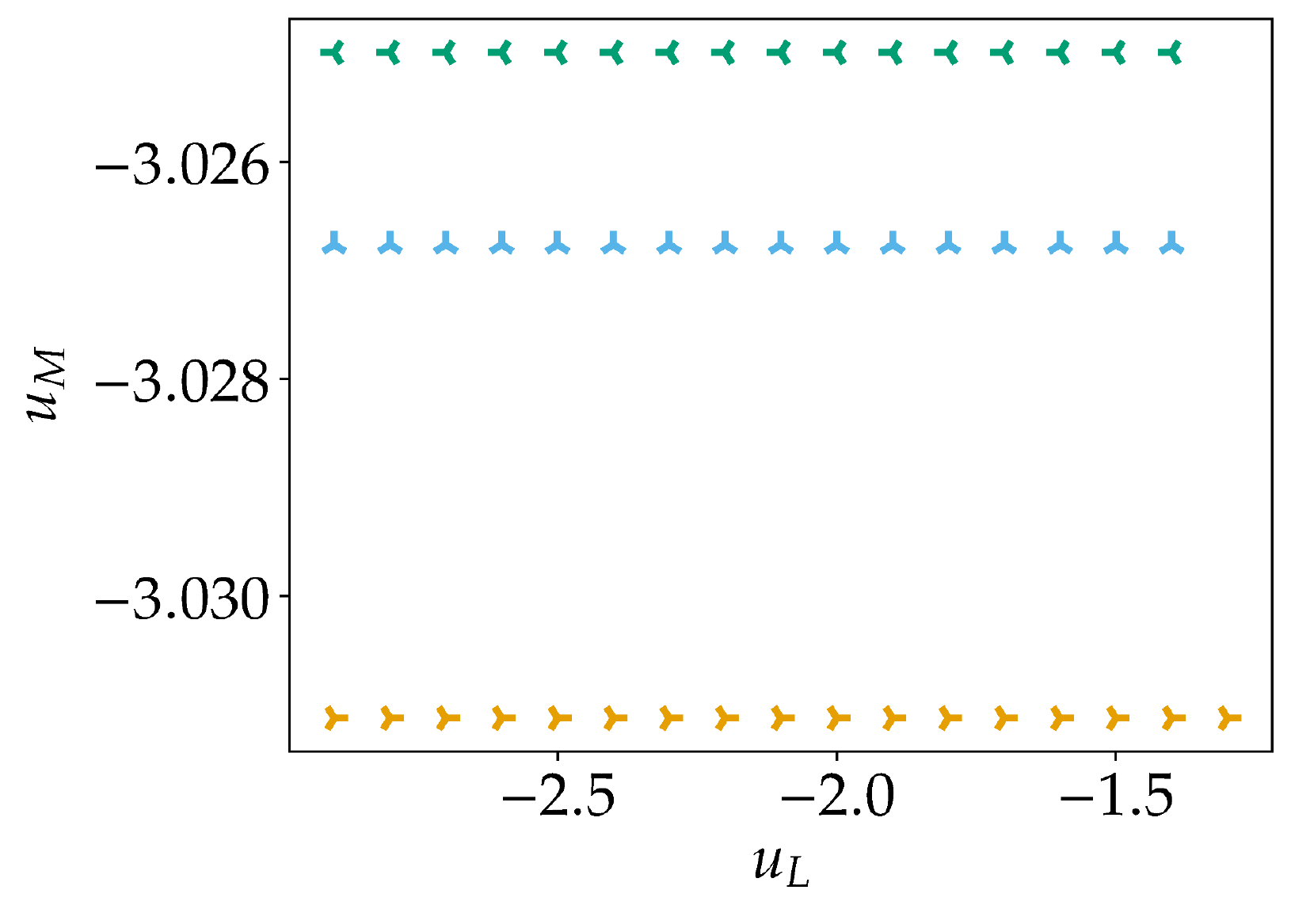}
      \caption{$\epsilon_2 = 0$, $\epsilon_4 = 300$, $\epsilon_6 = 0$.}
    \end{subfigure}%
    \hspace*{\fill}
    \begin{subfigure}[b]{0.45\textwidth}
    \centering
      \includegraphics[width=\textwidth]{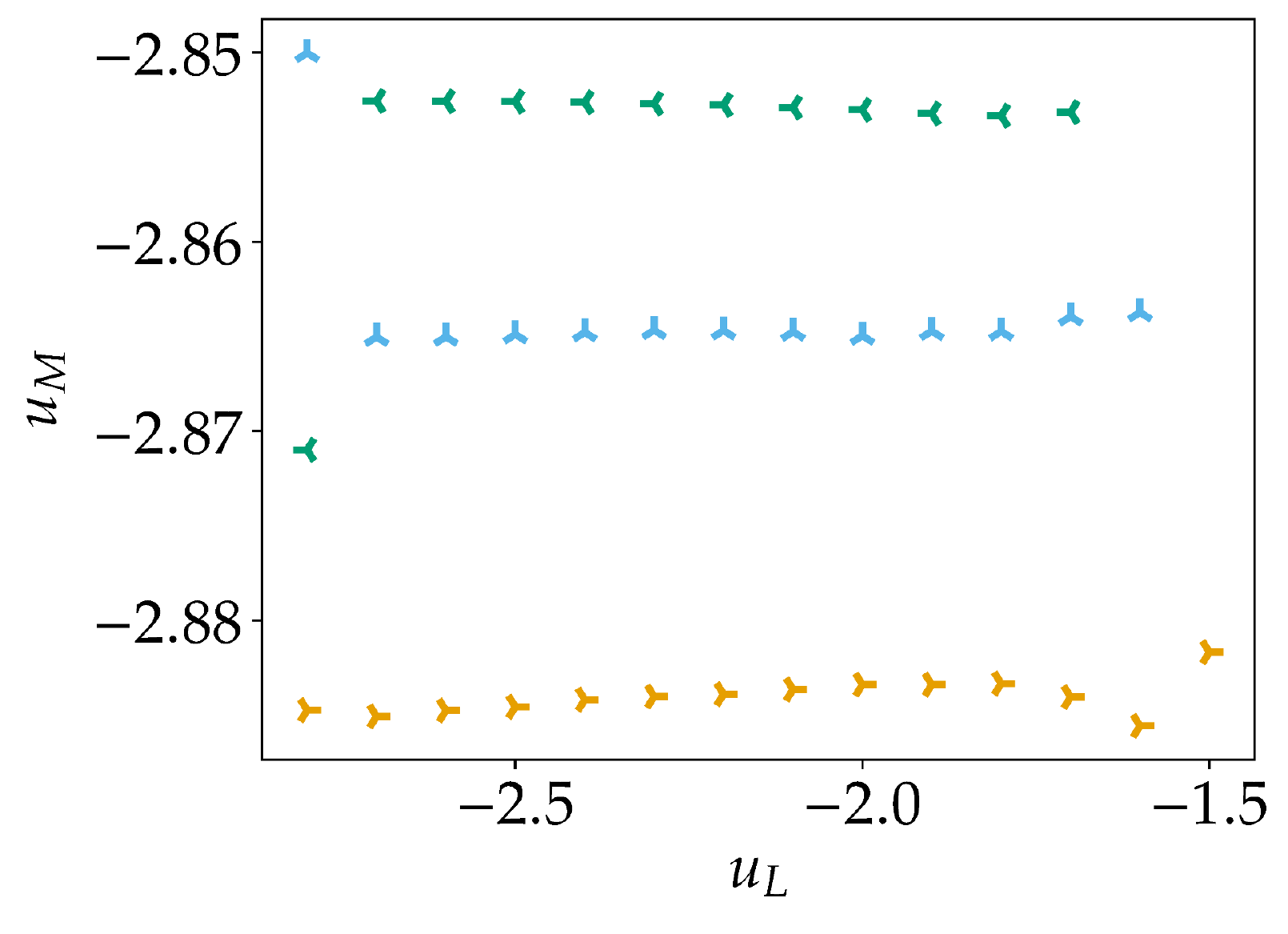}
      \caption{$\epsilon_2 = 0$, $\epsilon_4 = 0$, $\epsilon_6 = 400$.}
    \end{subfigure}%
    \caption{Kinetic function for FD methods with different orders
             of accuracy $p$, $N = \num{16384}$ grid nodes, and different
             strengths of the artificial dissipation.}
    \label{fig:Quartic_FD_Strength2_0_Strength4_0300_Strength6_0400_N_16384__kinetic_functions}
  \end{figure}
\end{enumerate}

\subsection{Fourier collocation methods}
\label{sec:kinetic-function-quartic-Fourier}

For Fourier methods, the viscosity strengths $\epsilon \in \{10 / N,
50 / N, 100 / N\}$ for the standard and convergent choices of
\cite{tadmor2012adaptive} and $N \in \{2^{10}, 2^{11}, \dots, 2^{14}\}$
grid nodes have been used.
The following observations have been made.
\begin{enumerate}
  \item
  As for DG and FD methods, nonclassical intermediate states
  $u_M < u_L$ occured only if $u_L \in [-4, -0.5]$ for the
  investigated range of parameters.

  \item
  In contrast to DG and FD methods, not all numerically obtained
  kinetic functions $\phi^\flat$ are approximately constant.
  If the strength of the standard spectral viscosity is low,
  there is a zig-zag behavior of the kinetic function that is
  stable under grid refinement.
  If the strength of the spectral viscosity is bigger, the kinetic
  functions are approximately constant and vary slightly under
  grid refinement, similarly to DG and FD methods. This is visualized in
  Figure~\ref{fig:Quartic_Fourier_TadmorWaagan2012Standard_Strength_1001000_Split_1__kinetic_function}.
  \begin{figure}[!htb]
  \centering
    \begin{subfigure}{\textwidth}
    \centering
      \includegraphics[width=0.9\textwidth]{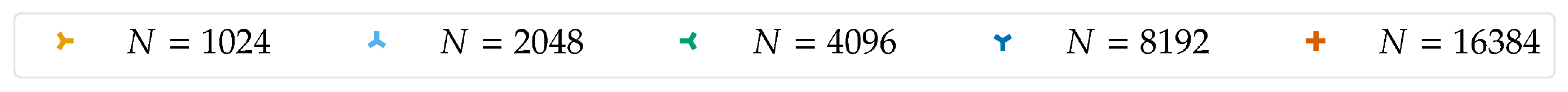}
    \end{subfigure}%
    \\
    \begin{subfigure}[b]{0.45\textwidth}
    \centering
      \includegraphics[width=\textwidth]{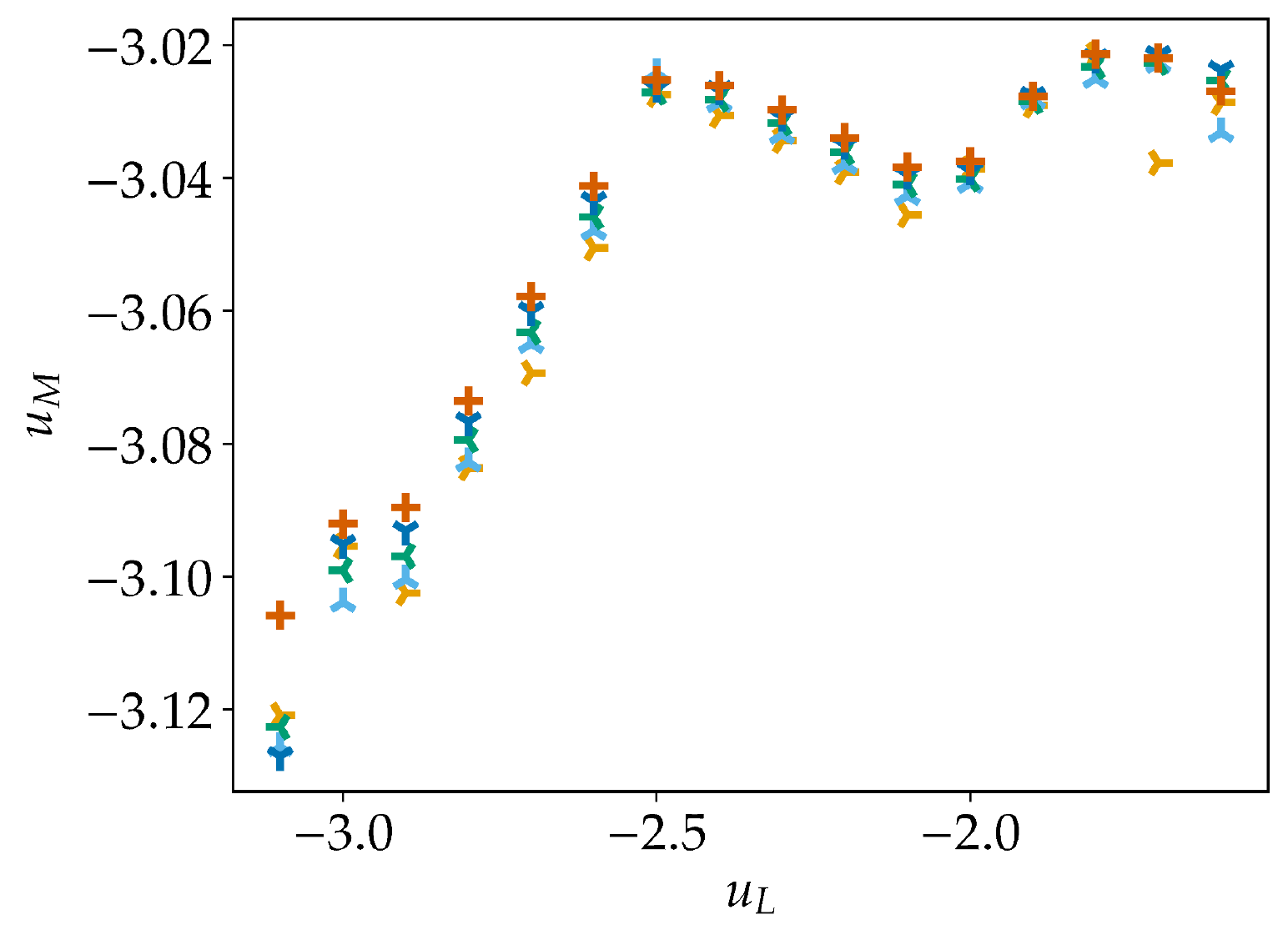}
      \caption{$\epsilon = 10 / N$.}
    \end{subfigure}%
    \hspace*{\fill}
    \begin{subfigure}[b]{0.45\textwidth}
    \centering
      \includegraphics[width=\textwidth]{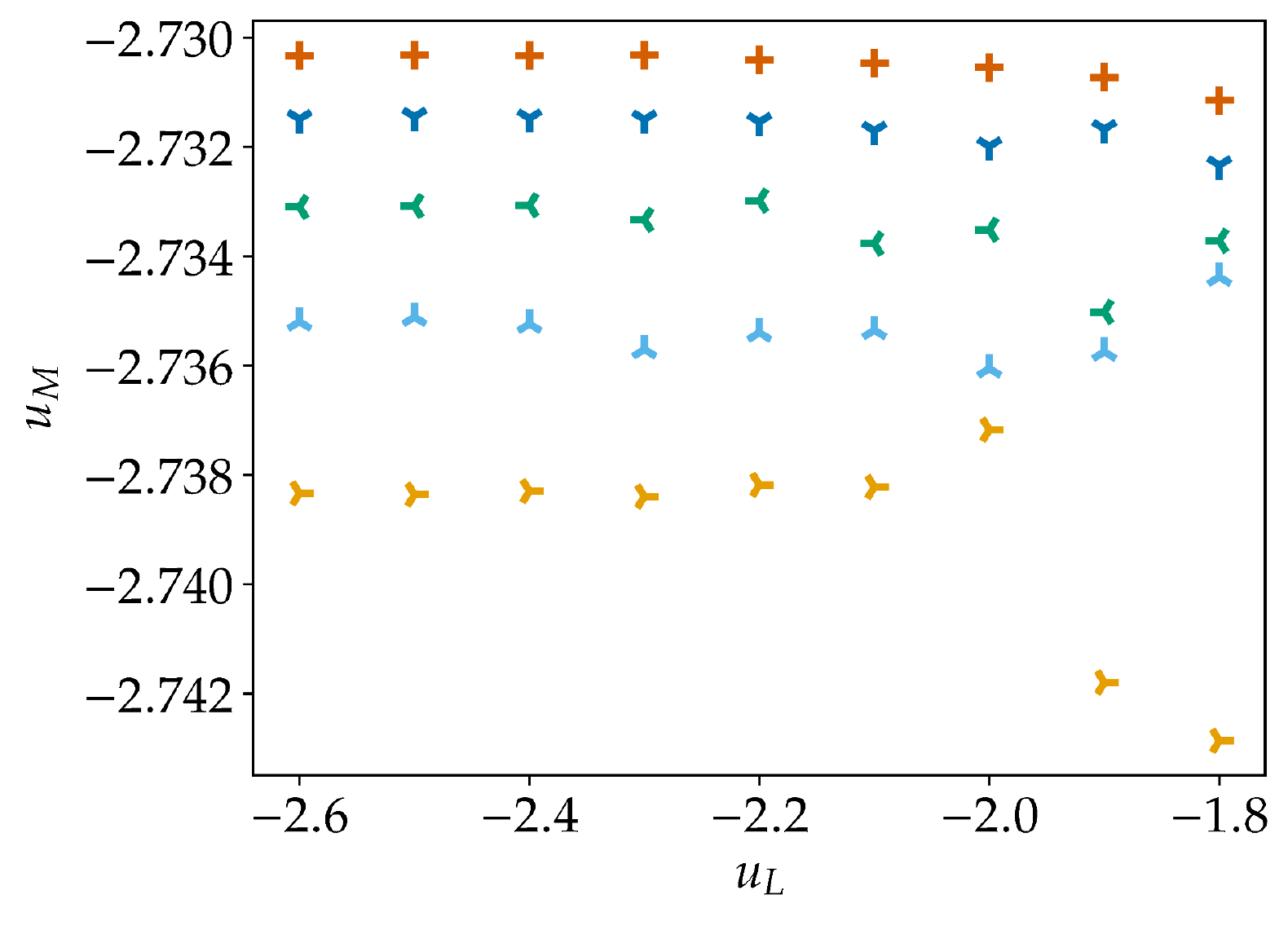}
      \caption{$\epsilon = 100 / N$.}
    \end{subfigure}%
    \caption{Kinetic function for Fourier methods using the
             standard choice of \cite{tadmor2012adaptive} with different
             strengths of the spectral viscosity and numbers $N$
             of grid nodes.}
    \label{fig:Quartic_Fourier_TadmorWaagan2012Standard_Strength_1001000_Split_1__kinetic_function}
  \end{figure}

  \item
  The kinetic functions depend on the strength of the spectral viscosity.
  For the convergent choice of \cite{tadmor2012adaptive}, nonclassical
  shocks can occur or not, depending on the strength $\epsilon$.
  For example, many nonclassical shocks occured for $\epsilon = 10 / N$,
  none for $\epsilon = 50 / N$ and only a few for $\epsilon = 100 / N$.
  In contrast, nonclassical solutions occured for all of these strengths
  if the classical spectral viscosity is applied.
  This can be seen in Figure~\ref{fig:Quartic_Fourier_N_16384_Split_1__kinetic_function}.
  \begin{figure}[!htb]
  \centering
    \begin{subfigure}[b]{0.4\textwidth}
    \centering
      \includegraphics[width=\textwidth]{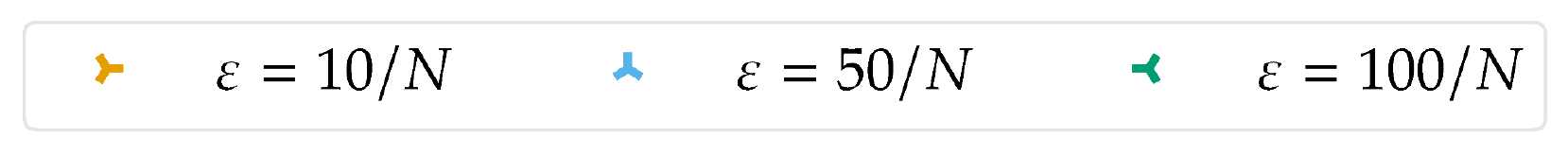}
    \end{subfigure}%
    \\
    \begin{subfigure}{0.45\textwidth}
    \centering
      \includegraphics[width=\textwidth]{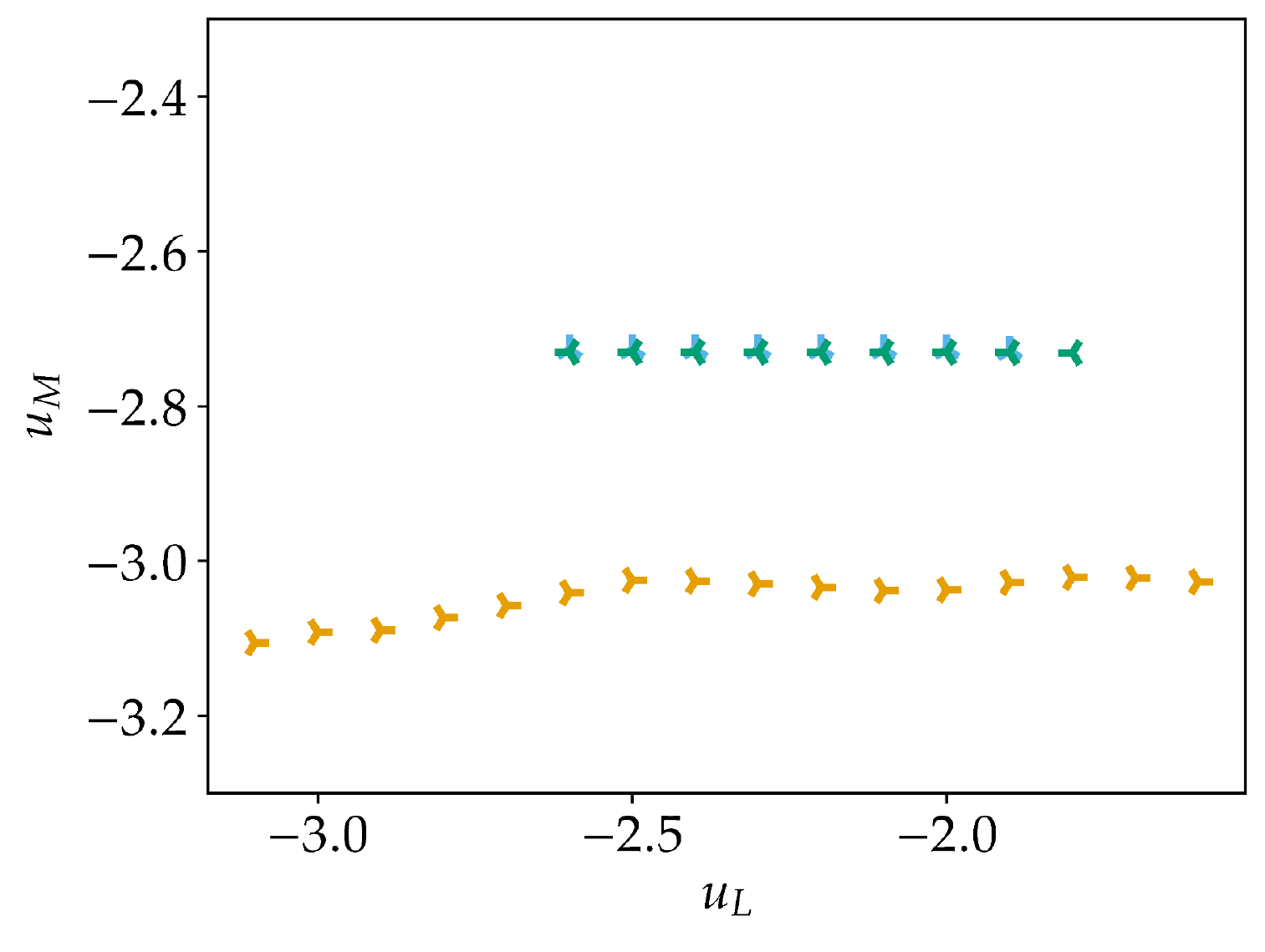}
      \caption{Standard choice of \cite{tadmor2012adaptive}.}
    \end{subfigure}%
    \hspace*{\fill}
    \begin{subfigure}{0.45\textwidth}
    \centering
      \includegraphics[width=\textwidth]{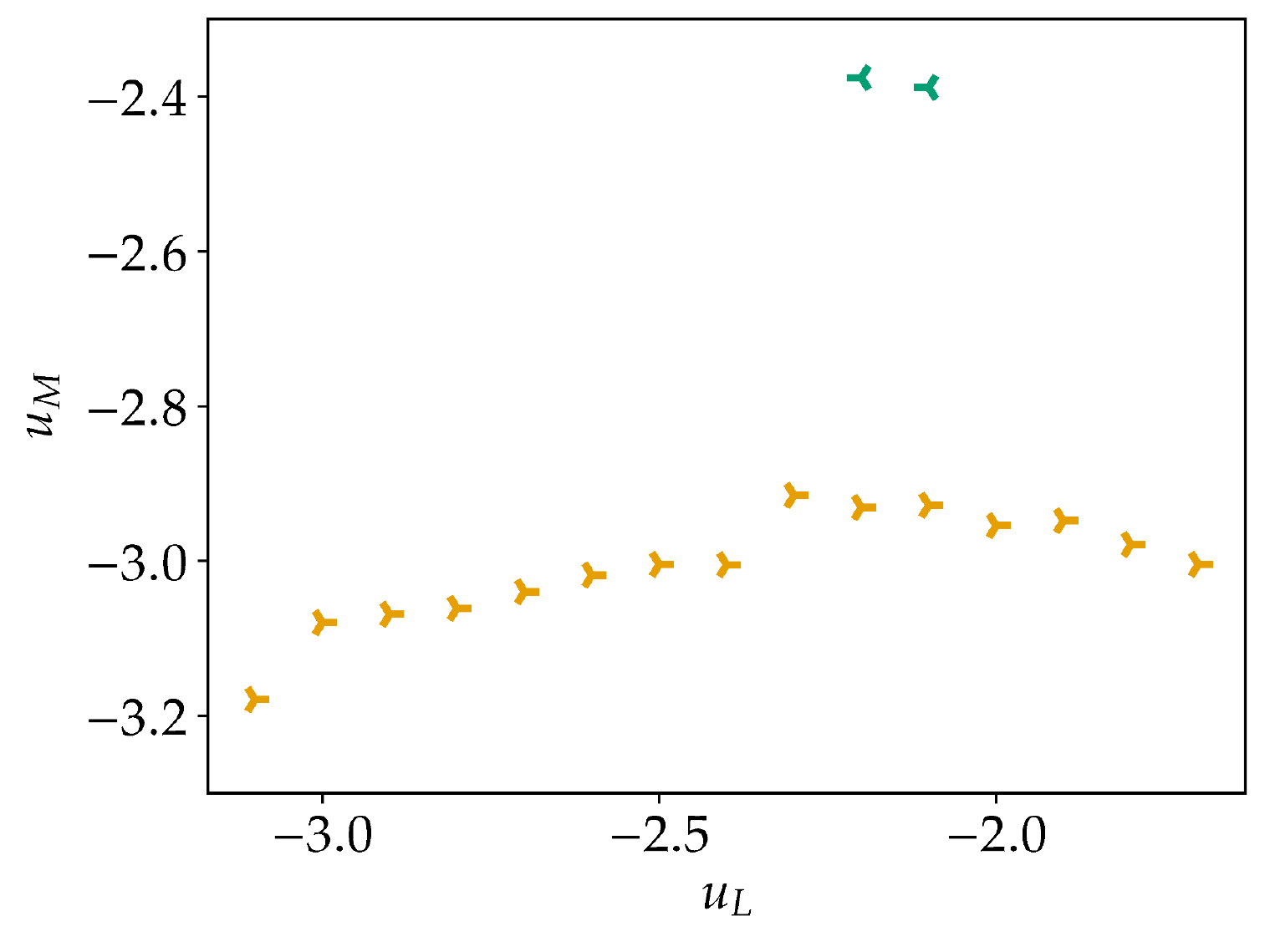}
      \caption{Convergent choice of \cite{tadmor2012adaptive}.}
    \end{subfigure}%
    \caption{Kinetic function for Fourier methods with different choices
             of the spectral viscosity and $N = \num{16384}$ grid nodes.}
    \label{fig:Quartic_Fourier_N_16384_Split_1__kinetic_function}
  \end{figure}
\end{enumerate}

\begin{remark}
\label{rem:quartic}
  Numerical experiments indicate that TeCNO methods based on the $L^2$ and
  $L^4$ entropy approximate the classical entropy solution for the setup
  considered here. However, it is unclear whether this is the case for all
  possible entropies and initial conditions.
\end{remark}

\section{Boundedness property for the Keyfitz-Kranzer system}
\label{section7}

\subsection{Preliminaries}

The Keyfitz-Kranzer system introduced in \cite{keyfitz1995spaces}
\begin{equation}
\label{eq:keyfitz-kranzer}
  \partial_t u_1 + \partial_x \bigl( u_1^2 - u_2 \bigr) = 0,
\qquad
  \partial_t u_2 + \partial_x \bigl( \nicefrac{u_1^3}{3} - u_1 \bigr) = 0,
\end{equation}
is strictly hyperbolic with eigenvalues $\lambda_\pm = u_1 \pm 1$ and right
eigenvectors $r_\pm = \big(1, u_1 \mp 1\big)^T$. Thus, it is genuinely nonlinear.
Moreover, it is straightforward to check that the convex function
$
  U(u) = \exp\bigl( \nicefrac{u_1^2}{2} - u_2 \bigr)
$
is an entropy with associated entropy flux $F(u) = u_1 U(u)$. Indeed, the entropy variables are
\begin{equation}
\label{eq:keyfitz-kranzer-w}
  w(u) = U'(u)
  =
  \begin{pmatrix}
    u_1 U(u)
    \\
    - U(u)
  \end{pmatrix}
  \quad \implies \quad
  U''(u)
  =
  \begin{pmatrix}
    1 + u_1^2 & -u_1
    \\
    -u_1 & 1
  \end{pmatrix} U(u).
\end{equation}
Thus, the flux potential is
$  \psi(u)
  =
  w(u) \cdot f(u) - F(u)
  =
  \left( \frac{2}{3} u_1^3 - u_1 u_2 \right) \exp\bigl( \nicefrac{u_1^2}{2} - u_2 \bigr).
$
Although there is a convex entropy and the system \eqref{eq:keyfitz-kranzer} is
strictly hyperbolic and genuinely nonlinear, the Riemann problem can be solved in
general only if measures are allowed in the solution \cite{keyfitz1995spaces}. These solutions
are measures in space for fixed time. Especially, singular shock waves appear in
the solution of certain Riemann problems, cf.\ \cite{lefloch-IMA-1990,keyfitz2011singular}.

In order to construct entropy-stable semi-discretizations as entropy-conservative
ones with additional dissipation, entropy-conservative numerical fluxes are sought.
Using the procedure to derive entropy-conservative fluxes described in
\cite[Procedure~4.1]{ranocha2018comparison}, the following steps have to be performed:
\begin{enumerate}
  \item
  Choose a set of variables, e.g.\ conservative variables, entropy variables,
  or something else.

  \item
  Apply scalar differential mean values for $w$, $\psi$ to get an entropy
  conservative flux fulfilling $\jump{w} \cdot \fnum = \jump{\psi}$.
\end{enumerate}

Due to the analytical form of the entropy variables $w$ \eqref{eq:keyfitz-kranzer-w}
and the flux potential $\psi$ given earlier,
the variables
$u_1$ and $U(u)$ are used in the following. Of course, other choices are also
possible.
The jumps of the entropy variables can be written using the discrete product rule
$
  \jump{a b} = \mean{a} \jump{b} + \mean{b} \jump{a},
$
with
$
  \jump{a} := a_+ - a_-$
and
  $
  \mean{a} := \frac{a_+ + a_-}{2},
$
as
\begin{equation}
  \jump{w_1} = \jump{u_1 U} = \mean{u_1} \jump{U} + \mean{U} \jump{u_1},
  \qquad
  \jump{w_2} = - \jump{U}.
\end{equation}
Using the discrete chain rule
\begin{equation}
  \jump{\log a} = \frac{1}{\logmean{a}} \jump{a},
  \qquad
  \logmean{a} := \frac{\jump{a}}{\jump{\log a}},
\end{equation}
for the logarithmic mean $\logmean{a}$ \cite{ismail2009affordable}
and $-u_2 = \log(U) - \frac{1}{2} u_1^2$, the jump of the
flux potential $\psi = \left(  \frac{1}{6} u_1^3 + u_1 \log(U) \right) U$ can be
written as
\begin{equation}
\begin{aligned}
  \jump{\psi}
  &=
  \mean{\frac{1}{6} u_1^3 + u_1 \log(U)} \jump{U}
  + \mean{U} \jump{\frac{1}{6} u_1^3 + u_1 \log(U)}
  \\
  &=
  \mean{\frac{1}{6} u_1^3 + u_1 \log(U)} \jump{U}
  + \frac{\mean{U}}{\logmean{U}} \mean{u_1} \jump{U}
  \\
  &\quad
  + \frac{1}{6} \mean{U} \frac{\jump{u_1^3}}{\jump{u_1}} \jump{u_1}
  + \mean{U} \mean{\log(U)} \jump{u_1}.
\end{aligned}
\end{equation}
Here, we have
$
  \frac{\jump{a^3}}{\jump{a}}
  =
  \frac{a_+^3 - a_-^3}{a_+ - a_-}
  =
  a_+^2 + a_+ a_- + a_-^2.
$
Therefore, the condition $\jump{w} \cdot \fnum - \jump{\psi} = 0$ \eqref{eq:fnum-EC}
for an entropy-conservative flux becomes
\begin{multline}
  0
  =
  \left(
    \mean{u_1} \fnum_1
    - \fnum_2
    - \mean{\frac{1}{6} u_1^3 + u_1 \log(U)}
    - \frac{\mean{U}}{\logmean{U}} \mean{u_1}
  \right) \jump{U}
  \\
  + \left(
    \mean{U} \fnum_1
    - \frac{1}{6} \mean{U} \frac{\jump{u_1^3}}{\jump{u_1}}
    - \mean{U} \mean{\log(U)}
  \right) \jump{u_1}.
\end{multline}
Hence,
\begin{equation}
\label{eq:keyfitz-kranzer-fnum-EC}
\begin{aligned}
  \fnum_1
  &=
  \frac{1}{6} \frac{\jump{u_1^3}}{\jump{u_1}} + \mean{\log(U)},
  \\
  \fnum_2
  &=
  \mean{u_1} \fnum_1 - \mean{\frac{1}{6} u_1^3 + u_1 \log(U)}
  - \frac{\mean{U}}{\logmean{U}} \mean{u_1},
\end{aligned}
\end{equation}
can be seen to yield an entropy-conservative numerical flux for the Keyfitz-Kranzer
system \eqref{eq:keyfitz-kranzer} with the entropy  $U(u) = \exp( \nicefrac{u_1^2}{2} - u_2 )$.
In order to create an entropy-stable numerical flux, the local Lax-Friedrichs/Rusanov
type dissipation $-\frac{\lambda}{2} \jump{u}$ will be added to the numerical
flux \eqref{eq:keyfitz-kranzer-fnum-EC}, where $\lambda = \max\set{\abs{u_{1,-}},
\abs{u_{1,+}}} + 1$ is the maximal eigenvalue of both arguments $u_\pm$.

In the following, the Riemann problem at $x=0$ with left and right initial data
\begin{equation}
  u_L = \vect{1.5 \\ 0}, \qquad u_R = \vect{-2.065426 \\ 1.410639},
\end{equation}
given in \cite[Section~4]{sanders2003numerical} will be considered in the domain
$(\xmin, \xmax) = (-3/4, 1/4)$ up to the final time
$T = 2$. Since the solution does not interact with the boundary during this time,
constant boundary values are assumed.
The continuous rate of change of the entropy is
$
  \od{}{t} \int_\xmin^\xmax U
  =
  - F \big|_\xmin^\xmax
  =
  F(g_L) - F(g_R),
$
if the left and right boundary data are $g_L, g_R$. Thus, the entropy rate is bounded
by data if $F_L \geq 0$ and $F_R \leq 0$, which is fulfilled for the given Riemann
problem, since $F(u) = u_1 \exp(u_1^2/2 - u_2)$ and $u_1 > 0$ on the left-hand side
and $u_1 < 0$ on the right-hand side.
The semi-discrete entropy rate is
\begin{equation}
  \vec{w}^T \mat{M} \partial_t \vec{u}
  =
  \left( w_L \cdot \fnum_L - \psi_L \right) - \left( w_R \cdot \fnum_R - \psi_R \right).
\end{equation}
In order to bound the semi-discrete entropy rate by the continuous one, only the
left-hand side is considered in the following. The other side can be handled
similarly. Here, the condition
\begin{equation}
\label{eq:bc-condition-necessary}
  \left( w_L \cdot \fnum_L - \psi_L \right) \leq F(g_L) = w(g_L) \cdot f(g_L) - \psi(g_L)
\end{equation}
is required for the semi-discretization. If the numerical flux $\fnum$ is entropy
stable,
$\bigl( w_L - w(g_L) \bigr) \cdot \fnum_L - \bigl( \psi_L - \psi(g_L) \bigr) \leq 0$.
Thus, \eqref{eq:bc-condition-necessary} is fulfilled if
$
  w(g_L) \cdot \fnum_L - \psi(g_L)
  \leq
  w(g_L) \cdot f(g_L) - \psi(g_L),
$
which is equivalent to
\begin{equation}
\label{eq:bc-condition-sufficient}
  w(g_L) \cdot \left( \fnum(g_L, u_L) - f(g_L) \right) \leq 0.
\end{equation}
It seems to be technically difficult to check this condition for the entropy
stable flux using the local Lax-Friedrichs-Rusanov type dissipation. However,
it can be checked during the numerical experiments. There, it is fulfilled up
to machine accuracy in the Riemann problems under consideration.

\subsection{Numerical results}

\begin{figure}[!htb]
\centering
  \begin{subfigure}{0.45\textwidth}
    \centering
    \includegraphics[width=\textwidth]{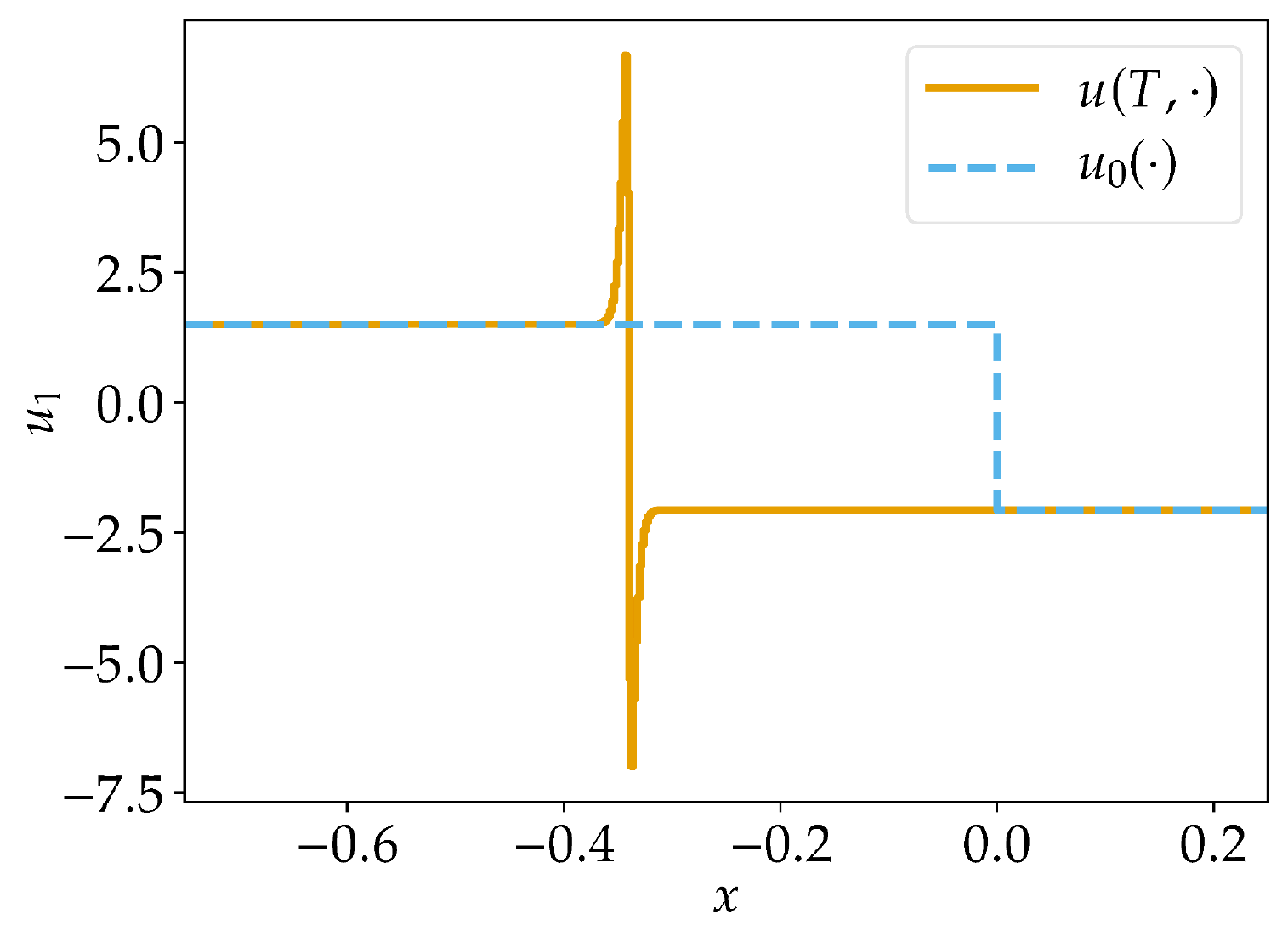}
    \caption{First component $u_1$.}
  \end{subfigure}%
  \hspace*{\fill}
  \begin{subfigure}{0.45\textwidth}
    \centering
    \includegraphics[width=\textwidth]{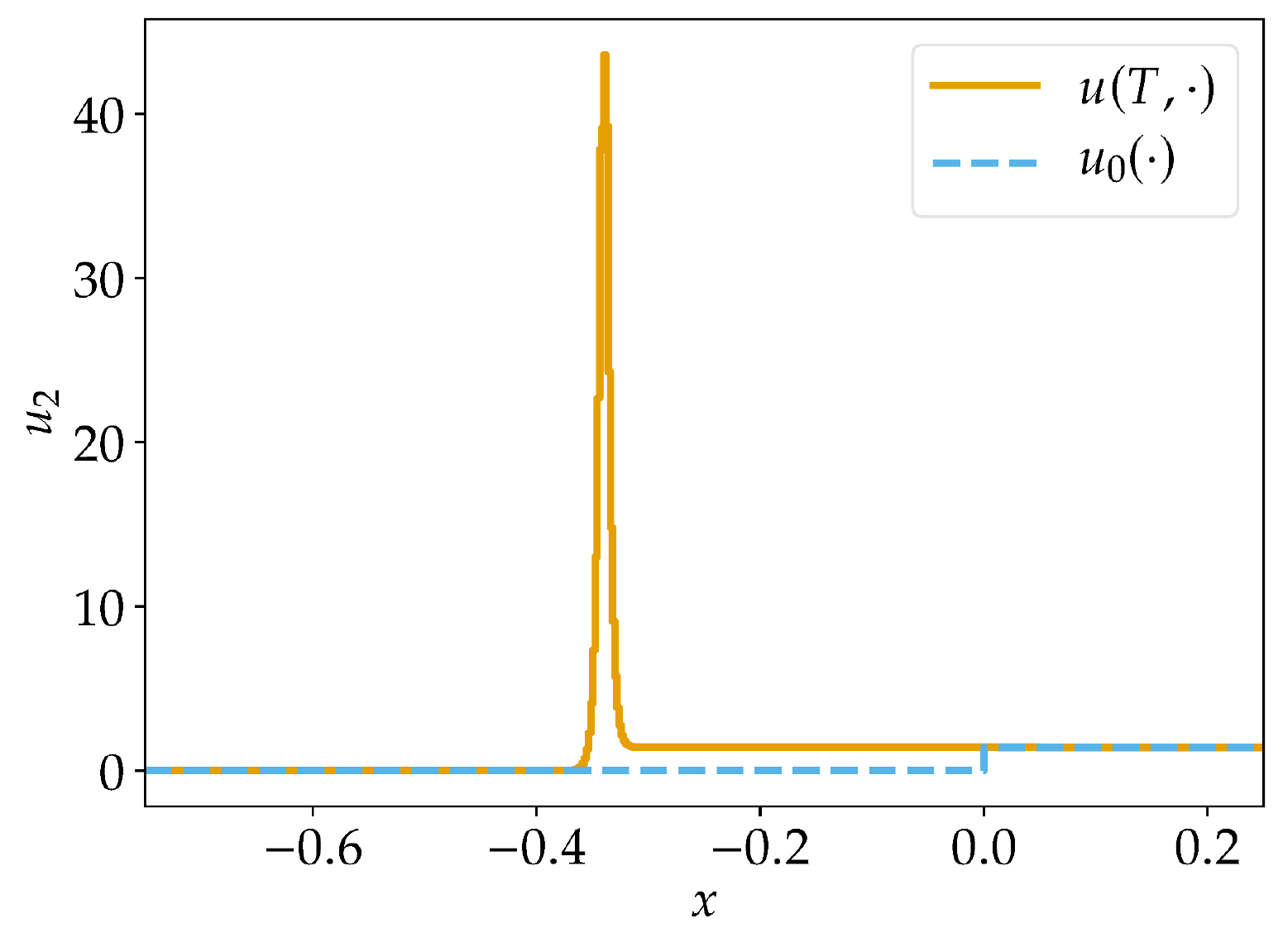}
    \caption{Second component $u_2$.}
  \end{subfigure}%
  \caption{Numerical solutions using the first order finite volume method with
           $N = 512$ cells at the final time $t=T$.}
  \label{fig:KeyfitzKranzer_solution}
\end{figure}

\begin{figure}[!htb]
\centering
  \includegraphics[width=0.45\textwidth]{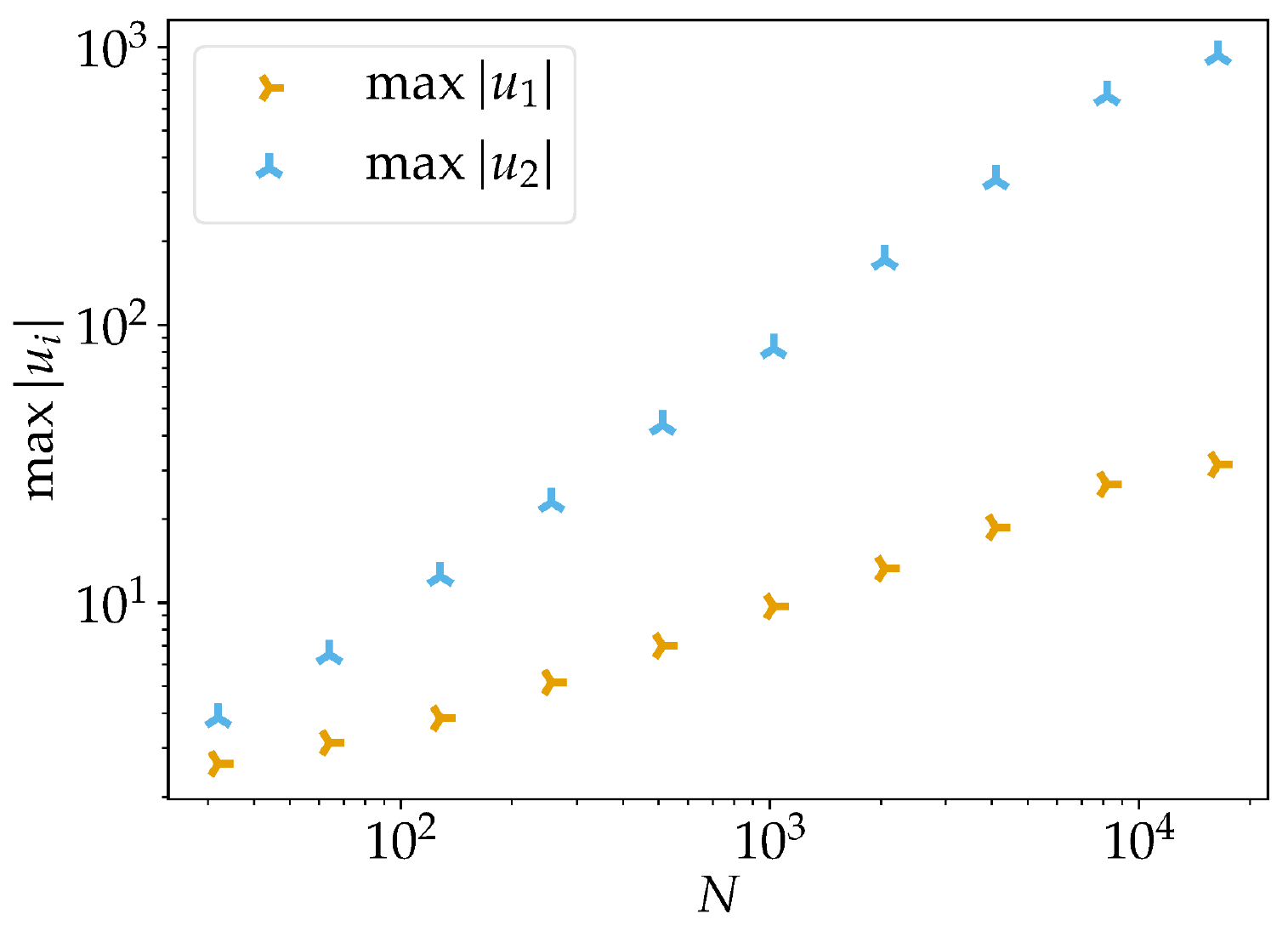}
  \caption{Extreme values of the numerical solutions using the first order finite
           volume method with $N$ cells at the final time $t=T$.}
  \label{fig:KeyfitzKranzer_extreme_values}
\end{figure}

First-order finite volume discretizations based on the entropy-stable flux above are now applied.
An explicit Euler method is used in order to advance the
numerical solution in time, up to $T=2$, and the time step is chosen adaptively as
$\Delta t = \frac{\Delta x}{2\lambda}$, where $\lambda = \max\set{\abs{u_{i,1}} + 1}$
is the maximal advection speed of the numerical solution.
A typical solution at $t=T$ is displayed in \autoref{fig:KeyfitzKranzer_solution}.
A singular shock occurs and moves to the left, cf.\ \cite[Section~4]{sanders2003numerical}.
The maximal values of each component increase with
increasing resolution, since the analytical solution contains a Dirac mass moving with the shock.
This behavior is visualized in \autoref{fig:KeyfitzKranzer_extreme_values}.
The maximal values of the numerical solutions increase over time, in agreement
with the increasing mass of the delta measure. Clearly, in
general, entropy stability does not imply boundedness of numerical solutions.

\section{Summary and conclusions}
\label{sec:summary}

We studied a variety of entropy-dissipative numerical methods for nonlinear scalar
conservation laws with non-convex flux in one space dimension, including finite difference methods
with artificial dissipation, discontinuous Galerkin methods with or without
filtering, and Fourier methods with (super-) spectral viscosity. We demonstrated
experimentally that all these numerical methods can converge to different weak solutions,
either the classical entropy solution or weak solutions involving nonclassical
shock waves. For the same class of methods, the convergence depends also on
associated parameters; changing a parameter such as the order or strength of
dissipation changes the limiting solution in general.
Moreover, we demonstrated that numerical solutions can also depend on the
specific choice of the entropy function, e.g.\ for TeCNO schemes.
To distinguish the different convergence behavior of the numerical methods,
we also computed the associated kinetic functions for a variety of schemes,
including also high-order WENO methods.
Finally, we developed entropy-dissipative numerical methods for the Keyfitz-Kranzer
system, demonstrating that entropy-dissipative methods may generate numerical
solutions that do not remain bounded under grid refinement.

Our results provide important contributions to the theory of nonclassical shocks
by comparing a variety of numerical schemes and their kinetic functions.
On the other hand, these results also demonstrate limitations of modern high-order
entropy-dissipative methods, complementing the recent investigations
\cite{gassner2020stability,ranocha2020preventing}. Hence, we stress the importance
of choosing appropriate regularizations to compute numerical solutions of
conservation laws, in particular if they support nonclassical shock waves
as weak solutions.


\appendix

\section*{Acknowledgments}

The authors would like to thank David Ketcheson for very interesting discussions.
This research work was supported by the King Abdullah University of Science and
Technology (KAUST). The first author (PLF) was partially supported by the
Innovative Training Network (ITN) grant 642768 (ModCompShocks).
The second author (HR) was partially supported by the German Research Foundation
(DFG, Deutsche Forschungsgemeinschaft) under Grant SO~363/14-1.

\printbibliography

\end{document}